\documentclass[12pt]{article}%
\usepackage{graphicx}
\usepackage{amsmath}
\usepackage{amsfonts}
\usepackage{amssymb}%
\usepackage{amsthm}

\providecommand{\U}[1]{\protect\rule{.1in}{.1in}}
\newtheorem{theorem}{Theorem}[section]

\newtheorem{corollary}[theorem]{Corollary}

\newtheorem{definition}[theorem]{Definition}
\newtheorem{example}[theorem]{Example}

\newtheorem{lemma}[theorem]{Lemma}

\newtheorem{problem}[theorem]{Problem}
\newtheorem{proposition}[theorem]{Proposition}
\newtheorem{remark}[theorem]{Remark}

\newtheorem*{quote1}{Theorem \ref{mainresult}}
\begin{document}

\author{Peter Scott\thanks{Partially supported by NSF grants DMS 9626537 and 0203883}\\Mathematics Department\\University of Michigan\\Ann Arbor, Michigan 48109, USA.\\email:pscott@umich.edu
\and Gadde A. Swarup\\718 High Street Road\\Glen Waverley\\Victoria 3150, Australia.\\email: anandaswarupg@gmail.com }
\title{Canonical decompositions for Poincar\'{e} duality pairs}
\maketitle

\begin{abstract}
The authors previously described an algebraic analogue of the
JSJ--decomposition of a $3$--manifold. This analogue is defined for any
finitely presented, one-ended group. We study this analogue in the special
case of Poincar\'{e} duality pairs.

\end{abstract}
\date{}
\tableofcontents

\begin{center}
\textit{Dedicated to Terry Wall}
\end{center}

\textbf{{\Large Introduction}}

In \cite{SS2}, as corrected in \cite{SS2errata}, we obtained canonical
decompositions for almost finitely presented groups analogous to the
$JSJ$--decomposition of a $3$--manifold. In particular, for many almost
finitely presented groups $G$, and any integer $n\geq1$, we defined a
decomposition $\Gamma_{n,n+1}(G)$, and we showed that when $G$ is the
fundamental group of an orientable Haken $3$--manifold $M$ with incompressible
boundary, then $\Gamma_{1,2}(G)$ essentially yields the $JSJ$--decomposition
of $M$. Further details are discussed in \cite{SS4}. We recall that the
$JSJ$--decomposition of $M$ is given by a possibly disconnected compact
submanifold $V(M)$ which is called the characteristic submanifold of $M$, such
that the frontier of $V(M)$ consists of essential annuli and tori, and each
component of $V(M)$ is an $I$--bundle or a Seifert fibre space. Further each
component of the closure of $M-V(M)$ is simple. The frontier of $V(M)$
determines a graph of groups structure for $G=\pi_{1}(M)$, in which all edge
groups are free abelian of rank $1$ or $2$, and this is essentially the same
as $\Gamma_{1,2}(G)$.

In this paper, we consider the structure of $\Gamma_{n,n+1}(G)$, in the case
of Poincar\'{e} duality pairs $(G,\partial G)$ of dimension $n+2$, where
$n\geq1$. The results we obtain are very closely analogous to the above
description for $3$--manifolds. This greatly generalises work of Kropholler in
\cite{K} and of Castel in \cite{Castel}. Kropholler described a canonical
decomposition of such Poincar\'{e} duality pairs in any dimension at least
three, and Castel described a canonical decomposition of such Poincar\'{e}
duality pairs in dimension three only, but their decompositions are analogous
to the Torus Decomposition of an orientable Haken $3$--manifold. See section
\ref{torusdecompofPDgroup} for a discussion.

A Poincar\'{e} duality pair is the algebraic analogue of an aspherical
manifold with aspherical boundary components whose fundamental groups inject.
It consists of a group $G$ which corresponds to the fundamental group of the
manifold, and a family $\partial G$ of subgroups which correspond to the
fundamental groups of the boundary components, and the whole setup satisfies
an appropriate version of Poincar\'{e} duality. The main difficulty in
establishing the results in this paper is that if one considers any of the
decompositions of $G$ described in \cite{SS2}, then a priori there is no
connection between the decomposition and the boundary groups of the pair. In
the topological case, if one considers the full characteristic submanifold
$V(M)$ of an orientable Haken manifold $M$ with incompressible boundary, one
can double $M$ along its boundary to obtain a closed Haken $3$--manifold $DM$,
and there is a natural submanifold $DV$ of $DM$ which is the double of $V(M)$.
Further, in most cases, $DV$ is the characteristic submanifold of $DM.$ In the
algebraic context, doubling a $(n+2)$--dimensional Poincar\'{e} duality pair
$(G,\partial G)$ along its boundary yields a $(n+2)$--dimensional Poincar\'{e}
duality group $DG$, but in general there is no natural way to double an
algebraic decomposition of $G$ to obtain a corresponding decomposition of
$DG$. However, after establishing all the properties of the decomposition
$\Gamma_{n,n+1}(G)$, when $(G,\partial G)$ is a Poincar\'{e} duality pair, we
will show in section \ref{comparingdecompositions} that the algebraic
situation is very similar to the topological one. In the topological setting,
one can also reverse this process and construct the full characteristic
submanifold of an orientable Haken manifold $M$ with incompressible boundary
by starting with the characteristic submanifold $V(DM)$ of $DM$ and
\textquotedblleft undoubling\textquotedblright\ to obtain the required
submanifold $V(M)$ of $M$. This greatly simplifies the construction of the
characteristic submanifold of $M$. If we start with the Poincar\'{e} duality
group $DG$, and the decomposition $\Gamma_{n+1}(DG)$, then the natural
algebraic analogue of \textquotedblleft undoubling\textquotedblright\ is
simply to restrict this decomposition to $G$ using the Subgroup Theorem. This
determines a graph of groups structure on $G$, and it follows from our results
in section \ref{comparingdecompositions} that this decomposition of $G$ is
$\Gamma_{n,n+1}(G)$. However this does not simplify the proof that
$\Gamma_{n,n+1}(G)$ has the properties we require. In fact, the proof of this
"undoubling" result in section \ref{comparingdecompositions} depends on first
establishing all the properties of the decomposition $\Gamma_{n,n+1}(G)$. The
difficulty is that we are unable to show directly that the decomposition of
$G$ obtained by restricting $\Gamma_{n+1}(DG)$ to $G$ has any of the enclosing
properties required of $\Gamma_{n,n+1}(G)$. This is partly because our idea of
enclosing is much stronger than simply requiring certain subgroups of $G$ to
be conjugate into certain vertex groups of a given decomposition of $G$. For
example if a splitting of $G$ is enclosed by a vertex $v$ of a graph of groups
decomposition $\Gamma$ of $G$, then we can split $\Gamma$ at $v$ to obtain a
refined graph of groups structure for $G$ which has an extra edge associated
to the given splitting.

An important point about our ideas in \cite{SS2} and in this paper is that we
consider all almost invariant subsets of a group $G$ rather than just those
which correspond to splittings. In the topological setting, this corresponds
to considering essential maps of codimension--$1$ manifolds rather than just
essential embeddings. We do not know how to carry out the program in this
paper using only splittings. The main goal in \cite{SS2} was to enclose the
algebraic analogues of immersed annuli and tori (the analogues were almost
invariant sets over virtually polycyclic groups), which is the natural
generalization of the approaches in \cite{NS} and \cite{Sc2}. It turns out
that the analogy is stronger in the case of Poincar\'{e} duality pairs. The
decomposition $\Gamma_{1,2}(G)$ which we obtained in \cite{SS2} was
constructed to enclose the analogues of immersions of annuli and tori whereas
in $JSJ$ theory \cite{JS, JO, Waldhausen}, the aim was to enclose essential
Seifert pairs. When $G$ is the fundamental group of a $3$--manifold, the
difference turns out to be minor (consisting of small Seifert fibre spaces)
and one can easily go from one decomposition to the other. In the case of
groups it seems more natural to enclose almost invariant sets over virtually
polycyclic subgroups, that is the analogues of immersions of annuli and tori,
rather than Seifert pairs and to make the distinction clear, we will call the
decompositions that we obtain Annulus--Torus decompositions. It should be
pointed out that when doubling a manifold, the $JSJ$--decompositions behave
better than the Annulus--Torus decompositions, and we will switch from one to
the other when it is convenient. Similar comments apply to our algebraic
decompositions of Poincar\'{e} duality pairs. This is made precise in Theorem
\ref{doublingtheorem} and Remark \ref{doublingremark}.

Our main result, Theorem \ref{mainresult}, is a description of the
decomposition $\Gamma_{n,n+1}(G)$, and its completion $\Gamma_{n,n+1}^{c}(G)$,
for an orientable $PD(n+2)$ pair $(G,\partial G)$, when $n\geq1$. We leave the
precise statement till later because it requires the introduction of a
substantial amount of terminology. Some of the result is simply a restatement
of our results from \cite{SS2}. But an important part of the result is that
all the edge splittings of $\Gamma_{n,n+1}(G)$ are induced by essential annuli
and tori in $(G,\partial G)$, in a sense which we define in section
\ref{essentialannuli}. This means that, for each $n\geq1$, the decomposition
$\Gamma_{n,n+1}(G)$ of $G$ is closely analogous to the $JSJ$--decomposition of
an orientable $3$--manifold.

In section \ref{preliminaries}, we recall the main definitions and results
from our work in \cite{SS2}, as corrected in \cite{SS2errata}, and from our
more recent paper with Guirardel \cite{GSS}. In section \ref{essentialannuli},
we discuss the definition of a Poincar\'{e} duality group and pair, and then
we discuss essential annuli and tori (the terms `annulus' and `torus' are used
in a generalised sense) in orientable Poincar\'{e} duality pairs. We show that
each essential annulus and torus has a naturally associated almost invariant
set, which we call its dual. We also show that all almost invariant sets over
$VPCn$ groups in an orientable $PD(n+2)$ pair are `generated' by duals of
essential annuli. In section \ref{themaintheorem}, we give several more
definitions which finally allow us to state our main theorem, Theorem
\ref{mainresult}. In section \ref{torusdecompofPDgroup}, we discuss the
analogues of torus decompositions for orientable $PD(n+2)$ groups and pairs.
These decompositions were already obtained by Kropholler \cite{K} under the
extra hypothesis that any $VPC$ subgroup has finitely generated centraliser, a
condition which he called Max-c. In \cite{K2} Kropholler showed that the Max-c
condition holds in dimension three. However an example due to Mess \cite{Mess}
shows that the Max-c condition is not always satisfied. In dimension three,
this decomposition was also obtained by Castel \cite{Castel}. The comparison
between our results and Kropholler's results in the case of orientable
$PD(n+2)$ groups is discussed briefly in \cite{SS2}, but we discuss this in
more detail here. In section \ref{furtherpropertiesoftorusdecomp}, we analyse
further our torus decomposition of Poincar\'{e} duality pairs. In section
\ref{enclosingpropertiesofATdecomps}, we continue studying orientable
$PD(n+2)$ pairs, and consider the crossing of almost invariant subsets over
$VPCn$ groups with almost invariant subsets over $VPC(n+1)$ groups. This is a
new feature of our arguments in this paper, which could not be handled in the
more general setting of \cite{SS2}. In section \ref{proofofmaintheorem}, we
bring together the various pieces and prove our main result, Theorem
\ref{mainresult}. In section \ref{comparingdecompositions}, we are able to
prove Theorem \ref{doublingtheorem} which shows that, for an orientable
$PD(n+2)$ pair $(G,\partial G)$, one can double the decomposition
$\Gamma_{n,n+1}^{c}(G)$ to obtain $\Gamma_{n+1}^{c}(DG)$. Finally in section
\ref{concludingremarks}, we discuss some natural further questions.

This paper is a revised version of \cite{SS-pdnarxiv}. The main changes are
that sections 2 and 3 of that paper have been removed, as the theory therein
has now been developed more thoroughly and generally in \cite{GSS}, and a new
section, numbered \ref{comparingdecompositions}, has been added. In addition,
there are several minor corrections and improvements in the exposition.

\section{Preliminaries\label{preliminaries}}

In this section we recall the main definitions and results from \cite{SS2}, as
corrected in \cite{SS2errata}, which we will use. But we will start by briefly
discussing some $3$--manifold theory which motivates and guides all our work.
Let $M$ be an orientable Haken $3$--manifold with incompressible boundary.
Jaco and Shalen \cite{JS} and Johannson \cite{JO} proved the existence and
uniqueness of the characteristic submanifold of $M$. We will denote this
submanifold by $JSJ(M)$. Its frontier consists of essential annuli and tori in
$M$, and each component of $JSJ(M)$ is a Seifert fibre space or $I$--bundle.
Further any essential map of a Seifert fibre space into $M$ can be properly
homotoped to lie in $JSJ(M)$, and this condition characterizes $JSJ(M)$. In
order to compare this with algebraic generalisations, we note that, in
particular, any essential map of the annulus or torus into $M$ can be properly
homotoped to lie in $JSJ(M)$. This weaker condition does not characterise
$JSJ(M)$, but does characterise a submanifold of $M$ which we denote by
$AT(M)$. The letters $AT$ stand for Annulus--Torus. This is discussed in
detail in chapter 1 of \cite{SS2}, but the notation $AT(M)$ is not used. Any
essential map of the annulus or torus into $M$ can be properly homotoped to
lie in $AT(M)$, and $AT(M)$ is minimal, up to isotopy, among all essential
submanifolds of $M$ with this property. (A compact submanifold of $M$ is
essential if its frontier consists of essential embedded surfaces.) We will
say that the family of all essential annuli and tori in $M$ \textit{fills}
$AT(M)$, and we regard $AT(M)$ as a kind of regular neighbourhood of this
family. The connection between $AT(M)$ and $JSJ(M)$ can be described as
follows. The submanifold $JSJ(M)$ has certain exceptional components. These
are of two types. One type is a solid torus $W$ whose frontier consists of
three annuli each of degree $1$ in $W$, or of one annulus of degree $2$ in
$W$, or of one annulus of degree $3$ in $W$. The other type lies in the
interior of $M$ and is homeomorphic to the twisted $I$--bundle over the Klein
bottle. (Note that as $M$ is orientable, only one such bundle can occur.) Then
$AT(M)$ can be obtained from $JSJ(M)$ by discarding all these exceptional
components, replacing each of them by a regular neighbourhood of its frontier,
and finally discarding any redundant product components from the resulting submanifold.

Recall from the previous paragraph that the family of all essential annuli and
tori in $M$ fills $AT(M)$, and that we regard $AT(M)$ as a kind of regular
neighbourhood of this family. However, $AT(M)$ is filled by a smaller family
of essential annuli and tori in $M$, and this turns out to be crucial for the
algebraic analogues we are discussing in this paper. Let $AT_{\partial}(M)$
denote the union of those components of $AT(M)$ which meet $\partial M$, and
let $AT_{int}(M)$ denote the union of the remaining components of $AT(M)$.
Then any essential annulus in $M$ can be properly homotoped to lie in
$AT_{\partial}(M)$. Thus it is clear that any essential torus in $M$ which is
homotopic into $AT_{int}(M)$ cannot cross any such annulus, and that
$AT_{int}(M)$ must be filled by tori which are homotopic into $AT_{int}(M)$.
Further, it is easy to show that $AT_{\partial}(M)$ is filled by the family of
all essential annuli in $M$. We conclude that $AT(M)$ is filled by the family
of all essential annuli in $M$ together with those essential tori in $M$ which
do not cross any essential annulus in $M$.\label{AT(M)}

For future reference, we will also need to discuss the Torus Decomposition of
$M$ and its relationship with $AT(M)$. As above one can characterise a
submanifold $T(M)$ of $M$ by the property that any essential map of a torus
into $M$ can be homotoped into $T(M)$ and that $T(M)$ is minimal, up to
isotopy, among all essential submanifolds of $M$ with this property. Of course
if $M$ admits no essential annulus, then $AT(M)$ and $T(M)$ are equal. In
general, $T(M)$ is obtained from $AT(M)$ as follows. Any component of
$AT_{int}(M)$ is left unchanged. Now $AT_{\partial}(M)$ has three types of
component. The first type is a Seifert fibre space which is not a solid torus
and such that each boundary torus lies in the interior of $M$, or is contained
in $\partial M$, or meets $\partial M$ in vertical annuli. The second type is
a solid torus which meets $\partial M$ in annuli, and the third type is an
$I$--bundle over a surface $F$ which meets $\partial M$ in the associated
$\partial I$--bundle over $F$. As no essential torus in $M$ can be homotopic
into a component of $AT_{\partial}(M)$ of the second or third type, all such
components are omitted when we form $T(M)$. Finally let $W$ denote a component
of $AT_{\partial}(M)$ which is of the first type. Thus $W$ is a Seifert fibre
space which is not a solid torus and $W\cap\partial M$ consists of tori and
vertical annuli in $\partial W$. By pushing into the interior of $W$ each
torus component of $\partial W$ which meets $\partial M$ in annuli, we obtain
a Seifert fibre space $W^{\prime}$ which is contained in and homeomorphic to
$W$. Note that the components of the closure of $W-W^{\prime}$ are
homeomorphic to $T\times I$, and $W^{\prime}\cap\partial M$ consists only of
tori. Replacing $W$ by $W^{\prime}$ for each such component of $AT_{\partial
}(M)$ finally yields $T(M)$. Note that it is clear that $T(M)\subset AT(M)$
and that an essential torus in $M$ is homotopic into $AT(M)$ if and only if it
is homotopic into $T(M)$.

Next we recall the cohomological formulation of the theory of ends and of
almost invariant subsets of a group. Let $G$ be a group and let $E$ be a set
on which $G$ acts on the right. Let $PE$ denote the power set of $E$. Under
Boolean addition (\textquotedblleft symmetric difference\textquotedblright)
this is an additive group of exponent $2$. Write $FE$ for the additive
subgroup of finite subsets. We refer to two sets $A$ and $B$ whose symmetric
difference lies in $FE$ as \textit{almost equal}, and write $A\overset{a}{=}%
B$. This amounts to equality in the quotient group $PE/FE$. Now define
\[
QE=\{A\subset E:\forall g\in G,\;A\overset{a}{=}Ag\}.
\]
The action of $G$ on $PE$ by right translation preserves the subgroups
$QE\ $and $FE$, and $QE/FE$ is the subgroup of elements of $PE/FE$ fixed under
the induced action. Elements of $QE$ are said to be \textit{almost invariant}.
If we take $E$ to be $G$ with the action of $G$ being right multiplication,
then the number of ends of $G$ is
\[
e(G)=\dim_{\mathbb{Z}_{2}}\;(QG/FG).
\]
If $G$ is finite, all subsets are finite and clearly $e(G)=0$. Otherwise, $G$
is an infinite set which is invariant (not merely \textquotedblleft
almost\textquotedblright), so $e(G)\geq1$.

If $H$ is a subgroup of $G$, and we take $E$ to be the coset space
$H\backslash G$ of all cosets $Hg$, still with the action of $G$ being right
multiplication, then the number of ends of the pair $(G,H)$ is
\[
e(G,H)=\dim_{\mathbb{Z}_{2}}\;\left(  \frac{Q(H\backslash G)}{F(H\backslash
G)}\right)  .
\]

When $H$ is trivial, so that $e(G,H)=e(G)$, this can be formulated in terms of
group cohomology as follows. The abelian group $PG$ is naturally a (right)
$\mathbb{Z}_{2}G$--module, and the submodule $FG$ can be identified with the
group ring $\mathbb{Z}_{2}G$. Thus the invariant subgroup $QG/FG$ equals
$H^{0}(G;PG/\mathbb{Z}_{2}G)$. Now the short exact sequence of coefficients%
\[
0\rightarrow\mathbb{Z}_{2}G\rightarrow PG\rightarrow PG/\mathbb{Z}%
_{2}G\rightarrow0
\]
yields the following long exact cohomology sequence.%
\[
H^{0}(G;\mathbb{Z}_{2}G)\rightarrow H^{0}(G;PG)\rightarrow H^{0}%
(G;PG/\mathbb{Z}_{2}G)\overset{\delta}{\rightarrow}H^{1}(G;\mathbb{Z}%
_{2}G)\rightarrow H^{1}(G;PG)\rightarrow
\]

For any group $G$, the group$\,H^{n}(G;PG)$ is zero if $n\neq0$, and
isomorphic to $\mathbb{Z}_{2}\ $when$\ n=0$. And if $G$ is infinite, then
$H^{0}(G;\mathbb{Z}_{2}G)=0$. Also when $G$ is infinite, the non-zero element
of $H^{0}(G;PG)$ maps to the element of $H^{0}(G;PG/\mathbb{Z}_{2}G)$ which
corresponds to the equivalence class of $G$ in $QG/FG$, under the
identification of these two groups. It follows that when $G$ is infinite, the
group $H^{1}(G;\mathbb{Z}_{2}G)$ can be identified with the collection of all
almost invariant subsets of $G$ modulo almost equality and complementation.

If $H\ $is nontrivial, let $E$ denote the coset space $H\backslash G$. The
abelian group $PE$ is naturally a (right) $\mathbb{Z}_{2}G$--module, and we
denote the submodule $FE$ by $\mathbb{Z}_{2}E$. Thus the invariant subgroup
$QE/FE$ equals $H^{0}(G;PE/\mathbb{Z}_{2}E)$. Now the short exact sequence of
coefficients%
\[
0\rightarrow\mathbb{Z}_{2}E\rightarrow PE\rightarrow PE/\mathbb{Z}%
_{2}E\rightarrow0
\]
yields the following long exact cohomology sequence.%
\[
H^{0}(G;\mathbb{Z}_{2}E)\rightarrow H^{0}(G;PE)\rightarrow H^{0}%
(G;PE/\mathbb{Z}_{2}E)\overset{\delta}{\rightarrow}H^{1}(G;\mathbb{Z}%
_{2}E)\rightarrow H^{1}(G;PE)\rightarrow\label{coboundarymap}%
\]

For any group $G$ and any subgroup $H$, the group$\,H^{n}(G;PE)$ is isomorphic
to $H^{n}(H;\mathbb{Z}_{2})$, and if the index of $H$ in $G$ is infinite, then
$H^{0}(G;\mathbb{Z}_{2}E)=0$. Thus, as above, when $H$ has infinite index in
$G$, the image of the coboundary map $H^{0}(G;PE/\mathbb{Z}_{2}E)\rightarrow
H^{1}(G;\mathbb{Z}_{2}E)$ can be identified with the collection of all almost
invariant subsets of $H\backslash G$ modulo almost equality and complementation.

If $G$ is finitely presented, this description connects very nicely with
topology. Let $X$ be an Eilenberg-MacLane space $K(G,1)$ with finite
$2$--skeleton, let $\widetilde{X}$ denote its universal cover, and let $X_{H}$
denote the cover of $X$ with fundamental group $H$. Thus $X_{H}$ is the
quotient of $\widetilde{X}$ by the (left) action of $H$ acting as a covering
group. Because $X\ $has finite $2$--skeleton, the part of the long exact
cohomology sequence shown above is isomorphic to the corresponding part of the
following long exact cohomology sequence for $X_{H}$,%
\[
H_{f}^{0}(X_{H};\mathbb{Z}_{2})\rightarrow H^{0}(X_{H};\mathbb{Z}%
_{2})\rightarrow H_{e}^{0}(X_{H};\mathbb{Z}_{2})\overset{\delta}{\rightarrow
}H_{f}^{1}(X_{H};\mathbb{Z}_{2})\rightarrow H^{1}(X_{H};\mathbb{Z}%
_{2})\rightarrow
\]
where $H_{f}^{i}(X_{H};\mathbb{Z}_{2})$ denotes cellular cohomology with
finite supports. Thus when $H$ has infinite index in $G$, the image of the
coboundary map\linebreak$H_{e}^{0}(X_{H};\mathbb{Z}_{2})\rightarrow H_{f}%
^{1}(X_{H};\mathbb{Z}_{2})$ can be identified with the collection of all
almost invariant subsets of $H\backslash G$ modulo almost equality and
complementation. If $G$ is finitely generated but not finitely presented, we
can take $X$ to be a $K(G,1)$ with finite $1$--skeleton but it is no longer
correct to identify $H^{1}(G,\mathbb{Z}_{2}E)$ with $H_{f}^{1}(X_{H}%
;\mathbb{Z}_{2})$. In fact, cellular cohomology with finite supports for a
cell complex which is not locally finite is an unreasonable idea, as the
coboundary of a finite cochain need not be finite. However it is easy to
define a modified version of this theory in our particular setting. For future
reference we set this out as a remark.

\begin{remark}
\label{defnofH1fwhenGisnotfp} Let $G$ be a group which is finitely generated
but need not be finitely presented, and let $X$ be a $K(G,1)$ with finite
$1$--skeleton. Let $H$ be a subgroup of $G$, and let $X_{H}$ denote the cover
of $X$ with fundamental group $H$. We replace the use of cochains on $X_{H}$
with finite support by cochains whose support consists of only finitely many
cells above each cell of $X$. Note that as $X$ has only finitely many
$1$--cells, such $1$--cochains are finite. In this paper, $H_{f}^{1}%
(X_{H};\mathbb{Z}_{2})$ will denote the appropriate cohomology group of this
cochain complex. This enables us to identify $H^{1}(G,\mathbb{Z}_{2}E)$ with
$H_{f}^{1}(X_{H};\mathbb{Z}_{2})$.

This is the topological formulation of Lemma 7.4 in \cite{Brown}.
\end{remark}

Recall that the invariant subgroup $QE/FE$ equals $H^{0}(G;PE/\mathbb{Z}%
_{2}E)$. Thus the elements of $H^{0}(G;PE/\mathbb{Z}_{2}E)$ are equivalence
classes of almost invariant subsets of $H\backslash G$ under the equivalence
relation of almost equality. Also the elements of $H_{e}^{0}(X_{H}%
;\mathbb{Z}_{2})$ are equivalence classes of cellular $0$--cochains on $X_{H}$
which have finite coboundary. The support of such a cochain is a subset of the
vertex set of $X_{H}$. Thus whether or not $G\ $is finitely presented, the
isomorphism between $H^{0}(G;PE/\mathbb{Z}_{2}E)$ and $H_{e}^{0}%
(X_{H};\mathbb{Z}_{2})$ associates to an almost invariant subset $Y$ of
$H\backslash G$ a subset $Z$ of the vertex set of $X_{H}$ with finite
coboundary, where $Z$ is unique up to almost equality. This is a convenient
fact which we will use on several occasions. Note that $Y$ is trivial, i.e.
finite or co-finite, if and only if $Z$ is finite or co-finite. Also if
$Y\ $and $Y^{\prime}$ are almost invariant subsets of $H\backslash G$ with
corresponding subsets $Z$ and $Z^{\prime}$ of the vertex set of $X_{H}$, then
the intersections $Y\cap Y^{\prime}$ and $Z\cap Z^{\prime}$ also correspond.

Next we recall some more basic facts about almost invariant sets. If $X$ and
$Y$ are subsets of $G$, the four sets $X\cap Y$, $X\cap Y^{\ast}$, $X^{\ast
}\cap Y$ and $X^{\ast}\cap Y^{\ast}$ are called the \textit{corners} of the
pair $(X,Y)$. If $X$ is $H$--almost invariant and $Y$ is $K$--almost
invariant, a corner of the pair $(X,Y)$ is \textit{small} if it is $H$--finite
or $K$--finite. (These two conditions are equivalent so long as $X$ and $Y$
are both nontrivial and $G$ is finitely generated.) We say that $X$
\textit{crosses} $Y$ if all four corners of the pair $(X,Y)$ are
$K$--infinite. The preceding parenthetical comment shows that if $X$ and $Y$
are nontrivial and $G$ is finitely generated, then $X$ crosses $Y$ if and only
$Y$ crosses $X$. In \cite{SS1} and \cite{SS2}, we defined a partial order
$\leq$ on certain families of almost invariant subsets of a finitely generated
group $G$ as follows. The idea of the definition is that $Y\leq X$ means that
$Y$ is \textquotedblleft almost\textquotedblright\ contained in $X$. If $E$ is
a family of almost invariant subsets of $G$, we say that the elements of $E$
are in \textit{good position} if whenever $U$ and $V$ are elements of $E$ such
that two of the corners of the pair are small, then one corner is empty. If
the elements of $E$ are in good position, then we defined $Y\leq X$ to mean
that either $Y\cap X^{\ast}$ is empty or it is the only small corner of the
pair $(X,Y)$. We showed that $\leq$ is a partial order on $E$. Note that if
$Y\subset X$, then automatically $Y\leq X$. Note also that the requirement of
good position was needed to avoid the possibility of having distinct sets $X$
and $Y$ such that $Y\leq X$ and $X\leq Y$. In \cite{NibloSageevScottSwarup},
we defined an even stronger condition. If $E$ is a family of almost invariant
subsets of $G$, we say that the elements of $E$ are in \textit{very good
position} if whenever $U$ and $V$ are elements of $E$, either none of the four
corners of the pair is small or one is empty. This is equivalent to the
partial orders on $E$ induced by inclusion and by $\leq$ being the same. We
also showed that one can often arrange that families of almost invariant sets
are in very good position by replacing the given sets by equivalent ones.

Next we recall the theory of algebraic regular neighbourhoods. In \cite{SS2},
we defined an algebraic regular neighbourhood (Definition 6.1) and a reduced
algebraic regular neighbourhood (Definition 6.18) of a family of almost
invariant subsets of a finitely generated group $G$. See also Definition 9.1
in \cite{GSS}. We discuss the difference between these objects immediately
after Definition \ref{defnofisolated} below. Each is a bipartite graph of
groups structure $\Gamma$ for $G$ with certain properties. The basic property
of $\Gamma$ is that the $V_{0}$--vertices enclose the given almost invariant
sets. See chapters 4 and 5 of \cite{SS2}, or section 3 of \cite{GSS}, for the
definition and basic properties of enclosing. Algebraic regular neighbourhoods
need not exist, but we showed that when they exist they are unique up to
isomorphism of bipartite graphs of groups. We also showed that if one has a
finite family of almost invariant sets each over a finitely generated subgroup
of $G$, then it always has an algebraic regular neighbourhood and a reduced
algebraic regular neighbourhood. The main results of \cite{SS2} were existence
results for algebraic regular neighbourhoods of infinite families in several
special cases. In this paper, we will use the existence results for reduced
algebraic regular neighbourhoods.

Before stating these existence results, we briefly discuss how the topological
and algebraic situations are related. Groups which are virtually polycyclic
($VPC$) play an important role in this paper. The Hirsch length of such a
group will simply be called the length for brevity. A group which is $VPC$ of
length $n$ will be called $VPCn$. We will often need to refer to a group which
is $VPC$ of length at most $n$. Such a group will be called $VPC(\leq n)$. We
will also use the notation $VPC(<n)$ in a similar way. Almost invariant sets
which do not cross other almost invariant sets play a special role. We will
need the following definition from \cite{SS2}.

\begin{definition}
\label{definecanonical}Let $G$ be a one-ended finitely generated group and let
$X$ be a nontrivial almost invariant subset over a subgroup $H$ of $G$.

For $n\geq1$, we will say that $X$ is $n$\textsl{--canonical} if $X$ crosses
no nontrivial $K$--almost invariant subset of $G$, for which $K$ is $VPC(\leq
n)$.
\end{definition}

Let $G$ denote any one-ended almost finitely presented group. The natural
algebraic analogue of an essential annulus in a $3$--manifold $M$ is a
nontrivial almost invariant subset of $G$ over a $VPC1$ subgroup, and the
natural analogue of an essential torus in $M$ is a nontrivial almost invariant
subset of $G$ over a $VPC2$ subgroup. As $AT(M)$ is a kind of regular
neighbourhood of the family of all essential annuli and tori in $M$, it would
seem natural to consider the algebraic regular neighbourhood in $G$ of the
family of all nontrivial almost invariant subsets of $G$ which are over $VPC$
subgroups of length $1$ or $2$. However, we showed in Example 11.7 of
\cite{SS2} that, even when $G$ is the fundamental group of a $3$--manifold,
such a family need not possess an algebraic regular neighbourhood with the
right properties. But recall from our discussion on page \pageref{AT(M)} that
$AT(M)$ is filled by the family of all essential annuli in $M$ together with
those essential tori in $M$ which do not cross any essential annulus in $M$.
The algebraic analogue of this family is the family $\mathcal{F}_{1,2}$ of
equivalence classes of all nontrivial almost invariant subsets of $G$ which
are over $VPC1$ subgroups and of equivalence classes of all $1$--canonical
almost invariant subsets of $G$ which are over $VPC2$ subgroups. In
\cite{SS2}, as corrected in \cite{SS2errata}, we showed that $\mathcal{F}%
_{1,2}$ has an algebraic regular neighbourhood which is precisely analogous to
the decomposition given by $AT(M)$. We also described analogous constructions
for $VPC$ subgroups of $G$ of higher length. This is contained in Theorem
\ref{JSJforVPCoftwolengths} below.

Next we need to introduce some more definitions which we used in \cite{SS2}.

\begin{definition}
If $E$ is a $G$--invariant family of nontrivial almost invariant subsets of a
group $G$, we will say that an element of $E$ which crosses no element of $E$
is \textsl{isolated}\textit{ in }$E$.
\end{definition}

When forming an algebraic regular neighbourhood of a family $E$ of almost
invariant sets, isolated elements yield special vertices which we also call isolated.

\begin{definition}
\label{defnofisolated}A vertex of a graph of groups $\Gamma$ is
\textsl{isolated} if it has exactly two incident edges for each of which the
inclusion of the associated edge group into the vertex group is an isomorphism.
\end{definition}

\begin{remark}
If $\Gamma$ consists of a single vertex $v$ and a single edge, then $v$ is
\emph{not} isolated, as only one edge is incident to $v$.
\end{remark}

Note that the two edges incident to an isolated vertex have the same
associated edge splitting. Conversely if two distinct edges $e$ and
$e^{\prime}$ of a minimal graph $\Gamma$ of groups have associated edge
splittings which are conjugate, there is an edge path $\lambda$ in $\Gamma$
which starts with $e$ and ends with $e^{\prime}$ such that all the interior
vertices of $\lambda$ are isolated, and all the edges in $\lambda$ have
associated edge splittings which are conjugate. When one forms the algebraic
regular neighbourhood $\Gamma$ of a family $E$ of nontrivial almost invariant
subsets of a group $G$, it may happen that $\Gamma$ has such edge paths with
more than two edges. One could reduce $\Gamma$ by simply collapsing each
maximal such edge path in $\Gamma$ to a single edge. However $\Gamma$ is
bipartite and one wants to preserve this property, so one may instead need to
collapse such a maximal edge path to two edges. The resulting bipartite graph
of groups is the reduced algebraic regular neighbourhood of $E$. It never has
three distinct edges such that the associated splittings of $G$ are all
conjugate. We formalise this in the following definition.

\begin{definition}
A minimal bipartite graph of groups $\Gamma$ is called \textsl{reduced
bipartite} if it does not have three distinct edges such that the associated
splittings of $G$ are all conjugate.
\end{definition}

Some other special types of vertices may occur when one forms an algebraic
regular neighbourhood.

\begin{definition}
\label{defnofVPC-by-Fuchsiantype}Let $\Gamma$ be a minimal graph of groups
decomposition of a group $G$. A vertex $v$ of $\Gamma$ is of \textsl{$VPC$%
--by--Fuchsian type} if $G(v)$ is a $VPC$--by--Fuchsian group, where the
Fuchsian group is finitely generated and is not finite nor two-ended, and
there is exactly one edge of $\Gamma$ which is incident to $v$ for each
peripheral subgroup $K$ of $G(v)$, and this edge carries $K$.

If the length of the normal $VPC$ subgroup of $G(v)$ is $n$, we will say that
$v$ is of \textsl{$VPCn$--by--Fuchsian type}.
\end{definition}

\begin{remark}
It is possible that a single edge of $\Gamma$ can have both ends incident to
$v$. In this case, the two inclusions of the associated edge group into $G(v)$
must have images which are distinct peripheral subgroups of $G(v)$ up to conjugacy.
\end{remark}

Note that if $G=G(v)$, then $\Gamma$ must consist of $v$ alone, and the
Fuchsian quotient group of $G$ corresponds to a closed orbifold. Conversely if
the Fuchsian quotient group of $G(v)$ corresponds to a closed orbifold, then
$\Gamma$ must consist of $v$ alone, and $G=G(v)$. Note also that if $v$ is of
$VPCn$--by--Fuchsian type, then each peripheral subgroup of $G(v)$ is
$VPC(n+1)$.

The assumption in Definition \ref{defnofVPC-by-Fuchsiantype} that the Fuchsian
quotient of $G(v)$ not be finite nor two-ended is made to ensure the
uniqueness of the $VPCn$ normal subgroup of $G(v)$ with Fuchsian quotient.
This is immediate from the following result.

\begin{lemma}
\label{VPC-by-Fuchsiangrouphasuniquefibre}Let $G$ be a group with a normal
$VPCk$\ subgroup $L$ with Fuchsian quotient $\Phi$. Suppose that $\Phi$ is not
finite nor two-ended. If $L^{\prime}$ is a $VPCk$ normal subgroup of $G$ with
Fuchsian quotient, then $L^{\prime}$ must equal $L$.
\end{lemma}

\begin{proof}
If $L^{\prime}$ is not contained in $L$, the image of $L^{\prime}$ in $\Phi$
is a nontrivial normal $VPC$ subgroup, which we denote by $N$. As no Fuchsian
group can be $VPC2$, it follows that $N$ must be $VPC0$ or $VPC1$. If $N$ is
$VPC0$, i.e. finite, then $\Phi$ must also be finite. If $N$ is $VPC1$, then
it must be of finite index in $\Phi$. As we are assuming that $\Phi$ is not
finite nor two-ended, it follows that $L^{\prime}$ must be contained in $L$.
Similarly $L$ must be contained in $L^{\prime}$, so that $L\ $and $L^{\prime}$
are equal, as required.
\end{proof}

We note that in \cite{SS2} we used the word Fuchsian to include discrete
groups of isometries of the Euclidean plane, as well as the hyperbolic plane.
The additional groups were all virtually $\mathbb{Z}\times\mathbb{Z}$. The
reason for this abuse of language was that we wanted to include the case of
$VPC$ groups in the statements of our results. However in all of the main
results of this paper, it will be convenient to exclude the case of $VPC$ groups.

Next we prove two simple results about $VPC$ groups which will be needed on
several occasions.

\begin{lemma}
\label{splittingsofVPCgroups}Let $G$ be a $VPC(n+1)$ group which splits over a
subgroup $L$. Then $L$ is $VPCn$, and is normal in $G$ with quotient which is
isomorphic to $\mathbb{Z}$ or $\mathbb{Z}_{2}\ast\mathbb{Z}_{2}$.
\end{lemma}

\begin{proof}
The result is equivalent to asserting that, for a $VPC(n+1)$ group $G$, any
minimal $G$--tree must be a point or a line. We will prove this by induction
on the length of $G$. The induction starts when $G$ has length $1$. Then $G$
has two ends and the result is standard. Now suppose that $G$ has length
$n+1\geq2$, and that the result is known for $VPC(\leq n)$ groups. Let $T$ be
a minimal $G$--tree. There is a subgroup $G^{\prime}$ of finite index in $G$
which normalises some $VPCn$ subgroup $L^{\prime}$. By our induction
assumption, the action of $L^{\prime}$ on $T$ must fix a point or have a
minimal subtree $T^{\prime}$ which is a line. In the second case, the minimal
subtree of $T$ left invariant by $G^{\prime}$ must also be $T^{\prime}$. In
the first case, we let $T^{\prime}$ denote the fixed subtree of $L^{\prime}$,
i.e. $T^{\prime}$ consists of all vertices and edges fixed by $L^{\prime}$.
The action of $G^{\prime}$ must preserve $T^{\prime}$, so that the quotient
group $L^{\prime}\backslash G^{\prime}$ acts on $T^{\prime}$. As this quotient
group has two ends, it follows that the minimal subtree of $T^{\prime}$ left
invariant by $G^{\prime}$ is a point or a line. Thus in either case, the
minimal subtree of $T$ left invariant by $G^{\prime}$ is a point or a line. As
$G^{\prime}$ has finite index in $G$, this minimal subtree must equal $T$, so
that $T$ itself is a point or a line as required.
\end{proof}

\begin{lemma}
\label{normalsubgroupsofVPCgroupsofcolength1} Let $G$ be a $VPC(n+1)$ group
with a normal $VPCn$ subgroup $L$ such that $L\backslash G$ is isomorphic to
$\mathbb{Z}$ or to $\mathbb{Z}_{2}\ast\mathbb{Z}_{2}$. Let $K$ be a normal
$VPCn$ subgroup of $G$ which is commensurable with $L$. Then the following hold:

\begin{enumerate}
\item $K$ is contained in $L$.

\item If $K\backslash G$ is isomorphic to $\mathbb{Z}$ or to $\mathbb{Z}%
_{2}\ast\mathbb{Z}_{2}$, then $K=L$.
\end{enumerate}
\end{lemma}

\begin{proof}
1) If $K$ is not contained in $L$, then the image of $K$ in the quotient
$L\backslash G$ is a nontrivial finite normal subgroup. As neither
$\mathbb{Z}$ nor $\mathbb{Z}_{2}\ast\mathbb{Z}_{2}$ possesses such a subgroup,
it follows that $K$ must be contained in $L$ as required.

2) If $K\backslash G$ is isomorphic to $\mathbb{Z}$ or to $\mathbb{Z}_{2}%
\ast\mathbb{Z}_{2}$, we can apply the first part with the roles of $K$ and $L$
reversed. We deduce that $L$ is contained in $K$, so that $K=L$ as required.
\end{proof}

We also prove the following useful technical results about $VPC$--by--Fuchsian groups.

\begin{lemma}
\label{torusinVPC-by-Fuchsiangroup} Let $G$ be a group with a normal
$VPCn$\ subgroup $L$ with Fuchsian quotient $\Phi$, and let $K$ be a
$VPC(n+1)$ subgroup of $G$.

\begin{enumerate}
\item Then $L\cap K$ is a normal $VPCn$ subgroup of $K$ with quotient
isomorphic to $\mathbb{Z}$ or to $\mathbb{Z}_{2}\ast\mathbb{Z}_{2}$.

\item If $H$ is a normal $VPCn$ subgroup of $K$ with quotient isomorphic to
$\mathbb{Z}$ or to $\mathbb{Z}_{2}\ast\mathbb{Z}_{2}$, and if $H$ is
commensurable with $L$, then $H$ equals $L\cap K$.
\end{enumerate}
\end{lemma}

\begin{proof}
1) As $L$ is $VPCn$ and normal in $G$, the intersection $L\cap K$ must be
$VPC(\leq n)$, and normal in $K$. Hence the quotient of $K$ by $L\cap K$ is a
$VPC$ subgroup of $\Phi$ of length at least $1$. As a Fuchsian group can have
no $VPC2$ subgroups, it follows that $k$ must equal $n$ and the quotient of
$K$ by $L\cap K$ must be $VPC1$. As the only $VPC1$ subgroups of a Fuchsian
group are isomorphic to $\mathbb{Z}$ or to $\mathbb{Z}_{2}\ast\mathbb{Z}_{2}$,
the result follows.

2) This follows from Lemma \ref{normalsubgroupsofVPCgroupsofcolength1}.
\end{proof}

Next we have a uniqueness result for a $VPC$--by--Fuchsian structure on a group.

\begin{lemma}
\label{uniqueVPC-by-Fuchsianstructure}Let $G$ be a group with a normal
$VPCn$\ subgroup $L$ with Fuchsian quotient $\Phi$ which is not virtually
cyclic, and suppose $G$ also has a normal $VPCm$\ subgroup $L^{\prime}$ with
Fuchsian quotient $\Phi^{\prime}$ which is not virtually cyclic. Then
$L=L^{\prime}$.
\end{lemma}

\begin{proof}
Suppose $L^{\prime}$ is not contained in $L$. Then the image of $L^{\prime}$
in $\Phi$ is a nontrivial normal $VPC$ subgroup. But a Fuchsian group which is
not virtually cyclic cannot contain such a subgroup. It follows that
$L^{\prime}\subset L$. Similarly, it follows that $L\subset L^{\prime}$ so
that $L=L^{\prime}$, as required.
\end{proof}

Now we can state the results from \cite{SS2} which play a basic role in this
paper. Recall that if $H$ is a subgroup of a group $G$, the
\textit{commensuriser}, $Comm_{G}(H)$, of $H$ in $G\ $is the subgroup of $G$
consisting of all elements $g$ such that the conjugate of $H$ by $g$ is
commensurable with $H$. Trivially, $Comm_{G}(H)$ contains $H$. We will say
that $Comm_{G}(H)$ is \textit{large} if it contains $H$ with infinite index,
and is \textit{small} otherwise.

The following existence result is essentially the statement of Theorem 12.3 of
\cite{SS2}. We have made one slight modification in the last sentence of part
3), where we refer to the number of coends of a subgroup in a group. See page
33 of \cite{SS2} for a brief discussion of this concept. It was introduced
independently by Bowditch \cite{B2} under the name of coends, by Geoghegan
\cite{Geoghegan} under the name of filtered coends, and by Kropholler and
Roller \cite{KR} under the name of relative ends.

\begin{theorem}
\label{JSJforVPCofgivenlength}Let $n\geq1$, and let $G$ be a one-ended, almost
finitely presented group which is not $VPC$ and does not admit any nontrivial
almost invariant subsets over $VPC(<n)$ subgroups, and let $\mathcal{F}_{n}$
denote the collection of equivalence classes of all nontrivial almost
invariant subsets of $G$ which are over $VPCn$ subgroups.

Then $\mathcal{F}_{n}$ has an unreduced and a reduced algebraic regular
neighbourhood in $G$. Let $\Gamma_{n}=\Gamma(\mathcal{F}_{n}:G)$ denote the
reduced algebraic regular neighbourhood of $\mathcal{F}_{n}$ in $G$.

Then $\Gamma_{n}$ is a minimal, reduced bipartite, graph of groups
decomposition of $G$. Each $V_{0}$--vertex $v$ of $\Gamma_{n}$ satisfies one
of the following conditions:

\begin{enumerate}
\item $v$ is isolated, and $G(v)$ is $VPCn$.

\item $v$ is of $VPC(n-1)$--by--Fuchsian type, and elements of $\mathcal{F}%
_{n}$ enclosed by $v$ cross strongly if at all.

\item $G(v)$ is the full commensuriser $Comm_{G}(H)$ for some $VPCn$ subgroup
$H$, such that $e(G,H)\geq2$, and elements of $\mathcal{F}_{n}$ enclosed by
$v$ cross weakly if at all.

Further, if $H$ is a $VPCn$ subgroup of $G$ such that $e(G,H)\geq2$, then
$\Gamma_{n}$ will have a non-isolated $V_{0}$--vertex $v$ such that
$G(v)=Comm_{G}(H)$ if and only if $H$ has at least $4$ coends in $G$.
\end{enumerate}

$\Gamma_{n}$ consists of a single vertex if and only if $\mathcal{F}_{n}$ is
empty, or $G$ itself satisfies condition 2) or 3). In the first case,
$\Gamma_{n}$ consists of a single $V_{1}$--vertex. In the second case,
$\Gamma_{n}$ consists of a single $V_{0}$--vertex.
\end{theorem}

\begin{remark}
\label{remarksonJSJ}The assumption that $G$ is not $VPC$ is made to simplify
the statement of this result. If $G\ $is $VPC(n+1)$, then $\Gamma_{n}$
consists of a single $V_{0}$--vertex $v$ with associated group $G$. In this
case, $v$ may not satisfy any of the conditions in the above theorem.

We will say that a $V_{0}$--vertex in case 3) is \textsl{of commensuriser
type} if $v$ is not isolated, and is \textsl{of large commensuriser type} if,
in addition, $H$ has large commensuriser.

Note that we showed in Example 11.1 of \cite{SS2} that, even if $G$ is
finitely presented, the group associated to a $V_{0}$--vertex of commensuriser
type need not be finitely generated.

The statement at the end of part 3) with the assumption that $H$ has at least
$4$ coends in $G$ comes from the proofs of Propositions 7.16, 7.17, 8.1 and
8.6 of \cite{SS2}. Note that if $H$ has large commensuriser in $G$, the proof
of Proposition 8.1 of \cite{SS2} shows that $H$ has infinitely many coends in
$G$. This uses the fact that we have excluded the case when $G\ $is $VPC$.
\end{remark}

When $n=2$, this result is the algebraic analogue of the torus decomposition
$T(M)$ of a closed orientable Haken $3$--manifold $M$. See the start of this
section for a discussion of $T(M)$. As $M$ is compact, its fundamental group
$G$ is finitely presented. As $M$ is irreducible, it follows that $G$ is also
one-ended. As $M$ is closed, our discussion in \cite{SS4} implies that $G$
does not admit any nontrivial almost invariant subsets over $VPC1$ subgroups.
Now the above result asserts that $\Gamma_{2}(G)$ exists, and we showed in
\cite{SS2} that it is the graph of groups structure for $G$ determined by the
frontier of $T(M)$ in $M$. The $V_{0}$--vertices of $\Gamma_{2}(G)$ correspond
to the components of $T(M)$. In this case, $\Gamma_{2}(G)$ has no $V_{0}%
$--vertices of commensuriser type. An isolated $V_{0}$--vertex of $\Gamma
_{2}(G)$ corresponds to a component of $T(M)$ homeomorphic to $T\times I$. A
$V_{0}$--vertex of $\Gamma_{2}(G)$ which is of $VPC1$--by--Fuchsian type
corresponds to a component of $T(M)$ which is homeomorphic to a Seifert fibre space.

Next we come to the following more general result which is essentially the
statement of Theorem 13.12 of \cite{SS2}. Again we assume that $G$ is not
$VPC$ in order to simplify the statement.

\begin{theorem}
\label{JSJforVPCoftwolengths}Let $n\geq1$, and let $G$ be a one-ended, almost
finitely presented group which is not $VPC$ and does not admit any nontrivial
almost invariant subsets over $VPC(<n)$ subgroups, and let $\mathcal{F}%
_{n,n+1}$ denote the collection of equivalence classes of all nontrivial
almost invariant subsets of $G$ which are over a $VPCn$ subgroup, together
with the equivalence classes of all $n$--canonical almost invariant subsets of
$G$ which are over a $VPC(n+1)$ subgroup.

Then $\mathcal{F}_{n,n+1}$ has an unreduced and a reduced algebraic regular
neighbourhood in $G$. Let $\Gamma_{n,n+1}=\Gamma(\mathcal{F}_{n,n+1}:G)$
denote the reduced algebraic regular neighbourhood of $\mathcal{F}_{n,n+1}$ in
$G$.

Then $\Gamma_{n,n+1}$ is a minimal, reduced bipartite, graph of groups
decomposition of $G$. Each $V_{0}$--vertex $v$ of $\Gamma_{n,n+1}$ satisfies
one of the following conditions:

\begin{enumerate}
\item $v$ is isolated, and $G(v)$ is $VPC$ of length $n$ or $n+1$.

\item $v$ is of $VPCk$--by--Fuchsian type, where $k$ equals $n-1$ or $n$, and
elements of $\mathcal{F}_{n,n+1}$ enclosed by $v$ cross strongly if at all.

\item $G(v)$ is the full commensuriser $Comm_{G}(H)$ for some $VPC$ subgroup
$H$ of length $n$ or $n+1$, such that $e(G,H)\geq2$, and elements of
$\mathcal{F}_{n,n+1}$ which are enclosed by $v$ and are over groups
commensurable with $H$ cross weakly if at all.

Further, if $H$ is a $VPC$ subgroup of $G$ of length $n$ or $n+1$, such that
$e(G,H)\geq2$, then $\Gamma_{n,n+1}$ will have a non-isolated $V_{0}$--vertex
$v$ such that $G(v)=Comm_{G}(H)$ if and only if $H$ has at least $4$ coends in
$G$.
\end{enumerate}

$\Gamma_{n,n+1}$ consists of a single vertex if and only if $\mathcal{F}%
_{n,n+1}$ is empty, or $G$ itself satisfies condition 2) or 3). In the first
case, $\Gamma_{n,n+1}$ consists of a single $V_{1}$--vertex. In the second
case, $\Gamma_{n,n+1}$ consists of a single $V_{0}$--vertex.
\end{theorem}

When $n=1$, this result is the algebraic analogue of the Annulus-Torus
decomposition $AT(M)$ of an orientable Haken $3$--manifold $M$ with
incompressible boundary. See the start of this section for a discussion of
$AT(M)$. As $M$ is compact, its fundamental group $G$ is finitely presented.
As $M$ is irreducible and has incompressible boundary, it follows that $G$ is
also one-ended. Now the above result asserts that $\Gamma_{1,2}(G)$ exists,
and we showed in \cite{SS2} that it is the graph of groups structure for $G$
determined by the frontier of $AT(M)$ in $M$. The $V_{0}$--vertices of
$\Gamma_{1,2}(G)$ correspond to the components of $AT(M)$. An isolated $V_{0}%
$--vertex of $\Gamma_{1,2}(G)$ corresponds to a component of $AT(M)$
homeomorphic to $A\times I$ or $T\times I$. A $V_{0}$--vertex of
$VPC1$--by--Fuchsian type corresponds to a component of $AT(M)$ which is a
Seifert fibre space not meeting $\partial M$. A $V_{0}$--vertex of
$VPC0$--by--Fuchsian type corresponds to a component of $AT(M)$ which is an
$I$--bundle and meets $\partial M$ in the associated $S^{0}$--bundle. A
$V_{0}$--vertex $v$ of commensuriser type, where $G(v)$ is the full
commensuriser $Comm_{G}(H)$ for some $VPC1$ subgroup $H$, corresponds to a
component of $AT(M)$ which is a Seifert fibre space and meets $\partial M$.
Finally, in this case, $\Gamma_{1,2}(G)$ has no $V_{0}$--vertex $v$ of
commensuriser type, where $G(v)$ is the full commensuriser of a $VPC2$
subgroup of $G$.

In the special case when $M$ is closed or admits no essential annuli, then $G$
has no nontrivial almost invariant subsets over any $VPC1$ subgroup, so that
$\Gamma_{2}(G)$ is defined and equals $\Gamma_{1,2}(G)$.

Recall that $JSJ(M)$ can be obtained from $AT(M)$ by adding certain
exceptional submanifolds of $M$. We will describe analogous algebraic
constructions for any $n$ and any group $G$ for which $\Gamma_{n+1}(G)$ or
$\Gamma_{n,n+1}(G)$ exist. The results of these constructions are graphs of
groups structures $\Gamma_{n+1}^{c}(G)$ and $\Gamma_{n,n+1}^{c}(G)$ for $G$,
which we call the completions of $\Gamma_{n+1}(G)$ and $\Gamma_{n,n+1}(G)$ respectively.

\begin{definition}
\label{defnofcompletion}Let $G$ be a group for which the decompositions
$\Gamma_{n+1}(G)$ or $\Gamma_{n,n+1}(G)$ exist. The \textsl{completions} of
these decompositions, denoted $\Gamma_{n+1}^{c}(G)$ and $\Gamma_{n,n+1}%
^{c}(G)$ respectively are graphs of groups structures for $G$ obtained as follows:

If $\Gamma_{n+1}(G)$ or $\Gamma_{n,n+1}(G)$ has a $V_{1}$--vertex $w$ such
that $G(w)$ is $VPC(n+1)$, and if $w$ has a single incident edge $e$ with
$G(e)$ of index $2$ in $G(w)$, then we subdivide $e$ into two edges. The new
vertex is a $V_{1}$--vertex and $w$ becomes a $V_{0}$--vertex. If the original
$V_{0}$--vertex of $e$ is isolated, then in addition we collapse $e$ to a
point, which becomes a new $V_{0}$--vertex.

If $\Gamma_{n,n+1}(G)$ has a $V_{1}$--vertex $w$ such that $G(w)$ is $VPCn$,
and if $w$ has a single incident edge $e$ with $G(e)$ of index $2$ or $3$ in
$G(w)$, or if $w$ has exactly three incident edges each carrying $G(w)$, then
we subdivide each of the incident edges into two edges. The new vertices are
$V_{1}$--vertices and $w$ becomes a $V_{0}$--vertex. If the original $V_{0}%
$--vertex of any of the edges incident to $w$ is isolated, then in addition we
collapse that edge.

Making all these changes for every such $V_{1}$--vertex of $\Gamma_{n+1}(G)$
and $\Gamma_{n,n+1}(G)$ yields $\Gamma_{n+1}^{c}(G)$ and $\Gamma_{n,n+1}%
^{c}(G)$ respectively.
\end{definition}

\begin{remark}
If $n=1$, and $G$ is the fundamental group of an orientable Haken
$3$--manifold $M$ with incompressible boundary, then the discussion in
\cite{SS2} shows that $\Gamma_{1,2}^{c}(G)$ is the graph of groups determined
by the frontier of $JSJ(M)$ in $M$. The $V_{0}$--vertices of $\Gamma_{1,2}%
^{c}(G)$ correspond to the components of $JSJ(M)$. Those $V_{0}$--vertices of
$\Gamma_{1,2}^{c}(G)$ which are obtained from $V_{1}$--vertices of
$\Gamma_{1,2}(G)$ correspond to the exceptional components of $JSJ(M)$. If $M$
is closed or admits no essential annuli, then $\Gamma_{2}^{c}(G)$ is defined
and equals $\Gamma_{1,2}^{c}(G)$. Also those $V_{0}$--vertices of $\Gamma
_{2}^{c}(G)$ which are obtained from $V_{1}$--vertices of $\Gamma_{2}(G)$
correspond to the exceptional components of $JSJ(M)$.
\end{remark}

\section{Poincar\'{e} duality pairs and essential annuli and
tori\label{essentialannuli}}

We refer to Brown \cite{Brown}, Bieri and Eckmann \cite{B-E}, Kapovich and
Kleiner \cite{K-K}, and Wall \cite{Wall3} for various definitions of
Poincar\'{e} duality groups and pairs. The definition of Bieri and Eckmann is
that $G$ is a Poincar\'{e} duality group of dimension $n+2$ if $H^{i}%
(G;\mathbb{Z}G)$ is $0$, when $i\neq n+2$, and is isomorphic to $\mathbb{Z}$
when $i=n+2$. Further $G$ is orientable if the action of $G$ on $H^{n+2}%
(G;\mathbb{Z}G)$ is trivial. In the following, we will be mostly concerned
with orientable Poincar\'{e} duality groups and pairs. A Poincar\'{e} duality
pair is a pair $(G,\partial G)$, where $\partial G=\{S_{1},...,S_{m}\}$ is a
system of subgroups of $G$, such that the double of $G$ along $\partial G$ is
a Poincar\'{e} duality group. Theorem 8.1 of \cite{B-E} shows that each
$S_{i}$ must be a $PD(n+1)$ group. Note that the order of the $S_{i}$'s is
irrelevant, and that repetitions are allowed. However, if any repetition
occurs, or even if two distinct $S_{i}$'s are conjugate, there is a $PD(n+1)$
group $H$, such that the pair $(G,\partial G)$ equals $(H,\{H,H\})$, which is
a trivial $PD(n+2)$ pair analogous to the product of a closed $(n+1)$%
--manifold with the unit interval. We will usually assume that $n\geq1$, so
that our Poincar\'{e} duality groups and pairs are at least $3$--dimensional.
Bieri and Eckmann \cite{B-E} show that their definition implies that $G$ is
almost finitely presented, which suffices for the accessibility results that
we use. If $G$ is finitely presented then the corresponding $K(G,1)$ space is
dominated by a finite complex (see Theorem 7.1 in Chapter 8 of \cite{Brown})
and then $G$ is a Poincar\'{e} duality group in the sense of Wall \cite{Wall3}.

A Poincar\'{e} duality pair is a special case of what the authors of
\cite{GSS} call a group system, but it seems natural to use the language of
pairs in the setting of this paper. In \cite{GSS}, many of our results require
that the group system be of finite type, but this condition is automatic for
Poincar\'{e} duality pairs. An important idea which we first introduced in
\cite{SS-pdnarxiv}, and is worked out in more detail in \cite{GSS}, is that of
an almost invariant subset of a group $G$ being adapted to a family of
subgroups. As this idea will play an important role in this paper, we
reproduce the definition from \cite{GSS}, and some of the following remarks.
Lemma \ref{propertiesofadapted} summarises the properties of adapted almost
invariant sets which we will need in this paper.

\begin{definition}
\label{defnofadapted} (Definition 5.1 of \cite{GSS}) Let $G$ be a group and
let\ $H$ and $S$ be subgroups. Let $\mathcal{S}=\{S_{i}\}_{i\in I}$ be a
family of subgroups of $G$, with repetitions allowed.

A $H$--almost invariant subset $X$ of $G$ is \textsl{strictly adapted to the
subgroup} $S$ if, for all $g\in G$, the coset $gS$ is contained in $X$ or in
$X^{\ast}$.

A $H$--almost invariant subset $X$ of $G$ is \textsl{adapted to} $S$, or is
$S$--\textsl{adapted}, if it is equivalent to a $H^{\prime}$--almost invariant
subset $X^{\prime}$ of $G$ such that $X^{\prime}$ is strictly adapted to $S$.

A $H$--almost invariant subset $X$ of $G$ is \textsl{adapted to the family}
$\mathcal{S}$, or is $\mathcal{S}$--\textsl{adapted}, if it is adapted to each
$S_{i}$.
\end{definition}

\begin{remark}
If $X$ is $\mathcal{S}$--adapted, then so is any almost invariant subset of
$G$ which is equivalent to $X$.

If $X$ is adapted to the family $\mathcal{S}$, and we replace each $S_{i}$ by
some conjugate, then $X$ is also adapted to the new family. For if $X$ is
strictly adapted to a subgroup $S$, and $k\in G$, then $Xk$ is $H$--almost
invariant and equivalent to $X$, and is strictly adapted to $k^{-1}Sk$.

Note that if $X$ is adapted to the family $\mathcal{S}$, then, for each $i$,
there is a $K_{i}$--almost invariant subset $X_{i}$ of $G$ which is equivalent
to $X$ and is strictly adapted to $S_{i}$, but the $X_{i}$'s may all be
different. In general, it is difficult for an almost invariant set to be
strictly adapted to more than one subgroup of $G$.
\end{remark}

The following definition, due to M\"{u}ller \cite{Muller}, is natural when one
considers splittings of a group.

\begin{definition}
\label{defnofadaptedsplitting}Let $K$ be a group with a splitting $\sigma$
over a subgroup $H$, and let $\mathcal{S}=\{S_{i}\}$ be a family of subgroups
of $K$. The splitting $\sigma$ of $K$ is \textsl{adapted to} $\mathcal{S}$, or
is $\mathcal{S}$--\textsl{adapted}, if each $S_{i}$ is conjugate into a vertex
group of $\sigma$.
\end{definition}

There is also a natural generalisation of Definition
\ref{defnofadaptedsplitting} to graphs of groups.

\begin{definition}
\label{defnofadaptedgraphofgroups}Let $K$ be a group with a graph of groups
structure $\Gamma$, and let $\mathcal{S}=\{S_{i}\}$ be a family of subgroups
of $K$. Then $\Gamma$ is \textsl{adapted to} $\mathcal{S}$, or is
$\mathcal{S}$--\textsl{adapted}, if each $S_{i}$ is conjugate into a vertex
group of $\Gamma$.
\end{definition}

In \cite{GSS}, it is shown that these terminologies are all compatible.

In topological terms, this concept seems very natural when one considers
manifolds with boundary, and this is how Muller's definition arose. Consider a
manifold $M$, and a codimension--$1$ embedded submanifold $F$ in the interior
of $M$, such that $F$ is two-sided, closed and $\pi_{1}$--injective. Then
clearly the splitting of $\pi_{1}(M)$ over $\pi_{1}(F)$ determined by $F$ is
adapted to the family $\mathcal{S}$ of subgroups of $M$ carried by the
components of $\partial M$. On the other hand, if $F$ has boundary and is
properly embedded, then the corresponding splitting of $\pi_{1}(M)$ over
$\pi_{1}(F)$ is likely not to be adapted to $\mathcal{S}$.

For the purposes of this paper we need to recall from \cite{GSS} some
properties of adapted almost invariant subsets of a group, which we summarize
in the following lemma.

\begin{lemma}
\label{propertiesofadapted}Let $\overline{G}$ be a group with subgroups $H$
and $G$, such that $H\subset G$, and suppose that $\overline{G}$ has a minimal
graph of groups decomposition $\Gamma$ with a vertex $V$ whose associated
group is $G$. Let $\mathcal{S}$ denote the family of subgroups of $G$
associated to the edges of $\Gamma$ which are incident to $V$.

\begin{enumerate}
\item If $\overline{X}$ is a $H$--almost invariant subset of $\overline{G}$
which is enclosed by $V$, and if $X$ denotes the $H$--almost invariant subset
$\overline{X}\cap G$ of $G$, then $X$ is adapted to the family $\mathcal{S}$.

\item Suppose that the group system $(G,\mathcal{S})$ is of finite type. If
$X$ is a nontrivial $\mathcal{S}$--adapted $H$--almost invariant subset of
$G$, then $X$ has an extension $\overline{X}$ to $\overline{G}$, i.e. there is
a $H$--almost invariant subset $\overline{X}$ of $\overline{G}$ which is
enclosed by $V$, such that $\overline{X}\cap G=X$. Further if $X$ is
associated to a splitting of $G$ over $H$, then $\overline{X}$ is associated
to a splitting of $\overline{G}$ over $H$.
\end{enumerate}
\end{lemma}

\begin{remark}
As Poincar\'{e} duality pairs are automatically group systems of finite type,
this condition in part 2) of the lemma can be ignored in this paper.
\end{remark}

We will collect here some useful basic results about Poincar\'{e} duality
groups and pairs. The following result due to Kropholler and Roller is Lemma
2.2 of \cite{KR2}. It will be needed at several points in this paper. Similar
results are well known in the topology of $3$--manifolds. In this setting,
cases 1) and 2) of the lemma below occur when one has an $I$--bundle over a
closed surface.

\begin{lemma}
\label{boundarygroupismaximal} (Kropholler and Roller) Let $(G,\partial G)$ be
a $PD(n+2)$ pair with $\partial G$ non-empty. Then one of the following holds:

\begin{enumerate}
\item $G$ is a $PD(n+1)$ group and $\partial G$ consists of a single group $S$
which has index $2$ in $G$.

\item $G$ is a $PD(n+1)$ group and the pair $(G,\partial G)$ is the trivial
pair $(G,\{G,G\})$.

\item For each group $S$ in $\partial G$, the index of $S$ in $G$ is infinite,
and $Comm_{G}(S)=S$. Further if $S$ and $S^{\prime}$ are distinct groups in
$\partial G$, they are not conjugate commensurable.
\end{enumerate}
\end{lemma}

It will also be convenient to state separately the following easy consequences.

\begin{corollary}
\label{PDgroupsplitsoverPDHimpliesHismaximal} Let $(G,\partial G)$ be an
orientable $PD(n+2)$ pair. Then

\begin{enumerate}
\item If $\partial G$ consists of a single group $S$ which has index $2$ in
$G$, then $G$ itself must be a non-orientable $PD(n+1)$ group.

\item If $S$ is a group in $\partial G$, and $K$ is an orientable $PD(n+1)$
subgroup of $G$ commensurable with $S$, then $K$ is contained in $S$.

\item If $(G,\partial G)$ splits, adapted to $\partial G$, over a $PD(n+1)$
subgroup $H$, and if $H$ is commensurable with an orientable $PD(n+1)$
subgroup $K$ of $G$, then $K$ is contained in $H$.
\end{enumerate}
\end{corollary}

\begin{proof}
1) As $(G,\partial G)$ is an orientable $PD(n+2)$ pair, $S$ must be an
orientable $PD(n+1)$ group. Thus $G$ is also a $PD(n+1)$ group. If $G\ $were
orientable, the fact that $S$ has index $2$ in $G$ would imply that the
induced map $\mathbb{Z}\cong H_{n+1}(S)\rightarrow H_{n+1}(G)\cong\mathbb{Z}$
would be multiplication by $2$. But this map must be zero as $S=\partial G$.
This contradiction shows that $G$ must be non-orientable.

2) We apply Lemma \ref{boundarygroupismaximal}. If case 3) of that lemma
holds, the fact that $K$ must commensurise $S$ implies that $K$ is contained
in $S$. If case 2) of that lemma holds, the result is trivial. If case 1) of
that lemma holds, then part 1) of this lemma tells us that $G$ is a
non-orientable $PD(n+1)$ group. Hence $S$ is the orientation subgroup of $G$,
so that any orientable $PD(n+1)$ subgroup of $G$ must be contained in $S$.

3) First suppose that $\partial G$ is empty so that $G$ is an orientable
$PD(n+2)$ group. As $G$ splits over $H$ we have $G$ equal to $A\ast_{H}B$ or
to $A\ast_{H}$ for some subgroups $A$ or $B$. In the first case, Theorem 8.1
of \cite{B-E} tells us that each of the pairs $(A,H)$ and $(B,H)$ is an
orientable $PD(n+2)$ pair. In the second case, there are two inclusions of $H$
into $A$ whose images we denote by $H_{1}$ and $H_{2}$, and Theorem 8.1 of
\cite{B-E} tells us that the pair $(A,\{H_{1},H_{2}\})$ is an orientable
$PD(n+2)$ pair. As $K$ is commensurable with $H$, it must be conjugate into
$A$ or $B$. Now part 2) shows that $K$ is contained in $H$, as required.

If $\partial G$ is not empty, we recall that the given splitting of $G$ over
$H$ is adapted to $\partial G$. Thus if we consider $DG$, the double of $G$
over $\partial G$, and apply part 2) of Lemma \ref{propertiesofadapted}, the
given splitting of $G\ $over $H$ induces a splitting of $DG$ over $H$. Now we
can apply the above argument to the orientable $PD(n+2)$ group $DG$ to deduce
that $K$ is contained in $H$, as required.
\end{proof}

The following observation will also be useful at several points in this paper.
In the setting of $3$--manifolds, the corresponding result is that if a Haken
manifold $M$ with non-empty incompressible boundary has $VPC$ fundamental
group, than $M$ must be an $I$--bundle over the torus or Klein bottle.

\begin{corollary}
\label{PD(n+2)andVPCgroup} Let $(G,\partial G)$ be a $PD(n+2)$ pair, and
suppose that $G$ is $VPC$. Then one of the following holds:

\begin{enumerate}
\item $\partial G$ is empty and $G$ is $VPC(n+2)$.

\item $\partial G$ is non-empty, $G$ is $VPC(n+1)$, and either $(G,${$\partial
G$}$)$ is the trivial pair $(G,\{G,G\})$, or $\partial G$ is a single group
$S$, and $G$ contains $S$ with index $2$.
\end{enumerate}
\end{corollary}

\begin{proof}
If $\partial G$ is empty, then $G$ has cohomological dimension $n+2$, and so
$G$ must be $VPC(n+2)$. Thus we have case 1) of the corollary.

If $\partial G$ is non-empty, each group in $\partial G$ is $PD(n+1)$, and is
also $VPC.$ Thus each group in $\partial G$ is $VPC(n+1)$. Now $G$ must have
cohomological dimension $n+1$, and so $G$ is also $VPC(n+1)$. As this implies
that any group in $\partial G$ has finite index in $G$, Lemma
\ref{boundarygroupismaximal}\ implies that we must have case 2) of the corollary.
\end{proof}

Now we will begin our discussion of almost invariant sets in Poincar\'{e}
duality groups and pairs. But first here are two important facts about such
sets which follow immediately from Lemma 4.3 of Kropholler in \cite{K}.

\begin{lemma}
\label{PD(n+2)groupshavenoa.i.setsoverVPCk<n}(Kropholler and Roller) Let
$n\geq1$, and let $(G,\partial G)$ be a $PD(n+2)$ pair. Then

\begin{enumerate}
\item $G$ has no nontrivial almost invariant subset over a $VPC(<n)$ subgroup.

\item If $\partial G$ is empty, so that $G$ is a $PD(n+2)$ group, then $G$ has
no nontrivial almost invariant subset over a $VPC(\leq n)$ subgroup.
\end{enumerate}
\end{lemma}

\begin{remark}
In the $3$--manifold setting, so that $n=1$, part 1) corresponds to the fact
that a Haken manifold $M$ with incompressible boundary must have a one-ended
fundamental group, and part 2) corresponds to the additional fact that when
$M$ is closed it cannot admit essential annuli.
\end{remark}

The analogy between $\pi_{1}$--injective maps of surfaces into $3$--manifolds
and almost invariant subsets of groups was one of the guiding principles in
\cite{SS1} and \cite{SS2}. In particular, essential maps of annuli and tori
into orientable $3$--manifolds have corresponding nontrivial almost invariant sets.

We start by describing the analogous correspondence when one considers an
orientable $PD(n+1)$ subgroup $H$ of an orientable $PD(n+2)$ group $G$. As $G$
and $H$ are orientable, it follows that $e(G,H)=2$. Thus $G\ $has a nontrivial
$H$--almost invariant subset $X_{H}$ which is unique up to equivalence and
complementation. This is the almost invariant subset of $G$ which we associate
to $H$. We call it the \textit{dual} of $H$. The restriction that $H$ and $G$
be orientable is not crucial to ensure that $e(G,H)=2$, but it does simplify
the statements somewhat. What is crucial is that when we consider the
inclusion of $H$ into $G$, it should commute with the orientation
homomorphisms. For otherwise $e(G,H)=1$, so that $G$ has no nontrivial
$H$--almost invariant subsets. In topological terms, we want our
codimension--$1$ manifolds to have trivial normal bundle or equivalently to be
two-sided. Our higher dimensional algebraic analogue of a torus is an
orientable $PD(n+1)$ group $H$ which is also $VPC(n+1)$, and then a torus in
an orientable $PD(n+2)$ group $G$ is an injective homomorphism $\Phi
:H\rightarrow G$. Note that a $PD(n+2)$ group $G$ is torsion free, so that a
$VPC(n+1)$ subgroup of $G$ is always $PD(n+1)$.

Next we consider $PD(n+2)$ pairs. If $(G,\partial G)\ $is an orientable
$PD(n+2)$ pair with nonempty boundary, and if $H$ is an orientable $PD(n+1)$
subgroup of $G$, then $e(G,H)$ may not equal $2$. Thus it is no longer clear
how to naturally associate a $H$--almost invariant subset of $G$ to $H$.
However in the topological context, a map of a codimension--$1$ closed
orientable manifold into an orientable manifold determines a corresponding
almost invariant set in a natural way. For simplicity, consider the case when
$n=1$, and $G$ is the fundamental group of an orientable $3$--manifold $M$
with incompressible boundary, and $H$ is isomorphic to the fundamental group
of a closed orientable surface $F$. Pick a map of $F$ into the interior of $M$
so that $\pi_{1}(F)$ maps to $H$ by the given isomorphism, and consider its
lift to the cover $M_{F}$ of $M$ with fundamental group $H$. Then $M_{F}$ need
not have two ends, so that the group $H$ does not determine a unique
$H$--almost invariant subset of $G$. But, as $F$ is two-sided in $M$, the lift
of $F$ into $M_{F}$ separates $M_{F}$ into two pieces, and so $F$ \emph{does}
determine a unique $H$--almost invariant subset $Y$ of $G$. Note that $Y$ will
be trivial if and only if $F$ is homotopic into a component of $\partial M$.
Note also that as $F$ does not meet $\partial M$, the associated almost
invariant set $Y$ is adapted to $\partial G$, the family of subgroups of
$G=\pi_{1}(M)$ carried by the components of $\partial M$.

In the algebraic setting, where $(G,\partial G)\ $is an orientable $PD(n+2)$
pair with nonempty boundary and $H$ is an orientable $PD(n+1)$ subgroup of
$G$, we will also associate to $H$ a $H$--almost invariant subset of $G$ which
is adapted to $\partial G$. We do this by considering the double $DG$ of $G$
along $\partial G$. Since $DG$ is orientable, $H$ determines a nontrivial
$H$--almost invariant subset $X_{H}$ of $DG$. The $H$--almost invariant subset
of $G$ which we associate to $H$ is the intersection $Y=X_{H}\cap G$. It is
clear that $Y$ is a trivial $H$--almost invariant subset of $G$ if $H$ is
conjugate to a subgroup of one of the $S_{i}$'s. The following result shows
that, as in the topological situation, this is the only way in which $Y$ can
be trivial.

\begin{lemma}
\label{naturalrestriction} Let $(G,\partial G)\ $be an orientable $PD(n+2)$
pair with nonempty boundary, and let $H$ be an orientable $PD(n+1)$ subgroup
of $G$. Let $X_{H}$ denote the $H$--almost invariant subset of $DG$ determined
by $H$, and let $Y$ denote the $H$--almost invariant subset $X_{H}\cap G$ of
$G$.

Then $Y$ is adapted to $\partial G$. Further, if $Y$ is trivial, then $H$ is
conjugate to a subgroup of one of the $S_{i}$'s.
\end{lemma}

\begin{proof}
To prove this result, we need to use some techniques from \cite{SS2}. Let
$\Delta$ denote the graph of groups structure for $DG$ corresponding to its
construction by doubling. Thus $\Delta$ has two vertices $w$ and $\overline
{w}$, and edges corresponding to the groups in $\partial G$, each joining the
two vertices. We identify $G$ with the vertex group $G(w)$. As $DG$ is an
orientable $PD(n+2)$ group, each $S_{i}$ determines, up to complementation and
equivalence, a unique $S_{i}$--almost invariant subset $X_{i}$ of $DG$. Thus
$X_{i}$ is associated to the edge splitting of $\Delta$ which is over $S_{i}$.
Since each of $H$ and $S_{i}$ has two coends in $DG$, Proposition 7.4 of
\cite{SS2} shows that the associated almost invariant sets $X_{H}$ and $X_{i}$
must cross strongly if they cross at all. As $H\subset G=G(w)$, it is clear
that $X_{H}$ cannot cross any $X_{i}$ strongly, and so does not cross any
$X_{i}$ at all. Thus $X_{H}$ must be enclosed by one of the vertices of
$\Delta$.

If $X_{H}$ is enclosed by $\overline{w}$, then $H$ must be a subgroup of
$G(\overline{w})$. As $H$ is also a subgroup of $G(w)$, this implies that $H$
is conjugate into some $S_{i}$, as required. It also implies that $Y$ is
trivial, and so is automatically adapted to $\partial G$, completing the proof
in this case.

For the rest of this proof, we will assume that $X_{H}$ is enclosed by $w$.
Thus Lemma \ref{propertiesofadapted} shows that $Y=X_{H}\cap G$ is adapted to
$\partial G$, which proves the first part of the lemma.

Now suppose that $Y$ is trivial, so that one of $Y$ or $Y^{\ast}$ is
$H$--finite. We consider the action of $DG$ on the universal covering
$DG$--tree $T$ of $\Delta$. There is a vertex $v$ of $T$ with stabiliser
$G=G(v)$ which lies above $w$ such that $X_{H}$ is enclosed by $v$. Corollary
4.16 of \cite{SS2} tells us that $X_{H}$ determines a nontrivial partition of
the edges of $T$ which are incident to $v$. Suppose that $Y$ is $H$--finite
and let $e$ be an edge of $T$ which is incident to $v$ and on the $X_{H}%
$--side of $v$. Then $X_{H}$ is equivalent to an almost invariant subset $W$
of $G$ which contains $G(e)$. The fact that $Y=X_{H}\cap G(v)$ is $H$--finite
implies that $W\cap G(v)$ is also $H$--finite, and hence in particular that
$G(e)$ itself is $H$--finite. Thus $H\cap G(e)$ has finite index in $G(e)$.
Now $G(e)$ is a conjugate of some $S_{i}$. As $H$ and $S_{i}$ are both
$PD(n+1)$, a subgroup of $G(e)$ of finite index is also $PD(n+1)$ and hence of
finite index in $H$. It follows that $H$ is commensurable with a conjugate of
$S_{i}$. As $H$ is orientable, part 2) of Corollary
\ref{PDgroupsplitsoverPDHimpliesHismaximal} shows that $H$ is conjugate to a
subgroup of $S_{i}$. This completes the proof of the lemma.
\end{proof}

The conclusion of the above discussion is that if $(G,\partial G)\ $is an
orientable $PD(n+2)$ pair with non-empty boundary, and $H$ is an orientable
$PD(n+1)$ subgroup of $G$, we can associate to $H$ a $H$--almost invariant
subset $Y$ of $G$, and $Y$ will be nontrivial so long as $H$ is not conjugate
into some $S_{i}$. Further $Y$ is adapted to $\partial G$. We will say that
$Y$ is \textit{dual} to $H$. When $H$ is $VPC(n+1)$ and $Y$ is nontrivial, we
will say that $H$ is an \textit{essential} torus in $(G,\partial G)$. If $G$
is an orientable $PD(n+2)$ group, then, as above, we associate $X_{H}$ itself
to $H$. In this case, we may also refer to $H$ as an essential torus in $G$,
though the word `essential' is redundant in this case.

For later reference, we briefly consider the situation where $(G,\partial
G)\ $is an orientable $PD(n+2)$ pair with non-empty boundary, and $H$ is a
non-orientable $PD(n+1)$ subgroup of $G$. Unlike the case when $\partial G$ is
empty, $G$ may possess nontrivial $H$--almost invariant subsets. However the
following result says that no such subset of $G$ can be adapted to $\partial
G$.

\begin{lemma}
\label{aisetovernon-orblesubgroupisnotadaptedtodG}Let $(G,\partial G)\ $be an
orientable $PD(n+2)$ pair with non-empty boundary, let $H$ be a non-orientable
$PD(n+1)$ subgroup of $G$, and let $X$ be a nontrivial $H$--almost invariant
subset of $G$. Then $X$ is not adapted to $\partial G$.
\end{lemma}

\begin{proof}
As usual we let $DG$ denote the double of $G$ along $\partial G$, so that $DG$
is an orientable $PD(n+2)$ group. We recall that if $K$ is a $PD(n+1)$
subgroup of $DG$, then $e(DG,K)$ equals $2$ if $K$ is orientable and equals
$1$ otherwise. In particular, $DG$ has no nontrivial $H$--almost invariant subset.

Suppose that $X$ is adapted to $\partial G$. Then Lemma
\ref{propertiesofadapted} tells us that there is a $H$--almost invariant
subset $\overline{X}$ of $DG$ such that $\overline{X}\cap G$ equals $X$. As
$X$ is a nontrivial almost invariant set, so is $\overline{X}$. This
contradiction shows that $X$ cannot be adapted to $\partial G$, as required.
\end{proof}

For the rest of this section, we will discuss essential maps of higher
dimensional `annuli' into orientable $PD(n+2)$ pairs. It turns out that we
need two types of higher dimensional analogue of an annulus. The first and
most obvious type is a trivial orientable $PD(n+1)$ pair $A_{H}=(H;\{H,H\})$,
where $n\geq1$, and $H$ is an orientable $PDn$ group which is also $VPCn$. We
call this an \textit{untwisted annulus}. The second type is an orientable
$PD(n+1)$ pair $\Lambda_{H}=(H,H_{0})$, where $H$ is a non-orientable $PDn$
group which is also $VPCn$, and $H_{0}$ is the orientation subgroup of $H$. We
call this a \textit{twisted annulus}. When $n=2$, an example of a
$3$--manifold of this type is the twisted $I$--bundle over the Klein bottle
with orientable total space. Note that there are no twisted annuli when $n=1$,
as any $PD1$ group is orientable. In particular, there are no twisted annuli
in a $3$--manifold. Algebraically, an annulus in a $PD(n+2)$ pair $(G,\partial
G)$ is an injective homomorphism of group pairs $\Theta:A_{H}\rightarrow
(G,\partial G)$, or $\Theta:\Lambda_{H}\rightarrow(G,\partial G)$. This means
that $\Theta$ maps $H$ to $G$ and also maps each group in $\partial\Lambda
_{H}$ to a conjugate of some group in $\partial G$. Again the above
orientation restrictions are not crucial, but they do simplify the statements.
What is crucial is that $\Theta$ should commute with the orientation
homomorphisms, so that, in topological terms, our annuli have trivial normal
bundle or equivalently are two-sided. The reason for this is that we need to
associate a nontrivial $H$--almost invariant subset of $G$ to each essential
annulus. We will do this in Definition \ref{defnofa.i.setassociatedtoannulus}.

It will be very helpful to consider a map $\theta$ of aspherical spaces such
that the induced map on fundamental groups is $\Theta$. For this we need to
choose a $K(H,1)$, and use a mapping cylinder construction to make a $K(G,1)$
with the $K(S_{i},1)$'s as disjoint subcomplexes, for $i\geq1$. To simplify
the notation for an untwisted annulus, we write $A$ for $K(H,1)\times I$, and
$\partial A$ for $K(H,1)\times\partial I$. The two components of $\partial A$
will be denoted by $\partial_{0}A$ and $\partial_{1}A$. If $n=1$, then $A$ can
be chosen to be the usual annulus $S^{1}\times I$. For a twisted annulus, we
write $A$ for the twisted $I$--bundle over $K(H,1)$ determined by the
orientation homomorphism of $H$, and write $\partial A$ for the induced
$S^{0}$--bundle. Finally we write $M$ for $K(G,1)$, and $\partial M$ for the
union of the $K(S_{i},1)$'s, for $i\geq1$. Then $\Theta$ is induced by a map
$\theta:(A,\partial A)\rightarrow(M,\partial M)$. Note that in the untwisted
case, such a map $\theta$ is determined up to homotopy by choosing a copy of
$H$ in two conjugates of groups in $\partial G$, such that the two copies of
$H$ are conjugate in $G$. And in the twisted case, $\theta$ is determined up
to homotopy by choosing a copy of $H$ in $G$ and a conjugate of some group in
$\partial G$ such that the intersection of $H$ with this conjugate contains
$H_{0}$. Thus an annulus can be thought of purely algebraically. We will say
that $\Theta$ is \textit{essential} if $\theta$ cannot be homotoped relative
to $\partial A$ into $\partial M$. It is clear that the essentiality of an
annulus is also a purely algebraic property. An untwisted annulus is essential
if and only if the images of the two boundary groups are not conjugate in a
group in $\partial G$. And a twisted annulus is essential if and only if
$H_{0}$ lies in a boundary group $K$ in $\partial G$, and $H\cap K=H_{0}$.
Note that as $G\ $is finitely generated, we can choose $M$ to have finite
$1$--skeleton. If $G\ $is finitely presented, we can also choose $M$ to have
finite $2$--skeleton, as each $S_{i}$ is finitely generated.

Now suppose that $\theta$ is an essential map of an untwisted annulus $A$ into
$M$, where $\pi_{1}(A)$ is equal to $H$, and identify $H$ with its image in
$G$ under $\theta_{\ast}$. Let $M_{H}$ denote the cover of $M$ with
fundamental group $H$, and let $\widetilde{M}$ denote the universal cover of
$M$. Let $\theta_{H}:(A,\partial A)\rightarrow(M_{H},\partial M_{H})$ be the
lift of $\theta$. The induced map on homology sends the fundamental class
$[A]\in H_{n+1}(A,\partial A;\mathbb{Z})$ to an element $\alpha$ of
$H_{n+1}(M_{H},\partial M_{H};\mathbb{Z})$. Let $\Sigma$ denote the component
of $\partial M_{H}$ which contains $\theta_{H}(\partial_{0}A)$. We will assume
that the base point of $A$ lies in $\partial_{0}A$. As $M_{H}$ and $\Sigma$
are aspherical and the inclusion of $\Sigma$ in $M_{H}$ induces an isomorphism
of fundamental groups, $\Sigma$ is a deformation retract of $M_{H}$. As
$\theta$ is essential, it follows that $\theta_{H}(\partial_{1}A)$ must lie in
a different component $\Sigma^{\prime}$ of $\partial M_{H}$. This implies that
$\alpha$ is a nontrivial element of $H_{n+1}(M_{H},\partial M_{H};\mathbb{Z})$.

If $\theta$ is an essential map into $M$ of a twisted annulus $A$ with
$\pi_{1}(A)$ equal to $H$, we again identify $H$ with its image in $G$ under
$\theta_{\ast}$. Recall that $\pi_{1}(\partial A)$ equals the orientation
subgroup $H_{0}$ of $H$. Let $M_{H}$ denote the cover of $M$ with fundamental
group $H$, and let $\theta_{H}:(A,\partial A)\rightarrow(M_{H},\partial
M_{H})$ be the lift of $\theta$. Again we let $\alpha$ denote the image in
$H_{n+1}(M_{H},\partial M_{H};\mathbb{Z})$ of the fundamental class $[A]\in
H_{n+1}(A,\partial A;\mathbb{Z})$. Let $\Sigma$ denote the component of
$\partial M_{H}$ which contains $\theta_{H}(\partial A)$, and consider the
inclusions $H_{0}=\pi_{1}(\partial A)\subset\pi_{1}(\Sigma)\subset\pi
_{1}(M_{H})=H$. If $\pi_{1}(\Sigma)$ equals $H$, then $M_{H}$ deformation
retracts to $\Sigma$, which contradicts the hypothesis that $\theta$ is
essential. It follows that we must have $\pi_{1}(\Sigma)=H_{0}$. In turn this
implies that $\alpha$ is a nontrivial element of $H_{n+1}(M_{H},\partial
M_{H};\mathbb{Z})$. We note that the double cover $A_{0}$ of $A$ with $\pi
_{1}(A_{0})$ equal to $H_{0}$ is an untwisted annulus, and the induced map
$\theta_{0}:A_{0}\rightarrow M_{0}$ is also essential, where $M_{0}$ is the
cover of $M$ with fundamental group $H_{0}$.

\begin{remark}
The above discussion shows that if $A$ is an annulus, twisted or untwisted,
and $\theta:(A,\partial A)\rightarrow(M,\partial M)$ is an essential map, then
in either case the image of the fundamental cycle $[A]$ with $\mathbb{Z}_{2}%
$-coefficients is also nontrivial in $H_{n+1}(M_{H},\partial M_{H}%
;\mathbb{Z}_{2})$ and is the specialization of $\alpha$.
\end{remark}

We also note that, whether or not $A$ is twisted, the induced action of $H$ on
the universal cover $\widetilde{M}$ of $M$ preserves the union of two distinct
components of $\partial\widetilde{M}$. If some element of $H$ interchanges
these two components, then $A$ is twisted. Otherwise, $A$ is untwisted.

Conversely, suppose that two distinct components $\Sigma$ and $T$ of
$\partial\widetilde{M}$ are each stabilised by an orientable $PDn$ subgroup
$H$ of $G$. Then there is an essential map $\theta$ of an untwisted annulus
$A$ to $M$ with $\pi_{1}(A)$ equal to $H$, and a lift $\theta_{H}:A\rightarrow
M_{H}$ which maps $\partial_{0}A$ to $H\backslash\Sigma$ and maps
$\partial_{1}A$ to $H\backslash T$. We denote this essential untwisted annulus
in $(M,\partial M)$ by $H_{\Sigma,T}$.

If the union of $\Sigma$ and $T$ is stabilised by a $PDn$ subgroup $H$ of $G$
and if some element of $H$ interchanges $\Sigma$ and $T$, let $H_{0}$ denote
the subgroup of $H$ of index $2$ which stabilises both $\Sigma$ and $T$. If
$H_{0}$ is the orientation subgroup of $H$, there is an essential map $\theta$
of a twisted annulus $A$ to $M$ with $\pi_{1}(A)$ equal to $H$, and the
induced essential map $\theta_{0}$ of the double cover $A_{0}$ of $A$ with
$\pi_{1}(A_{0})$ equal to $H_{0}$ has a lift to $M_{H_{0}}$ which maps
$\partial_{0}A_{0}$ to $H\backslash\Sigma$ and maps $\partial_{1}A_{0}$ to
$H\backslash T$. We denote this essential twisted annulus in $(M,\partial M)$
by $H_{\Sigma,T}$.

In order to associate an almost invariant subset of $G$ to an essential
annulus, we will need the following result, which was proved by Swarup in
\cite{Swarup3}. See \cite{KR3} for a purely algebraic proof. The statement we
give here is equivalent to Theorem 2 of \cite{Swarup3}.

\begin{lemma}
\label{swarup}(Swarup) Let $G$ be a finitely generated group, and let $H$ be a
subgroup of infinite index in $G$. If $\varphi$ is a homomorphism from $H$ to
$\mathbb{Z}$, denote the kernel of $\varphi$ by $N$. Suppose that whenever
$\varphi$ is non-zero, we have $e(G,N)=1$. Then the restriction map
$r:H^{1}(G;\mathbb{Z}[H\backslash G])\rightarrow H^{1}(H;\mathbb{Z})$ is trivial.
\end{lemma}

\begin{proof}
As before, we let $M$ denote a $K(G,1)$ with the $K(S_{i},1)$'s as disjoint
subcomplexes, for $i\geq1$. Recall from Remark \ref{defnofH1fwhenGisnotfp},
that we can identify $H^{1}(G,\mathbb{Z}[H\backslash G])$ with $H_{f}%
^{1}(M_{H};\mathbb{Z})$. If $G$ is finitely generated but not finitely
presented, we need to modify the usual definition of $H_{f}^{1}(M_{H}%
;\mathbb{Z})$ as discussed there. In any case, $M$ need not be locally finite,
so any reference we make to the number of ends of a cover of $M$ really refers
to the number of ends of the $1$--skeleton of the cover. We also let $r$
denote the natural map $H_{f}^{1}(M_{H};\mathbb{Z})\rightarrow H^{1}%
(M_{H};\mathbb{Z})$.

Let $\beta$ denote an element of $H_{f}^{1}(M_{H};\mathbb{Z})$, and let
$\gamma\in H^{1}(M_{H};\mathbb{Z})$ denote $r(\beta)$. Thus $\gamma$ can be
represented by a map $g:M_{H}\rightarrow S^{1}$. Represent $\beta$ by a finite
cocycle $c:M_{H}^{(1)}\rightarrow\mathbb{Z}$ on the $1$--skeleton $M_{H}%
^{(1)}$ of $M_{H}$, and let $\Sigma$ denote the support of $c$. Thus $\Sigma$
is a finite subcomplex of $M_{H}^{(1)}$, and $c$ restricted to any simplex of
$M_{H}^{(1)}-\Sigma$ is trivial. As $c$ also represents $\gamma$, for any
component $L$ of $M_{H}^{(1)}-\Sigma$, the map $\pi_{1}(L)\rightarrow\pi
_{1}(S^{1})$ induced by $g$ is trivial. Since $\Sigma$ is finite, $M_{H}%
^{(1)}-\Sigma$ has at least one unbounded component $L$ whose coboundary
$\delta L$ must be finite, as $\delta L\subset\delta\Sigma$.

Now suppose that $\gamma$ is non-zero. We consider the induced map $g_{\ast
}:H\rightarrow\mathbb{Z}$, let $N$ denote the kernel of $g_{\ast}$, and
consider the cover $q_{N}:M_{N}\rightarrow M_{H}$. As the map $\pi
_{1}(L)\rightarrow\pi_{1}(S^{1})$ induced by $g$ is trivial, $L$ lifts to
$M_{N}$. As the infinite quotient group $H/N$ acts on $M_{N}$, it follows that
$L$ has infinitely many disjoint lifts to $M_{N}$. This implies that $M_{N}$
has infinitely many ends, which contradicts the hypothesis that $e(G,N)=1$. It
follows that $\gamma$ must be zero which completes the proof of the lemma.
\end{proof}

We will want to apply this result to Poincar\'{e} duality pairs. The result we
obtain is the following.

\begin{lemma}
\label{restrictionmapistrivial}Let $(G,\partial G)\ $be an orientable
$PD(n+2)$ pair, and let $H$ be a $PDn$ subgroup of $G$. Then the restriction
map $r:H^{1}(G;\mathbb{Z}[H\backslash G])\rightarrow H^{1}(H;\mathbb{Z})$ is trivial.
\end{lemma}

\begin{remark}
The $PDn$ subgroup $H$ need not be orientable.
\end{remark}

\begin{proof}
As before, we let $M$ denote a $K(G,1)$ with the $K(S_{i},1)$'s as disjoint
subcomplexes, for $i\geq1$. We need to check all the hypotheses of Lemma
\ref{swarup}. Certainly $G$ is finitely generated and $H\ $has infinite index
in $G$. Now suppose that $N$ is the kernel of a non-zero homomorphism from $H$
to $\mathbb{Z}$. We need to show that $e(G,N)=1$. As $N$ has infinite index in
the $PDn$ group $H$, a theorem of Strebel \cite{Strebel} tells us that $N$ has
cohomological dimension $\leq n-1$. Hence $H_{k}(N;\mathbb{Z})=0$, for any
$k\geq n$, and the same holds for any subgroup of $N$. In particular,
$H_{n+1}(M_{N};\mathbb{Z})$ and $H_{n}(\partial M_{N};\mathbb{Z})$ are both
zero. It follows from the exact sequence of the pair $(M_{N},\partial M_{N})$
that $H_{n+1}(M_{N},\partial M_{N};\mathbb{Z})$ is zero. By Poincar\'{e}
duality, this implies that $H_{f}^{1}(M_{N};\mathbb{Z})$ is zero, so that
$M_{N}$ has only one end. Hence $e(G,N)=1$ as required.
\end{proof}

Our real interest lies in the corresponding restriction map with
$\mathbb{Z}_{2}$ coefficients in place of $\mathbb{Z}$. This map need not be
trivial, but the fact that $r$ is trivial yields enough information about the
case of $\mathbb{Z}_{2}$ coefficients for our purposes.

\begin{corollary}
\label{imageofrhoequalsimageofdelta}Let $(G,\partial G)\ $be an orientable
$PD(n+2)$ pair, and let $H$ be a $PDn$ subgroup of $G$. Let $\rho$ denote the
map of cohomology groups given by reduction of the coefficients modulo $2$.
Then the image of $\rho:H^{1}(G;\mathbb{Z}[H\backslash G])\rightarrow
H^{1}(G;\mathbb{Z}_{2}[H\backslash G])$ is contained in the image of the
coboundary map $\delta:H^{0}(G;P[H\backslash G]/\mathbb{Z}_{2}[H\backslash
G])\rightarrow H^{1}(G;\mathbb{Z}_{2}[H\backslash G]),$ given on page
\pageref{coboundarymap}.
\end{corollary}

\begin{proof}
As before, we let $M$ denote a $K(G,1)$ with finite $1$--skeleton, and let
$M_{H}$ denote the cover of $M$ with fundamental group $H$.

Consider the diagram%
\[%
\begin{array}
[c]{ccccc}
&  & H^{1}(G;\mathbb{Z}[H\backslash G]) & \overset{r}{\rightarrow} &
H^{1}(H;\mathbb{Z})\\
&  & \downarrow\rho &  & \downarrow\rho\\
H^{0}(G;P[H\backslash G]/\mathbb{Z}_{2}[H\backslash G]) & \overset{\delta
}{\rightarrow} & H^{1}(G;\mathbb{Z}_{2}[H\backslash G]) & \overset{\overline
{r}}{\rightarrow} & H^{1}(H;\mathbb{Z}_{2})\\
\cong\downarrow &  & \cong\downarrow &  & \cong\downarrow\\
H_{e}^{0}(M_{H};\mathbb{Z}_{2}) & \overset{\delta}{\rightarrow} & H_{f}%
^{1}(M_{H};\mathbb{Z}_{2}) & \overset{\overline{r}}{\rightarrow} & H^{1}%
(M_{H};\mathbb{Z}_{2})
\end{array}
\]
where the bottom two rows come from the long exact cohomology sequences given
before Remark \ref{defnofH1fwhenGisnotfp}. Lemma \ref{restrictionmapistrivial}
tells us that the map $r$ is zero. As the image of $\delta$ equals the kernel
of the map $\overline{r}$, it follows immediately that the image of $\rho$ is
contained in the image of $\delta$.
\end{proof}

Now we can associate an almost invariant set to an essential annulus as promised.

\begin{definition}
\label{defnofa.i.setassociatedtoannulus} Let $(G,\partial G)\ $be an
orientable $PD(n+2)$ pair, and let $\theta$ be an essential annulus in
$(M,\partial M)$ with fundamental group $H$.

As discussed before Lemma \ref{swarup}, the essential annulus $\theta$
determines a non-zero element $\alpha\in H_{n+1}(M_{H},\partial M_{H}%
;\mathbb{Z})\cong H_{n+1}(G,\partial G;\mathbb{Z}[H\backslash G])$. The
Poincar\'{e} dual of $\alpha$, regarded as an element of this second group, is
a nontrivial element $\beta$ of $H^{1}(G;\mathbb{Z}[H\backslash G])$.
Corollary \ref{imageofrhoequalsimageofdelta} shows that $\rho(\beta)$ is
contained in the image of $\delta$, and so, as discussed on page
\pageref{coboundarymap}, determines an almost invariant subset of $H\backslash
G$ modulo almost equality and complementation. The pre-image in $G$ of such a
set is a $H$--almost invariant subset $X_{\theta}$ of $G$.

We will say that $X_{\theta}$ is \textsl{dual} to the essential annulus
$\theta$. On occasion, it will also be convenient to say that the almost
invariant subset $H\backslash X$ of $H\backslash G$ is \textsl{dual} to
$\theta$.
\end{definition}

If $(G,\partial G)\ $is an orientable $PD(n+2)$ pair, and the almost invariant
subset of $G$ associated to a splitting is dual to an essential annulus or
torus, we will say that the splitting itself is dual to an essential annulus
or torus, as appropriate.

If $(G,\partial G)\ $is an orientable $PD(n+2)$ pair which admits an essential
annulus, the dual almost invariant subset of $G$ is nontrivial and is over a
$VPCn$ subgroup. The converse is not true. In general $G$ will have many
nontrivial such subsets which are not dual to any essential annulus. In the
following two results, we consider this in more detail. In particular, we show
that if $G$ has a nontrivial almost invariant subset over a $VPCn$ subgroup,
then $(G,\partial G)$ admits an essential annulus.

\begin{proposition}
\label{aisetimpliesannulus}Let $(G,\partial G)$ be an orientable $PD(n+2)$
pair and let $H$ be a $VPCn$ subgroup of $G$. Let $\alpha$ be an element of
$H_{n+1}(M_{H},\partial M_{H};\mathbb{Z})$, so that $\partial\alpha\in
H_{n}(\partial M_{H};\mathbb{Z})$ is supported by some finite number $k$ of
components of $\partial M_{H}$. Then the following statements hold:

\begin{enumerate}
\item If $\alpha$ is non-zero, then $k$ is non-zero.

\item If $\alpha$ is non-zero, then each of the $k$ components of $\partial
M_{H}$ which support $\partial\alpha$ carries a subgroup of finite index in
$H$. If $H$ is orientable then $k\geq2$.

\item If $G$ has a nontrivial $H$--almost invariant subset, then $(G,\partial
G)$ admits an essential map of an annulus whose fundamental group is a
subgroup of finite index in $H$.
\end{enumerate}
\end{proposition}

\begin{proof}
1) Consider the exact sequence.
\[
H_{n+1}(M_{H};\mathbb{Z})\rightarrow H_{n+1}(M_{H},\partial M_{H}%
;\mathbb{Z})\overset{\partial}{\rightarrow}H_{n}(\partial M_{H};\mathbb{Z}%
)\rightarrow H_{n}(M_{H};\mathbb{Z})
\]
Since $H$ is a $PDn$ group, it follows that $H_{n+1}(M_{H};\mathbb{Z})$ is
zero. Thus the map $\partial$ in this sequence is injective. If $\alpha$ is
non-zero, it follows that $\partial\alpha$ is non-zero, and hence that $k$ is non-zero.

2) If $\Sigma$ is one of the components of $\partial M_{H}$ which supports
$\partial\alpha$, then $H_{n}(\Sigma;\mathbb{Z})$ must be nontrivial. Note
that the fundamental group of any component of $\partial M_{H}$ is a subgroup
of $H$. As $H$ is $PDn$, Strebel's result in \cite{Strebel} implies that a
subgroup of $H$ of infinite index has cohomological dimension less than $n$.
Hence $\pi_{1}(\Sigma)$ must have finite index $d$ in $H$, so that
$H_{n}(\Sigma;\mathbb{Z})$ is infinite cyclic. Now suppose that $H$ is
orientable. Then $H_{n}(M_{H};\mathbb{Z})$ must also be infinite cyclic.
Further the map $H_{n}(\Sigma;\mathbb{Z})\rightarrow H_{n}(M_{H};\mathbb{Z})$
is multiplication by $d$, and so is injective. We know that the map
$H_{n}(\partial M_{H};\mathbb{Z})\rightarrow H_{n}(M_{H};\mathbb{Z})$ has
nontrivial kernel as it contains $\partial\alpha$. It follows that $\partial
M_{H}$ has a second boundary component $\Sigma^{\prime}$ which carries
$\partial\alpha$. Thus when $H$ is orientable, we must have $k\geq2$.

3) If $G$ has a nontrivial $H$--almost invariant subset, Corollary
\ref{imageofrhoequalsimageofdelta} implies that $H^{1}(G;\mathbb{Z}%
[H\backslash G])$, and hence $H_{n+1}(M_{H},\partial M_{H};\mathbb{Z})$, is
nontrivial. Suppose first that $H$ is orientable. Then part 2) shows that
there are two distinct components $\Sigma$ and $\Sigma^{\prime}$ of $\partial
M_{H}$ such that each of $\pi_{1}(\Sigma)$ and $\pi_{1}(\Sigma^{\prime})$ is a
subgroup of finite index in $H$. It follows that $(G,\partial G)$ admits an
essential map of an untwisted annulus whose fundamental group is a subgroup of
finite index in $H$. If $H$ is non-orientable, we consider the orientation
subgroup $H_{0}$ of $H$ which is of index $2$. As $H_{0}$ is a subgroup of
$H\ $of finite index, $G$ has a nontrivial $H_{0}$--almost invariant subset.
Now we apply the above discussion to $H_{0}$ in place of $H$ and obtain an
essential map of an untwisted annulus whose fundamental group is a subgroup of
finite index in $H_{0}$, and hence of finite index in $H$, as required. This
completes the proof of the proposition.
\end{proof}

In the next result, we prove more.

\begin{proposition}
\label{annuligenerate}Let $(G,\partial G)$ be an orientable $PD(n+2)$ pair and
let $H$ be a $VPCn$ subgroup of $G$. Then we have the following results.

\begin{enumerate}
\item No nontrivial $H$--almost invariant subset of $G$ can be adapted to
$\partial G$.

\item If $G$ has a nontrivial $H$--almost invariant subset $X$, there is an
orientable subgroup $H^{\prime}$ of finite index in $H$, such that $X$ is
equivalent to a sum of $H^{\prime}$--almost invariant subsets of $G$ each dual
to an untwisted annulus.

\item If $\partial M_{H}$ has $k$ components each of which carries a subgroup
of finite index in $H$, then the number of coends of $H$ in $G$ is at least
$k$.

\item If the number of coends of $H$ in $G$ is at least $k$, then $H$ has a
subgroup $L$ of finite index such that $\partial M_{L}$ has $k$ components
each of which carries $L$.

\item If $G$ has a nontrivial $H$--almost invariant subset $X$ which crosses
no nontrivial almost invariant subset of $G$ over any finite index subgroup of
$H$, then $X$ is dual to an annulus.
\end{enumerate}
\end{proposition}

\begin{proof}
1) Suppose there is a nontrivial $H$--almost invariant subset $Y$ of $G$ which
is adapted to $\partial G$. We will consider the $PD(n+2)$ group $DG$ obtained
by doubling $G\ $along $\partial G$, and the corresponding graph of groups
decomposition $\Delta$ of $DG$, which has a vertex $w$ with associated group
$G$. Lemma \ref{propertiesofadapted} tells us that there is a nontrivial
$H$--almost invariant subset $X$ of $G$ which is enclosed by $w$ such that
$X\cap G$ equals $Y$. But part 2) of Lemma
\ref{PD(n+2)groupshavenoa.i.setsoverVPCk<n} tells us that $DG$ has no
nontrivial almost invariant subsets over $VPCn$ subgroups. This contradiction
shows that no nontrivial $H$--almost invariant subset of $G$ can be adapted to
$\partial G$, as required.

2) Let $X$ be a nontrivial $H$--almost invariant subset of $G$, and let $Y$
denote the almost invariant subset $H\backslash X$ of $H\backslash G$. As
discussed before Remark \ref{defnofH1fwhenGisnotfp}, the equivalence class of
$Y$ under almost equality is an element of $H^{0}(G;P[H\backslash
G]/\mathbb{Z}_{2}[H\backslash G])$. We let $[Y]$ denote the image of this
equivalence class in $H^{1}(G;\mathbb{Z}_{2}[H\backslash G])$ under the
coboundary map $\delta$ given on page \pageref{coboundarymap}. Thus $[Y]$ is
represented by any almost invariant subset of $H\backslash G$ which is almost
equal to $Y$ or to $Y^{\ast}$. As $H$ is $VPCn$, and torsion free, it is
$PDn$. Thus Corollary \ref{imageofrhoequalsimageofdelta} tells us that there
is an element $\beta$ of $H^{1}(G;\mathbb{Z}[H\backslash G])$ such that
$\rho(\beta)=[Y]$. Let $\alpha$ denote the element of $H_{n+1}(G,\partial
G;\mathbb{Z}[H\backslash G])$ which is Poincar\'{e} dual to $\beta$. Regard
$\alpha$ as an element of $H_{n+1}(M_{H},\partial M_{H};\mathbb{Z})$, and
consider $\partial\alpha\in H_{n}(\partial M_{H};\mathbb{Z})$. From
Proposition \ref{aisetimpliesannulus}, $\partial\alpha$ is non-zero, so there
is at least one component of $\partial M_{H}$ which lies in its support. Each
component $\Sigma$ of $\partial M_{H}$ which supports $\partial\alpha$ carries
a subgroup of $H$ of finite index. Let $H^{\prime\prime}$ denote the
intersection of all conjugates in $H$ of these subgroups, let $H_{0}$ denote
the maximal orientable subgroup of $H$ of index at most $2$, and let
$H^{\prime}$ denote the intersection $H^{\prime\prime}\cap H_{0}$. Thus
$H^{\prime}$ is an orientable normal subgroup of $H$ of finite index. By
replacing $H$ by $H^{\prime}$, we can assume that $H$ is orientable and that
the finitely many components $\Sigma_{1},\ldots,\Sigma_{k}$ of $\partial
M_{H}$ which support $\partial\alpha$ all carry $H$. As $H$ is orientable, we
must have $k\geq2$.

If $i$ and $j$ are distinct integers, there is an essential untwisted annulus
$A_{ij}$ in $M_{H}$ with fundamental group $H$ and whose boundary components
lie in $\Sigma_{i}$ and $\Sigma_{j}$. Let $\alpha_{ij}$ denote the image of
the fundamental class of $A_{ij}$ in $H_{n+1}(M_{H},\partial M_{H}%
;\mathbb{Z}).$ Then $\partial\alpha_{ij}\in H_{n}(\partial M_{H};\mathbb{Z})$
is supported by $\Sigma_{i}$ and $\Sigma_{j}$. Regard $\alpha_{ij}$ as an
element of $H_{n+1}(G,\partial G;\mathbb{Z}[H\backslash G])$, and let
$\beta_{ij}\in H^{1}(G;\mathbb{Z}[H\backslash G])$ be Poincar\'{e} dual to
$\alpha_{ij}$. Recall that $\rho$ denotes reduction of coefficients mod $2$,
and that Corollary \ref{imageofrhoequalsimageofdelta} tells us that
$\rho(\beta_{ij})\in H^{1}(G;\mathbb{Z}_{2}[H\backslash G])$ equals $[Y_{ij}%
]$, for some almost invariant subset $Y_{ij}$ of $H\backslash G$ which is said
to be dual to $A_{ij}$.

Let $\overline{Y}$ denote $\rho\partial\alpha\in H_{n}(\partial M_{H}%
;\mathbb{Z}_{2})$, and let $\overline{Y_{ij}}$ denote $\rho\partial\alpha
_{ij}\in H_{n}(\partial M_{H};\mathbb{Z}_{2})$. Then $\overline{Y_{ij}}$ is
supported by $\Sigma_{i}$ and $\Sigma_{j}$. Let $\gamma$ denote any nontrivial
element of $H_{n+1}(M_{H},\partial M_{H};\mathbb{Z}_{2})$ whose image in
$H_{n}(\partial M_{H};\mathbb{Z}_{2})$ is supported by some subset of the
$\Sigma_{i}$'s. As each $\Sigma_{i}$ carries $H$, this image must be supported
by at least two components of $\partial M_{H}$. Thus a simple induction
argument on $k$ shows that $\overline{Y}$ must be equal to a sum of
$\overline{Y_{ij}}$'s.

Now consider the following commutative diagram. The vertical maps are
Poincar\'{e} duality isomorphisms, the top horizontal map is the boundary map
in the long exact homology sequence of the pair $(M_{H},\partial M_{H})$, and
the bottom horizontal map is induced by the inclusion of $\partial M_{H}$ into
$M_{H}$.%
\[%
\begin{array}
[c]{ccc}%
H_{n+1}(M_{H},\partial M_{H};\mathbb{Z}_{2}) & \overset{\partial}{\rightarrow}
& H_{n}(\partial M_{H};\mathbb{Z}_{2})\\
\downarrow\cong &  & \downarrow\cong\\
H_{f}^{1}(M_{H};\mathbb{Z}_{2}) & \overset{i^{\ast}}{\rightarrow} & H_{f}%
^{1}(\partial M_{H};\mathbb{Z}_{2})
\end{array}
\]
As $H$ is $PDn$, it follows that $H_{n+1}(M_{H};\mathbb{Z}_{2})$ is zero. Thus
the map $\partial$ in this diagram must be injective, so that $i^{\ast}$ is
also injective.

Recall that $H^{1}(G;\mathbb{Z}_{2}[H\backslash G])$ and $H_{f}^{1}%
(M_{H};\mathbb{Z}_{2})$ are naturally isomorphic, so that we can identify
$[Y]$ and $[Y_{ij}]$ with elements of $H_{f}^{1}(M_{H};\mathbb{Z}_{2})$. Now,
under Poincar\'{e} duality, $i^{\ast}[Y]$ corresponds to $\overline{Y}\in
H_{n}(\partial M_{H};\mathbb{Z}_{2})$, and $i^{\ast}[Y_{ij}]$ corresponds to
$\overline{Y_{ij}}$. As $\overline{Y}$ is equal to a sum of $\overline{Y_{ij}%
}$'s, and $i^{\ast}$ is injective, it follows that $[Y]$ is equal to a sum of
$[Y_{ij}]$'s. Thus $Y$ is equivalent to a sum of $Y_{ij}$'s and their
complements. Hence $X$ is equivalent to a sum of $H$--almost invariant subsets
of $G$ each dual to an untwisted annulus in $M_{H}$, as required.

3) Suppose that $\partial M_{H}$ has $k$ components each of which carries a
subgroup of finite index in $H$. In order to show that the number of coends of
$H$ in $G$ is at least $k$, it suffices to show there is a subgroup
$H^{\prime}$ of $H$ of finite index such that $e(G,H^{\prime})$ is at least
$k$. As in part 2), by replacing $H$ by a suitable subgroup of finite index,
we can suppose that $\partial M_{H}$ has components $\Sigma_{1},\ldots
,\Sigma_{k}$ each carrying $H$. Consider the composite map%
\[
H_{e}^{0}(M_{H};\mathbb{Z}_{2})\overset{\delta}{\rightarrow}H_{f}^{1}%
(M_{H};\mathbb{Z}_{2})\overset{i^{\ast}}{\rightarrow}H_{f}^{1}(\partial
M_{H};\mathbb{Z}_{2})\overset{\cong}{\rightarrow}H_{n}(\partial M_{H}%
;\mathbb{Z}_{2}).
\]
If $1\leq i<j\leq k$, the almost invariant subset $Y_{ij}$ of $H\backslash G$
in part 2) determines an element of $H_{e}^{0}(M_{H};\mathbb{Z}_{2})$ whose
image $\overline{Y_{ij}}$ in $H_{n}(\partial M_{H};\mathbb{Z}_{2})$ under this
composite map is supported by $\Sigma_{i}$ and $\Sigma_{j}$. As the
$\overline{Y_{ij}}$'s span a $(k-1)$--dimensional subgroup of $H_{n}(\partial
M_{H};\mathbb{Z}_{2})$, and the kernel of $\delta$ is nontrivial, it follows
that $H_{e}^{0}(M_{H};\mathbb{Z}_{2})$ has dimension at least $k$, so that
$e(G,H)$ is at least $k$, as required.

4) If the number of coends of $H$ in $G$ is at least $k$, then $H$ has a
subgroup $H_{1}$ of finite index such that $e(G,H_{1})\geq k$. It follows from
part 2) of this proposition that $H_{1}$ has a subgroup $L$ of finite index
such that the space of $L$--almost invariant subsets of $G$ spanned by such
sets which are dual to an untwisted annulus has dimension at least $k$. Pick a
finite family of $L$--almost invariant subsets of $G$, each dual to an
untwisted annulus, which together span a space with dimension at least $k$.
The corresponding annuli in $M_{L}$ have boundaries in a finite family of
components of $\partial M_{L}$. As in part 2), we can replace $L$ by a
subgroup of finite index so that all these boundary components carry $L$. It
follows immediately that there must be at least $k$ such components of
$\partial M_{L}$, as required.

5) Let $X$ be a nontrivial $H$--almost invariant subset of $G$ which crosses
no nontrivial almost invariant subset of $G$ over any finite index subgroup of
$H$. In part 2) of this lemma, we showed that, after replacing $H$ by a
suitable subgroup of finite index, the almost invariant subset $Y=H\backslash
X$ of $H\backslash G$ is equivalent to a sum of almost invariant subsets
$Y_{ij}$, $i\neq j$, of $H\backslash G$ where $Y_{ij}$ is dual to an untwisted
annulus in $M_{H}$ with boundary in $\Sigma_{i}\cup\Sigma_{j}$. By
re-labelling the $\Sigma_{i}$'s if needed, we can assume that $1\leq i,j\leq
m$ and that each index between $1$ and $m$ occurs.

If $m=2$, then $Y$ is equivalent to $Y_{12}$, and we are done. We will show
that no other case is possible.

If $m=3$, then, after renumbering, $Y$ must be equivalent to $Y_{12}+Y_{23}$.
But this implies that $\overline{Y}$ is supported on $\Sigma_{1}$ and
$\Sigma_{3}$ which contradicts our assumption that $m=3$.

Now suppose that $m\geq4$. Recall that an almost invariant subset $Y$ of
$H\backslash G$ determines an element $[Y]$ of $H_{e}^{0}(M_{H};\mathbb{Z}%
_{2})$ which can be represented by a $0$--cochain with finite coboundary. The
support $Z$ of this cochain is an infinite subset of vertices in $M_{H}$, with
infinite complement $Z^{\ast}$. Similarly each $Y_{ij}$ yields a corresponding
infinite subset of vertices $Z_{ij}$ in $M_{H}$, with infinite complement
$Z_{ij}^{\ast}$. It will be convenient to consider $Z$ and the $Z_{ij}$'s
rather than $Y$ and the $Y_{ij}$'s. Lemma \ref{thereexisttwocrossingannuli}
below shows that there are distinct integers $i$, $j$, $k$ and $l$ such that
the vertex set of $\Sigma_{k}$ is almost contained in $Z_{ij}$, and the vertex
set of $\Sigma_{l}$ is almost contained in $Z_{ij}^{\ast}$. This implies that
the four corners of the pair $(Z,Z_{ij})$ are infinite, as $Z$ and $Z^{\ast}$
each meet both $\Sigma_{k}$ and $\Sigma_{l}$ in an infinite set of vertices.
Hence the four corners of the pair $(Y,Y_{ij})$ are infinite, so that $Y$
crosses $Y_{ij}$. This contradicts our hypothesis that $X$ crosses no
nontrivial almost invariant subset of $G$ over any finite index subgroup of
$H$, which completes the proof of the lemma.
\end{proof}

In the next lemma we consider how two untwisted annuli can cross, by which we
mean that the dual almost invariant sets cross. The corresponding picture in
the $3$--manifold setting is very simple. Start with a $2$--disc $D$ with four
disjoint open intervals in its boundary. Then remove the rest of $\partial D$.
The resulting surface has four boundary components, and it is trivial that of
the six arcs which join pairs of distinct boundary components, there are two
which cross. The product of this manifold with $S^{1}$ is a $3$--manifold $M$
with four annulus boundary components, and there are two annuli in $M$ which cross.

\begin{lemma}
\label{thereexisttwocrossingannuli} Let $(G,\partial G)$ be an orientable
$PD(n+2)$ pair and let $H$ be an orientable $VPCn$ subgroup of $G$. Suppose
that $\Sigma_{1},\ldots,\Sigma_{4}$ are distinct components of $\partial
M_{H}$ each with fundamental group $H$. Let $Y_{ij}$ denote the almost
invariant subset of $H\backslash G$ dual to the annulus $A_{ij}$ in $M_{H}$
which has fundamental group $H$ and has boundary in $\Sigma_{i}$ and
$\Sigma_{j}$. Then there are distinct integers $i$, $j$, $k$ and $l$ such that
$Y_{ij}$ crosses $Y_{kl}$.
\end{lemma}

\begin{proof}
As at the end of the previous lemma, it will be convenient to consider
$Z_{ij}$ rather than $Y_{ij}$, where $Z_{ij}$ is the support of a $0$--cochain
on $M_{H}$ with finite coboundary which represents the element of $H_{e}%
^{0}(M_{H};\mathbb{Z}_{2})$ determined by $Y_{ij}$.

For three distinct integers $i$, $j$ and $k$, the vertex set of $\Sigma_{k}$
must be almost contained in $Z_{ij}$ or $Z_{ij}^{\ast}$. For simplicity we
will say that $\Sigma_{k}$ is almost contained in $Z_{ij}$ to mean that the
vertex set of $\Sigma_{k}$ is almost contained in $Z_{ij}$. Note that
$\Sigma_{i}$ and $\Sigma_{j}$ are not almost contained in $Z_{ij}$ or in
$Z_{ij}^{\ast}$.

Now let $i$, $j$, $k$ and $l$ be distinct integers. We will say that $Z_{ij}$
\textit{separates} $\Sigma_{k}$ and $\Sigma_{l}$ if $\Sigma_{k}$ is almost
contained in $Z_{ij}$, and $\Sigma_{l}$ is almost contained in $Z_{ij}^{\ast}%
$, or vice versa.

\medskip

\textbf{Claim:} $Y_{ij}$\textit{ crosses }$Y_{kl}$\textit{ if and only if
}$Z_{ij}$\textit{ separates }$\Sigma_{k}$\textit{ and }$\Sigma_{l}$\textit{,
and }$Z_{kl}$\textit{ separates }$\Sigma_{i}$\textit{ and }$\Sigma_{j}%
$\textit{.}

\medskip

Suppose first that $Z_{ij}$ separates $\Sigma_{k}$ and $\Sigma_{l}$. Then
$Y_{ij}$ must cross $Y_{kl}$, because all four corners of the pair
$(Z_{ij},Z_{kl})$ will be infinite as each has infinite intersection with
$\Sigma_{k}$ or $\Sigma_{l}$.

Next suppose that $Z_{ij}$ does not separate $\Sigma_{k}$ and $\Sigma_{l}$.
Without loss of generality, we can assume that $\Sigma_{k}$ and $\Sigma_{l}$
are both almost contained in $Z_{ij}$. Thus each of the corners $Z_{ij}^{\ast
}\cap Z_{kl}$ and $Z_{ij}^{\ast}\cap Z_{kl}^{\ast}$ intersects $\Sigma_{k}$
and $\Sigma_{l}$ in a finite set. As each of $\Sigma_{i}$ and $\Sigma_{j}$ is
almost contained in one of $Z_{kl}$ or $Z_{kl}^{\ast}$, it follows that one of
these two corners intersects at most one of $\Sigma_{i}$ and $\Sigma_{j}$ in
an infinite set. Without loss of generality, we can suppose that this corner
is $Z_{ij}^{\ast}\cap Z_{kl}$, which we denote by $V$. Let $W$ denote
$Y_{ij}^{\ast}\cap Y_{kl}$. Thus $W$ is an almost invariant subset of
$H\backslash G$ and the corresponding element $[W]$ of $H_{f}^{1}%
(M_{H};\mathbb{Z}_{2})$ is represented by a $1$--cocycle equal to the
coboundary $\delta V$. Our choice of $V$ means that at most one component of
$\partial M_{H}$ is not almost contained in $V^{\ast}$. This implies that the
element $\overline{W}$ of $H_{n}(\partial M_{H};\mathbb{Z}_{2})$ is supported
on at most one component of $\partial M_{H}$. As $H$ is orientable, part 2) of
Proposition \ref{aisetimpliesannulus} implies that $W$ must be trivial. It
follows that $Y_{ij}$ and $Y_{kl}$ do not cross, which completes the proof of
the claim.

\medskip

Now suppose that the lemma is false. Then the above claim shows that, for any
four distinct integers $i$, $j$, $k$ and $l$, the boundary components
$\Sigma_{k}$ and $\Sigma_{l}$ must both be almost contained in $Z_{ij}$ or
both in $Z_{ij}^{\ast}$. By replacing each of $Z_{12}$ and $Z_{23}$ by its
complement if needed, we can arrange that $Z_{12}$ meets each of $\Sigma_{3}$
and $\Sigma_{4}$ in a finite set and that $Z_{23}$ meets each of $\Sigma_{1}$
and $\Sigma_{4}$ in a finite set. Note that we cannot have $Y_{12}\leq Y_{23}$
as $Z_{12}$ meets $\Sigma_{1}$ in an infinite set, and $Z_{23}$ meets
$\Sigma_{1}$ in a finite set. And we cannot have $Y_{23}\leq Y_{12}$ as
$Z_{23}$ meets $\Sigma_{3}$ in an infinite set and $Z_{12}$ meets $\Sigma_{3}$
in a finite set. It follows that $Y_{12}\cup Y_{23}$ is equivalent to $Y_{13}$
or its complement. Thus $Z_{12}\cup Z_{23}$ is almost equal to $Z_{13}$ or
$Z_{13}^{\ast}$. Now $Z_{12}\cup Z_{23}$ meets $\Sigma_{4}$ in a finite set,
and has infinite intersection with $\Sigma_{2}$. Hence it is not the case that
$\Sigma_{2}$ and $\Sigma_{4}$ are both almost contained in $Z_{13}$ or in
$Z_{13}^{\ast}$, so that $Z_{13}$ separates $\Sigma_{2}$ and $\Sigma_{4}$. Now
the above claim shows that $Y_{13}$ and $Y_{24}$ must cross. This contradicts
our supposition, which completes the proof of the lemma.
\end{proof}

Let $(G,\partial G)$ be an orientable $PD(n+2)$ pair and, as usual, let $M$ be
an aspherical space with fundamental group $G$ and with aspherical subspaces
corresponding to $\partial G$ whose union is denoted $\partial M$. Let $DM$
denote the space obtained by doubling $M$ along $\partial M$, and let $DG$
denote the fundamental group of $DM$. If $M$ is a $3$--manifold, an annulus in
$M$ can be doubled to yield a torus in $DM$. If the annulus is inessential, it
can be homotoped to have image a loop in $\partial M$. Thus the torus can be
homotoped to have the same image and so is not $\pi_{1}$--injective. However,
if the annulus is essential, then the torus will be $\pi_{1}$--injective. The
same construction works in the more general setting of this section. Doubling
a topological annulus $f:(A,\partial A)\rightarrow(M,\partial M)$ yields a map
$Df:DA\rightarrow DM$, where $DA$ denotes the double of $A$ along its
boundary. Thus $Df$ is a map of a torus into $DM$. Again, if $f$ is
inessential, $Df$ will not be $\pi_{1}$--injective, but if $f$ is essential,
then $Df$ will be $\pi_{1}$--injective. To see this, first suppose that $A$ is
untwisted, so that $(A,\partial A)$ is of the form $(C\times I,C\times\partial
I)$, let $H=\pi_{1}(A)$, and consider the lift $f_{H}$ of $f$ to the cover
$M_{H}$ of $M$ such that $\pi_{1}(M_{H})=H$. As $f$ is essential in
$(M,\partial M)$, the images $f_{H}(C\times\{0\})$ and $f_{H}(C\times\{1\})$
must lie in distinct components of $\partial M_{H}$. It follows that a
component of the pre-image of $Df(DA)$ in $(DM)_{H}$ is homeomorphic to
$C\times\mathbb{R}$, so that $Df$ must be $\pi_{1}$--injective. If $A$ is
twisted, we simply apply the above argument to the untwisted double cover of
$A$. Now part 3) of Proposition \ref{aisetimpliesannulus} implies that if $G$
has any nontrivial almost invariant subset over a $VPCn$ subgroup, then there
is an essential torus in $DG$.

In terms of almost invariant sets, the preceding discussion shows that if $Y$
is a nontrivial $H$--almost invariant set dual to an essential annulus in an
orientable $PD(n+2)$ pair $(G,\partial G)$, then there is a natural way to
double $Y$. One obtains a nontrivial almost invariant subset $X$ of $DG$, such
that $X\cap G$ equals $Y$, and $X$ is over the double of $H$, i.e. an
essential torus in $DG$. At first sight, this result sounds somewhat similar
to that in Lemma \ref{propertiesofadapted}, as both results are about
constructing an almost invariant subset $X$ of $DG$ from an almost invariant
subset $Y$ of $G$. However they are completely different as Lemma
\ref{propertiesofadapted} requires that $Y$ be adapted to $\partial G$,
whereas a nontrivial $H$--almost invariant set over a $VPCn$ subgroup of $G$
is never adapted to $\partial G$, by part 1) of Proposition
\ref{annuligenerate}.

\section{The main theorem\label{themaintheorem}}

In the previous section, we discussed the analogues in an orientable $PD(n+2)$
pair of annuli and tori in a $3$--manifold. In this section we will finally
state our main theorem, but first we need to discuss the analogues of the
various types of component of the characteristic submanifold of a
$3$--manifold. Recall that if $(G,\partial G)$ is a Poincar\'{e} duality pair
our aim is to produce a bipartite graph of groups structure for $G$ in which
$V_{0}$--vertices are analogous to components of the characteristic
submanifold. In particular if $G$ is the fundamental group of a Haken
$3$--manifold $M$, this graph of groups structure is dual to the frontier of
the characteristic submanifold of $M$. In our earlier discussion in section
\ref{preliminaries}, we described only two types of such component, namely
$I$--bundles and Seifert fibre spaces. But in order to describe the algebraic
analogues correctly, we will need to subdivide into several cases. There are
special cases when an $I$--bundle has infinite cyclic fundamental group, and
we also need to distinguish between Seifert fibre spaces depending on how they
meet $\partial M$.

We start by considering a component $W$ of the characteristic submanifold of
an orientable Haken $3$--manifold $M$ such that $W$ is an $I$--bundle over a
surface $F$, and $F$ is not an annulus or Moebius band. Thus $\pi_{1}(W)$
equals $\pi_{1}(F)$, so is not finite nor two-ended, and the frontier of $W$
in $M$ consists of the restriction of the $I$--bundle to $\partial F$ and so
consists of essential annuli. In addition, if $F$ is orientable, the
$I$--bundle is trivial and $W$ meets $\partial M$ in two copies of $F$, and if
$F$ is non-orientable, then the $I$--bundle is nontrivial and $W$ meets
$\partial M$ in one copy of the orientable double cover of $F$. Let
$\widetilde{F}$, $\widetilde{W}$ and $\widetilde{M}$ denote the universal
covers of $F$, $W$ and $M$ respectively. Thus $\widetilde{W}$ is homeomorphic
to $\widetilde{F}\times I$. Further each component of the pre-image in
$\widetilde{M}$ of $W$ consists of a copy of $\widetilde{W}$ such that
$\widetilde{W}\cap\partial\widetilde{M}=\widetilde{F}\times\{0,1\}$, and
$\widetilde{F}\times\{0\}$ and $\widetilde{F}\times\{1\}$ lie in distinct
components of $\partial\widetilde{M}$. Thus the induced action of $\pi_{1}(W)$
on $\widetilde{M}$ preserves the union of two distinct components of
$\partial\widetilde{M}$. If some element of $\pi_{1}(W)$ interchanges these
two components, then $W$ is a twisted $I$--bundle. Otherwise, $W$ is
untwisted. This leads to the following definition.

\begin{definition}
\label{defnofI-bundletype} Let $(G,\partial G)$ be an orientable $PD(n+2)$
pair, and let $M$ be an aspherical space with fundamental group $G$ and with
aspherical subspaces corresponding to $\partial G$ whose union is denoted
$\partial M$. Let $\Gamma$ be a minimal graph of groups decomposition of $G$,
and let $v$ be a vertex of $\Gamma$ which is of $VPC(n-1)$--by--Fuchsian type.
(See Definition \ref{defnofVPC-by-Fuchsiantype}. Note that each peripheral
subgroup of $G(v)$ is $VPCn$.)

Then $v$ is \textsl{of }$I$\textsl{--bundle type} if there are two distinct
components $\Sigma$ and $T$ of $\partial\widetilde{M}$ such that

\begin{enumerate}
\item the induced action of $G(v)$ on $\widetilde{M}$ preserves the union of
$\Sigma$ and $T$, and

\item for each peripheral subgroup $K$ of $G(v)$, if $e_{K}$ denotes the edge
of $\Gamma$ which is incident to $v$ and carries $K$, then the edge splitting
associated to $e_{K}$ is given by the essential annulus $K_{\Sigma,T}$. (See
the discussion just before Lemma \ref{swarup}.)
\end{enumerate}
\end{definition}

Next we consider a component $W$ of the characteristic submanifold of $M$ such
that $W$ is a Seifert fibre space, and the orbifold fundamental group of the
base orbifold of $W$ is not finite nor two-ended. Thus $\pi_{1}(W)$ is
$VPC1$--by--Fuchsian, the frontier of $W$ in $M$ consists of boundary torus
components or of vertical annuli in its boundary, and $W$ meets $\partial M$
in boundary torus components or vertical annuli in its boundary. We
distinguish three types of such components $W$.

We will say that $W$ is an \textit{interior} component if it lies in the
interior of $M$. Otherwise $W$ is a \textit{peripheral} component. Our first
definition is the algebraic analogue of an interior component of the
characteristic submanifold of a $3$--manifold.

\begin{definition}
\label{defnofinteriorSeiferttype}Let $(G,${$\partial G$}$)$ be an orientable
$PD(n+2)$ pair, and let $\Gamma$ be a minimal graph of groups decomposition of
$G$. Let $v$ be a vertex of $\Gamma$ which is of $VPCn$--by--Fuchsian type.
(See Definition \ref{defnofVPC-by-Fuchsiantype}.)

Then $v$ is of \textsl{interior Seifert type} if each edge of $\Gamma$ which
is incident to $v$ determines a splitting of $G$ over an essential torus.
\end{definition}

Recall that Lemma \ref{VPC-by-Fuchsiangrouphasuniquefibre} shows that the
$VPCn$ normal subgroup $L$ of $G(v)$ with Fuchsian quotient is unique. Note
that as $G$ is torsion free, so is $L$, so that $L$ is $PDn$. The following
little result applied to any edge group $K$ of $v$ tells us that if $v$ is of
interior Seifert type, then $L$ must be orientable. Note that $L$ is normal in
$K$ with quotient which must be isomorphic to $\mathbb{Z}$ or to
$\mathbb{Z}_{2}\ast\mathbb{Z}_{2}$, as these are the only possible peripheral
subgroups of a finitely generated Fuchsian group. Of course this question did
not arise in the case of a $3$--manifold as then $L$ was the fundamental group
of a closed $1$--manifold, and the only such manifold is orientable.

\begin{lemma}
\label{Lisorientable} Let $K$ be an orientable $PD(n+1)$ group, and let $L$ be
a $VPCn$ normal subgroup of $K$ with quotient isomorphic to $\mathbb{Z}$ or to
$\mathbb{Z}_{2}\ast\mathbb{Z}_{2}$. Then $L$ is an orientable $PDn$ group.
\end{lemma}

\begin{proof}
As $K\ $is torsion free, so is $L$, and hence $L$ is $PDn$. As $\mathbb{Z}%
_{2}\ast\mathbb{Z}_{2}$ has an infinite cyclic subgroup of index $2$, there is
a subgroup $K_{0}$ of $K$, of index at most $2$, such that $K_{0}$ contains
$L$, and $L$ is normal in $K_{0}$ with infinite cyclic quotient. As $K$ is
orientable, so is $K_{0}$. Now Theorem 7.3 of \cite{B-E} shows that $L$ is orientable.
\end{proof}

We will say that a Seifert fibre space component $W$ of the characteristic
submanifold of $M$ is \textit{adapted to} $\partial M$ if there are no annuli
in its frontier. In this case each boundary torus of $W$ is either a component
of $\partial M$ or lies in the interior of $M$ and so is a component of the
frontier of $W$ in $M$. Our next definition is the algebraic analogue of such
a component of the characteristic submanifold of a $3$--manifold.

\begin{definition}
\label{defnofSeiferttypeadaptedtodG}Let $(G,${$\partial G$}$)$ be an
orientable $PD(n+2)$ pair, and let $\Gamma$ be a minimal graph of groups
decomposition of $G$. Let $v$ be a vertex of $\Gamma$ such that $G(v)$ is a
$VPCn$--by--Fuchsian group, where the Fuchsian group is not finite nor two-ended.

Then $v$ is of \textsl{Seifert type adapted to} $\partial G$ if the following
conditions hold:

\begin{enumerate}
\item If $K$ is a peripheral subgroup of $G(v)$, then either $K$ is a
conjugate of a group in $\partial G$, or $K$ is carried by an edge of $\Gamma$
which is incident to $v$.

\item For each peripheral subgroup $K$ of $G(v)$, there is at most one edge
which is incident to $v$ and carries $K$.

\item Each edge of $\Gamma$ which is incident to $v$ carries a peripheral
subgroup of $G(v)$ and determines a splitting of $G$ over an essential torus
in $G$.
\end{enumerate}
\end{definition}

\begin{remark}
If $v$ is of Seifert type adapted to $\partial G$, the two possibilities in 1)
for each peripheral subgroup $K$ of $G(v)$ are mutually exclusive. For an
essential torus in $G${ }cannot be a conjugate of a group in $\partial G$.

Note that if $v$ is of interior Seifert type, it is automatically of Seifert
type adapted to $\partial G$.
\end{remark}

If $W$ is not adapted to $\partial M$, we can push into the interior of $W$
those components of $\partial W$ which meet $\partial M$ in annuli to obtain a
Seifert fibre space $W^{\prime}$ which is homeomorphic to $W$ and adapted to
$\partial M$. Note that each component of the closure of $W-W^{\prime}$ is
homeomorphic to $T\times I$. Recall that the annuli in which $W$ meets
$\partial M$ must be vertical in $W$. Our next definition is the algebraic
analogue of such a component of the characteristic submanifold of a
$3$--manifold $M$.

\begin{definition}
\label{defnofSeiferttype} Let $(G,${$\partial G$}$)$ be an orientable
$PD(n+2)$ pair, and let $\Gamma$ be a minimal graph of groups decomposition of
$G$. Let $v$ be a vertex of $\Gamma$ such that $G(v)$ is a $VPCn$%
--by--Fuchsian group, where the Fuchsian group is not finite nor two-ended.
Let $L$ denote the $VPCn$ normal subgroup of $G(v)$ with Fuchsian quotient.

Then $v$ is of \textsl{Seifert type} if $\Gamma$ can be refined by splitting
at $v$ to a graph of groups structure $\Gamma^{\prime}$ of $G$ with the
following properties:

\begin{enumerate}
\item There is a vertex $v^{\prime}$ of $\Gamma^{\prime}$ with $G(v^{\prime
})=G(v)$ such that $v^{\prime}$ is of Seifert type adapted to $\partial G$.
Thus each edge of $\Gamma^{\prime}$ which is incident to $v^{\prime}$ carries
a peripheral subgroup of $G(v^{\prime})$ and determines a splitting of $G$
over an essential torus in $G$.

\item The projection map $\Gamma^{\prime}\rightarrow\Gamma$ sends $v^{\prime}$
to $v$ and is an isomorphism apart from the fact that certain edges incident
to $v^{\prime}$ are collapsed to $v$.

\item Let $e$ denote an edge of $\Gamma^{\prime}$ which is incident to
$v^{\prime}$ and collapsed to $v$. Thus $G(e)$ is a peripheral subgroup $K$ of
$G(v^{\prime})$. Let $w$ denote the other vertex of $e$. Then $G(w)=K$, and
there is at least one other edge incident to $w$. Further for each such edge
the associated edge splitting is dual to an essential annulus, and the
boundary of each such annulus carries $L$.
\end{enumerate}
\end{definition}

\begin{remark}
When comparing this definition with the topological situation, think of
$\Gamma$ as being dual to the frontier of $W$, and the refinement
$\Gamma^{\prime}$ as being dual to the union of the frontiers of $W$ and of
$W^{\prime}$. Note that if $v$ is of Seifert type adapted to $\partial G$,
then $v$ is trivially of Seifert type. One simply takes $\Gamma^{\prime}$
equal to $\Gamma$ in the above definition.

Part 3) of the definition corresponds to the facts that each component of the
closure of $W-W^{\prime}$ is homeomorphic to $T\times I$, and that the
frontier annuli of $W$ must be vertical in $W$. The reason for the formulation
involving the boundary of each annulus is that some of the annuli involved may
be twisted, a phenomenon with no analogue in $3$--manifold theory.
\end{remark}

There are some special cases which are not covered by the above definitions.
These occur when the Fuchsian quotient group of a vertex group is finite or
two-ended. In the case of the characteristic submanifold of an orientable
Haken $3$--manifold, such a vertex corresponds to a component which is
homeomorphic to one of $S^{1}\times D^{2}$, $T\times I$, or a twisted
$I$--bundle over the Klein bottle.

\begin{definition}
\label{defnofsolidtorustype}Let $(G,${$\partial G$}$)$ be an orientable
$PD(n+2)$ pair, and let $\Gamma$ be a minimal graph of groups decomposition of
$G$. Let $v$ be a vertex of $\Gamma$ such that $G(v)$ is a $VPCn$ group.

Then $v$ is of \textsl{solid torus type} if it is not isolated, and for each
edge of $\Gamma$ which is incident to $v$ the associated edge splitting is
dual to an essential annulus, and there is a $VPCn$ subgroup $H$ of finite
index in $G(v)$ such that the boundary of each such annulus carries $H$.

A vertex of solid torus type is of \textsl{special solid torus type} if either
$v$ has valence $3$ and $H=G(v)$, or if $v$ has valence $1$ and $H$ has index
$2$ or $3$ in $G(v)$.
\end{definition}

\begin{remark}
As mentioned after the previous definition, the annuli involved in this
definition may be twisted. However, in the special case when $v$ is of special
solid torus type, the conditions imply that each of the annuli must be
untwisted. This is proved during the proof of Theorem \ref{mainresult} in
section \ref{proofofmaintheorem}.
\end{remark}

\begin{definition}
\label{defnofspecialSeiferttype}Let $(G,${$\partial G$}$)$ be an orientable
$PD(n+2)$ pair, and let $\Gamma$ be a minimal graph of groups decomposition of
$G$.

A vertex $v$ of $\Gamma$ is of \textsl{special Seifert type} if $v$ has only
one incident edge $e$, the splitting of $G$ associated to $e$ is dual to an
essential torus, and $G(e)$ is a subgroup of index $2$ in $G(v)$.
\end{definition}

\begin{remark}
As $G(v)$ is a finite, torsion free extension of $G(e)$, it follows that
$G(v)$ is also $PD(n+1)$. Also part 3) of Corollary
\ref{PDgroupsplitsoverPDHimpliesHismaximal} implies that $G(v)$ must be non-orientable.
\end{remark}

Our last definition is the algebraic analogue of a component $W$ of the
characteristic submanifold of a $3$--manifold which is homeomorphic to
$T\times I$ or to $K\widetilde{\times}I$, but has some annuli in its frontier.
There are subcases here, depending on whether or not $W$ has a component of
its frontier which is a torus, and whether or not $W$ contains a torus
component of $\partial M$. It is not possible to have both. Note that in all
cases, there is a Seifert fibration of $W$ for which all the annuli in its
frontier are vertical. As $W$ is $T\times I$ or $K\widetilde{\times}I$, this
is equivalent to the condition that all the frontier annuli carry the same
subgroup of $\pi_{1}(W)$ and that $\pi_{1}(W)$ splits over this subgroup. This
is what we generalise in the definition below. Recall that as $G$ is torsion
free, a $VPC(n+1)$ subgroup is automatically $PD(n+1)$.

\begin{definition}
\label{defnoftorustype} Let $(G,${$\partial G$}$)$ be an orientable $PD(n+2)$
pair, and let $\Gamma$ be a minimal graph of groups decomposition of $G$. A
vertex $v$ of $\Gamma$ is of \textsl{torus type} if $G(v)$ is $VPC(n+1)$ and
one of the following cases hold:

\begin{enumerate}
\item $G(v)$ is orientable and is one of the groups in {$\partial G$, }and for
each edge of $\Gamma$ which is incident to $v$ the associated edge splitting
is dual to an essential annulus. Further there is a $VPCn$ subgroup $H$ of
$G(v)$ such that the boundary of each such annulus carries $H$, and $G(v)$
splits over $H$.

\item $G(v)$ is orientable, one of the edges of $\Gamma$ incident to $v$
carries $G(v)$, and the associated edge splitting is dual to an essential
torus in $G$. For each of the remaining edges of $\Gamma$ incident to $v$, the
associated edge splitting is dual to an essential annulus. Further there is a
$VPCn$ subgroup $H$ of $G(v)$ such that the boundary of each such annulus
carries $H$, and $G(v)$ splits over $H$.

\item $G(v)$ is orientable, and for each edge of $\Gamma$ which is incident to
$v$ the associated edge splitting is dual to an essential annulus. Further
there is a $VPCn$ subgroup $H$ of $G(v)$ such that the boundary of each such
annulus carries $H$, and $G(v)$ splits over $H$. In addition $\Gamma$ can be
refined by splitting at $v$ to a graph of groups structure $\Gamma^{\prime}$
of $G$ such that the projection map $\Gamma^{\prime}\rightarrow\Gamma$ sends
an edge $e$ to $v$ and otherwise induces a bijection of edges and vertices.
The group $G(e)$ associated to $e$ is equal to $G(v)$, and the edge splitting
associated to $e$ is dual to an essential torus in $G$.

\item $G(v)$ is non-orientable, and we denote the orientable subgroup of index
$2$ by $G(v)_{0}$. For each edge of $\Gamma$ which is incident to $v$ the
associated edge splitting is dual to an essential annulus. Further there is a
$VPCn$ subgroup $H$ of $G(v)_{0}$ such that the boundary of each such annulus
carries $H$, and $G(v)$ splits over $H$. In addition $\Gamma$ can be refined
by splitting at $v$ to a graph of groups structure $\Gamma^{\prime}$ of $G$
such that the projection map $\Gamma^{\prime}\rightarrow\Gamma$ sends an edge
$e$ to $v$ and otherwise induces a bijection of edges and vertices. The group
associated to $e$ is equal to $G(v)_{0}$, and the associated edge splitting is
dual to an essential torus in $G$. Finally one vertex of $e$ has valence $1$,
and associated group $G(v)$.
\end{enumerate}
\end{definition}

\begin{remark}
In part 3) of this definition, the comparable $3$--manifold situation occurs
when $W$ is $T\times I$, and the frontier of $W$ consists of annuli. Think of
$\Gamma$ as being dual to the frontier of $W$, and the refinement
$\Gamma^{\prime}$ as being dual to the union of the frontier of $W$ with the
torus $T\times\{\frac{1}{2}\}$.

In part 4) of this definition, the comparable $3$--manifold situation occurs
when $W$ is $K\widetilde{\times}I$, and the frontier of $W$ consists of
annuli. Think of $\Gamma$ as being dual to the frontier of $W$, and the
refinement $\Gamma^{\prime}$ as being dual to the union of the frontier of $W$
with a torus in the interior of $W$ which is parallel to $\partial W$.

In all cases, the reason for the formulation involving the boundary of each
annulus is that some of the annuli involved may be twisted, a phenomenon with
no analogue in $3$--manifold theory.

Note that Lemma \ref{splittingsofVPCgroups} tells us that if a $VPC(n+1)$
group $G$ splits over a subgroup $H$, then $H$ is $VPCn$, and is normal in $G$
with quotient which is isomorphic to $\mathbb{Z}$ or $\mathbb{Z}_{2}%
\ast\mathbb{Z}_{2}$.
\end{remark}

Now we are ready to state the main result of this paper. Lemma
\ref{PD(n+2)groupshavenoa.i.setsoverVPCk<n} implies that if $(G,\partial G)$
is an orientable $PD(n+2)$ pair, then the decomposition $\Gamma_{n,n+1}(G)$ of
Theorem \ref{JSJforVPCoftwolengths} exists. That theorem tells us that for any
group $G$, the $V_{0}$--vertices of $\Gamma_{n,n+1}(G)$ are of four types,
namely they are isolated, of $VPCk$--by--Fuchsian type, where $k$ is $n-1$ or
$n$, or of commensuriser type. Our main result is the following theorem which
asserts that $\Gamma_{n,n+1}(G)$ and its completion $\Gamma_{n,n+1}^{c}(G)$
(see Definition \ref{defnofcompletion}) have properties analogous to the
topological picture in dimension $3$. If $n=1$, and $G$ is the fundamental
group of an orientable Haken $3$--manifold $M$, then $\Gamma_{1,2}(G)$ is dual
to the frontier of $AT(M)$, and $\Gamma_{1,2}^{c}(G)$ is dual to the frontier
of $JSJ(M)$, where $AT(M)$ and $JSJ(M)$ are the submanifolds of $M$ discussed
at the start of section \ref{preliminaries}.

\begin{theorem}
\label{mainresult} (Main Result) Let $n\geq1$, and let $(G,\partial G)$ be an
orientable $PD(n+2)$ pair such that $G$ is not $VPC$. Let $\mathcal{F}%
_{n,n+1}$ denote the family of equivalence classes of all nontrivial almost
invariant subsets of $G$ which are over a $VPCn$ subgroup, together with the
equivalence classes of all $n$--canonical almost invariant subsets of $G$
which are over a $VPC(n+1)$ subgroup. Finally let $\Gamma_{n,n+1}$ denote the
reduced algebraic regular neighbourhood of $\mathcal{F}_{n,n+1}$ in $G$, and
let $\Gamma_{n,n+1}^{c}$ denote the completion of $\Gamma_{n,n+1}$. Thus
$\Gamma_{n,n+1}$ and $\Gamma_{n,n+1}^{c}$ are bipartite graphs of groups
structures for $G$, with vertices of $V_{0}$--type and of $V_{1}$--type.

Then $\Gamma_{n,n+1}$ and $\Gamma_{n,n+1}^{c}$ have the following properties:

\begin{enumerate}
\item Each $V_{0}$--vertex $v$ of $\Gamma_{n,n+1}$ satisfies one of the
following conditions:

\begin{enumerate}
\item $v$ is isolated, and $G(v)$ is $VPC$ of length $n$ or $n+1$, and the
edge splittings associated to the two edges incident to $v$ are dual to
essential annuli or tori in $G$.

\item $v$ is of $VPC(n-1)$--by--Fuchsian type, and is of $I$--bundle type.
(See Definition \ref{defnofI-bundletype}.)

\item $v$ is of $VPCn$--by--Fuchsian type, and is of interior Seifert type.
(See Definition \ref{defnofinteriorSeiferttype}.)

\item $v$ is of commensuriser type. Further $v$ is of Seifert type (see
Definition \ref{defnofSeiferttype}), or of torus type (see Definition
\ref{defnoftorustype}) or of solid torus type (see Definition
\ref{defnofsolidtorustype}).
\end{enumerate}

\item The $V_{0}$--vertices of $\Gamma_{n,n+1}^{c}$ obtained by the completion
process are of special Seifert type (see Definition
\ref{defnofspecialSeiferttype}) or of special solid torus type (see Definition
\ref{defnofsolidtorustype}).

\item Each edge splitting of $\Gamma_{n,n+1}$ and of $\Gamma_{n,n+1}^{c}$ is
dual to an essential annulus or torus in $G$.

\item Any nontrivial almost invariant subset of $G$ over a $VPC(n+1)$ group
and adapted to $\partial G$ is enclosed by some $V_{0}$--vertex of
$\Gamma_{n,n+1}$, and also by some $V_{0}$--vertex of $\Gamma_{n,n+1}^{c}$.

\item If $H$ is a $VPC(n+1)$ subgroup of $G$ which is not conjugate into
$\partial G$, then $H$ is conjugate into a $V_{0}$--vertex group of
$\Gamma_{n,n+1}^{c}$.
\end{enumerate}
\end{theorem}

\begin{remark}
Part 3) follows immediately from parts 1) and 2), as the definitions of the
various types of $V_{0}$--vertex in the statements of parts 1) and 2) all
contain the requirement that the edge splittings be dual to an essential
annulus or torus.

Part 4) does not follow from the properties of an algebraic regular
neighbourhood as an almost invariant subset of $G$ over a $VPC(n+1)$ group
which is adapted to $\partial G$ need not be $n$--canonical, and so need not
lie in the family $\mathcal{F}_{n,n+1}$. Note that, from \cite{SS4}, we know
that there may be almost invariant subsets of $G$ over $VPC(n+1)$ subgroups
which are not adapted to $\partial G$.

Part 5) also does not follow from the properties of an algebraic regular
neighbourhood as a $VPC(n+1)$ subgroup $H$ of $G$ may be non-orientable.
\end{remark}

\section{Torus Decompositions for $PD(n+2)$ groups and
pairs\label{torusdecompofPDgroup}}

Before embarking on the proof of Theorem \ref{mainresult}, we will consider a
simpler graph of groups structure analogous to the torus decomposition $T(M)$
of a $3$--manifold $M$, discussed in section \ref{preliminaries}.

We will start this section by considering an orientable $PD(n+2)$ group $G$,
where $n\geq1$, and the graph of groups structure $\Gamma_{n+1}(G)=\Gamma
(\mathcal{F}_{n+1}:G)$ of Theorem \ref{JSJforVPCofgivenlength}$.$ Recall from
Lemma \ref{PD(n+2)groupshavenoa.i.setsoverVPCk<n} that such a group cannot
admit any nontrivial almost invariant subset over a $VPC(\leq n)$ subgroup, so
that $\Gamma_{n+1}(G)$ does exist. It is a reduced algebraic regular
neighbourhood of $\mathcal{F}_{n+1}$ in $G$, where $\mathcal{F}_{n+1}$ is the
collection of equivalence classes of all nontrivial almost invariant subsets
of $G$ which are over $VPC(n+1)$ subgroups. Since $G$ is orientable, the
nontrivial almost invariant sets in $\mathcal{F}_{n+1}$ are automatically over
orientable $VPC(n+1)$ groups. Recall that if $H$ is an orientable $VPC(n+1)$
subgroup of $G$, we call $H$ an essential torus in $G$, and $G$ possesses a
unique nontrivial $H$--almost invariant subset, up to equivalence and
complementation. Thus $\mathcal{F}_{n+1}$ can be thought of as the collection
of all essential tori in $G$. It will be convenient to denote $\Gamma
_{n+1}(G)$ by $T_{n+1}(G)$, and we will call $T_{n+1}(G)$ the \textit{torus
decomposition} of $G$. We will also use $T_{n+1}^{c}(G)$ to denote the
completion $\Gamma_{n+1}^{c}(G)$ of $\Gamma_{n+1}(G)$. See Definition
\ref{defnofcompletion}.

In chapter 12 of \cite{SS2}, we briefly discussed the connection between
$T_{n+1}(G)$ and Kropholler's decomposition in \cite{K} in the case of
orientable $PD(n+2)$ groups. Using the present notation, it follows from
Theorem \ref{torusdecompofPDgroupproperties} below that $T_{n+1}^{c}(G)$ is
the same as Kropholler's decomposition. We will show shortly that the
analogous statement holds for the case of a $PD(n+2)$ pair. In dimension $3$,
it is clear that $T_{2}(G)$ is the same as Castel's decomposition in
\cite{Castel} as his decomposition is defined to be $\Gamma_{2}(G)$. As
discussed in section \ref{preliminaries}, if $M\ $is a closed orientable Haken
$3$--manifold and $G$ equals $\pi_{1}(M)$, then the topological torus
decomposition $T(M)$ determines the decomposition $T_{2}^{c}(G)$ of $G$. As in
the statement of Theorem \ref{JSJforVPCofgivenlength}, it will be convenient
to state Theorem \ref{torusdecompofPDgroupproperties} excluding the case when
$G$ is $VPC$. For if $G$ is a $VPC(n+2)$ group, then $\Gamma_{n+1}(G)$
consists of a single $V_{0}$--vertex.

\begin{theorem}
\label{torusdecompofPDgroupproperties} Let $n\geq1$. If $G$ is an orientable
$PD(n+2)$ group which is not $VPC$, then $T_{n+1}(G)$ and $T_{n+1}^{c}(G)$
have the following properties:

\begin{enumerate}
\item Each $V_{0}$--vertex $v$ of $T_{n+1}(G)$ satisfies one of the following conditions:

\begin{enumerate}
\item $v$ is isolated, and $G(v)$ is a torus in $G$.

\item $v$ is of interior Seifert type. (See Definition
\ref{defnofinteriorSeiferttype}.)
\end{enumerate}

\item The $V_{0}$--vertices of $T_{n+1}^{c}(G)$ obtained by the completion
process are of special Seifert type. (See Definition
\ref{defnofspecialSeiferttype}.)

\item Each edge splitting of $T_{n+1}(G)$ and of $T_{n+1}^{c}(G)$ is dual to
an essential torus in $G$.

\item If $H$ is a $VPC(n+1)$ subgroup of $G$, then $H$ is conjugate into a
$V_{0}$--vertex group of $T_{n+1}^{c}(G)$.
\end{enumerate}
\end{theorem}

\begin{remark}
If $G$ has no nontrivial almost invariant subsets which are over $VPC(n+1)$
subgroups, then $\mathcal{F}_{n+1}$ is empty and $T_{n+1}(G)$, and hence
$T_{n+1}^{c}(G)$, consists of a single $V_{1}$--vertex.
\end{remark}

\begin{proof}
1) Theorem \ref{JSJforVPCofgivenlength} tells us that $T_{n+1}(G)$ is a
minimal, reduced bipartite, graph of groups decomposition of $G$, and that
each $V_{0}$--vertex of $T_{n+1}(G)$ is isolated, of $VPCn$--by--Fuchsian
type, or of commensuriser type. If $v$ is a $V_{0}$--vertex of commensuriser
type, then $v$ encloses an element $X$ of $\mathcal{F}_{n+1}$ which is over
some $VPC(n+1)$ group $H$, and $G(v)$ is of the form $Comm_{G}(H)$. Further
$X$ crosses weakly some of its translates by $Comm_{G}(H)$. In the present
situation, $G$ is $PD(n+2)$ and $H$ is $VPC(n+1)$, so that the number of
coends of $H$ in $G$ is $2$. Now Proposition 7.4 of \cite{SS2} implies that no
almost invariant set can cross $X$ or any of its translates weakly, so that
$V_{0}$--vertices of commensuriser type cannot occur. It follows that each
$V_{0}$--vertex of $T_{n+1}(G)$ is isolated or of $VPCn$--by--Fuchsian type.
As each edge splitting is over a $VPC(n+1)$ subgroup of $G$, it follows that
each edge splitting of $T_{n+1}(G)$ is dual to an essential torus in $G$.
Hence any $V_{0}$--vertex of $T_{n+1}(G)$ which is of $VPCn$--by--Fuchsian
type must be of interior Seifert type. This completes the proof of part 1) of
the theorem.

2) The construction of $T_{n+1}^{c}(G)$ from $T_{n+1}(G)$ described in section
\ref{preliminaries} can only introduce $V_{0}$--vertices of special Seifert
type, so part 2) of the theorem holds.

3) Part 3) of the theorem follows immediately from parts 1)\ and 2).

4) First note that as $G$ is torsion free, so is $H$. Thus $H$ must be
$PD(n+1)$.

If $H$ is orientable, the pair $(G,H)$ has two ends, so there is a nontrivial
$H$--almost invariant subset $X$ of $G$. As $X$ is enclosed by some $V_{0}%
$--vertex $v$ of $T_{n+1}(G)$, it follows that $H$ is conjugate into $G(v)$.
Hence $H$ is also conjugate into a $V_{0}$--vertex group of $T_{n+1}^{c}(G)$,
as required.

If $H$ is non-orientable, let $H_{0}$ denote its orientable subgroup of index
$2$. As the pair $(G,H_{0})$ has two ends, there is a nontrivial $H_{0}%
$--almost invariant subset of $G$. As $T_{n+1}(G)$ has no $V_{0}$--vertices of
commensuriser type, it follows from Theorem \ref{JSJforVPCofgivenlength} and
Remark \ref{remarksonJSJ} that $H_{0}$ must have small commensuriser, which we
denote by $K$. This means that $K$ contains $H_{0}$ with finite index, so that
$K$ is itself $VPC(n+1)$ and $PD(n+1)$. Note that $K$ must contain $H$, so
that $K$ is non-orientable. We let $K_{0}$ denote its orientable subgroup of
index $2$. Note that $K_{0}$ is a maximal torus subgroup of $G$. The preceding
paragraph shows that there is a $V_{0}$--vertex $v$ of $T_{n+1}(G)$ such that
$K_{0}$ is conjugate into $G(v)$. As $K$ contains $K_{0}$ with finite index,
it follows that there is a vertex $w$ of $T_{n+1}(G)$ such that $K$ is
conjugate into $G(w)$. If $w$ is a $V_{0}$--vertex, then $K$,$\ $and hence
$H$, is conjugate into a $V_{0}$--vertex group of $T_{n+1}^{c}(G)$, as
required. So we now consider the case when $w$ is a $V_{1}$--vertex. In
particular, $v$ and $w$ must be distinct. Thus there is an edge $e$ of
$T_{n+1}(G)$ which is incident to $w$ such that, after a conjugation, $G(e)$
contains $K_{0}$. As all the edge groups of $T_{n+1}(G)$ are torus groups, the
group $G(e)$ must equal $K_{0}$. If $E(w)$ denotes the family of subgroups of
$G(w)$ which are edge groups for the edges incident to $w$, then Theorem 8.1
of \cite{B-E} tells us that the pair $(G(w),E(w))$ is $PD(n+2)$. The
commensuriser in $G(w)$ of $G(e)=K_{0}$ contains a conjugate of $K$ and so is
not equal to $K_{0}$. Thus Lemma \ref{boundarygroupismaximal} shows that
$K_{0}$ is the only element of the family $E(w)$, and $G(w)$ contains $K_{0}$
with index $2$. Thus $w$ has valence $1$, and $G(w)$ equals a conjugate of
$K$. Now it follows that $w$ becomes a $V_{0}$--vertex in the completion
$T_{n+1}^{c}(G)$, so that $K$, and hence $H$, is conjugate into a $V_{0}%
$--vertex group of $T_{n+1}^{c}(G)$, as required.
\end{proof}

Next we discuss the torus decomposition of an orientable $PD(n+2)$ pair
$(G,\partial G)$ with non-empty boundary. Recall from section
\ref{essentialannuli} the discussion of an essential torus in $G$. In
particular, we let $DG$ denote the orientable $PD(n+2)$ group obtained by
doubling $G$ along its boundary. Then given an orientable $PD(n+1)$ subgroup
$H$ of $G$, there is a $H$--almost invariant subset $X_{H}$ of $DG$ associated
to $H$, and the intersection $Y_{H}=X_{H}\cap G$ is a $H$--almost invariant
subset of $G$ which is nontrivial unless $H$ is conjugate into $\partial G$.
Let $\mathcal{T}_{n+1}$ denote the family of equivalence classes of all such
nontrivial subsets $Y_{H}$ of $G$, where $H$ is $VPC(n+1)$. Then the torus
decomposition $T_{n+1}(G,\partial G)$ of $(G,\partial G)$ will be the reduced
algebraic regular neighbourhood in $G$ of $\mathcal{T}_{n+1}$. This is the
natural definition, but it is not obvious that this algebraic regular
neighbourhood exists. One immediate problem is that $G$ may have nontrivial
almost invariant subsets over $VPCn$ subgroups, so this decomposition is
different from any of those proved to exist in \cite{SS2}. In order to show
that $T_{n+1}(G,\partial G)$ exists, we will use the fact that $DG$ does not
have nontrivial almost invariant subsets over $VPCn$ subgroups, so we can
apply results from \cite{SS2}.

\begin{theorem}
\label{torusdecompofapairexists}Let $(G,\partial G)$ be an orientable
$PD(n+2)$ pair, such that $\partial G$ is non-empty, and let $\mathcal{T}%
_{n+1}$ denote the family of equivalence classes of almost invariant subsets
$Y_{H}$ of $G$ described above. Then $\mathcal{T}_{n+1}$ has a reduced
algebraic regular neighbourhood $T_{n+1}(G,\partial G)$ in $G$. Further
$T_{n+1}(G,${$\partial G$}$)$ is adapted to $\partial G$.
\end{theorem}

\begin{proof}
If $G$ is $VPC$, part 2) of Corollary \ref{PD(n+2)andVPCgroup} tells us that
either $(G,${$\partial G$}$)$ is the trivial pair $(G,\{G,G\})$, or that
$\partial G$ is a single group $S$, and $G$ contains $S$ with index $2$. In
either case, the pair $(G,\partial G)$ admits no essential tori, so that
$\mathcal{T}_{n+1}$ is empty. Thus $\mathcal{T}_{n+1}$ has a reduced algebraic
regular neighbourhood which consists of a single $V_{1}$--vertex, and this is
trivially adapted to $\partial G$. For the rest of this proof we will assume
that $G$ is not $VPC$.

Recall our discussion at the start of section \ref{essentialannuli}. The
natural graph of groups structure $\Delta$ for the orientable $PD(n+2)$ group
$DG$ has two vertices $w$ and $\overline{w}$ and edges joining them which
correspond to the groups of $\partial G$. Given an orientable $PD(n+1)$
subgroup $H$ of $G$, there is a $H$--almost invariant subset $X_{H}$ of $DG$
associated to $H$, and the intersection $Y_{H}=X_{H}\cap G$ is a $H$--almost
invariant subset of $G$. An important point about $X_{H}$ is that it is
enclosed by the vertex $w$ of $\Delta$, where we identify $G$ with $G(w)$, so
that $Y_{H}$ is adapted to $\partial G$.

Now consider $\mathcal{F}_{n+1}(DG)$, which is the collection of equivalence
classes of all nontrivial almost invariant subsets of $DG$ which are over
$VPC(n+1)$ subgroups. Recall that as $DG$ is $PD(n+2)$, it is torsion free so
that a $VPC(n+1)$ subgroup $K$ of $DG$ must be $PD(n+1)$. Further, as $DG$ is
orientable, if there is a nontrivial almost invariant subset of $DG$ which is
over $K$, then $K$ must be orientable, and all such almost invariant subsets
of $DG$ are equivalent up to complementation. Also recall from Theorem
\ref{torusdecompofPDgroupproperties} that the reduced algebraic regular
neighbourhood $\Gamma(\mathcal{F}_{n+1}(DG):DG)$ exists, and is denoted
$T_{n+1}(DG)$, and its $V_{0}$--vertices are either isolated or of
$VPCn$--by--Fuchsian type. As we will need to use our construction of
unreduced algebraic regular neighbourhoods in \cite{SS2}, as corrected in
\cite{SS2errata}, we note that Theorem \ref{JSJforVPCofgivenlength} tells us
that the unreduced algebraic regular neighbourhood of $\mathcal{F}_{n+1}(DG)$
in $DG$ also exists. We will use the notation $\overline{\Gamma}%
(\mathcal{F}_{n+1}(DG):DG)$ for this unreduced algebraic regular
neighbourhood. As in the case of the reduced algebraic regular neighbourhood,
its $V_{0}$--vertices are either isolated or of $VPCn$--by--Fuchsian type.

Let $\mathcal{E}_{n+1}(DG)$ denote the subfamily of $\mathcal{F}_{n+1}(DG)$
which consists of nontrivial almost invariant subsets of $DG$ which are over
subgroups of $G$ which are non-peripheral in $G$. Note that $\mathcal{E}%
_{n+1}(DG)$ is $G$--invariant but is not $DG$--invariant. We claim that
$\mathcal{E}_{n+1}(DG)$ possesses an unreduced algebraic regular neighbourhood
$\overline{\Gamma}(\mathcal{E}_{n+1}(DG):DG)$ in $DG$, and hence a reduced
algebraic regular neighbourhood $\Gamma(\mathcal{E}_{n+1}(DG):DG)$. We recall
from \cite{SS2} that any finite subset of $\mathcal{F}_{n+1}(DG)$ possesses an
unreduced algebraic regular neighbourhood in $DG$. As $\overline{\Gamma
}(\mathcal{F}_{n+1}(DG):DG)$ has no $V_{0}$--vertices of commensuriser type,
the same holds for an algebraic regular neighbourhood of any subset of
$\mathcal{F}_{n+1}(DG)$. Thus the proof of Theorem
\ref{JSJforVPCofgivenlength} which we gave in chapter 12 of \cite{SS2} shows
that the $V_{0}$--vertices of the unreduced algebraic regular neighbourhood of
a finite subset of $\mathcal{F}_{n+1}(DG)$ must be either isolated or of
$VPCn$--by--Fuchsian type. Now we consider the construction in the proof of
Theorem 9.2 of \cite{SS2}. This constructs $\overline{\Gamma}(\mathcal{E}%
_{n+1}(DG):DG)$ by expressing $\mathcal{E}_{n+1}(DG)$ as an ascending union of
finite subsets $\mathcal{E}_{n+1}^{k}$. We obtain a sequence of bipartite
graphs of groups decompositions $\overline{\Gamma}^{k}$ of $DG$ for each of
which, every $V_{0}$--vertex must be either isolated or of $VPCn$%
--by--Fuchsian type. The accessibility results of section 13 of \cite{GSS}
imply that the sequence must stabilise so that $\overline{\Gamma}%
(\mathcal{E}_{n+1}(DG):DG)$ exists as required. Further each $V_{0}$--vertex
of $\overline{\Gamma}(\mathcal{E}_{n+1}(DG):DG)$ is either isolated or of
$VPCn$--by--Fuchsian type.

Let $T_{DG}$ denote the universal covering $DG$--tree of $\overline{\Gamma
}(\mathcal{E}_{n+1}(DG):DG)$, and recall that the $V_{0}$--vertices of
$T_{DG}$ are the CCC's of all the translates by $DG$ of elements of
$\mathcal{E}_{n+1}(DG)$. Now we consider the construction of the unreduced
algebraic regular neighbourhood of $\mathcal{T}_{n+1}$ described in chapter 6
of \cite{SS2}. One starts by choosing one element of $\mathcal{T}_{n+1}$ from
each equivalence class, and then considers the CCC's of these elements. There
is a natural map $\varphi$ from the equivalence classes of $\mathcal{T}_{n+1}$
to the equivalence classes of $\mathcal{E}_{n+1}(DG)$, given by sending the
class of $Y_{H}$ to the class of $X_{H}$. Recall that $Y_{H}=X_{H}\cap G$.
Clearly if two elements of $\mathcal{T}_{n+1}$ cross, then the corresponding
elements of $\mathcal{E}_{n+1}(DG)$ also cross. Now recall that elements of
$\mathcal{E}_{n+1}(DG)$ which cross must do so strongly. It follows that if
two elements of $\mathcal{E}_{n+1}(DG)$ cross, then the corresponding elements
of $\mathcal{T}_{n+1}$ also cross. Hence $\varphi$ induces a $G$--equivariant
bijection between the collection $P$ of all CCC's of $\mathcal{T}_{n+1}$ and a
subset $Q$ of the $V_{0}$--vertices of $T_{DG}$. It also follows that the
pretree structures on $P$ and $Q$ are the same. In particular, the pretree
structure on the collection $P$ of all CCC's of $\mathcal{T}_{n+1}$\ is
discrete. Now the proof of Theorem 3.8 of \cite{SS2} shows that there is a
bipartite $G$--tree $T_{G}$ whose quotient by $G$ is the unreduced algebraic
regular neighbourhood of $\mathcal{T}_{n+1}$ in $G$. Thus $\mathcal{T}_{n+1}$
also has a reduced algebraic regular neighbourhood in $G$, so the torus
decomposition $T_{n+1}(G,\partial G)$ of $G$ exists.

Recall that each of the edge splittings of $\Delta$ crosses no element of
$\mathcal{E}_{n+1}(DG)$, and so must be enclosed by some $V_{1}$--vertex of
the unreduced algebraic regular neighbourhood $\overline{\Gamma}%
(\mathcal{E}_{n+1}(DG):DG)$. Thus we can refine $\overline{\Gamma}%
(\mathcal{E}_{n+1}(DG):DG)$ by splitting at $V_{1}$--vertices to obtain a
graph of groups decomposition $\overline{\Gamma}\Delta$ of $DG$ which also
refines $\Delta$. Recall that each element of $\mathcal{E}_{n+1}(DG)$ is
enclosed by the vertex $w$ of $\Delta$. Thus if we remove the interiors of the
edges of $\overline{\Gamma}\Delta$ which correspond to the edges of $\Delta$,
we will be left with a connected graph of groups $\overline{\Gamma}_{w}$, and
the single vertex $\overline{w}$. The fundamental group of $\overline{\Gamma
}_{w}$ is $G(w)$, which we continue to identify with $G$.

\medskip

\textbf{Claim:} $\overline{\Gamma}_{w}$\textit{ is isomorphic to the unreduced
algebraic regular neighbourhood of }$\mathcal{T}_{n+1}$\textit{ in }%
$G$\textit{.}

\medskip

This claim immediately implies that the graph of groups $\Gamma_{w}$ obtained
by reducing $\overline{\Gamma}_{w}$ is isomorphic to the torus decomposition
$T_{n+1}(G,\partial G)$. It also implies that each group in $\partial G$ is
conjugate into some vertex group of $T_{n+1}(G,\partial G)$, so that
$T_{n+1}(G,${$\partial G$}$)$ is adapted to {$\partial G$.}

To prove our claim, recall from two paragraphs previously that the map
$\varphi$ induces a $G$--equivariant injection from the $V_{0}$--vertices of
$T_{G}$ to the $V_{0}$--vertices of $T_{DG}$. It follows that $\varphi$
induces a $G$--equivariant injection from the $V_{0}$--vertices of $T_{G}$ to
the $V_{0}$--vertices of $T_{\overline{\Gamma}\Delta}$. Let $T$ denote the
subtree of $T_{\overline{\Gamma}\Delta}$ spanned by the $V_{0}$--vertices in
the image of $\varphi$. As $T_{G}$ is a minimal $G$--tree, $T$ must be the
minimal $G$--invariant subtree of $T_{\overline{\Gamma}\Delta}$. It follows
that $T$ is the universal covering $G$--tree of $\overline{\Gamma}_{w}$,
proving that $\overline{\Gamma}_{w}$ is isomorphic to the unreduced algebraic
regular neighbourhood of $\mathcal{T}_{n+1}$ in $G$, as claimed.{ Note that it
also follows that }$T_{n+1}(G,\partial G)$ is the decomposition of $G$ induced
from $\Gamma(\mathcal{E}_{n+1}(DG):DG)$.
\end{proof}

The properties of $T_{n+1}(G,${$\partial G$}$)$ when $\partial G$ is non-empty
are similar to those in the case when $\partial G$ is empty. The fact that
$T_{n+1}(G,${$\partial G$}$)$ is adapted to {$\partial G$ plays an important
}role. Before listing these properties, it will be convenient to introduce the
completion $T_{n+1}^{c}(G,${$\partial G$}$)$ of $T_{n+1}(G,${$\partial G$}$)$
which is defined in the same way as we defined the completions of
$\Gamma_{n+1}(G)$ and $\Gamma_{n,n+1}(G)$ in Definition \ref{defnofcompletion}.

The following result lists the properties of $T_{n+1}(G,${$\partial G$}$)$ and
its completion $T_{n+1}^{c}(G,${$\partial G$}$)$ when $\partial G$ is
non-empty. As usual, it will be convenient to exclude the case when $G$ is
$VPC$. For in that case $T_{n+1}(G,${$\partial G$}$)$ and $T_{n+1}^{c}%
(G,${$\partial G$}$)$ are equal and consist of a single $V_{1}$--vertex.

\begin{theorem}
\label{torusdecompofapair}Let $n\geq1$, and let $(G,${$\partial G$}$)$ be an
orientable $PD(n+2)$ pair, such that $\partial G$ is non-empty. If $G$ is not
$VPC$, then $T_{n+1}(G,${$\partial G$}$)$ and $T_{n+1}^{c}(G,${$\partial G$%
}$)$ have the following properties:

\begin{enumerate}
\item Each edge splitting of $T_{n+1}(G,${$\partial G$}$)$ and of $T_{n+1}%
^{c}(G,\partial G)$ is dual to an essential torus in $(G,${$\partial G$}$)$.

\item Each $V_{0}$--vertex $v$ of $T_{n+1}(G,\partial G)$ satisfies one of the
following conditions:

\begin{enumerate}
\item $v$ is isolated, and $G(v)$ is an essential torus in $G$.

\item $v$ is of Seifert type adapted to $\partial G$. (See Definition
\ref{defnofSeiferttypeadaptedtodG}. Note that this includes the possibility
that $v$ is of interior Seifert type.)
\end{enumerate}

\item The $V_{0}$--vertices of $T_{n+1}^{c}(G,\partial G)$ obtained by the
completion process are of special Seifert type. (See Definition
\ref{defnofspecialSeiferttype}.)

\item If $H$ is a $VPC(n+1)$ subgroup of $G$ which is not conjugate into
$\partial G$, then $H$ is conjugate into a $V_{0}$--vertex group of
$T_{n+1}^{c}(G,${$\partial G$}$)$.
\end{enumerate}
\end{theorem}

\begin{remark}
It follows from property 4) that $T_{n+1}^{c}(G,${$\partial G)$ is the same as
Kropholler's decomposition in \cite{K}. In dimension }$3$, i{t is also easy to
see that} $T_{2}(G,${$\partial G)$ is the same as Castel's decomposition in
\cite{Castel}. Finally }if $M\ $is an orientable Haken $3$--manifold and
$(G,${$\partial G)$ is the corresponding Poincar\'{e} duality pair}, then the
topological torus decomposition $T(M)$ determines the decomposition $T_{2}%
^{c}(G,\partial G)$ of $G$.
\end{remark}

\begin{proof}
It follows from the proof of Theorem \ref{torusdecompofapairexists} that the
reduced graphs of groups decompositions $\Gamma(\mathcal{E}_{n+1}(DG):DG)$ and
$\Delta$ have a common refinement $\Gamma\Delta$, obtained from $\Gamma
(\mathcal{E}_{n+1}(DG):DG)$ by splitting at $V_{1}$--vertices, which consists
of $T_{n+1}(G,\partial G)$ and a single extra vertex $\overline{w}$ which is
joined to $T_{n+1}(G,\partial G)$ by edges $e_{1},\ldots,e_{m}$ whose
associated splittings are those of $\Delta$. Note that the $V_{0}$--vertices
of $\Gamma(\mathcal{E}_{n+1}(DG):DG)$ and their incident edges are unaffected
by the refinement process. Thus we will refer to the vertices of $\Gamma
\Delta$ which are obtained from $V_{0}$--vertices of $\Gamma(\mathcal{E}%
_{n+1}(DG):DG)$ as $V_{0}$--vertices of $\Gamma\Delta$. As the $V_{0}%
$--vertices of $\Gamma(\mathcal{E}_{n+1}(DG):DG)$ are isolated or of
$VPCn$--by--Fuchsian type, the same holds for the $V_{0}$--vertices of
$\Gamma\Delta$.

1) Each edge of $T_{n+1}(G,\partial G)$, regarded as an edge of $\Gamma\Delta
$, determines a splitting of the orientable $PD(n+2)$ group $DG$ over a
$VPC(n+1)$ subgroup. It follows that this splitting is dual to an essential
torus in $DG$. Hence each edge splitting of $T_{n+1}(G,${$\partial G$}$)$ is
dual to an essential torus in $(G,${$\partial G$}$)$, proving part 1) of the theorem.

2) Let $v$ be an isolated $V_{0}$--vertex of $\Gamma\Delta$. If one of the
$e_{i}$'s is incident to $v$, then $v$ has only one incident edge $e$ in
$T_{n+1}(G,\partial G)$ and the inclusion of $G(e)$ into $G(v)$ is an
isomorphism. But this contradicts the minimality of $T_{n+1}(G,\partial G)$.
Thus no $e_{i}$ is incident to $v$, and $v$ must be an isolated $V_{0}%
$--vertex of $T_{n+1}(G,\partial G)$.

Now let $v$ be a $V_{0}$--vertex of $\Gamma\Delta$ of $VPCn$--by--Fuchsian
type. If no $e_{i}$ is incident to $v$, then $v$ is a $V_{0}$--vertex of
$T_{n+1}(G,\partial G)$ of $VPCn$--by--Fuchsian type. As each edge splitting
of $T_{n+1}(G,${$\partial G$}$)$ is dual to an essential torus in
$(G,${$\partial G$}$)$, it follows that $v$ is of interior Seifert type. If
some $e_{i}$ is incident to $v$, the associated edge group is a group in
$\partial G$. It follows that $v$ is a $V_{0}$--vertex of $T_{n+1}(G,\partial
G)$ which is of Seifert type adapted to $\partial G$.

Thus the $V_{0}$--vertices of $T_{n+1}(G,\partial G)$ are isolated, or of
Seifert type adapted to $\partial G$, which completes the proof of part 2) of
the theorem.

3) The construction of $T_{n+1}^{c}(G,\partial G)$ from $T_{n+1}(G,\partial
G)$ described in section \ref{preliminaries} can only introduce $V_{0}%
$--vertices of special Seifert type, so part 3) of the theorem holds.

4) First note that as $G$ is torsion free, so is $H$. Thus $H$ must be
$PD(n+1)$.

Suppose that $H$ is orientable. The hypothesis that $H$ is not conjugate into
$\partial G$ implies that $H$ is an essential torus in $(G,\partial G)$, so
that there is a nontrivial $H$--almost invariant subset $X$ of $G$ dual to
$H$. As $X$ is enclosed by some $V_{0}$--vertex $v$ of $T_{n+1}(G,\partial
G)$, it follows that $H$ is conjugate into $G(v)$. Hence $H$ is also conjugate
into a $V_{0}$--vertex group of $T_{n+1}^{c}(G)$, as required.

Now suppose that $H$ is non-orientable, and let $H_{0}$ denote its orientable
subgroup of index $2$. If $H_{0}$ is conjugate into a group $S$ in $\partial
G$, this conjugate of $H_{0}$ will be a $PD(n+1)$ subgroup of the $PD(n+1)$
group $S$ and so will be of finite index. Thus $S$ must itself be $VPC(n+1)$
and be conjugate commensurable with $H_{0}$. As $H$ is not conjugate into $S$,
it follows that $Comm_{G}(S)\neq S$, so Lemma \ref{boundarygroupismaximal}
implies that $G$ contains $S$ with index $2$. But this implies that $G$ is
$VPC(n+1)$ which contradicts the hypothesis that $G$ is not $VPC$. This
contradiction shows that $H_{0}$ is not conjugate into $\partial G$.

It follows from the discussion in section \ref{essentialannuli} that $H_{0}$
is an essential torus in $(G,\partial G)$, so that there is a nontrivial
$H_{0}$--almost invariant subset of $G$ which is dual to $H_{0}$. As
$T_{n+1}(G,\partial G)$ has no $V_{0}$--vertices of commensuriser type, it
follows from Theorem \ref{JSJforVPCofgivenlength} and Remark
\ref{remarksonJSJ} that $H_{0}$ must have small commensuriser, which we denote
by $K$. This means that $K$ contains $H_{0}$ with finite index, so that $K$ is
itself $VPC(n+1)$ and $PD(n+1)$. Note that $K$ must contain $H$, so that $K$
is non-orientable. We let $K_{0}$ denote its orientable subgroup of index $2$.
Note that $K_{0}$ is a maximal torus subgroup of $G$, and is not conjugate
into $\partial G$. Thus there is a nontrivial $K_{0}$--almost invariant subset
of $G$ which must be enclosed by some $V_{0}$--vertex of $T_{n+1}(G,\partial
G)$. It follows that there is a $V_{0}$--vertex $v$ of $T_{n+1}(G,\partial G)$
such that $K_{0}$ is conjugate into $G(v)$. As $K$ contains $K_{0}$ with
finite index, it follows that there is a vertex $w$ of $T_{n+1}(G,\partial G)$
such that $K$ is conjugate into $G(w)$. If $w$ is a $V_{0}$--vertex, then
$K$,$\ $and hence $H$, is conjugate into a $V_{0}$--vertex group of
$T_{n+1}^{c}(G,\partial G)$, as required. So we now consider the case when $w$
is a $V_{1}$--vertex. In particular, $v$ and $w$ must be distinct. Thus there
is an edge $e$ of $T_{n+1}(G,\partial G)$ which is incident to $w$ such that
$G(e)$ contains $K_{0}$. As all the edge groups of $T_{n+1}(G,\partial G)$ are
torus groups, the group $G(e)$ must equal $K_{0}$. Let $E(w)$ denote the
family of subgroups of $G(w)$ which are edge groups for the edges incident to
$w$.

Recall that each edge splitting of $T_{n+1}(G,\partial G)$ is over a $PD(n+1)$
group. If $\partial G$ were empty so that $G$ was a $PD(n+2)$ group, then
Theorem 8.1 of \cite{B-E} would tell us that the pair $(G(w),E(w))$ is
$PD(n+2)$. As $\partial G$ is not empty this need not be the case, but instead
we apply Theorem 8.1 of \cite{B-E} to the graph of groups structure
$\Gamma\Delta$ of $DG$. This shows that the pair $(G(w),E(w))$ becomes
$PD(n+2)$ when $E(w)$ is augmented by suitable groups in $\partial G$. Let
$\overline{E}(w)$ denote this augmented family of subgroups of $G(w)$, so that
the pair $(G(w),\overline{E}(w))$ is $PD(n+2)$. As the commensuriser in $G(w)$
of the group $G(e)=K_{0}$ is not equal to $K_{0}$, Lemma
\ref{boundarygroupismaximal} shows that $K_{0}$ is the only element of the
family $\overline{E}(w)$, and $G(w)$ contains $K_{0}$ with index $2$. Thus $w$
has valence $1$, and $G(w)$ equals a conjugate of $K$. Now it follows that $w$
becomes a $V_{0}$--vertex in the completion $T_{n+1}^{c}(G,\partial G)$, so
that $K$, and hence $H$, is conjugate into a $V_{0}$--vertex group of
$T_{n+1}^{c}(G,\partial G)$, as required.
\end{proof}

\section{Further properties of Torus
Decompositions\label{furtherpropertiesoftorusdecomp}}

In the previous section, we showed that any $PD(n+2)$ pair has a torus
decomposition, and established the basic properties of this decomposition. In
this section, we will establish more detailed information about the vertices
of this decomposition.

Let $(G,\partial G)$ be an orientable $PD(n+2)$ pair, and let $w$ be a $V_{1}%
$--vertex of $T_{n+1}(G,\partial G)$. Let $E(w)$ denote the family of
subgroups of $G(w)$ which are edge groups for the edges incident to $w$.
Recall that each edge splitting of $T_{n+1}(G,\partial G)$ is over a $PD(n+1)$
group. If $\partial G$ is empty so that $G$ is a $PD(n+2)$ group, then Theorem
8.1 of \cite{B-E} tells us that the pair $(G(w),E(w))$ is an orientable
$PD(n+2)$ pair. In general, the pair $(G(w),E(w))$ need not be $PD(n+2)$, but
instead we apply Theorem 8.1 of \cite{B-E} to the graph of groups structure
$\Gamma\Delta$ of $DG$ described in the proof of Theorem
\ref{torusdecompofapair}. This shows that the pair $(G(w),E(w))$ becomes
$PD(n+2)$ when $E(w)$ is augmented by suitable groups in $\partial G$. As any
essential torus in $(G,\partial G)$ is enclosed by some $V_{0}$--vertex of
$T_{n+1}(G,\partial G)$, it follows that any orientable $VPC(n+1)$ subgroup of
$G(w)$ is conjugate into one of the groups in $E(w)$. It will be convenient to
give a name to this property of a $PD(n+2)$ pair.

\begin{definition}
\label{defnofatoroidal}An orientable $PD(n+2)$ pair $(G,\partial G)$ is
\textsl{atoroidal} if any orientable $VPC(n+1)$ subgroup of $G$ is conjugate
into one of the groups in $\partial G$.
\end{definition}

\begin{remark}
As $G$ is torsion free, a $VPC(n+1)$ subgroup of $G$ is also torsion free and
hence is $PD(n+1)$. Thus it makes sense to say that such a subgroup is
orientable. Recall from the preceding paragraph that if $w$ is a $V_{1}%
$--vertex of $T_{n+1}(G,\partial G)$, then the pair $(G(w),E(w))$ becomes
$PD(n+2)$ when $E(w)$ is augmented by suitable groups in $\partial G$. The
resulting $PD(n+2)$ pair is atoroidal.
\end{remark}

This is precisely analogous to the definition of atoroidal for a
$3$--manifold. In the case when an orientable atoroidal $3$--manifold $M$ has
incompressible boundary, it is easy to show that $M$ admits no essential
annulus, unless $M$ is homeomorphic to $T\times I$ or to a twisted $I$--bundle
over the Klein bottle. We will now prove the algebraic analogue of this fact.

\begin{proposition}
\label{atoroidalimpliesannulusfree}Let $(G,\partial G)$ be an orientable
atoroidal $PD(n+2)$ pair, where $n\geq1$. Let $A$ and $B$ be $VPC(n+1)$ groups
in $\partial G$, possibly $A=B$. Let $S$ and $T$ be $VPCn$ subgroups of $A$
and $B$ respectively, and let $g$ be an element of $G$ such that $gSg^{-1}=T$.
Then one of the following holds:

\begin{enumerate}
\item $A$ and $B$ are the same element of $\partial G$, and $g\in A$.

\item $A$ and $B$ are distinct elements of $\partial G$, are the only groups
in $\partial G$, and $A=G=B$. Thus $(G,\partial G)$ is the trivial pair
$(G,\{G,G\})$.

\item $A$ and $B$ are the same element of $\partial G$. Further $A$ is the
only group in $\partial G$, and has index $2$ in $G$.
\end{enumerate}
\end{proposition}

\begin{remark}
The hypothesis that there is $g$ in $G$ such that $gSg^{-1}=T$ means that the
pair $(G,\partial G)$ admits an annulus. The conclusion of the proposition is
that either this annulus is inessential (case 1) or that we have the special
cases in 2) or 3).
\end{remark}

\begin{proof}
As $(G,\partial G)$ is atoroidal, any orientable $VPC(n+1)$ subgroup of $G$ is
conjugate into one of the groups in $\partial G$. Suppose that $G$ contains a
non-orientable $VPC(n+1)$ subgroup $K$. Then $K$ has a subgroup of index $2$
which must be conjugate into a group $H$ of $\partial G$. As $H$ and $K$ are
orientable $PD(n+1)$ groups, it follows that $H$ and $K$ are conjugate
commensurable. As $H$ cannot contain a non-orientable $VPC(n+1)$ subgroup,
this implies that $Comm_{G}(H)\neq H$. Now Lemma \ref{boundarygroupismaximal}
shows that $\partial G$ consists of a single group $H$ which has index $2$ in
$G$, so that we have case 3) of this proposition.

Thus in what follows we will assume that every $VPC(n+1)$ subgroup of $G$ is
orientable, and show that we have case 1) or case 2) of the proposition. We
will consider separately the cases when $A$ and $B$ are the same or distinct.

\medskip

\textbf{Case:} $A$\textit{ and }$B$\textit{ are distinct elements of
}$\partial G$.

\medskip

In this case, we will show that $A$ and $B$ are conjugate commensurable. Then
Lemma \ref{boundarygroupismaximal} shows that we must have case 2) of the proposition.

We suppose that $A$ and $B$ are not conjugate commensurable, and will obtain a
contradiction. After a suitable conjugation, we can arrange that $A\cap B$ is
$VPCn$. Thus $A$, $B$ and $A\cap B$ contain respectively finite index
subgroups $A^{\prime}$, $B^{\prime}$ and $L$ such that $L$ is normal in each
of $A^{\prime}$ and $B^{\prime}$, and $L\backslash A^{\prime}$ and
$L\backslash B^{\prime}$ are both infinite cyclic (see Lemma 13.2 of
\cite{SS2}). Thus $A^{\prime}$ and $B^{\prime}$ are orientable $PD(n+1)$
groups, and $L$ is $PDn$. Now let $K$ denote the amalgamated free product
$A^{\prime}\ast_{L}B^{\prime}$, so that $L$ is normal in $K$ with quotient a
free group $F$ of rank $2$. We identify $F$ with the fundamental group of a
surface $M$ which is a disc with two holes, in such a way that two of the
boundary components of $M$ carry the groups $L\backslash A^{\prime}$ and
$L\backslash B^{\prime}$. Thus $F$ together with the subgroups $L\backslash
A^{\prime}$ and $L\backslash B^{\prime}$ and a third infinite cyclic subgroup
forms a $PD2$ pair. The pre-images in $K$ of the three boundary groups of this
pair yield three $PD(n+1)$ subgroups of $K$. Theorem 7.3 of \cite{B-E} implies
that $K$ together with these subgroups forms a $PD(n+2)$ pair. Two of these
three boundary subgroups of $K$ are equal to $A^{\prime}$ and $B^{\prime}$. As
$A^{\prime}$ and $B^{\prime}$ are orientable and together generate $K$, it
follows that the pair $(K,\partial K)$ is orientable. The inclusions of
$A^{\prime}$ and $B^{\prime}$ into $G$ determine a homomorphism of $K$ into
$G$, and we consider the image $H$ in $G$ of the third boundary subgroup
$\partial_{3}K$ of $K$. As $\partial_{3}K$ is an extension of $L$ by an
infinite cyclic group, $H$ is an extension of the $VPCn$ group $L$ by a cyclic
group. Thus $H$ is $VPC(n+1)$ or $VPCn$.

If $H$ is $VPC(n+1)$, it must be orientable as we are assuming that every
$VPC(n+1)$ subgroup of $G$ is orientable. As $(G,\partial G)$ is atoroidal,
this implies that $H$ is conjugate into a group in $\partial G$. Thus the map
from $K$ to $G$ can be regarded as a map of $PD(n+2)$ pairs. Recall that the
maps from the boundary subgroups $A^{\prime}$ and $B^{\prime}$ of $K$ to the
boundary subgroups $A$ and $B$ of $G$ each have non-zero degree. As $A$ and
$B$ are distinct, and $K$ has only one other boundary group, it follows that
the map from $K$ to $G$ must also have non-zero degree. Hence the image of $K$
in $G$ is a subgroup $G^{\prime}$ of finite index. Further $L$ must be normal
in $G^{\prime}$. Recall that we are considering the case where $(G,\partial
G)$ admits an annulus with fundamental group $L$ whose boundary lies in the
groups $A$ and $B$ of $\partial G$. As we are also assuming that $A$ and $B$
are distinct elements of $\partial G$, this annulus is automatically
essential. As in Definition \ref{defnofa.i.setassociatedtoannulus}, an
essential annulus with fundamental group $L$ determines a nontrivial
$L$--almost invariant subset of $G$. In particular, it follows that
$e(G,L)>1$, so that $e(G^{\prime},L)>1.$ As $L$ is normal in $G^{\prime}$, it
follows that $e(L\backslash G^{\prime})>1.$ Now we apply Stallings' structure
theorem \cite{Stallings1}\cite{Stallings2} for groups with more than one end.
If $e(L\backslash G^{\prime})=2$, then $L\backslash G^{\prime}$ is virtually
infinite cyclic, so that $G^{\prime}$, and hence $G$, must be $VPC(n+1)$. But
then Corollary \ref{PD(n+2)andVPCgroup} implies that $G$, $A$ and $B$ are all
equal, which contradicts our assumption that $A$ and $B$ are not
commensurable. If $e(L\backslash G^{\prime})=\infty$, then either $L\backslash
G^{\prime}$ is of the form $P\ast_{R}Q$, where $R$ is finite, $P\neq R\neq Q$
and one of $P\ $and $Q$ contains $R$ with index at least $3$, or $L\backslash
G^{\prime}$ is of the form $P\ast_{R}$, where $R$ is finite and at least one
of the inclusions of $R$ into $P$ is not an isomorphism. In either case, it is
easy to see that $L\backslash G^{\prime}$ contains infinitely many conjugacy
classes of maximal infinite cyclic subgroups. It follows that $L\backslash
G^{\prime}$ contains infinitely many conjugacy classes of maximal two-ended
subgroups. As a group is two-ended if and only if it is $VPC1$, the pre-images
of these subgroups in $G^{\prime}$ form an infinite collection of conjugacy
classes of maximal $VPC(n+1)$ subgroups of $G^{\prime}$. Recall that we are
assuming that every $VPC(n+1)$ subgroup of $G$ is orientable. As $G^{\prime}$
is of finite index in $G$, there is a finite family $\partial G^{\prime}$ of
$VPC(n+1)$ subgroups of $G^{\prime}$ such that $(G^{\prime},\partial
G^{\prime})$ is an orientable atoroidal $PD(n+2)$ pair. Thus any maximal
$VPC(n+1)$ subgroup of $G^{\prime}$ must be conjugate to one of the groups in
$\partial G^{\prime}$. As $\partial G^{\prime}$ is a finite family, this is a
contradiction, which completes the proof that $H$ cannot be $VPC(n+1)$.

Now consider the case when $H$ is $VPCn$. Recall that $L$ is normal in
$\partial_{3}K$ with infinite cyclic quotient. Thus $L$ is normal in $H\ $with
finite cyclic quotient of some order $d$. There is a $d$--fold regular cover
$M_{d}$ of the surface $M$ in which the pre-image of the third boundary
component of $M$ is a single circle $C$. Now $M_{d}$ determines a subgroup
$K_{d}$ of $K$ of index $d$, and $L$ is normal in $K_{d}$ with quotient
$\pi_{1}(M_{d})$. The boundary component $C$ of $M_{d}$ determines a boundary
subgroup $\partial_{C}$ of $K_{d}$. By construction the image of $\partial
_{C}$ in $G$ is equal to $L$. We let $\overline{M}$ be obtained from $M_{d}$
by gluing a disc onto $C$, and let $\overline{K}$ denote the corresponding
quotient of $K_{d}$. Thus $L$ is normal in $\overline{K}$ with quotient
$\pi_{1}(\overline{M})$, and $\overline{K}$ yields an orientable $PD(n+2)$
pair with one less boundary subgroup than $K_{d}$. Further the homomorphism
from $K_{d}$ to $G$ factors through $\overline{K}$. Each boundary subgroup of
$\overline{K}$ maps to a conjugate of $A$ or $B$, so that the map from
$\overline{K}$ to $G$ is a map of $PD(n+2)$ pairs. Again we have a map of
non-zero degree, as it is of non-zero degree on the boundary. Now we argue
exactly as in the preceding paragraph to obtain a contradiction. This
completes the proof that if $A$ and $B$ are not conjugate commensurable, we
have a contradiction, thus completing the proof that if $A$ and $B$ are
distinct elements of $\partial G$, then we have case 2) of the proposition.

\medskip

\textbf{Case:} $A$\textit{ and }$B$\textit{ are the same element of }$\partial
G$.

\medskip

Recall that we are assuming that every $VPC(n+1)$ subgroup of $G$ is
orientable. Thus case 3) cannot occur. Also recall that, as $A=B$, there are
$VPCn$ subgroups $S$ and $T$ of $A$, and an element $g$ of $G$ such that
$gSg^{-1}=T$. If $g$ lies in $A$, we have case 1) of the proposition. Thus for
the rest of this proof, we will suppose that $g$ does not lie in $A$, and will
obtain a contradiction.

As $g$ does not lie in $A$, it follows that $(G,\partial G)$ admits an
essential annulus with both ends in $A$. Let $K$ denote the double of $G$
along $A$. As $(G,\partial G)$ is an orientable $PD(n+2)$ pair, and $A$ is one
of the groups in the family $\partial G$, there is a natural structure of an
orientable $PD(n+2)$ pair on $K$. As discussed at the end of section
\ref{essentialannuli}, we can double this essential annulus in $(G,\partial
G)$ to obtain an essential torus in $(K,\partial K)$ which clearly crosses the
torus in $K$ represented by the subgroup $A$. Now we consider the uncompleted
torus decomposition $T_{n+1}(K,\partial K)$ of the orientable $PD(n+2)$ pair
$(K,\partial K)$. As $K$ admits an essential torus, either this decomposition
consists of a single $V_{0}$--vertex or it has at least one edge.

If $T_{n+1}(K,\partial K)$ has at least one edge, then the associated
splitting $\sigma$ of $K$ along an essential torus cannot cross any torus in
$K$. In particular, it cannot cross the torus $A$, nor can it equal this
torus. It follows that $\sigma$ determines an essential torus in $(G,\partial
G)$. But this contradicts the fact that $(G,\partial G)$ is atoroidal. It
follows that $T_{n+1}(K,\partial K)$ must consist of a single $V_{0}$--vertex,
so that either $K$ is $VPC$ or the pair $(K,\partial K)$ is $VPCn$--by--Fuchsian.

If $K$ is $VPC$, then $G$ must also be $VPC$, so Corollary
\ref{PD(n+2)andVPCgroup} tells us that either $G$ has two boundary groups each
equal to $G$, or $G$ has one boundary group which is a subgroup of $G$ of
index $2$. The first case is not possible as we assumed $g$ does not lie in
$A$, and the second case is not possible, as part 1) of Corollary
\ref{PDgroupsplitsoverPDHimpliesHismaximal} shows that $G$ would be $VPC(n+1)$
and non-orientable which contradicts our assumption that every $VPC(n+1)$
subgroup of $G\ $is orientable.

If the pair $(K,\partial K)$ is $VPCn$--by--Fuchsian, we let $L$ denote the
$VPCn$ normal subgroup and let $\Phi$ denote the quotient Fuchsian group. We
can assume that $K$ is not $VPC$, so Lemma
\ref{VPC-by-Fuchsiangrouphasuniquefibre} tells us that $L\ $is unique. Also
part 1) of Lemma \ref{torusinVPC-by-Fuchsiangroup} tells us that as $A$ is
$VPC(n+1)$ the intersection $A\cap L$ must be $VPCn$ and hence of finite index
in $L$. As $K$ is the double of $G$ along $A$, it follows that $L$ must be
conjugate into a vertex group of this splitting. As $L$ is normal in $K$, it
now follows that $A$ must contain $L$. Thus $G$ is itself isomorphic to a
$VPCn$--by--Fuchsian group, where the normal $VPCn$ subgroup is $L$. We denote
the quotient group by $\Theta$. As $A$ is a $VPC(n+1)$ subgroup of $G$, the
group $\Theta$ must be infinite. If $\Theta$ is two-ended, then $G$ is
$VPC(n+1)$, and we have a contradiction by the preceding paragraph. Thus we
can assume that $\Theta$ is not two-ended. This implies that there are
elements $\delta$ and $\varepsilon$ in $\Theta$ of infinite order such that
$\delta$ and $\varepsilon$ have non-zero geometric intersection number. The
pre-images in $G$ of the infinite cyclic subgroups of $\Theta$ generated by
$\delta$ and $\varepsilon$ are $VPC(n+1)$ subgroups $D$ and $E$ of $G$. Note
that $D\cap E=L$. By replacing $\delta$ and $\varepsilon$ by their squares if
needed, we can ensure that they are orientable elements of $\Phi$, so that $D$
and $E$ will be orientable. As $\delta$ and $\varepsilon$ have non-zero
geometric intersection number, it follows that $D$ and $E$ are tori in $K$
which cross. As $D$ and $E$ are subgroups of $G$, it follows that $D$ and $E$
are tori in $(G,\partial G)$ which cross. But this contradicts the hypothesis
that $(G,\partial G)$ is atoroidal. This contradiction completes the proof of
the proposition.
\end{proof}

We now apply Proposition \ref{atoroidalimpliesannulusfree} to get information
about the $V_{1}$--vertices of the torus decomposition of a Poincar\'{e}
duality pair.

\begin{proposition}
\label{V1verticesareannulusfree}Let $(G,\partial G)$ be an orientable
$PD(n+2)$ pair, and let $w$ be a $V_{1}$--vertex of the uncompleted torus
decomposition $T_{n+1}(G,\partial G)$. Let $s$ and $t$ be edges of
$T_{n+1}(G,\partial G)$ which are incident to $w$, where possibly $s=t$. Let
$S$ and $T$ be $VPCn$ subgroups of $G(s)$ and $G(t)$ respectively, and let $g$
be an element of $G(w)$ such that $gSg^{-1}=T$. Then one of the following holds:

\begin{enumerate}
\item $s=t$, and $g\in G(s)$.

\item $s$ and $t$ are distinct and $v$ is isolated, so that $s$ and $t$ are
the only edges incident to $w$, and $G(s)=G(w)=G(t)$.

\item $s=t$, the vertex $w$ has valence $1$, and $G(s)$ has index $2$ in
$G(w)$.
\end{enumerate}
\end{proposition}

\begin{remark}
If $n=1$ and $M$ is an orientable Haken $3$--manifold, and a component $W$ of
$M-T(M)$ admits a $\pi_{1}$--injective annulus with boundary in $fr(W)$, then
either this annulus can be homotoped into $fr(W)$, or $W$ is homeomorphic to
$T\times I$ or to a twisted $I$--bundle over the Klein bottle.
\end{remark}

\begin{proof}
Let $E(w)$ denote the family of subgroups of $G(w)$ which are edge groups for
the edges incident to $w$. If $\partial G$ is empty, then Theorem 8.1 of
\cite{B-E} shows that the pair $(G(w),E(w))$ is $PD(n+2)$. In general, as
discussed just before Definition \ref{defnofatoroidal}, the pair becomes
$PD(n+2)$ when $E(w)$ is augmented by some groups in $\partial G$, and the
$PD(n+2)$ pair obtained this way is atoroidal. Applying Proposition
\ref{atoroidalimpliesannulusfree} to this pair yields three cases, which yield
the three cases of this proposition.
\end{proof}

\begin{remark}
If we consider the completed torus decomposition $T_{n+1}^{c}(G,\partial G)$
of the $PD(n+2)$ pair, then the third case in Proposition
\ref{V1verticesareannulusfree} cannot occur. For such $V_{1}$--vertices of
$T_{n+1}(G,\partial G)$ become $V_{0}$--vertices when $T_{n+1}(G,\partial G)$
is completed to $T_{n+1}^{c}(G,\partial G)$.
\end{remark}

We will also need information about the $V_{0}$--vertices of the torus
decomposition of a Poincar\'{e} duality pair. Part 1) of Theorem
\ref{torusdecompofapair} states that a $V_{0}$--vertex of the uncompleted
torus decomposition $T_{n+1}(G,\partial G)$ must be isolated or of Seifert
type adapted to $\partial G$. In the next result, we consider the second type
of $V_{0}$--vertex.

\begin{lemma}
\label{essentialannulusinFuchsiangroup}Let $(G,\partial G)$ be an orientable
$PD(n+2)$ pair, and let $v$ be a $V_{0}$--vertex of the uncompleted torus
decomposition $T_{n+1}(G,\partial G)$ which is of Seifert type adapted to
$\partial G$. (See Definition \ref{defnofSeiferttypeadaptedtodG}.) Let $L$
denote the $VPCn$ normal subgroup of $G(v)$ with Fuchsian quotient $\Phi$. Let
$s$ and $t$ be edges of $T_{n+1}(G,\partial G)$ which are incident to $v$,
where possibly $s=t$. Let $S$ and $T$ be $VPCn$ subgroups of $G(s)$ and $G(t)$
respectively, and let $g$ be an element of $G(v)$ such that $gSg^{-1}=T$. Then
one of the following holds:

\begin{enumerate}
\item $s=t$, and $g\in G(s)$.

\item $S$ is commensurable with $L$.
\end{enumerate}
\end{lemma}

\begin{remark}
If $n=1$ and $M$ is an orientable Haken $3$--manifold, then any component of
$T(M)$ is a Seifert fibre space. The corresponding result is a standard result
about Seifert fibre spaces. Namely that if $W$ is a Seifert fibre space with
incompressible boundary which admits a $\pi_{1}$--injective annulus with
boundary in $\partial W$, then either this annulus can be homotoped into
$\partial W$, or this annulus is vertical in $W$.
\end{remark}

\begin{proof}
Let $E(v)$ denote the family of subgroups of $G(v)$ which are edge groups for
the edges incident to $v$. If $\partial G$ is empty, then as discussed just
before Definition \ref{defnofatoroidal}, the pair $(G(v),E(v))$ is orientable
$PD(n+2)$. In general, the pair becomes $PD(n+2)$ when $E(v)$ is augmented by
some groups in $\partial G$. Let $\overline{E}(v)$ denote this augmented
family of groups. As $v$ is of Seifert type adapted to $\partial G$, the
normal subgroup $L$ of $G(v)$ is contained in each group in $\overline{E}(v)$.
Let $\partial\Phi$ denote the family of subgroups of $\Phi$ obtained by taking
the quotient by $L$ of each group in $\overline{E}(v)$. Thus the pair
$(\Phi,\partial\Phi)$ is the orbifold fundamental group of a compact
$2$--dimensional orbifold $(F,\partial F)$. Note that as $v$ is of Seifert
type adapted to $\partial G$, the group $\Phi$ is not finite nor two-ended.

Let $S^{\prime}$, $T^{\prime}$ and $g^{\prime}$ denote the images of $S$, $T$
and $g$ in $\Phi$, so that we have the equation $g^{\prime}S^{\prime}%
g^{\prime-1}=T^{\prime}$ in $\Phi$. If $S^{\prime}$ is finite, then $S\cap L$
has finite index in $S$, so that $S$ is commensurable with $L$, and we have
case 2) of the lemma. Otherwise $S^{\prime}$ is an infinite subgroup of a
group $H$ in $\partial\Phi$. Let $A$ denote an infinite cyclic subgroup of
$S^{\prime}$, and let $F_{A}$ denote the orbifold cover of $F$ with
fundamental group $A$. Thus $F_{A}$ has a boundary component which carries
$A$. Also, as $A$ is torsion free, $F_{A}$ is a surface. If $g^{\prime}$ does
not lie in $H$, then this surface admits an essential annulus, and so must be
an annulus. In particular, $F_{A}$ is compact and hence a finite cover of $F$.
As $\Phi$ is not two-ended, this is impossible so that $g^{\prime}$ must lie
in $H$, and hence so does $T^{\prime}$. Thus $S^{\prime}$, $T^{\prime}$ and
$g^{\prime}$ all lie in the same group $H$ in $\partial\Phi$, which means that
we have case 1) of the lemma.
\end{proof}

An important consequence of Proposition \ref{V1verticesareannulusfree} is the following.

\begin{lemma}
\label{normaliserequalscommensuriser} Let $(G,\partial G)$ be an orientable
$PD(n+2)$ pair, let $V$ be a non-isolated $V_{0}$--vertex of the uncompleted
torus decomposition $T_{n+1}(G,\partial G)$ which is of Seifert type adapted
to $\partial G$, and let $L$ denote the $VPCn$ normal subgroup of $G(V)$ with
Fuchsian quotient. Then
\[
G(V)=N_{G}(L)=Comm_{G}(L),
\]
where $N_{G}(L)$ denotes the normalizer of $L$ in $G$.
\end{lemma}

\begin{remark}
If $n=1$ and $M$ is an orientable Haken $3$--manifold, then any component of
$T(M)$ is a Seifert fibre space, and so its fundamental group has an infinite
cyclic normal subgroup. Let $W$ be a component of $T(M)$, and $L$ this normal
subgroup. Then the corresponding result is that $\pi_{1}(W)$ is equal to the
commensuriser of $L$ in $\pi_{1}(M)$.
\end{remark}

\begin{proof}
The inclusions $G(V)\subset N_{G}(L)\subset Comm_{G}(L)$ are all clear. Thus
it remains to prove that $Comm_{G}(L)\subset G(V)$.

Let $T$ denote the universal covering $G$--tree of $T_{n+1}(G,\partial G)$ and
let $v$ denote a vertex of $T$ above $V$ with stabiliser $G(v)$ equal to
$G(V)$. Suppose that $Comm_{G}(L)$ does not equal $G(V)$. Then there is an
element $g$ of $Comm_{G}(L)$ which does not fix $v$, and we let $L^{\prime}$
denote the intersection $L\cap gLg^{-1}$. Thus $L^{\prime}$ is a $VPCn$
subgroup of $G$ which fixes both $v$ and $gv$, and hence fixes every edge on
the path $\lambda$ joining $v$ to $gv$.

If $w$ is a $V_{1}$--vertex of $\lambda$, then $L^{\prime}$ fixes $w$ and two
distinct incident edges. Now Proposition \ref{V1verticesareannulusfree} shows
that $w$ has valence $2$ in $T$ and that the two incident edges each have the
same stabiliser. It also implies that $G(w)$ contains this stabiliser with
index $1$ or $2$, but we will not need to distinguish these cases. Recall that
any $V_{0}$--vertex of $T_{n+1}(G,\partial G)$ is isolated or of Seifert type
adapted to $\partial G$. If $v_{0}$ is an interior $V_{0}$--vertex of
$\lambda$ which is of Seifert type adapted to $\partial G$, then Lemma
\ref{essentialannulusinFuchsiangroup} shows that $L^{\prime}$ must be
commensurable with the $VPCn$ subgroup $L_{0}$ of $G(v_{0})$ which is normal
in $G(v_{0})$ with Fuchsian quotient.

Now let $\mu$ denote the subinterval of $\lambda$ between $v$ and the first
$V_{0}$--vertex $w$ of $\lambda$ which is of Seifert type adapted to $\partial
G$. Possibly $\mu$ equals $\lambda$. The preceding discussion shows that any
pair of adjacent edges of $\mu$ have the same stabiliser. Thus the stabiliser
of $\mu$ equals an edge group $H$ of $G(v)$, which is also an edge group of
$G(w)$. Hence the subgroup $K$ of $G$ generated by $G(v)$ and $G(w)$ equals
$G(v)\ast_{H}G(w)$, and $K$ has the natural structure of an orientable
$PD(n+2)$ pair. The splitting of $K$ over $H$ is over an essential torus in
$K$. Recall that $L$ is the $VPCn$ normal subgroup of $G(v)$ with Fuchsian
quotient, and let $L^{\prime\prime}$ denote the corresponding $VPCn$ normal
subgroup of $G(w)$. If $w$ equals $gv$, then $L^{\prime\prime}=gLg^{-1}$, so
that $L^{\prime}$ is a subgroup of $L^{\prime\prime}$. If $w$ is not equal to
$gv$, then the discussion in the preceding paragraph shows that $L^{\prime}$
is commensurable with $L^{\prime\prime}$. As $L^{\prime}$ is contained in $L$,
it follows that in either case $L$ and $L^{\prime\prime}$ are commensurable.
As each is a normal subgroup of $H$ with quotient isomorphic to $\mathbb{Z}$
or to $\mathbb{Z}_{2}\ast\mathbb{Z}_{2}$, part 2) of Lemma
\ref{normalsubgroupsofVPCgroupsofcolength1} shows that $L$ and $L^{\prime
\prime}$ must be equal. Thus the pair $(K,\partial K)$ is $VPCn$%
--by--Fuchsian. It follows that there is an essential torus in $K$ which
crosses the essential torus determined by $H$, and hence an essential torus in
$G$ which crosses the torus determined by $H$. But this contradicts the fact
that the edge splittings of $T_{n+1}(G,\partial G)$ do not cross any essential
torus in $G$.

We conclude that every element of $Comm_{G}(L)$ lies in $G(V)$. Thus
$G(V)=N_{G}(L)=Comm_{G}(L)$ as required.
\end{proof}

There is another more technical consequence of Proposition
\ref{V1verticesareannulusfree} which will only be needed in the case when $G$
is a $PD(n+2)$ group, i.e. when $\partial G$ is empty, but the general result
is no more difficult. Before stating this result it will be convenient to have
the following definition.

\begin{definition}
\label{defnofperipheral} Let $\Gamma$ be a minimal graph of groups structure
for a group $G$, let $v$ be a vertex of $\Gamma$, and let $X$ be a $H$--almost
invariant subset of $G$ which is enclosed by $v$.

Then $X$ is \textsl{peripheral in }$v$ if $X$ is equivalent to the almost
invariant subset of $G$ associated to some edge of $\Gamma$ which is incident
to $v$.
\end{definition}

The word peripheral here is very natural, but the reader should note that when
$(G,\partial G)$ is a Poincar\'{e} duality pair, this idea may have nothing to
do with $\partial G$.

\begin{lemma}
\label{torienclosedbydistinctverticesmustlieinone}Let $(G,\partial G)$ be an
orientable $PD(n+2)$ pair, and let $H\ $and $H^{\prime}$ be essential tori in
$(G,\partial G)$ whose dual almost invariant sets $X$ and $X^{\prime}$ are
enclosed by $V_{0}$--vertices $V$ and $V^{\prime}$ of $T_{n+1}(G,\partial G)$.
Suppose that $H\cap H^{\prime}$ is a $VPCn$ subgroup $K$ of $G$. Suppose
further that $X$ is not peripheral in $V$. Then $V$ must be of Seifert type
adapted to $\partial G$, each of $X$ and $X^{\prime}$ is enclosed by $V$, and
$K$ is commensurable with the $VPCn$ normal subgroup $L$ of $G(V)$ with
Fuchsian quotient.
\end{lemma}

\begin{proof}
As $X$ is not peripheral in $V$, the vertex $V$ of $T_{n+1}(G,\partial G)$
cannot be isolated. Thus Theorem \ref{torusdecompofapair} implies it must be
of Seifert type adapted to $\partial G$. Now part 1) of Lemma
\ref{torusinVPC-by-Fuchsiangroup} tells us that the intersection $H\cap L$
must be $VPCn$, and hence of finite index in $L$.

Suppose that $V$ equals $V^{\prime}$. Then part 1) of Lemma
\ref{torusinVPC-by-Fuchsiangroup} tells us that the intersection $H^{\prime
}\cap L$ must also be $VPCn$, and hence of finite index in $L$. Thus the
intersection $H\cap H^{\prime}\cap L$ is of finite index in $L$, and so is
also $VPCn$. As $K=H\cap H^{\prime}$ is assumed to be $VPCn$, it follows that
$K\ $and $L$ must be commensurable, as required.

In what follows, we will suppose that $V$ and $V^{\prime}$ are distinct. Let
$T$ denote the universal covering $G$--tree of $T_{n+1}(G,\partial G)$. Let
$v$ be a vertex above $V$ which is fixed by $H$, and let $v^{\prime}$ be a
vertex above $V^{\prime}$ which is fixed by $H^{\prime}$. Then $K$ fixes the
path $\lambda$ joining $v$ to $v^{\prime}$. Now we are in much the same
situation as in the proof of Lemma \ref{normaliserequalscommensuriser}, with
$v^{\prime}$ in place of $gv$. Thus if $w$ is a $V_{1}$--vertex of $\lambda$,
then $w$ has valence $2$ in $T$ and the two incident edges each have the same
stabiliser, and if $v_{0}$ is an interior $V_{0}$--vertex of $\lambda$ which
is of Seifert type adapted to $\partial G$, then $K$ must be commensurable
with the $VPCn$ normal subgroup $L_{0}$ of $G(v_{0})$ with Fuchsian quotient.

Now we consider the edge $e$ of $\lambda$ which is incident to $v$. As $K$
fixes $e$, it is a subgroup of the torus group $G(e)$. Our hypothesis that $X$
is not peripheral in $V$ implies that $H\cap G(e)$ is commensurable with the
normal $VPCn$ subgroup $L$ of $G(v)$. As $K$ is a $VPCn$ subgroup of $G(e)$
and of $H$, it follows that $K$ must be commensurable with $L$. It remains to
prove that $X^{\prime}$ must be enclosed by $V$.

As in the proof of Lemma \ref{normaliserequalscommensuriser}, we let $\mu$
denote the subinterval of $\lambda$ between $v$ and the first $V_{0}$--vertex
$v^{\prime\prime}$ of $\lambda$ which is of Seifert type adapted to $\partial
G$. As before any pair of adjacent edges of $\mu$ have the same stabiliser.
Thus the edge group $G(e)$ is equal to an edge subgroup of $G(v^{\prime\prime
})$.

If $v^{\prime\prime}$ is not equal to $v^{\prime}$, so that $\mu$ is properly
contained in $\lambda$, the above discussion shows that $K$ must be
commensurable with the normal $VPCn$ subgroup $L^{\prime\prime}$ of
$G(v^{\prime\prime})$. As $K$ is commensurable with $L$ and $L^{\prime\prime}%
$, it follows that $L$ and $L^{\prime\prime}$ are commensurable. As each is a
normal subgroup of $G(e)$ with quotient isomorphic to $\mathbb{Z}$ or to
$\mathbb{Z}_{2}\ast\mathbb{Z}_{2}$, Lemma
\ref{normalsubgroupsofVPCgroupsofcolength1} shows that $L$ and $L^{\prime
\prime}$ must be equal. Now Lemma \ref{normaliserequalscommensuriser} yields a
contradiction. We deduce that $\mu$ must equal $\lambda$.

Let $e^{\prime}$ denote the edge of $\lambda$ which is incident to $v^{\prime
}$. As all consecutive edges of $\lambda$ have the same stabiliser, it follows
that $G(e)=G(e^{\prime})$. If $X^{\prime}$ is not peripheral in $v^{\prime}$,
the preceding argument with the roles of $X$ and $X^{\prime}$ reversed shows
that $K$ must be commensurable with the normal subgroup $L^{\prime}$ of
$G(v^{\prime})$. If $X^{\prime}$ is peripheral in $v^{\prime}$ but is not
equivalent to the almost invariant set associated to the edge $e^{\prime}$,
then Lemma \ref{essentialannulusinFuchsiangroup} shows that $K$ must be
commensurable with $L^{\prime}$. Thus either $K$ is commensurable with
$L^{\prime}$, or $X^{\prime}$ is equivalent to the almost invariant subset of
$G$ associated to $e^{\prime}$. In the first case, we derive a contradiction
as before. In the second case, we recall that each interior vertex of
$\lambda$ is isolated or of special Seifert type. Thus consecutive edges of
$\lambda$ not only have the same stabilisers but they have equivalent
associated almost invariant subsets of $G$. Hence $X^{\prime}$ must be
equivalent to the almost invariant subset of $G$ associated to $e$, so that
$X^{\prime}$ is enclosed by $v$, and hence by $V$, as required. This completes
the proof of the lemma.
\end{proof}

\section{Enclosing properties of Annulus-Torus
Decompositions\label{enclosingpropertiesofATdecomps}}

Let $G$ denote any almost finitely presented group which has no nontrivial
almost invariant subsets over $VPC(<n)\ $subgroups. Recall that $\mathcal{F}%
_{n,n+1}$ denotes the family of equivalence classes of all nontrivial almost
invariant subsets of $G$ which are over $VPCn$ subgroups and of equivalence
classes of all $n$--canonical almost invariant subsets of $G$ which are over
$VPC(n+1)$ subgroups. The decomposition $\Gamma_{n,n+1}$ of $G$ is the reduced
algebraic regular neighbourhood of $\mathcal{F}_{n,n+1}$, and its $V_{0}%
$--vertices correspond to the cross connected components of $\mathcal{F}%
_{n,n+1}$. Some of its properties are described in Theorem
\ref{JSJforVPCoftwolengths}. For brevity, we will denote $\mathcal{F}_{n,n+1}$
by ${\mathcal{F}}$ in the rest of this section. In the case when $(G,\partial
G)$ is a $PD(n+2)$ pair, we will also consider the family ${\mathcal{F}%
^{\prime}}$ of equivalence classes of all nontrivial almost invariant subsets
of $G$ which are over $VPCn$ subgroups together with equivalence classes of
all nontrivial almost invariant subsets of $G$ which are over $VPC(n+1)$
subgroups and are adapted to {$\partial G$}. Note that, by Lemma
\ref{aisetovernon-orblesubgroupisnotadaptedtodG}, nontrivial almost invariant
subsets of $G$ which are over $VPC(n+1)$ subgroups and are adapted to
{$\partial G$ are dual to essential tori in }$(G,\partial G)$. Thus from the
topological point of view, it seems very natural to consider the family
$\mathcal{F}^{\prime}$. However, one obvious reason why our theory in
\cite{SS2} discusses algebraic regular neighbourhoods of $\mathcal{F}%
=\mathcal{F}_{n,n+1}$ rather than $\mathcal{F}^{\prime}$ is that, for general
groups, there is no analogue of $\mathcal{F}^{\prime}$. Another reason is that
elements of $\mathcal{F}^{\prime}-{\mathcal{F}}$ need not be $n$--canonical
and our methods in \cite{SS2} cannot handle this situation. In this section,
we will show that $\mathcal{F}$ is contained in $\mathcal{F}^{\prime}$
(Corollary \ref{FiscontainedinF'}) and that $\mathcal{F}$ and $\mathcal{F}%
^{\prime}$ have the same reduced algebraic regular neighbourhoods (Theorem
\ref{FandF'havesameregnbhds}). This technical result will be used in section
\ref{proofofmaintheorem} of this paper in a crucial way. It is the algebraic
analogue of the situation discussed in the second paragraph of section
\ref{preliminaries}, with $\Gamma_{n,n+1}$ being the algebraic analogue of
$AT(M)$, the family $\mathcal{F}^{\prime}$ being the algebraic analogue of the
family of all essential annuli and tori in $M$, and the family $\mathcal{F}$
being the algebraic analogue of the family of all essential annuli in $M$
together with those essential tori in $M$ which do not cross any essential
annulus in $M$.

First we will show that $\mathcal{F}$ is contained in $\mathcal{F}^{\prime}$.
To do this we need to consider almost invariant subsets of $G$ which are not
adapted to $\partial G$. We start with the following simple result.

\begin{lemma}
\label{HintersectH_iisVPCn}Let $n\geq1$, and let $(G,\partial G)$ be an
orientable $PD(n+2)$ pair. Suppose that $X$ is a nontrivial almost invariant
subset of $G$ over a $VPC(n+1)$ group $H$ and that $X$ is not adapted to
$\partial G$. Then there is a subgroup $S_{i}$ of $G$ with a conjugate in
$\partial G$ such that $H\cap S_{i}$ is $VPCn$.
\end{lemma}

\begin{proof}
As $X$ is not adapted to $\partial G$, there is $S_{i}$ in {$\partial G$}, and
$g\in G$, such that both $X\cap gS_{i}$ and $X^{\ast}\cap gS_{i}$ are not
$H$--finite. By replacing $S_{i}$ by a conjugate, we can arrange that both
$X\cap S_{i}$ and $X^{\ast}\cap S_{i}$ are not $H$--finite. Thus both of them
are not $(H\cap S_{i})$--finite. Let $K$ denote $H\cap S_{i}$. Then $X\cap
S_{i}$ is a nontrivial almost invariant subset of $S_{i}$ which is over $K$,
so that $e(S_{i},K)>1$. As $K$ is a subgroup of the $VPC(n+1)$ group $H$, it
must be $VPC(\leq n+1)$. As $K$ must have infinite index in the $PD(n+1)$
group $S_{i}$, Strebel's result \cite{Strebel} shows that $K$ has
cohomological dimension $\leq n$. Thus $K$ is $VPC(\leq n)$.

Recall the following long exact cohomology sequence from page
\pageref{coboundarymap}.%
\[
H^{0}(G;\mathbb{Z}_{2}E)\rightarrow H^{0}(G;PE)\rightarrow H^{0}%
(G;PE/\mathbb{Z}_{2}E)\rightarrow H^{1}(G;\mathbb{Z}_{2}E)\rightarrow
H^{1}(G;PE)\rightarrow
\]

Here $H\ $is a subgroup of a group $G$, and $E$ denotes $H\backslash G$. Also
$H^{0}(G;PE)\cong\mathbb{Z}_{2}$, and $e(G,H)$ equals the dimension over
$\mathbb{Z}_{2}$ of $H^{0}(G;PE/\mathbb{Z}_{2}E)$. Thus if $e(G,H)>1$, then
$H^{1}(G;\mathbb{Z}_{2}E)$ is non-zero. In the setting of the present lemma,
as $e(S_{i},K)>1$ we see that $H^{1}(S_{i};\mathbb{Z}_{2}(K\backslash S_{i}))$
must be non-zero. As $S_{i}$ is $PD(n+1)$, this last group is isomorphic to
$H_{n}(S_{i};\mathbb{Z}_{2}(K\backslash S_{i}))$ which is in turn isomorphic
to $H_{n}(K,\mathbb{Z}_{2})$. Thus $H_{n}(K,\mathbb{Z}_{2})$ is non-zero. As
$K$ is torsion free and $VPC(\leq n)$, it follows that $K$ must be $VPCn$, so
that $H\cap S_{i}$ is $VPCn$ as required.
\end{proof}

Now we can prove the following.

\begin{proposition}
\label{XnotrelativetodGimpliesXcrossessomeVPCn}Let $(G,\partial G)$ be an
orientable $PD(n+2)$ pair. Suppose that $X$ is a nontrivial almost invariant
subset of $G$ which is over a $VPC(n+1)$ subgroup $H$ and is not adapted to
{$\partial G$}. Then $X$ crosses some nontrivial almost invariant subset of
$G$ which is over some $VPCn$ subgroup of $G$ and is dual to an essential annulus.
\end{proposition}

\begin{remark}
In particular, if such $X$ exists, then there are essential annuli in the pair
$(G,\partial G)$ which are over $VPCn$ subgroups of $H$.
\end{remark}

\begin{proof}
The proof of this result is suggested by our arguments in \cite{SS1}. Since
$X$ is not adapted to {$\partial G$}, the proof of Lemma
\ref{HintersectH_iisVPCn} shows that, after conjugation, there is a group
$S_{i}$ in $\partial G$ such that $X\cap S_{i}$ is a nontrivial almost
invariant subset of $S_{i}$ over $H\cap S_{i}$, and $H\cap S_{i}$ is $VPCn$.
We denote $H\cap S_{i}$ by $K$. By replacing $H$ by a subgroup of finite index
if necessary, we can arrange that $H$ is orientable, that $K$ has a finite
index subgroup $L$ which is normal in $H$, and that the quotient $L\backslash
H$ is infinite cyclic. Theorem 7.3 of \cite{B-E} now implies that $L$ is orientable.

As in section \ref{essentialannuli}, we choose an aspherical space $M$ with
fundamental group $G$ and with aspherical subspaces corresponding to $\partial
G$ whose union is denoted $\partial M$. We denote the cover of $M$
corresponding to $L$ by $M_{L}$. As $S_{i}$ contains $L$, there is a component
$\Sigma$ of $\partial{M_{L}}$ with fundamental group $L$. Since $L$ is normal
in $H$ with infinite cyclic quotient, the quotient $L\backslash H$ acts
naturally on $M_{L}$, and we obtain infinitely many components of $\partial
M_{L}$ which are translates of $\Sigma$ and have fundamental group $L$. Let
$Z$ denote the support of a $0$--cochain on $M_{H}$ with finite coboundary
which represents the element of $H_{e}^{0}(M_{H};\mathbb{Z}_{2})$ determined
by $H\backslash X$. Let $p:M_{L}\rightarrow M_{H}$ denote the covering
projection. As $X\cap S_{i}$ is a nontrivial $L$--almost invariant subset of
$S_{i}$, the vertices of $p(\Sigma)$ in $M_{H}$ meet both $Z$ and $Z^{\ast}$
in infinite sets. As $X$ is $H$--invariant, the vertices of the image of each
translate of $\Sigma$ also meet both $Z$ and $Z^{\ast}$ in infinite sets.

As the number of these translates of $\Sigma$ is infinite, we can apply Lemma
\ref{thereexisttwocrossingannuli}. As in that lemma, let $\Sigma_{1}$,
$\Sigma_{2}$, $\Sigma_{3}$ and $\Sigma_{4}$ denote four distinct translates of
$\Sigma$ by elements of $L\backslash H$, and let $A_{ij}$ denote the annulus
in $M_{L}$ with fundamental group $L$ and joining $\Sigma_{i}$ and $\Sigma
_{j}$. Let $Y_{ij}$ denote the $L$--almost invariant subset of $G$ dual to
$A_{ij}$. Let $Z_{ij}$ denote the support of a $0$--cochain on $M_{L}$ with
finite coboundary which represents the element of $H_{e}^{0}(M_{L}%
;\mathbb{Z}_{2})$ determined by $L\backslash Y_{ij}$. The proof of Lemma
\ref{thereexisttwocrossingannuli} shows that there are distinct integers $i$,
$j$, $k$ and $l$ such that $Z_{ij}$ separates $\Sigma_{k}$ and $\Sigma_{l}$,
i.e. $\Sigma_{k}$ is almost contained in $Z_{ij}$ and $\Sigma_{l}$ is almost
contained in $Z_{ij}^{\ast}$, or vice versa. It follows that $X$ must cross
$Y_{ij}$ because each of the four corners of the pair $(Z,pZ_{ij})$ is
infinite, as it has infinite intersection with $p\Sigma_{k}$ or $p\Sigma_{l}$.
As $Y_{ij}$ is dual to an essential annulus, this completes the proof of the lemma.
\end{proof}

The point of this proposition is the following corollary.

\begin{corollary}
\label{FiscontainedinF'}Let $(G,\partial G)$ be an orientable $PD(n+2)$ pair,
and let $\mathcal{F}$ and $\mathcal{F}^{\prime}$ be defined as above. If $X$
is a $n$--canonical almost invariant subset of $G$ over a $VPC(n+1)$ subgroup
$H$, then $X$ is automatically adapted to {$\partial G$}. Thus $\mathcal{F}$
is contained in $\mathcal{F}^{\prime}$.
\end{corollary}

\begin{remark}
\label{adaptedtodGimpliesorientablegroup}As $X$ is adapted to $G$, it follows
from Lemma \ref{aisetovernon-orblesubgroupisnotadaptedtodG} that $H$ must be orientable.
\end{remark}

\begin{proof}
If $X$ is not adapted to {$\partial G$}, then Proposition
\ref{XnotrelativetodGimpliesXcrossessomeVPCn} tells us that $X$ crosses some
nontrivial almost invariant set over some $VPCn$ subgroup of $G$, which
contradicts the hypothesis that $X$ is $n$--canonical. Hence the definitions
of $\mathcal{F}$ and $\mathcal{F}^{\prime}$ show that $\mathcal{F}$ is
contained in $\mathcal{F}^{\prime}$, as claimed.
\end{proof}

Next we give several results about how elements of $\mathcal{F}^{\prime
}-\mathcal{F}$ can cross elements of $\mathcal{F}$.

\begin{lemma}
\label{XcrossesYimpliesitcrossesY'strongly} Let $(G,\partial G)$ be an
orientable $PD(n+2)$ pair, and let $\mathcal{F}$ and $\mathcal{F}^{\prime}$ be
defined as above. Let $X$ be an element of $\mathcal{F}^{\prime}-\mathcal{F}$
which crosses a nontrivial almost invariant set $Y$ over a $VPCn$ group $K$.
Then the following hold:

\begin{enumerate}
\item If $Y$ is dual to an essential annulus, then $X$ crosses $Y$ strongly.

\item There is an almost invariant set $Y^{\prime}$ over a group $K^{\prime}$
commensurable with $K$ such that $Y^{\prime}$ is dual to an essential annulus
and $X$ crosses $Y^{\prime}$ strongly.
\end{enumerate}
\end{lemma}

\begin{remark}
If $n=1$ and $M$ is an orientable Haken $3$--manifold, the corresponding
result holds. For $X$ corresponds to a torus which must be homotopic into a
component $W$ of $T(M)$ which meets $\partial M$. Such $W$ is a Seifert fibre
space whose intersection with $\partial M$ consists of tori and vertical
annuli. Thus $W$ is filled by vertical annuli, so $X$ must cross one of them.

Note that in this topological situation, an annulus can never cross a torus
strongly, for the intersection of their fundamental groups must be of finite
index in the fundamental group of the annulus.
\end{remark}

\begin{proof}
As $X$ lies in $\mathcal{F}^{\prime}-\mathcal{F}$, it is a nontrivial almost
invariant subset of $G$ which is over an orientable $VPC(n+1)$ group $H$ and
is adapted to $\partial G$.

1) Suppose that $Y$ is dual to an essential annulus $A$. As discussed
immediately after the proof of Lemma \ref{thereexisttwocrossingannuli}, the
double $DG$ of $G$ along $\partial G$ contains a torus which is the double of
the annulus $A$. We let $DY$ denote the $DK$--almost invariant subset of $DG$
associated to this torus. As $X$ is adapted to $\partial G$, Lemma
\ref{propertiesofadapted} yields a $H$--almost invariant subset $\overline{X}$
of $DG$ such that $\overline{X}\cap G$ equals $X$. As $X$ and $Y$ cross, it
follows that $\overline{X}\ $and $DY$ must also cross. As $H$ and $DK$ are
both orientable $VPC(n+1)$ subgroups of the orientable $PD(n+2)$ group $DG$,
each has two coends in $DG$. Now Proposition 7.4 of \cite{SS2} shows that
$\overline{X}\ $and $DY$ must cross each other strongly. This implies that $X$
crosses $Y$ strongly, thus completing the proof of part 1).

2) As in section \ref{essentialannuli}, we choose an aspherical space $M$ with
fundamental group $G$ and with aspherical subspaces corresponding to $\partial
G$ whose union is denoted $\partial M$. Consider the cover $M_{K}$, the
element of $H_{f}^{1}(M_{K};\mathbb{Z})$ corresponding to $K\backslash Y$ and
its dual class $\alpha$ in $H_{n+1}(M_{K},\partial{M_{K};}\mathbb{Z})$. The
boundary of $\alpha$ has support in only finitely many components $\Sigma
_{1},\ldots,\Sigma_{k}$ of $\partial M_{K}$. Note that Proposition
\ref{aisetimpliesannulus} tells us that there must be at least one such
component of $\partial M_{K}$, and each such component must carry a subgroup
of finite index in $K$. By replacing $K$ by a suitable subgroup of finite
index, we can arrange that $K$ is orientable and that each of $\Sigma
_{1},\ldots,\Sigma_{k}$ carries $K$ itself. As $K$ is orientable, we know that
$k\geq2$. For each pair of distinct integers $i$ and $j$ between $1$ and $k$,
there is an essential untwisted annulus $A_{ij}$ in $M_{K}$ whose boundary
lies in $\Sigma_{i}\cup\Sigma_{j}$. Let $Y_{ij}$ denote the dual $K$--almost
invariant subset of $G$. Let $Z_{ij}$ denote the support of a $0$--cochain on
$M_{K}$ with finite coboundary which represents the element of $H_{e}%
^{0}(M_{K};\mathbb{Z}_{2})$ determined by $K\backslash Y_{ij}$. We will show
that $X$ must cross one of these $Y_{ij}$'s. This will be the required
$Y^{\prime}$. Now part 1) will imply that $X$ crosses $Y^{\prime}$ strongly$.$

Let $Z$ denote the support of a $0$--cochain on $M_{H}$ with finite coboundary
which represents the element of $H_{e}^{0}(M_{H};\mathbb{Z}_{2})$ determined
by $H\backslash X$. As $X$ is adapted to $\partial G$, we know that for each
component $\Sigma$ of $\partial M_{H}$ the vertices of $\Sigma$ are almost all
in $Z$ or almost all in $Z^{\ast}$. Further the vertices of $\Sigma$ are
$H$--infinite, as $X$ is nontrivial, so that the vertices of $\Sigma$ cannot
be almost all in $Z$ and almost all in $Z^{\ast}$. Thus if $X$ crosses none of
the $Y_{ij}$'s, it follows that, after replacing $X$ by $X^{\ast}$ if needed,
the vertices of the images in $M_{H}$ of each $\Sigma_{i}$, $1\leq i\leq k$,
almost all lie in $Z$. As $X$ and $Y_{ij}$ do not cross, it now follows that
$Y_{ij}\leq X$ or $Y_{ij}^{\ast}\leq X$. By replacing $Y_{ij}$ by its
complement if needed, we can arrange that $Y_{ij}\leq X$, for all $i$ and $j$.
Part 2) of Proposition \ref{annuligenerate} tells us that $Y$ is equivalent to
a sum of some of the $Y_{ij}$'s and their complements. If there are no
complements in the sum, then we have $Y\leq X$, which contradicts the
assumption that $X$ crosses $Y$. The same inequality holds if the number of
complements in the sum is even. If the number of complements in the sum is
odd, then we have $X^{\ast}\leq Y$, which again contradicts the assumption
that $X$ crosses $Y$. This contradiction shows that $X$ must cross some
$Y_{ij}$ and so completes the proof of the lemma.
\end{proof}

Two easy consequences are the following results.

\begin{lemma}
\label{firstlemma}Let $(G,\partial G)$ be an orientable $PD(n+2)$ pair, and
let $\mathcal{F}$ and $\mathcal{F}^{\prime}$ be defined as above. Let $X$ be
an element of $\mathcal{F}^{\prime}-\mathcal{F}$ which crosses nontrivial
almost invariant sets $Y$ and $Y^{\prime}$ over $VPCn$ groups $L$ and
$L^{\prime}$ respectively. Then $L$ and $L^{\prime}$ are commensurable.
\end{lemma}

\begin{remark}
If $n=1$ and $M$ is an orientable Haken $3$--manifold, the corresponding
result holds. For $X$ corresponds to a torus which crosses an annulus in $M$,
and so must be homotopic into a component $W$ of $T(M)$ which meets $\partial
M$. Such $W$ is a Seifert fibre space whose intersection with $\partial M$
consists of tori and vertical annuli. If $W$ has a unique Seifert fibration,
up to isotopy, it follows that all essential annuli in $M$ which are homotopic
into $W$ carry commensurable groups as this group must be commensurable with
the fibre group of $W$. Otherwise, $W$ is homeomorphic to $T\times I$ or to
$K\widetilde{\times}I$, and in either case $W\cap\partial M$ must consist of
annuli alone. In this case, all essential annuli in $M$ which are homotopic
into $W$ carry commensurable groups as their boundaries lie in $W\cap\partial
M$.
\end{remark}

\begin{proof}
As $X$ lies in $\mathcal{F}^{\prime}-\mathcal{F}$, it is a nontrivial almost
invariant subset of $G$ which is over an orientable $VPC(n+1)$ group $H$ and
is adapted to $\partial G$. By part 2) of Lemma
\ref{XcrossesYimpliesitcrossesY'strongly}, we may assume that $X$ crosses $Y$
and $Y^{\prime}$ strongly. As $X$ crosses $Y$ strongly, it follows that
$e(H,H\cap L)\geq2$. As $H$ is $VPC(n+1)$, it follows that $H\cap L$ is $VPCn$
and hence of finite index in $L$. Thus by replacing $L$ by a subgroup of
finite index, we can assume that it is a subgroup of $H$. Similarly we can
assume that $L^{\prime}$ is also a subgroup of $H$. If $L$ and $L^{\prime}$
are not commensurable, it follows that some element of $L^{\prime}$ is
hyperbolic with respect to $Y$ which implies that $Y^{\prime}$ crosses $Y$
strongly. Now Corollary 7.10 of \cite{SS2} implies that $L$ has small
commensuriser which contradicts the fact that $H$ commensurises $L$. This
contradiction shows that $L$ and $L^{\prime}$ must be commensurable as required.
\end{proof}

\begin{lemma}
\label{XcrossesYandX'impliesX'crossesY'}Let $(G,\partial G)$ be an orientable
$PD(n+2)$ pair, and let $\mathcal{F}$ and $\mathcal{F}^{\prime}$ be defined as
above. Let $X$ be an element of $\mathcal{F}^{\prime}-\mathcal{F}$ and let
$X^{\prime}$ be an element of $\mathcal{F}^{\prime}$ dual to a torus in
$(G,\partial G)$. Suppose that $X$ crosses $X^{\prime}$, and also crosses a
nontrivial almost invariant set $Y$ over a $VPCn$ group $L$. Then $X^{\prime}$
crosses a nontrivial almost invariant set $Y^{\prime}$ over a $VPCn$ group
$L^{\prime}$ commensurable with $L$. In particular, $X^{\prime}$ lies in
$\mathcal{F}^{\prime}-\mathcal{F}$.
\end{lemma}

\begin{remark}
If $n=1$ and $M$ is an orientable Haken $3$--manifold, the corresponding
result is clear. For $X$ corresponds to a torus which crosses an annulus in
$M$, and so must be homotopic into a component $W$ of $T(M)$ which meets
$\partial M$. As $X$ crosses $X^{\prime}$, the corresponding tori cross, so
both must be homotopic into $W$. In addition, each of these tori cannot be
homotopic to a boundary component of $W$. Now any such torus in $W$ must cross
an annulus in $W$, so that $X^{\prime}$ lies in $\mathcal{F}^{\prime
}-\mathcal{F}$.
\end{remark}

\begin{proof}
The hypotheses imply that $X$ and $X^{\prime}$ are nontrivial almost invariant
subsets of $G$ which are over orientable $VPC(n+1)$ subgroups $H$ and
$H^{\prime}$ respectively, and are adapted to $\partial G$. As they are
adapted to $\partial G$, they can be extended to almost invariant sets in $DG$
over the same groups $H$ and $H^{\prime}$. As $DG$ is $PD(n+2)$, and $H$ and
$H^{\prime}$ are each $PD(n+1)$, it follows that $H$ and $H^{\prime}$ each has
two coends in $DG$, so that these extended almost invariant sets must cross
strongly if at all. As $X$ and $X^{\prime}$ cross, it follows that they must
cross each other strongly. Let $L^{\prime}$ denote the intersection $H\cap
H^{\prime}$, which is $VPCn$. By part 2) of Lemma
\ref{XcrossesYimpliesitcrossesY'strongly}, we may assume that $X$ crosses $Y$
strongly, and the proof of Lemma \ref{firstlemma} shows that we can assume
that $L$ is a subgroup of $H$.

If $L$ and $L^{\prime}$ are not commensurable subgroups of $H$, then some
element of $L^{\prime}$ is hyperbolic with respect to $Y$. This implies that
$X^{\prime}$ crosses $Y$ strongly, thus proving the lemma in this case.

Now suppose that $L$ and $L^{\prime}$ are commensurable. Lemma 13.1 of
\cite{SS2} implies that there are finite index subgroups $H_{1}$ of $H$ and
$L_{1}$ of $L$ such that $L_{1}$ is normal in $H_{1}$ with infinite cyclic
quotient. We let $\Gamma$ denote the Cayley graph of $G$ with respect to some
finite generating set, and consider the action of $L_{1}\backslash H_{1}$ on
the graph $L_{1}\backslash\Gamma$. Let $h$ be an element of $H_{1}-L_{1}$. As
$L_{1}\backslash\delta Y$ is finite, there is a finite connected subcomplex
$C$ of $L_{1}\backslash\Gamma$ which contains $L_{1}\backslash\delta Y$. Thus
for all but finitely many elements $g$ of $L_{1}\backslash H_{1}$, we have
$gC\cap C$ empty, so that $g(L_{1}\backslash Y)$ and $L_{1}\backslash Y$ are
nested. As $X$ crosses $Y$ strongly, the intersection $Y\cap H_{1}$ is a
nontrivial $L_{1}$--almost invariant subset of $H_{1}.$ As $L_{1}\backslash
H_{1}$ is infinite cyclic, there is a power $h^{m}$ of $h$ such that $Y\cap
H_{1}\subset h^{m}(Y\cap H_{1})$. Thus by replacing $h$ by a suitable power if
needed, we can arrange that $Y\subset hY$. Also $Y$ cannot be equivalent to
$hY$. For if this happens, then $L_{1}$ has infinite index in $\{g\in
H_{1}:gY$ is equivalent to $Y\}$, and we can apply Theorem 5.8 from
\cite{Scott-Wall:Topological} to the action of $L_{1}\backslash H_{1}$ on
$L_{1}\backslash\Gamma$. The proof of this result implies that $H_{1}$ must
have finite index in $G$, which is impossible as $X$ is a nontrivial almost
invariant subset of $G$ over $H$.

As $h$ normalises $L_{1}$, each translate $h^{n}Y$ of $Y$ is also $L_{1}%
$--almost invariant. If $X^{\prime}$ crosses $Y$ or any translate $h^{n}Y$, we
will have proved the lemma. Thus we can suppose that $X^{\prime}$ does not
cross $h^{n}Y$, for any integer $n$. By replacing $X^{\prime}$ and $Y$ by
their complements if needed, we can arrange that $X^{\prime}\leq Y$. (If we
replace $Y$ by its complement, we simultaneously replace $h$ by its inverse,
in order to preserve the inclusion $Y\subset hY$.) Now this implies that
$X^{\prime}\leq h^{k}Y$ for all $k\geq0$. We claim that the inequality
$X^{\prime}\leq h^{k}Y$ cannot hold for every integer $k$. To see this pick a
finite generating set for $G$ and let $C$ denote the corresponding Cayley
graph for $G$. Now recall that given a pair of nontrivial almost invariant
subsets $U\ $and $W$ of $G$, each over a finitely generated subgroup of $G$,
there is an integer $d$ such that if $U\leq gW$ then $U$ is contained in the
$d$--neighbourhood of $gW$, where distances are measured in $C$. As the
intersection of all the $h^{k}Y$ is empty, the inequalities $X^{\prime}\leq
h^{k}Y$, for every $k$, would imply that $X^{\prime}$ is empty, which is a
contradiction. This completes the proof of the claim. (Alternatively the claim
follows immediately using the fact, proved in \cite{NibloSageevScottSwarup},
that $X^{\prime}$ and $Y$ are equivalent to almost invariant sets in very good
position, which implies that the partial orders on $E$ induced by inclusion
and by $\leq$ are the same.) Thus there must be a least integer $k$ such that
$X^{\prime}\nleq h^{k}Y$. Now by replacing $Y$ by $h^{k}Y$ if needed, we can
suppose that $k=0$. Thus we have $X^{\prime}\leq Y$ and $X^{\prime}\nleq
h^{-1}Y$.

Recall that $L$ and $L^{\prime}=H\cap H^{\prime}$ are commensurable. This
implies that $h$ and its powers do not lie in $L^{\prime}$. Let $h_{+}$ denote
$\{h^{k}:k\geq0\}$, let $h_{-}$ denote $\{h^{k}:k\leq0\}$, and let $h_{\pm}$
denote $h_{+}\cup h_{-}$. As $X\ $and $X^{\prime}$ cross strongly, it follows
that $h_{\pm}\cap X^{\prime}$ and $h_{\pm}\cap X^{\prime\ast}$ each contain
points which are arbitrarily far from $\delta X^{\prime}$. In particular, both
sets are infinite. The fact that $X^{\prime}\leq Y\subset hY\subset
h^{2}Y\subset\ldots$ implies that $h_{+}\cap X^{\prime}$ is finite, so that
$h_{-}\cap X^{\prime}$ must be infinite.

As we are assuming that $X^{\prime}$ does not cross $h^{n}Y$, for any integer
$n$, we know that $X^{\prime}$ does not cross $h^{-1}Y$. Also recall that
$X^{\prime}\nleq h^{-1}Y$, so that we must have one of the inequalities
$X^{\prime}\leq h^{-1}Y^{\ast}$, $h^{-1}Y^{\ast}\leq X^{\prime}$ or
$h^{-1}Y\leq X^{\prime}$. If $X^{\prime}\leq h^{-1}Y^{\ast}$, then $X^{\prime
}\leq h^{-k}Y^{\ast}$, for every $k\geq0$, so that $h_{-}\cap X^{\prime}$ is
finite, which is a contradiction. If $h^{-1}Y^{\ast}\leq X^{\prime}$, the
inequality $X^{\prime}\leq Y$ implies that $h^{-1}Y^{\ast}\leq Y$, which is
impossible as $h^{-1}Y\subset Y$. Thus we must have $h^{-1}Y\leq X^{\prime}$.
As $X^{\prime}\leq Y$, we know $Y^{\ast}\leq X^{\prime\ast}$. Hence
$X^{\prime}$ crosses the nontrivial $L_{1}$--almost invariant set $Y^{\prime
}=Y^{\ast}\cup h^{-1}Y$, completing the proof of the lemma.
\end{proof}

\begin{remark}
Example 2.13 of \cite{SS4} shows that Lemmas \ref{firstlemma} and
\ref{XcrossesYandX'impliesX'crossesY'} are not valid if $X$ is not adapted to
$\partial G$. In that example, $G$ is the fundamental group of an orientable
Haken $3$--manifold $M$ constructed by gluing two Seifert fibre spaces $M_{1}$
and $M_{2}$ along a boundary torus $T$, so that the Seifert fibrations do not
match. Thus the given decomposition of $M$ is essentially its
JSJ\ decomposition. Denote $\pi_{1}(T)$ by $H$, and let $\sigma^{\prime}$
denote the splitting of $G$ over $H$ determined by $T$. Note that
$\sigma^{\prime}$ is adapted to $\partial G$. Let $X^{\prime}$ denote a
$H$--almost invariant subset of $G$ associated to $\sigma^{\prime}$.

In this example, each $M_{i}$ has at least one boundary component other than
$T$, and so admits essential annuli disjoint from $T$. In particular $M$
itself admits essential annuli, so that $G$ does possess nontrivial almost
invariant subsets over $VPC1$ subgroups. We described a splitting $\sigma$ of
$G$ over $H$ which crosses $\sigma^{\prime}$, and is not adapted to $\partial
G$. Let $X$ denote a $H$--almost invariant subset of $G$ associated to this
splitting. Thus $X$ crosses $X^{\prime}$. We also showed that $X$ crosses
annuli in $M_{1}$ and in $M_{2}$. As annuli in $M_{1}$ and $M_{2}$ carry
incommensurable subgroups of $G$, Lemma \ref{firstlemma} fails. It follows
from results we proved in \cite{SS4} that $X^{\prime}$ is $1$--canonical, i.e.
$X^{\prime}$ does not cross any nontrivial almost invariant subset of $G$ over
a $VPC1$ subgroup. Thus Lemma \ref{XcrossesYandX'impliesX'crossesY'} also fails.
\end{remark}

Now we can start on the proof that $\mathcal{F}$ and $\mathcal{F}^{\prime}$
have the same algebraic regular neighbourhoods. We want to show that the
reduced algebraic regular neighbourhood $\Gamma(\mathcal{F}^{\prime})$ of
$\mathcal{F}^{\prime}$ is isomorphic to the reduced algebraic regular
neighbourhood of $\mathcal{F}$, which is $\Gamma_{n,n+1}(G)$. It seems natural
to approach this by proving that each element of $\mathcal{F}^{\prime}$ is
enclosed by some $V_{0}$--vertex of $\Gamma_{n,n+1}$, but this seems hard to
do directly. Instead we will consider our construction of unreduced algebraic
regular neighbourhoods from \cite{SS2}. We will show that $\mathcal{F}$ and
$\mathcal{F}^{\prime}$ have the same unreduced algebraic regular
neighbourhoods which will immediately imply that they also have the same
reduced algebraic regular neighbourhoods. Briefly our construction in
\cite{SS2} goes as follows. Let $E$ denote a $G$--invariant family of
nontrivial almost invariant subsets of $G$, and assume that the elements of
$E$ are in good position. Thus we have a $G$--invariant partial order $\leq$
on $E$, essentially given by almost inclusion. A cross-connected component
(CCC) of $E$ is an equivalence class of the equivalence relation generated by
crossing. We let $P$ denote the collection of CCC's of $E$. We showed that
there is a natural idea of betweenness on the elements of $P$ which gives it a
pretree structure. If this pretree is discrete, we can construct a bipartite
$G$--tree $T$ with $P$ as its $V_{0}$--vertices, and the graph of groups
$\Gamma=G\backslash T$ is the unreduced algebraic regular neighbourhood of the
family $E$.

\begin{remark}
\label{canpretreesberemoved?}As $G$ is finitely generated, it may be that we
could simplify this approach using pretrees by one which uses cubings and very
good position, as discussed in section 9 of \cite{GSS}, but we have not
attempted to do this.
\end{remark}

Let $P(\mathcal{F})$ denote the collection of all CCC's of elements of
$\mathcal{F}$, and let $P(\mathcal{F}^{\prime})$ denote the collection of all
CCC's of elements of $\mathcal{F}^{\prime}$. Thus $P(\mathcal{F})$ and
$P(\mathcal{F}^{\prime})$ have natural pretree structures. As we know that
$\mathcal{F}$ has an unreduced algebraic regular neighbourhood, we know that
$P(\mathcal{F})$ is a discrete pretree. Recall that Corollary
\ref{FiscontainedinF'} shows that $\mathcal{F}$ is contained in $\mathcal{F}%
^{\prime}$. This inclusion yields a natural map $\varphi:P(\mathcal{F}%
)\rightarrow P(\mathcal{F}^{\prime})$, which is clearly $G$--equivariant. We
will show that $\varphi$ is a $G$--equivariant bijection such that $\varphi$
and $\varphi^{-1}$ preserve the pretree idea of betweenness. Thus $\varphi$
induces a $G$--equivariant isomorphism of the pretrees. This implies that the
pretree determined by $\mathcal{F}^{\prime}$ is discrete, so that
$\mathcal{F}^{\prime}$ has an unreduced algebraic regular neighbourhood in $G$
which is isomorphic to the unreduced algebraic regular neighbourhood of
$\mathcal{F}$. Thus the reduced algebraic regular neighbourhood of
$\mathcal{F}^{\prime}$ in $G$ is isomorphic to $\Gamma_{n,n+1}$ as required.

Recall that $\mathcal{F}^{\prime}-\mathcal{F}$ consists of all those
nontrivial almost invariant subsets of $G$ which are over orientable
$VPC(n+1)$ subgroups and are adapted to {$\partial G$}, but are not
$n$--canonical. In particular every element of $\mathcal{F}^{\prime
}-\mathcal{F}$ crosses some element of $\mathcal{F}$. It follows immediately
that $\varphi:P(\mathcal{F})\rightarrow P(\mathcal{F}^{\prime})$ is
surjective. We will show in the proof of Theorem \ref{FandF'havesameregnbhds}
that it is injective. Note also that $\varphi$ clearly preserves betweenness
in the sense that if three CCC's $A$, $B$ and $C$ of $P(\mathcal{F})$ satisfy
$ABC$, i.e. $B$ is between $A$ and $C$, and if the image vertices $A^{\prime}%
$, $B^{\prime}$ and $C^{\prime}$ of $P(\mathcal{F}^{\prime})$ are distinct,
then $A^{\prime}B^{\prime}C^{\prime}$. But it is not automatic that
$\varphi^{-1}$ preserves betweenness, and this needs proof in the special case
under consideration. We give here a simple example which makes clear that this
is a real and subtle problem.

\begin{example}
\label{examplerebetweenness} Let $M$ denote the compact surface obtained from
the $2$--disc by removing the interiors of two subdiscs. Thus $G=\pi_{1}(M)$
is free of rank $2$. Let $\lambda$ and $\mu$ denote two simple essential arcs
properly embedded in $M$ and intersecting transversely in a single point which
cannot be removed by an isotopy of $\lambda$ and $\mu$. Let $N(\lambda)$ and
$N(\lambda\cup\mu)$ denote regular neighbourhoods of $\lambda$ and
$\lambda\cup\mu$ respectively, and let $\Gamma(\lambda)$ and $\Gamma
(\lambda\cup\mu)$ denote the bipartite graphs of groups decompositions of $G$
which are determined by the frontiers in $M$ of $N(\lambda)$ and
$N(\lambda\cup\mu)$ respectively. Then $\Gamma(\lambda)$ and $\Gamma
(\lambda\cup\mu)$ will be non-isomorphic graphs of groups. Now $\lambda$ and
$\mu$ determine almost invariant subsets $X(\lambda)$ and $X(\mu)$ of $G$
obtained by lifting them to arcs $\overline{\lambda}$ and $\overline{\mu}$ in
the universal cover of $M$ and choosing one side of the lift. We can regard
$\Gamma(\lambda)$ as the algebraic regular neighbourhood of $X(\lambda)$ and
can regard $\Gamma(\lambda\cup\mu)$ as the algebraic regular neighbourhood of
the family $\{X(\lambda),X(\mu)\}$. As $\lambda$ is a simple arc on $M$, the
translates of $\overline{\lambda}$ are disjoint, so that the family
$E(\lambda)$ of translates of $X(\lambda)$ and its complement is nested. Thus
each CCC of $E(\lambda)$ consists of a single translate of $X(\lambda)$. If we
let $E(\lambda,\mu)$ denote the family of translates of $X(\lambda)$ and
$X(\mu)$ and their complements, then each CCC of $E(\lambda,\mu)$ will consist
of a single translate of $X(\lambda)$ and a single translate of $X(\mu)$. In
each case, the stabiliser of each CCC is trivial. Thus if $P(\lambda)$ and
$P(\lambda,\mu)$ denote the families of CCC's of $E(\lambda)$ and
$E(\lambda,\mu)$, the natural equivariant map $\varphi$ from $P(\lambda)$ to
$P(\lambda,\mu)$ is a bijection. Further it is clear that $\varphi$ preserves
betweenness. But this cannot be true for the inverse map $\varphi^{-1}$. For
the fact that $\Gamma(\lambda)$ and $\Gamma(\lambda\cup\mu)$ are not
isomorphic implies that the pretree structures on $P(\lambda)$ and
$P(\lambda,\mu)$ must be different.
\end{example}

The following technical result is the key to handling this difficulty.

\begin{lemma}
\label{FandF'havesamebetweenness}Let $(G,\partial G)$ be an orientable
$PD(n+2)$ pair, let $\mathcal{F}$ and $\mathcal{F}^{\prime}$ be defined as
above, and let $\varphi:P(\mathcal{F})\rightarrow P(\mathcal{F}^{\prime})$ be
the natural map. Let $X$ be an element of $\mathcal{F}^{\prime}-\mathcal{F}$
which crosses an almost invariant set $Y$ over a $VPCn$ group $K$, and let $B$
denote the CCC of $P(\mathcal{F})$ which contains $Y$ so that the CCC
$\varphi(B)$ of $P(\mathcal{F}^{\prime})$ contains $X$ and $Y$. Then

\begin{enumerate}
\item $Y$ cannot be isolated in $\mathcal{F}$.

\item If $U$ and $V$ belong to CCC's $A$ and $C$ of $P(\mathcal{F})$ such that
$A\neq B\neq C$, and if $U\leq X\leq V$, then there is an element $Z$ of $B$
such that $U\leq Z\leq V$.
\end{enumerate}
\end{lemma}

\begin{remark}
In the $3$--manifold setting, so that $n=1$, part 1) is easy to see. For let
$M$ be a Haken manifold with incompressible boundary, and suppose that an
essential torus $T$ in $M$ crosses an essential annulus $A$. Then both must be
homotopic into a component $W$ of $JSJ(M)$ which is a Seifert fibre space
which meets $\partial M$. As such a component $W$ is filled by essential
annuli, either $A$ crosses some essential annulus or it is homotopic into a
frontier component of $W$. The second case would imply that $A$ crosses no
essential torus of $M$, which is a contradiction, so that $A$ must cross some
essential annulus in $M$,$\ $as required.
\end{remark}

\begin{proof}
As $X$ lies in $\mathcal{F}^{\prime}$, it is a nontrivial almost invariant
subset of $G$ which is over an orientable $VPC(n+1)$ group $H$ and is adapted
to $\partial G$.

1) If $Y$ is isolated in $\mathcal{F}$, then part 4) of Proposition
\ref{annuligenerate} implies that $Y\ $is dual to an annulus. Now part 1) of
Lemma \ref{XcrossesYimpliesitcrossesY'strongly} implies that $X$ crosses $Y$
strongly. Thus $e(H,H\cap K)\geq2$, so that $H\cap K$ must be $VPCn$, and so
of finite index in $K$. Now Lemma 13.1 of \cite{SS2} shows that, by replacing
$H$ and $K$ by subgroups of finite index, we can suppose that $K\ $is a
subgroup of $H$ and is normal in $H$ with infinite cyclic quotient. Let $h$ be
an element of $H-K$. As $K$ is normal in $H$, the translate $hY$ of $Y$ is
also $K$--almost invariant. As $Y$ is assumed to be isolated in $\mathcal{F}$,
no translate $h^{n}Y$ of $Y$ can cross $Y$. As in the proof of Lemma
\ref{XcrossesYandX'impliesX'crossesY'}, by replacing $h$ by a suitable power
if needed, we can suppose that $Y\subset hY$. Also as in that proof, $Y$
cannot be equivalent to $hY$. Now it follows that $Y$ crosses the nontrivial
$K$--almost invariant set $hY^{\ast}\cup h^{-1}Y$, so that $Y$ cannot be
isolated. This contradiction completes the proof of part 1).

2) If $X$ crosses $Y$ weakly, part 2) of Lemma
\ref{XcrossesYimpliesitcrossesY'strongly} shows that there is a $K^{\prime}%
$--almost invariant set $Y^{\prime}$ such that $K^{\prime}$ is commensurable
with $K$, and $Y^{\prime}$ is dual to an annulus. It also shows that $X$
crosses $Y^{\prime}$ strongly and hence that $H$ commensurises $K^{\prime}$.
In particular $K^{\prime}$, and hence $K$, has large commensuriser. Now
Proposition 8.6 of \cite{SS2} shows that all nontrivial almost invariant
subsets of $G$ over subgroups commensurable with $K$ belong to a single CCC
apart from a finite number which are isolated. (Note that this proposition as
stated only applies to the case when $K$ is $VPC1$, but essentially the same
proof works in general. This is discussed in section 14 of \cite{GSS}.) As $X$
crosses $Y$ and $Y^{\prime}$, part 1) tells us that $Y$ and $Y^{\prime}$
cannot be isolated, so it follows that $Y$ and $Y^{\prime}$ belong to the same
CCC $B$.

It follows from the preceding paragraph that, by replacing $Y$ by $Y^{\prime}$
if needed, we can assume that $Y$ is dual to an annulus, so that $X$ crosses
$Y$ strongly.

As in part 1), we can further assume that $K\ $is a subgroup of $H$ and is
normal in $H$ with infinite cyclic quotient. In particular, as $H$
commensurises $K$, it follows that $H$ must preserve the CCC $B$, so that for
every $g\in H$, we have $gY$ belongs to $B$. As no element of the CCC $A$ can
cross any element of $B$, it follows that, for every $g\in H$, we have $U\leq
gY$ or $U\leq gY^{\ast}$. We will write $U\leq gY^{(\ast)}$ to indicate that
one of these two inequalities holds. Similarly, for every $g\in H$, we have
$gY^{(\ast)}\leq V$. By replacing $Y$ by its complement if needed, we can
suppose that $U\leq Y$.

Now we consider the action of $K\backslash H$ on $K\backslash\Gamma$, and let
$h$ be an element of $H-K$. As in the proof of Lemma
\ref{XcrossesYandX'impliesX'crossesY'}, by replacing $h$ by a suitable power,
we can suppose that $Y\subset hY$. Also as in that proof, $Y$ cannot be
equivalent to $h^{n}Y$, for any non-zero value of $n$. As $U\leq Y$, the
inclusion $Y\subset hY$ implies that $U\leq h^{k}Y$ for all $k\geq0$. As in
the proof of Lemma \ref{XcrossesYandX'impliesX'crossesY'}, the inequality
$U\leq h^{k}Y$ cannot hold for every integer $k$, and by replacing $Y$ by a
suitable translate $h^{k}Y$ if needed, we can suppose that $U\leq Y$ and
$U\nleq h^{-1}Y$. Thus we must have $U\leq h^{-1}Y^{\ast}\subset h^{-2}%
Y^{\ast}$, and hence $U\leq Y\cap h^{-2}Y^{\ast}$. Let $R$ denote the
intersection $Y\cap h^{-2}Y^{\ast}$. As $h^{-2}Y^{\ast}$ is $K$--almost
invariant, it follows that $R$ is also $K$--almost invariant. Now let $Z$
denote the intersection $R\cap X$. We claim that $Z$ is a $K$--almost
invariant subset of $G$. Certainly $KZ=Z$ as each of $R$ and $X$ is invariant
under the left action of $K$. The coboundary $\delta Z$ of $Z$ is the union of
subsets of $\delta R$ and $\delta X$. As $\delta R$ is $K$--finite, it remains
to show that $\delta Z\cap\delta X$ is $K$--finite. Now $\delta X$ is
$H$--finite, and contains the union of the translates $h^{n}(\delta
Z\cap\delta X)$, for every integer $n$. As the translates of $R$ by powers of
$h^{2}$ are all disjoint, it follows that the translates of $\delta
Z\cap\delta X$ by powers of $h^{2}$ are all disjoint. Now it follows that
$\delta Z\cap\delta X$ must be $K$--finite, so that $\delta Z$ itself is
$K$--finite as required. Recall that we have the inequalities $U\leq R$ and
$U\leq X$, so that we also have $U\leq Z$. Finally we claim that as $Z=Y\cap
h^{-2}Y^{\ast}\cap X$, it crosses $hZ=hY\cap h^{-1}Y^{\ast}\cap X$. This is
because none of the four corners of the pair $(Z,hZ)$ is small, which follows
from the inclusions $U\subset Z\cap hZ$, $hU\subset Z^{\ast}\cap hZ$,
$h^{-1}U\subset Z\cap hZ^{\ast}$ and $X^{\ast}\subset Z^{\ast}\cap hZ^{\ast}$.
Thus $Z$ is a nontrivial and non-isolated almost invariant subset of $G$ over
$K$ and hence must belong to the CCC $B$. Now the inequalities $U\leq Z\leq
X\leq V$ complete the proof of part 2) of the lemma.
\end{proof}

\begin{theorem}
\label{FandF'havesameregnbhds}Let $(G,\partial G)$ be an orientable $PD(n+2)$
pair. Let $\mathcal{F}$ denote the family of equivalence classes of all
nontrivial almost invariant subsets of $G$ which are over $VPCn$ subgroups and
of equivalence classes of all $n$--canonical almost invariant subsets of $G$
which are over $VPC(n+1)$ subgroups. Thus the reduced algebraic regular
neighbourhood of $\mathcal{F}$ exists and is $\Gamma_{n,n+1}$. Let
${\mathcal{F}^{\prime}}$ denote the family of equivalence classes of all
nontrivial almost invariant subsets of $G$ which are over $VPCn$ subgroups
together with all nontrivial almost invariant subsets of $G$ which are over
$VPC(n+1)$ subgroups and are adapted to {$\partial G$}. Then $\mathcal{F}%
^{\prime}$ has a reduced algebraic regular neighbourhood in $G$, and this is
naturally isomorphic to $\Gamma_{n,n+1}$.
\end{theorem}

\begin{proof}
Recall our discussion immediately before Example \ref{examplerebetweenness}.
There is a natural $G$--equivariant surjection $\varphi:P(\mathcal{F}%
)\rightarrow P(\mathcal{F}^{\prime})$. We need to show that $\varphi$ is an
injection and that $P(\mathcal{F})$ and $P(\mathcal{F}^{\prime})$ have the
same pretree idea of betweenness.

Suppose that $\varphi$ is not injective. Then there must be elements $Y_{0}$
and $Y_{1}$ of $\mathcal{F}$, belonging to distinct CCC's $v_{0}$ and $v_{1}$
of $P(\mathcal{F})$, and a sequence $X_{1},\ldots,X_{n}$ of elements of
$\mathcal{F}^{\prime}-\mathcal{F}$ such that $X_{1}$ crosses $Y_{0}$, and
$X_{n}$ crosses $Y_{1}$, and $X_{i}$ crosses both $X_{i-1}$ and $X_{i+1}$, for
each $i$ such that $2\leq i\leq n-1$. As $Y_{0}\in\mathcal{F}$, and $X_{1}%
\in\mathcal{F}^{\prime}-\mathcal{F}$, Lemma
\ref{XcrossesYandX'impliesX'crossesY'} shows that $Y_{0}$ cannot be dual to a
torus. Similarly $Y_{1}$ cannot be dual to a torus. Thus $Y_{0}$ and $Y_{1}$
are over $VPCn$ subgroups of $G$.

Let $X_{i}$ be almost invariant over $H_{i}$, and $Y_{j}$ over $K_{j}$, and
note that each $H_{i}$ is orientable. We claim that $H_{1}\cap K_{0}$ has
finite index in $K_{0}$. If $X_{1}$ crosses $Y_{0}$ strongly, then
$e(H_{1},H_{1}\cap K_{0})\geq2$. As $H_{1}$ is $VPC(n+1)$, it follows that
$H_{1}\cap K_{0}$ must be $VPCn$ and so the claim follows in this case. If
$X_{1}$ crosses $Y_{0}$ weakly, then part 2) of Lemma
\ref{XcrossesYimpliesitcrossesY'strongly} tells us that there is a nontrivial
almost invariant set $Y^{\prime}$ over a group $K^{\prime}$ commensurable with
$K_{0}$ such that $X_{1}$ crosses $Y^{\prime}$ strongly. As before this
implies that $H_{1}\cap K^{\prime}$ must be $VPCn$, so it follows that
$H_{1}\cap K_{0}$ must also be $VPCn$, and the claim follows in this case also.

Now Lemma 13.1 of \cite{SS2} implies that a subgroup of finite index in
$H_{1}$ commensurises $H_{1}\cap K_{0}$ and hence commensurises $K_{0}$
itself, so that $K_{0}$ has large commensuriser in $G$. Proposition 8.6 of
\cite{SS2} shows that the collection of equivalence classes of all nontrivial
almost invariant subsets of $G\ $which are over a subgroup commensurable with
$K_{0}$ form a single CCC apart from finitely many isolated elements. Part 1)
of Lemma \ref{FandF'havesamebetweenness} shows that $Y_{0}$ cannot be
isolated. Thus the CCC $v_{0}$ of $P(\mathcal{F})$ which contains $Y_{0}$
contains all the nontrivial and non-isolated almost invariant subsets of
$G\ $which are over a subgroup commensurable with $K_{0}$. As $X_{1}$ crosses
an element of $v_{0}$, and also crosses $X_{2}$, Lemma
\ref{XcrossesYandX'impliesX'crossesY'} shows that $X_{2}$ must cross some
nontrivial almost invariant subset $Y_{2}$ of $G\ $which is over a subgroup
commensurable with $K_{0}$. Again part 1) of Lemma
\ref{FandF'havesamebetweenness} shows that $Y_{2}$ cannot be isolated. It
follows that $Y_{2}$ lies in the CCC $v_{0}$. Now, by induction, it follows
that each $X_{i}$ must cross some element of $v_{0}$. Thus $X_{n}$ crosses an
element $Z_{0}$ of $v_{0}$ and the element $Y_{1}$ of $v_{1}$. Lemma
\ref{firstlemma} shows that $Z_{0}$ and $Y_{1}$ must have commensurable
stabilisers, and Lemma \ref{FandF'havesamebetweenness} shows that $Y_{1}$
cannot be isolated. But this implies that $Y_{1}$ also lies in the CCC $v_{0}%
$, which contradicts our assumption. This completes the proof that $\varphi$
is injective and hence is a bijection.

Finally Lemma \ref{FandF'havesamebetweenness} shows that the betweenness
relations on the two pretrees $P(\mathcal{F}^{\prime})$ and $P(\mathcal{F}%
^{\prime})$ are the same, which completes the proof of the theorem.
\end{proof}

\section{Proof of the Main Theorem\label{proofofmaintheorem}}

In this section we use our work from previous sections to show how our main
result, Theorem \ref{mainresult}, follows from Theorem
\ref{JSJforVPCoftwolengths}. One of the things we need to prove when
$(G,\partial G)$ is an orientable $PD(n+2)$ pair is that all the edge
splittings of $\Gamma_{n,n+1}(G)$ are dual to essential annuli and tori. If an
edge of $\Gamma_{n,n+1}(G)$ is incident to a $V_{0}$--vertex $v$ which is
isolated or of $VPCk$--by--Fuchsian type, with $k$ equal to $n-1$ or $n$, it
is trivial that the edge group is $VPCn$ or $VPC(n+1)$. But if $v$ is of
commensuriser type, then we do not even know that the edge group is finitely
generated. However we will show in this section that when $(G,\partial G)$ is
an orientable $PD(n+2)$ pair the edge groups of $\Gamma_{n,n+1}(G)$ are all
$VPCn$ or $VPC(n+1)$. Assuming this, the following result shows that the edge
splittings of $\Gamma_{n,n+1}(G)$ are all dual to essential annuli and tori.

\begin{lemma}
\label{edgesplittingsaredualtoannuliortori}Let $(G,\partial G)$ be an
orientable $PD(n+2)$ pair such that $G$ is not $VPC$. Let $\mathcal{F}%
_{n,n+1}$ denote the family of equivalence classes of all nontrivial almost
invariant subsets of $G$ which are over a $VPCn$ subgroup, together with
equivalence classes of all $n$--canonical almost invariant subsets of $G$
which are over a $VPC(n+1)$ subgroup. Let $\Gamma_{n,n+1}$ denote the reduced
algebraic regular neighbourhood of $\mathcal{F}_{n,n+1}$ in $G$. (See Theorem
\ref{JSJforVPCoftwolengths}.) Finally let $e$ be an edge of $\Gamma_{n,n+1}$
with associated edge splitting $\sigma$.

\begin{enumerate}
\item If $G(e)$ is $VPCn$, then $\sigma$ is dual to an essential annulus in
$G$.

\item If $G(e)$ is $VPC(n+1)$, then $\sigma$ is dual to an essential torus in
$G$.
\end{enumerate}
\end{lemma}

\begin{proof}
Let $X$ denote the $G(e)$--almost invariant subset of $G$ associated to the
edge splitting $\sigma$. As $\sigma$ is an edge splitting of $\Gamma_{n,n+1}$,
it follows that $X$ crosses no element of $\mathcal{F}_{n,n+1}$.

1) If $G(e)$ is $VPCn$, then $X$ is automatically an element of $\mathcal{F}%
_{n,n+1}$, and hence is an isolated element. Now part 5) of Proposition
\ref{annuligenerate} implies that $X$ is dual to an essential annulus in $G$,
as required.

2) If $G(e)$ is $VPC(n+1)$, the fact that $X$ crosses no nontrivial almost
invariant subset of $G$ over a $VPCn$ subgroup implies that $X$ is
$n$--canonical and so is again an element of $\mathcal{F}_{n,n+1}$. Now
Corollary \ref{FiscontainedinF'} and Remark
\ref{adaptedtodGimpliesorientablegroup} show that $X$ is adapted to $\partial
G$ and that $G(e)$ is orientable, so that $X$ is dual to an essential torus in
$G$, as required.
\end{proof}

Recall that Theorem \ref{JSJforVPCoftwolengths} states that each $V_{0}%
$--vertex $v$ of $\Gamma_{n,n+1}$ satisfies one of the following conditions:

1) $v$ is isolated, and $G(v)$ is $VPC$ of length $n$ or $n+1$.

2) $v$ is of $VPCk$--by--Fuchsian type, where $k$ equals $n-1$ or $n$.

3) $v$ is of commensuriser type, so that $G(v)$ is the full commensuriser
$Comm_{G}(H)$ for some $VPC$ subgroup $H$ of length $n$ or $n+1$, such that
$e(G,H)\geq2$.

We will consider cases 2) and 3) in the lemmas which follow. We start with
case 2), where $v$ is of $VPCk$--by--Fuchsian type.

\begin{lemma}
\label{V0vertexofVPCk-by-Fuchsiantype}Let $(G,\partial G)$ be an orientable
$PD(n+2)$ pair such that $G$ is not $VPC$, and let $\Gamma_{n,n+1}$ denote the
reduced algebraic regular neighbourhood of $\mathcal{F}_{n,n+1}$ in $G$. Let
$v$ be a $V_{0}$--vertex of $\Gamma_{n,n+1}$ which is of $VPCk$--by--Fuchsian
type, where $k$ equals $n-1$ or $n$.

\begin{enumerate}
\item If $k=n$, then $v$ is of interior Seifert type (see Definition
\ref{defnofinteriorSeiferttype}).

\item If $k=n-1$, then $v$ is of $I$--bundle type (see Definition
\ref{defnofI-bundletype}).
\end{enumerate}
\end{lemma}

\begin{proof}
1) As $v$ is of $VPCn$--by--Fuchsian type, the groups associated to the edges
of $\Gamma_{n,n+1}$ incident to $v$ are all $VPC(n+1)$. Now Lemma
\ref{edgesplittingsaredualtoannuliortori} implies that the edge splittings of
$G$ associated to these edges are dual to essential tori. It follows that $v$
is of interior Seifert type, as claimed.

2) In this case $G(v)$ is $VPC(n-1)$--by--Fuchsian where the Fuchsian quotient
group $\Theta$ is not finite nor two-ended, and there is exactly one edge of
$\Gamma_{n,n+1}$ which is incident to $v$ for each peripheral subgroup $K$ of
$G(v)$ and this edge carries $K$. Let $E(v)$ denote the collection of these
edge groups. Note that each group in $E(v)$ is $VPCn$. Thus Lemma
\ref{edgesplittingsaredualtoannuliortori} implies that the edge splittings of
$G$ associated to the edges of $\Gamma_{n,n+1}$ which are incident to $v$ are
dual to essential annuli. Note also that it is possible that the family $E(v)$
is empty. In this case, $\Gamma_{n,n+1}$ consists of a single $V_{0}$--vertex
$v$, and $G=G(v)$ is $VPC(n-1)$--by--Fuchsian.

As in section \ref{essentialannuli}, we choose an aspherical space $M$ with
fundamental group $G$ and with aspherical subspaces corresponding to $\partial
G$ whose union is denoted $\partial M$. Recall that in order to prove that $v$
is of $I$--bundle type, we need to show that there are two distinct components
$\Sigma$ and $T$ of $\partial\widetilde{M}$ such that the induced action of
$G(v)$ on $\widetilde{M}$ preserves the union of $\Sigma$ and $T$, and for
each peripheral subgroup $K$ of $G(v)$, if $e_{K}$ denotes the edge of
$\Gamma$ which is incident to $v$ and carries $K$, then the edge splitting
associated to $e_{K}$ is given by the essential annulus $K_{\Sigma,T}$.

\medskip

\textbf{Finding }$\Sigma$\textbf{ and }$T$

\medskip

We claim that the pair $(G(v),E(v))$ is $PD(n+1)$. Let $L$ denote the
$VPC(n-1)$ normal subgroup of $G(v)$ with Fuchsian quotient $\Theta$. Let
$\partial\Theta$ denote the family of subgroups of $\Theta$ which are the
images of the groups in $E(v)$, so that the pair $(\Theta,\partial\Theta)$ is
a Fuchsian pair. If $\Theta$ is torsion free, then the pair $(\Theta
,\partial\Theta)$ is $PD2$. In this case, as $L$ is $PD(n-1)$, Theorem 7.3 of
\cite{B-E} implies that the pair $(G(v),E(v))$ is $PD(n+1)$. In general
$\Theta$ has a torsion free subgroup of finite index, and the pre-image in
$G(v)$ of this subgroup yields a $PD(n+1)$ pair of finite index in $G(v)$. As
$G(v)$ is torsion free, it follows that the pair $(G(v),E(v))$ itself is
$PD(n+1)$, as claimed. Note that this pair need not be orientable.

Consider any orientation preserving non-peripheral element $\alpha$ in
$\Theta$ of infinite order. As $\Theta$ is not finite nor two-ended, there is
another orientation preserving non-peripheral element $\beta$ in $\Theta$ of
infinite order such that $\alpha$ and $\beta$ have non-zero geometric
intersection number. Thus $\alpha$ and $\beta$ determine nontrivial almost
invariant subsets of $\Theta$ such that each is over an infinite cyclic
subgroup, each is adapted to $\partial\Theta$, and they cross strongly. The
pre-images in $G(v)$ of these subsets of $\Theta$ are almost invariant subsets
of $G(v)$ such that each is over a $VPCn$ subgroup, each is adapted to $E(v)$,
and they cross strongly. We denote these $VPCn$ subgroups of $G(v)$ by $C$ and
$D$. Now Lemma \ref{propertiesofadapted} shows us that we can extend these
sets to almost invariant subsets $P$ and $Q$ of $G$ over $C$ and $D$
respectively such that $P$ and $Q$ are enclosed by $v$. As the almost
invariant subsets of $G(v)$ cross strongly, $P\ $and $Q$ must also cross
strongly. Now Proposition 7.2 of \cite{SS2} shows that $C$ and $D$ must each
have two coends in $G$. Note that $C\cap D$ equals the $VPC(n-1)$ normal
subgroup $L$ of $G(v)$. Now we consider the covers $M_{L}$, $M_{C}$ and
$M_{D}$ of $M$ with fundamental groups $L$, $C$ and $D$ respectively.

\medskip

\textbf{Case:} $C$\textit{ and }$D$\textit{ are orientable.}

\medskip

As $G$ possesses a nontrivial $C$--almost invariant subset, Proposition
\ref{aisetimpliesannulus} tells us that $M_{C}$ must have at least two
boundary components which carry a subgroup of $C$ of finite index. As $C$ has
two coends in $G$, part 3) of Proposition \ref{annuligenerate} shows that the
number of such boundary components of $M_{C}$ and of each finite cover of
$M_{C}$ is exactly $2$. Hence each of these two boundary components of $M_{C}$
must carry $C$ itself. We denote these two boundary components by
$\partial_{1}M_{C}$ and $\partial_{2}M_{C}$. Note that it follows that $C$
preserves precisely two boundary components of the universal cover
$\widetilde{M}$ of $M$, and that all other $C$--orbits of boundary components
of $\widetilde{M}$ are infinite. Now as in section \ref{essentialannuli},
there is an untwisted annulus $A$ and a map $\theta:(A,\partial A)\rightarrow
(M_{C},\partial M_{C})$ which is an isomorphism on fundamental groups and
sends $\partial_{1}A$ into $\partial_{1}M_{C}$ and $\partial_{2}A$ into
$\partial_{2}M_{C}$. Similarly $D$ preserves precisely two boundary components
of the universal cover $\widetilde{M}$ of $M$, and all other $D$--orbits of
boundary components of $\widetilde{M}$ are infinite. We denote by
$\partial_{1}M_{D}$ and $\partial_{2}M_{D}$ the images in $M_{D}$ of the two
boundary components of $\widetilde{M}$ which are preserved by $D$. Also there
is an untwisted annulus $B$ and a map $\phi:(B,\partial B)\rightarrow
(M_{D},\partial M_{D})$ which is an isomorphism on fundamental groups and
sends $\partial_{1}B$ into $\partial_{1}M_{D}$ and $\partial_{2}B$ into
$\partial_{2}M_{D}$.

Now we consider the pre-images of these annuli in the common cover $M_{L}$ of
$M_{C}$ and $M_{D}$. Above $\theta$ we have a map $\theta_{L}$ into $M_{L}$ of
a two-ended cover $A_{L}$ of $A$, and $\theta_{L}$ maps $\partial A_{L}$ into
the two components of $\partial M_{L}$ which lie above $\partial_{1}M_{C}$ and
$\partial_{2}M_{C}$. Similarly above $\phi$ we have a map $\phi_{L}$ into
$M_{L}$ of a two-ended cover $B_{L}$ of $B$, and $\phi_{L}$ maps $\partial
B_{L}$ into the two components of $\partial M_{L}$ which lie above
$\partial_{1}M_{D}$ and $\partial_{2}M_{D}$. The fact that $P$ crosses $Q$
strongly implies that $\theta_{L}(\partial_{1}A_{L})$ must cross the image of
$\phi_{L}$ strongly, in the natural sense that $\theta_{L}(\partial_{1}A_{L})$
contains points arbitrarily far from the image of $\phi_{L}$ and on each side.
In particular $\theta_{L}(\partial_{1}A_{L})$ must meet the image of $\phi
_{L}$, and hence must meet $\phi_{L}(\partial_{1}B_{L})$ or $\phi_{L}%
(\partial_{2}B_{L})$. As the same argument applies to $\theta_{L}(\partial
_{2}A_{L})$, it follows that the two components of $\partial M_{L}$ which
contain $\theta_{L}(\partial_{1}A_{L})$ and $\theta_{L}(\partial_{2}A_{L})$
are the same as the two components of $\partial M_{L}$ which contain $\phi
_{L}(\partial_{1}B_{L})$ and $\phi_{L}(\partial_{2}B_{L})$. Hence $C$ and $D$
preserve the same two boundary components $\Sigma$ and $T$ of the universal
cover $\widetilde{M}$ of $M$.

\medskip

\textbf{Case:} \textit{one or both of }$C$\textit{ and }$D$\textit{ is not
orientable.}

\medskip

We apply the above arguments to their orientable subgroups of index $2$. We
will see that the action of the group generated by $C$ and $D$ on the boundary
components of $\widetilde{M}$ has one orbit with the two elements $\Sigma$ and
$T$, and all other orbits are infinite.

\medskip

Let $\Theta_{0}$ denote the subgroup of $\Theta$ generated by all orientation
preserving non-peripheral elements, and let $G(v)_{0}$ denote the pre-image of
$\Theta_{0}$ in $G(v)$. Thus $\Theta_{0}$ has index at most $2$ in $\Theta$,
so that $G(v)_{0}$ has index at most $2$ in $G(v)$. Consider the action of
$G(v)_{0}$ on the boundary components of $\widetilde{M}$. As any two
orientation preserving non-peripheral elements of $\Theta$ of infinite order
belong to a finite sequence of such elements each crossing the next, and as
such elements generate $\Theta_{0}$, it follows that $\Sigma$ and $T$ form one
or two orbits under the action of $G(v)_{0}$ and that all other orbits are
infinite. It follows that the same statement holds for the action of $G(v)$ on
the boundary components of $\widetilde{M}$. Thus $\Sigma$ and $T$ are the
required components of $\partial\widetilde{M}$.

If the family $E(v)$ of edge splittings associated to edges of $\Gamma
_{n,n+1}$ incident to $v$ is empty, we have shown that $v$ is of $I$--bundle
type. In this case, either $\partial G$ consists of two copies of $G$, or
$\partial G$ consists of a single group which has index $2$ in $G$. It remains
to deal with the case when $E(v)$ is non-empty.

\medskip

\textbf{The edge splittings of }$v$

\medskip

Let $K$ denote a peripheral subgroup of $G(v)$, let $e_{K}$ denote the edge of
$\Gamma_{n,n+1}$ which is incident to $v$ and carries $K$, and let $\theta
_{K}:(A_{K},\partial A_{K})\rightarrow(M,\partial M)$ be an essential annulus
which is dual to the edge splitting associated to $e_{K}$. We need to show
that the annulus $\theta_{K}$ equals $K_{\Sigma,T}$ which was defined prior to
Lemma \ref{swarup}. In general there can be many distinct annuli all carrying
the same group, so there is something to be proved.

Let $S$ denote the stabiliser of $\Sigma$, so that $S$ is one of the groups in
$\partial G$. We first consider the case when no element of $G(v)$
interchanges $\Sigma$ and $T$, so that $G(v)$ is contained in $S$. As $v$ is a
vertex of $\Gamma_{n,n+1}$ with associated group $G(v)$, there is a natural
induced decomposition of $S$ as the fundamental group of a graph of groups
$\Gamma(S)$, with a vertex $V$ with associated group $S\cap G(v)=G(v)$.
Further the stars of $v$ in $\Gamma_{n,n+1}$ and of $V$ in $\Gamma(S)$ are
isomorphic, i.e. there is a bijection between the edges of $\Gamma_{n,n+1}$
incident to $v$ and the edges of $\Gamma(S)$ incident to $V$, and
corresponding edges have the same associated groups. Note that $\Gamma(S)$ is
probably not minimal, and may well be infinite. However, as $S$ is finitely
generated, there is a finite minimal subgraph $\Gamma_{\mu}(S)$ of $\Gamma(S)$
which carries $S$.

Recall from the start of this proof that the pair $(G(v),E(v))$ is $PD(n+1)$.
In particular, as we are now assuming that $E(v)$ is non-empty, $G(v)$ itself
is not $PD(n+1)$. Also each group in $E(v)$ is $VPCn$. As $S$ is $PD(n+1)$,
and $G(v)$ is not $PD(n+1)$ nor $VPCn$, it follows that the minimal graph
$\Gamma_{\mu}(S)$ must contain $V$, and at least one edge incident to $V$. As
$S$ is $PD(n+1)$, Theorem 8.1 of \cite{B-E} tells us that the pair
$(G(V),E_{\mu}(V))$ is $PD(n+1)$, where $E_{\mu}(V)$ denotes the family of
subgroups of $G(V)$ associated to edges of $\Gamma_{\mu}(S)$ incident to $V$.
As the pair $(G(v),E(v))$ is $PD(n+1)$, and $G(V)=G(v)$ and $E_{\mu}(V)$ is a
subfamily of $E(v)$, it follows that the families $E_{\mu}(V)$ and $E(v)$ must
be equal. In particular every edge of $\Gamma(S)$ which is incident to $V$
must also be an edge of $\Gamma_{\mu}(S)$. Recall that $K$ is a group in
$E(v)$, that $e_{K}$ denotes the edge of $\Gamma_{n,n+1}$ which is incident to
$v$ and carries $K$, and that $\theta_{K}:(A_{K},\partial A_{K})\rightarrow
(M,\partial M)$ denotes an essential annulus which is dual to the edge
splitting associated to $e_{K}$. Let $\varphi_{K}$ denote the lift of
$\theta_{K}$ into $M_{G(v)}$. As the edge of $\Gamma(S)$ corresponding to
$e_{K}$ is an edge of $\Gamma_{\mu}(S)$, it follows that the splitting of $G$
over $K$ associated to $e_{K}$ induces a splitting of $S$ over $K$. This
implies that $\varphi_{K}$ must have a boundary component on the image of
$\Sigma$ in $M_{G(v)}$. The same argument applied to the stabiliser of $T$
shows that $\varphi_{K}$ must also have a boundary component on the image of
$T$ in $M_{G(v)}$. Hence the annulus $\theta_{K}$ must be $K_{\Sigma,T}$,
which completes the proof of the lemma, on the assumption that no element of
$G(v)$ interchanges $\Sigma$ and $T$.

If there are elements of $G(v)$ which interchange $\Sigma$ and $T$, we let $S$
denote the stabiliser of $\Sigma$, and let $G(v)_{0}$ denote $S\cap G(v)$.
Note that in this case the images of $\Sigma$ and $T$ in $M_{G(v)}$ are the
same. Then we make essentially the same argument as in the preceding
paragraph. This time the vertex $V$ of $\Gamma(S)$ has associated group
$G(v)_{0}$. Given the edge $e_{K}$ of the star of $v$ in $\Gamma_{n,n+1}$
which has associated group $K$, the star of $V$ in $\Gamma(S)$ has either one
corresponding edge with associated group $K$ or it has two corresponding edges
each with associated group $K_{0}$, where $K_{0}=K\cap G(v)_{0}$. As the pair
$(G(v),E(v))$ is $PD(n+1)$, it follows that the pair $(G(V),E(V))$ is also
$PD(n+1)$. Hence as above, the minimal subgraph $\Gamma_{\mu}(S)$ of
$\Gamma(S)$ which carries $S$ must contain $V$ and each edge of $\Gamma(S)$
which is incident to $V$. Thus, as before, the lift $\varphi_{K}$ of the
annulus $\theta_{K}$ into $M_{G(v)}$ must have all of its boundary on the
image of $\Sigma$. Note that there are two cases here, depending on whether
the annulus $A_{K}$ is twisted or untwisted. In either case it follows that
the annulus $\theta_{K}$ must be $K_{\Sigma,T}$, which completes the proof of
the lemma.
\end{proof}

Now we consider the case when $v$ is a $V_{0}$--vertex of $\Gamma_{n,n+1}$
which is of commensuriser type. This is the most difficult case, and we will
need most of our previous work in this paper. Note that in Theorem
\ref{JSJforVPCoftwolengths}, the vertex group $G(v)=Comm_{G}(H)$ need not even
be finitely generated. Before we start, here is a preliminary result.

\begin{lemma}
\label{commensuriservertexenclosescrossingannuli}Let $(G,\partial G)$ be an
orientable $PD(n+2)$ pair such that $G$ is not $VPC$, and let $\Gamma_{n,n+1}$
denote the reduced algebraic regular neighbourhood of $\mathcal{F}_{n,n+1}$ in
$G$. Let $v$ be a $V_{0}$--vertex of $\Gamma_{n,n+1}$ of commensuriser type
such that $G(v)$ is the full commensuriser $Comm_{G}(H)$ for some $VPCn$
subgroup $H$ of $G$ with $e(G,H)\geq2$.

Then $v$ encloses two almost invariant subsets $X\ $and $X^{\prime}$ of $G$,
each over a subgroup of $H$ of finite index, and each dual to an annulus, such
that $X$ and $X^{\prime}$ cross.
\end{lemma}

\begin{proof}
As $\Gamma_{n,n+1}$ is the reduced algebraic regular neighbourhood of
$\mathcal{F}_{n,n+1}$ in $G$, any $V_{0}$--vertex arises from a
cross-connected component (CCC) of $\mathcal{F}_{n,n+1}$. As $v$ is of
commensuriser type, all the crossings in this CCC must be weak and all the
almost invariant sets in this CCC are over groups commensurable with $H$. As
in section \ref{essentialannuli}, we consider an aspherical space $M$ with
fundamental group $G$ and with aspherical subspaces corresponding to $\partial
G$ whose union is denoted $\partial M$. As $v$ is of commensuriser type, the
number of coends of $H$ in $G$ must be at least $4$. Now part 4) of
Proposition \ref{annuligenerate} tells us $H$ has a subgroup $K$ of finite
index such that $\partial M_{K}$ has $4$ (or more) components each of which
carries $K$. Now Lemma \ref{thereexisttwocrossingannuli} shows that we can
find almost invariant subsets $X\ $and $X^{\prime}$ of $G$, each over $K$, and
each dual to an annulus, such that $X$ and $X^{\prime}$ cross. When $n=1$,
Proposition 8.6 of \cite{SS2} shows that $X$ and $X^{\prime}$ must be enclosed
by the $V_{0}$--vertex $v$ of $\Gamma_{n,n+1}$, and essentially the same
argument applies for all values of $n$.
\end{proof}

Recall that we want to consider a $V_{0}$--vertex $v$ of $\Gamma_{n,n+1}$
which is of commensuriser type, and that we do not even know that
$G(v)=Comm_{G}(H)$ is finitely generated. We first deal with two special cases
when $Comm_{G}(H)$ is certainly finitely generated.

\begin{lemma}
\label{V0-vertexofsolidtorustype}Let $(G,\partial G)$ be an orientable
$PD(n+2)$ pair such that $G$ is not $VPC$, and let $\Gamma_{n,n+1}$ denote the
reduced algebraic regular neighbourhood of $\mathcal{F}_{n,n+1}$ in $G$. Let
$v$ be a $V_{0}$--vertex of $\Gamma_{n,n+1}$ which is of commensuriser type,
such that $G(v)$ is the full commensuriser $Comm_{G}(H)$ for some $VPCn$
subgroup $H$ of $G$ with $e(G,H)\geq2$.

If $H$ has finite index in $G(v)$, then $v$ is of solid torus type (see
Definition \ref{defnofsolidtorustype}).
\end{lemma}

\begin{proof}
Let $e_{1},\ldots,e_{m}$ denote the edges of $\Gamma_{n,n+1}$ which are
incident to $v$, and denote the associated subgroups of $G(v)$ by
$H_{1},\ldots,H_{m}$. The hypothesis implies that $G(v)=Comm_{G}(H)$ is itself
$VPCn$. Hence each $H_{i}$ is also $VPC(\leq n)$. Now Lemma
\ref{PD(n+2)groupshavenoa.i.setsoverVPCk<n} implies that $G$ cannot split over
a $VPC(<n)$ subgroup. It follows that each $H_{i}$ is $VPCn$, and now Lemma
\ref{edgesplittingsaredualtoannuliortori} tells us that the edge splitting of
$G$ associated to $e_{i}$ is dual to an essential annulus $A_{i}$. It remains
to show that the boundaries of these annuli all carry the same subgroup of
$G(v)$. We write $\partial H_{i}$ for the group carried by the boundary of
$A_{i}$. Note that as each $H_{i}$ is a $VPCn$ subgroup of the $VPCn$ group
$G(v)$, it must be of finite index in $G(v)$. Thus each $H_{i}$ is
commensurable with $H$.

Lemma \ref{commensuriservertexenclosescrossingannuli} shows that $v$ encloses
two almost invariant subsets $X\ $and $X^{\prime}$ of $G$, each over a
subgroup of $H$ of finite index, and each dual to an annulus, such that $X$
and $X^{\prime}$ cross. By replacing $H$ by a subgroup of finite index, we can
assume that $X$ and $X^{\prime}$ are both $H$--almost invariant. As discussed
in section \ref{essentialannuli}, the doubles of these annuli are tori $T$ and
$T^{\prime}$ in $DG$. These tori must cross and are therefore enclosed by the
same $V_{0}$--vertex $V$ of $T_{n+1}(DG)$. Hence $V$ is not an isolated
vertex, so that $V$ must be of $VPCn$--by--Fuchsian type. We let $L$ denote
the $VPCn$ normal subgroup of $G(V)$ with Fuchsian quotient. As $T$ and
$T^{\prime}$ are non-peripheral tori enclosed by $V$, Lemma
\ref{torienclosedbydistinctverticesmustlieinone} shows that the intersection
group $T\cap T^{\prime}$ must be commensurable with $L$. As $T\cap T^{\prime}$
contains the $VPCn$ group $H$, it follows that $H$ is commensurable with $L$,
and hence that each $H_{i}$ is also commensurable with $L$. Let $T_{i}$ denote
the torus in $DG$ obtained by doubling the annulus $A_{i}$. The tori $T_{i}$
and $T$ have the $VPCn$ subgroup $H\cap H_{i}$ in common. As $T$ is
non-peripheral in $V$, Lemma \ref{torienclosedbydistinctverticesmustlieinone}
shows that $T_{i}$ must also be enclosed by $V$. Now Lemma \ref{deltaHequalsL}
below shows that $\partial H_{i}$ must equal $L$. As this holds for all $i$,
it shows that the groups $\partial H_{i}$ are all equal, as required. This
completes the proof that $v$ is of solid torus type.
\end{proof}

\begin{lemma}
\label{deltaHequalsL}Let $(G,\partial G)$ be an orientable $PD(n+2)$ pair such
that $G$ is not $VPC$, and suppose that $G$ has a splitting over a $VPCn$
subgroup $H_{i}$ which is dual to an essential annulus $A_{i}$ in $(G,\partial
G)$. Let $\partial H_{i}$ denote the group carried by the boundary of $A_{i}$,
and let $T_{i}$ denote the torus in $DG$ obtained by doubling $A_{i}$. Let $V$
denote a vertex of $T_{n+1}(DG)$ of $VPCn$--by--Fuchsian type, and let $L$
denote the normal $VPCn$ subgroup of $G(V)$ with Fuchsian quotient. Suppose
that $T_{i}$ is enclosed by $V$, and that $H_{i}$ is commensurable with $L$.
Then $\partial H_{i}$ equals $L$.
\end{lemma}

\begin{proof}
Observe that $\partial H_{i}$ is a normal subgroup of $T_{i}$ with quotient
$Q$ isomorphic to $\mathbb{Z}$ or to $\mathbb{Z}_{2}\ast\mathbb{Z}_{2}$,
depending on whether the annulus $A_{i}$ is untwisted or twisted. As $\partial
H_{i}$ is assumed to be commensurable with $L$, part 2) of Lemma
\ref{torusinVPC-by-Fuchsiangroup} implies that $\partial H_{i}$ equals $L\cap
T_{i}$, so that $\partial H_{i}$ must be contained in $L$.

As in section \ref{essentialannuli}, it will again be convenient to consider
an aspherical space $M$ with fundamental group $G$ and with aspherical
subspaces which correspond to $\partial G$ whose union is denoted by $\partial
M$. Let $\Sigma$ be a component of $\partial M$ such that one end of $A_{i}$
is in $\Sigma$. The splitting of $G$ over $H_{i}$ determined by $A_{i}$
induces a splitting of $\pi_{1}(\Sigma)$ over $\partial H_{i}$. Thus part 3)
of Corollary \ref{PDgroupsplitsoverPDHimpliesHismaximal} shows that $\partial
H_{i}$ is a maximal orientable $VPCn$ subgroup of $\pi_{1}(\Sigma)$. Now $DG$
splits over $\pi_{1}(\Sigma)$ and this induces the splitting of $T_{i}$ over
$\partial H_{i}$. As $\partial H_{i}$ equals $L\cap T_{i}$, we can identify
the quotient $Q=T_{i}/\partial H_{i}$ with a subgroup of the Fuchsian quotient
group of $G(V)$ by $L$. Now we consider the full pre-image $\overline{T_{i}}$
of $Q$ in $G(V)$. Thus $L$ is normal in $\overline{T_{i}}$ with quotient $Q$.
In particular $\overline{T_{i}}$ is a $VPC(n+1)$ subgroup of $G(V)$ which
contains $T_{i}$ with finite index. As the splitting of $DG$ over $\pi
_{1}(\Sigma)$ induces a splitting of $T_{i}$ over $\partial H_{i}$, it follows
that it induces a splitting of $\overline{T_{i}}$ over some subgroup
$L^{\prime}$ of $\pi_{1}(\Sigma)$ which contains $\partial H_{i}$ with finite
index. Thus $L$ and $L^{\prime}$ must be commensurable. As $\overline{T_{i}}$
splits over $L^{\prime}$, Lemma \ref{splittingsofVPCgroups} shows that
$L^{\prime}$ must be a normal $VPCn$ subgroup of $\overline{T_{i}}$ with
quotient which is isomorphic to $\mathbb{Z}$ or $\mathbb{Z}_{2}\ast
\mathbb{Z}_{2}$. Now part 2) of Lemma
\ref{normalsubgroupsofVPCgroupsofcolength1} shows that $L^{\prime}$ must equal
$L$. In particular, it follows that $\partial H_{i}\subset L\subset\pi
_{1}(\Sigma)$. Recall from Lemma \ref{Lisorientable} that $L$ must be
orientable. Now the maximality of $\partial H_{i}$ among orientable $VPCn$
subgroups of $\pi_{1}(\Sigma)$ implies that $\partial H_{i}$ must equal $L$,
as required.
\end{proof}

Next we consider a case where $Comm_{G}(H)$ contains $H$ with infinite index,
but still must be finitely generated.

\begin{lemma}
\label{V0-vertexoftorustype}Let $(G,\partial G)$ be an orientable $PD(n+2)$
pair such that $G$ is not $VPC$, and let $\Gamma_{n,n+1}$ denote the reduced
algebraic regular neighbourhood of $\mathcal{F}_{n,n+1}$ in $G$. Let $v$ be a
$V_{0}$--vertex of $\Gamma_{n,n+1}$ which is of commensuriser type, such that
$G(v)$ is the full commensuriser $Comm_{G}(H)$ for some $VPCn$ subgroup $H$ of
$G$ with $e(G,H)\geq2$.

If $G(v)$ is $VPC(n+1)$, then $v$ is of torus type (see Definition
\ref{defnoftorustype}).
\end{lemma}

\begin{proof}
As $G(v)$ is $VPC(n+1)$, and $G$ cannot split over a $VPC(<n)$ subgroup, it
follows that the group associated to each edge of $\Gamma_{n,n+1}$ which is
incident to $v$ must be $VPCn$ or $VPC(n+1)$. Now Lemma
\ref{edgesplittingsaredualtoannuliortori} implies that the edge splittings of
$G$ associated to these edges are dual to essential annuli or tori. Let
$\Sigma(v)$ denote the collection of all these splittings. Recall that
Definition \ref{defnoftorustype} includes four cases, but in each case, at
most one of the splittings in $\Sigma(v)$ can be dual to an essential torus.

We start by showing that this condition holds.

\medskip

\textbf{The splittings in }$\Sigma(v)$

\medskip

Note that any splitting of $G$ dual to an essential torus must be over a
maximal orientable $VPC(n+1)$ subgroup of $G$, by part 3) of Corollary
\ref{PDgroupsplitsoverPDHimpliesHismaximal}. Thus if $G(v)$ is orientable, any
such edge splitting must be over $G(v)$. If $G(v)$ is non-orientable, then any
such edge splitting must be over $G(v)_{0}$, the orientable subgroup of $G(v)$
of index $2$. Let $T_{n,n+1}$ denote the universal covering $G$--tree of
$\Gamma_{n,n+1}$, and let $w$ denote a vertex of $T_{n,n+1}$ above $v$ and
with stabiliser $G(v)$. If $G(v)$ is orientable, suppose that two of the
splittings in $\Sigma(v)$ are dual to an essential torus, and if $G(v)$ is
non-orientable, suppose that one of the splittings in $\Sigma(v)$ is dual to
an essential torus. In either case, there are two distinct edges of
$T_{n,n+1}$ which are incident to $w$ and have the same stabiliser, which is
$G(v)$ or $G(v)_{0}$. Let $X$ and $Y$ denote the almost invariant subsets of
$G$ associated to these edges. As each is dual to the same essential torus in
$(G,\partial G)$ they must be equivalent. Now Corollary 4.16 of \cite{SS2}
implies that $w$ must have valence $2$, and hence is isolated. This implies
that either $v$ is isolated, or that $v$ has valence $1$ with edge group of
index $2$ in $G(v)$. Either case contradicts our assumption that $v$ is of
commensuriser type. We conclude that if $G(v)$ is orientable, then at most one
of the splittings in $\Sigma(v)$ can be dual to a torus, and if $G(v)$ is
non-orientable, then none of the splittings in $\Sigma(v)$ can be dual to a
torus. Further if $G(v)$ is orientable and one of the splittings in
$\Sigma(v)$ is dual to a torus, the minimality of $\Gamma_{n,n+1}$ shows that
there must also be a splitting in $\Sigma(v)$ which is dual to an essential
annulus. Thus in all cases, $\Sigma(v)$ must contain at least one splitting
which is dual to an essential annulus.

\medskip

\textbf{Notation}

\medskip

Let $e_{1},\ldots,e_{m}$ denote all those edges of $\Gamma_{n,n+1}$ which are
incident to $v$ and have associated edge splitting dual to an essential
annulus. Denote the associated $VPCn$ subgroups of $G(v)$ by $H_{1}%
,\ldots,H_{m}$. Let $T_{i}$ denote the torus in $DG$ obtained by doubling the
annulus $A_{i}$ associated to the edge splitting of $G$ over $H_{i}$.

Lemma \ref{commensuriservertexenclosescrossingannuli} implies that $v$
encloses two almost invariant subsets $X\ $and $X^{\prime}$ of $G$, each dual
to an annulus, such that $X$ and $X^{\prime}$ cross. Further we can assume
that $X$ and $X^{\prime}$ are both $H$--almost invariant. The doubles of these
annuli are tori in $DG$ which must cross, and are therefore both enclosed by a
$V_{0}$--vertex $V$ of $T_{n+1}(DG)$. Thus neither torus is peripheral in $V$,
and $V$ is of $VPCn$--by--Fuchsian type. If we let $L$ denote the $VPCn$
normal subgroup of $G(V)$ with Fuchsian quotient, then $H$ is commensurable
with $L$. Now Lemma \ref{normaliserequalscommensuriser} shows that
$G(V)=Comm_{DG}(L)$ which equals $Comm_{DG}(H)$. As $G(v)=Comm_{G}(H)$, it
follows that $G(v)$ is a subgroup of $G(V)$.

\medskip

\textbf{Case 1: }$G(v)$\textbf{ is orientable.}

\medskip

Thus $G(v)$ is an essential torus in $DG$, and we let $W$ denote the almost
invariant subset of $DG$ which is over $G(v)$. As $G(v)$ is a subgroup of
$G(V)$, it follows that $W$ is enclosed by $V$.

We will suppose that $W$ is not peripheral in $V$. We will not need to
consider here the case when $W$ is peripheral in $V$, although this can occur.

As $H_{i}$ is a subgroup of $G(v)$, the tori $T_{i}$ and $G(v)$ have the
$VPCn$ subgroup $H_{i}$ in common. Then Lemma
\ref{torienclosedbydistinctverticesmustlieinone} shows that the almost
invariant subset of $DG$ associated to the torus $T_{i}$ must be enclosed by
$V$, and that $H_{i}$ is commensurable with $L$. Now Lemma \ref{deltaHequalsL}
shows that $\partial H_{i}$ must equal $L$. As $\partial H_{i}$ is a subgroup
of $G(v)$ which in turn is a subgroup of $G(V)$, and as $L$ is normal in
$G(V)$, it follows that $\partial H_{i}$ is a normal $VPCn$ subgroup of $G(v)$
with quotient isomorphic to $\mathbb{Z}$ or to $\mathbb{Z}_{2}\ast
\mathbb{Z}_{2}$, so that $G(v)$ splits over $\partial H_{i}$. This proves part
of cases 1), 2) and 3) of the definition of torus type (Definition
\ref{defnoftorustype}).

We next consider how the $VPC(n+1)$ group $G(v)$ sits in $G$. There are three subcases.

\medskip

\textbf{Case 1a): }$G(v)$\textbf{ is conjugate to a subgroup of some group in
}$\partial G$\textbf{.}

\medskip

Note that this implies that no edge of $\Gamma_{n,n+1}$ incident to $v$ has
associated splitting dual to an essential torus. As in section
\ref{essentialannuli}, it will again be convenient to consider an aspherical
space $M$ with fundamental group $G$ and with aspherical subspaces which
correspond to $\partial G$ whose union is denoted by $\partial M$. Let
$\Sigma$ be a component of $\partial M$ such that $G(v)$ is a subgroup of
$\pi_{1}(\Sigma)$. As $\pi_{1}(\Sigma)$ must contain $G(v)$ with finite index,
it is conjugate into some vertex group of $\Gamma_{n,n+1}$. Thus either
$\pi_{1}(\Sigma)$ equals $G(v)$, or some edge incident to $v$ has associated
group equal to $G(v)$. As the second case is impossible, $G(v)$ must equal
$\pi_{1}(\Sigma)$. Consider the cover $M_{H}$ of $M$ with fundamental group
equal to $H$. There will be a component of $\partial M_{H}$ which covers
$\Sigma$ and has fundamental group $H$. As we already know that $M_{H}$ admits
essential annuli which carry $H$, it follows that there is an essential
annulus in $M$ which carries $H$ and has one end on $\Sigma$. Doubling such an
annulus yields a torus in $DG$ which crosses the torus $\Sigma$. It follows
that the $G(v)$--almost invariant subset $W$ of $DG$ associated to $G(v)$ is
not peripheral in the $V_{0}$--vertex $V$ of $T_{n+1}(DG)$ which encloses $W$,
so the argument in the preceding paragraph applies and shows that each
$\partial H_{i}$ equals $L$. As above, it also follows that $G(v)$ splits over
$\partial H_{i}$. Thus we have case 1) of Definition \ref{defnoftorustype}.

\medskip

\textbf{Case 1b): }$G(v)$\textbf{ is not conjugate to a subgroup of some group
in }$\partial G$\textbf{.}

\medskip

In this case, there is a nontrivial $G(v)$--almost invariant subset of $G$
which is adapted to $\partial G$, and we denote this set by $Y$. The
$G(v)$--almost invariant subset $W$ of $DG$ is enclosed by $G$, and $Y$ is
equal to $W\cap G$, and is adapted to $\partial G$. Now Theorem
\ref{FandF'havesameregnbhds} tells us that $Y$ must be enclosed by some
$V_{0}$--vertex $w$ of $\Gamma_{n,n+1}$. If $w$ is not equal to $v$, then
there must be an edge of $\Gamma_{n,n+1}$ which is incident to $v$, and
carries $G(v)$. Thus in any case, $Y$ is enclosed by $v$. Now we have two
subcases depending on whether or not $Y$ is peripheral in $v$.

If $Y$ is not peripheral in $v$, then $Y$ must cross some almost invariant
subset $Z$ of $G$ which is over a subgroup commensurable with $H$ and belongs
to the CCC of $\mathcal{F}_{n,n+1}$ which gives rise to $v$. Lemma
\ref{XcrossesYimpliesitcrossesY'strongly} shows that $Y$ must cross some
almost invariant subset $Z^{\prime}$ of $G$ also over a subgroup commensurable
with $H$ and dual to an annulus in $(G,\partial G)$. Now part 1) of Lemma
\ref{FandF'havesamebetweenness} shows that $Z^{\prime}$ cannot be isolated in
$\mathcal{F}_{n,n+1}$ and so must also belong to the CCC of $\mathcal{F}%
_{n,n+1}$ which gives rise to $v$. We conclude that $Y$ crosses some annulus
in $(G,\partial G)$ enclosed by $v$. It follows that $W$ crosses the torus in
$DG$ which is the double of this annulus. As above this implies that $W$ is
not peripheral in $V$, and hence that each $\partial H_{i}$ equals $L$, and
that $G(v)$ splits over $\partial H_{i}$. As $v$ is a $V_{0}$--vertex of
$\Gamma_{n,n+1}$, and $Y$ is enclosed by $v$ and is over $G(v)$, it follows
that $Y$ crosses no torus in $G$ so that $Y$ determines a splitting of $G$
over $G(v)$. As $Y$ is enclosed by $v$ and is not peripheral in $v$, this
shows that we have case 3) of Definition \ref{defnoftorustype}.

If $Y$ is peripheral in $v$, we let $e$ be the edge of $\Gamma_{n,n+1}$ to
which $Y$ is associated. We also let $X_{i}$ denote the $H_{i}$--almost
invariant subset of $G$ associated to the edge $e_{i}$. Now the definition of
betweenness in our construction of an algebraic regular neighbourhood in
\cite{SS2} implies that there is an element $Z$ of the CCC which gives rise to
$v$ such that $Z$ lies between $X_{i}$ and $Y$. More precisely, there is a
nontrivial almost invariant subset $Z$ of $G$ which is over a subgroup
commensurable with $H$ and which is enclosed by $v$ such that $Y^{(\ast)}\leq
Z\leq X_{i}^{(\ast)}$, where $Y^{(\ast)}$ denotes one of $Y$ or $Y^{\ast}$. As
$H_{i}$ stabilises both $X_{i}$ and $Y$, it follows that $H_{i}$ lies within
some bounded neighbourhood of $\delta Z$. In turn this implies that a subgroup
of finite index in $H_{i}$ must stabilise $Z$. Thus $H_{i}$ and $H$ must be
commensurable. As $H$ is commensurable with $L$, it follows that $H_{i}$ is
commensurable with $L$, and now Lemma \ref{deltaHequalsL} shows that $\partial
H_{i}$ must equal $L$. As above, it also follows that $G(v)$ splits over
$\partial H_{i}$. Thus we have case 2) of Definition \ref{defnoftorustype}.

This completes the proof that $v$ is of torus type, in the case when $G(v)$ is orientable.

\medskip

\textbf{Case 2: }$G(v)$\textbf{ is non-orientable.}

\medskip

Let $W$ denote the almost invariant subset of $DG$ which is over $G(v)_{0}$.
As before $W$ is enclosed by $V$. Recall that $V$ is the $V_{0}$--vertex of
$T_{n+1}(DG)$ with $G(V)=Comm_{DG}(L)=Comm_{DG}(H)$, and that $G(v)$ is a
subgroup of $G(V)$.

If $G(v)_{0}$ is conjugate to a subgroup of some group in $\partial G$, then
the non-orientable $PD(n+1)$ group $G(v)$ is commensurable with that group in
$\partial G$. Now part 1) of Corollary
\ref{PDgroupsplitsoverPDHimpliesHismaximal} shows that $G$ itself must be a
non-orientable $PD(n+1)$ group, and that $G$ must contain $G(v)_{0}$ with
finite index. In particular it would imply that $G$ is $VPC(n+1)$, which is
excluded in our hypotheses. Thus $G(v)_{0}$ cannot be conjugate to a subgroup
of any group in $\partial G$, so we let $Y=W\cap G$, a nontrivial $G(v)_{0}%
$--almost invariant subset of $G$ which is adapted to $\partial G$. As before
$Y$ must be enclosed by $v$. As we showed earlier that no edge incident to $v$
can have associated splitting dual to an essential torus, it follows that $Y$
cannot be peripheral in $v$. As in the case when $G(v)$ is orientable, this
implies that $W$ is not peripheral in $V$, and hence that each $\partial
H_{i}$ equals $L$, the $VPCn$ normal subgroup of $G(V)$ with Fuchsian
quotient. Thus $L$ is also a normal subgroup of $G(v)$, and the quotient
$G(v)/L$ is a $VPC1$ subgroup of the Fuchsian quotient $G(V)/L$. Hence
$G(v)/L$ is isomorphic to $\mathbb{Z}$ or to $\mathbb{Z}_{2}\ast\mathbb{Z}%
_{2}$, so that $G(v)$ splits over $L$. Thus $G(v)$ splits over $\partial
H_{i},$ for each $i$, as required by part 4) of Definition
\ref{defnoftorustype}.

As $v$ is a $V_{0}$--vertex of $\Gamma_{n,n+1}$, and $Y$ is enclosed by $v$
and is over a subgroup of index $2$ in $G(v)$, it follows that $Y$ determines
a splitting $\sigma$ of $G$ over $G(v)_{0}$. As $G(v)$ contains $G(v)_{0}$
with finite index, it must be conjugate into one of the vertex groups $K$ of
this splitting. Note that the splitting $\sigma$ is dual to an essential torus
in $(G,\partial G)$. Thus, if $\partial G$ is empty, Theorem 8.1 of \cite{B-E}
shows that the pair formed by $K$ and one or two copies of $G(v)_{0}$ is
$PD(n+2)$, where there will be two copies of $G(v)_{0}$ if $\sigma$ is a HNN
extension and only one copy otherwise. In general, as discussed just before
Definition \ref{defnofatoroidal}, the pair becomes $PD(n+2)$ when some family
of groups in $\partial G$ is added to the copies of $G(v)_{0}$. As $K$
contains a conjugate of $G(v)$, one of the copies of $G(v)_{0}$ in $\partial
K$ is not equal to its own commensuriser in $K$. Thus Lemma
\ref{boundarygroupismaximal} implies that $K$ contains this copy of $G(v)_{0}$
with index $2$, so that $K$ must be a conjugate of $G(v)$, and $\partial K$
consists only of $G(v)_{0}$. In particular it follows that $\sigma$ must be an
amalgamated free product and not a HNN extension.

As $Y$ is enclosed by $v$, we can refine $\Gamma_{n,n+1}$ by splitting at $v$
to obtain a graph of groups structure $\Gamma^{\prime}$ of $G$ such that the
projection map $\Gamma^{\prime}\rightarrow\Gamma_{n,n+1}$ sends an edge $e$ to
$v$ and otherwise induces a bijection of edges and vertices. The group
associated to $e$ is equal to $G(v)_{0}$, and the associated edge splitting is
$\sigma$. As $G(v)$ contains $G(v)_{0}$ with finite index, one vertex of $e$
must carry $G(v)$ and the other must carry $G(v)_{0}$. Let $w$ denote the
vertex of $e$ with $G(w)=G(v)$. Recall from the preceding paragraph that
$\sigma$ is an amalgamated free product and that one vertex group of $\sigma$
is conjugate to $G(v)$. It follows that the edge $e$ of $\Gamma^{\prime}$ is
separating, and that if we remove the interior of $e$ from $\Gamma^{\prime}$
then the component of the resulting subgraph which contains $w$ must carry the
group $G(v)$. As $G(w)$ equals $G(v)$, the minimality of $\Gamma_{n,n+1}$
implies that this subgraph must consist solely of $w$, so that $w$ has valence
$1$ in $\Gamma^{\prime}$. This shows that we have case 4) of Definition
\ref{defnoftorustype}, and so completes the proof that $v$ is of torus type in
all cases.
\end{proof}

Now we consider the general situation of a $V_{0}$--vertex $v$ of
$\Gamma_{n,n+1}$ which is of commensuriser type. Note that a priori, the group
$G(v)$ need not be finitely generated, but we show not only that $G(v)$ is
finitely generated, but also describe its structure.

\begin{lemma}
\label{V0-vertexofSeiferttype}Let $(G,\partial G)$ be an orientable $PD(n+2)$
pair such that $G$ is not $VPC$, and let $\Gamma_{n,n+1}$ denote the reduced
algebraic regular neighbourhood of $\mathcal{F}_{n,n+1}$ in $G$. Let $v$ be a
$V_{0}$--vertex of $\Gamma_{n,n+1}$ which is of commensuriser type, such that
$G(v)$ is the full commensuriser $Comm_{G}(H)$ for some $VPCn$ subgroup $H$ of
$G$ with $e(G,H)\geq2$.

If $H$ has infinite index in $G(v)$, and $G(v)$ is not $VPC(n+1)$, then $v$ is
of Seifert type (see Definition \ref{defnofSeiferttype}).
\end{lemma}

\begin{proof}
Lemma \ref{commensuriservertexenclosescrossingannuli} implies that $v$ must
enclose two almost invariant subsets of $G$, each dual to an annulus and
crossing each other. Further we can assume that each is over $H$. The doubles
of these annuli are tori in $DG$ which must cross, and are therefore both
enclosed by a $V_{0}$--vertex $V$ of $T_{n+1}(DG)$. Thus neither torus is
peripheral in $V$, and $V$ is of $VPCn$--by--Fuchsian type. If we let $L$
denote the $VPCn$ normal subgroup of $G(V)$ with Fuchsian quotient $\Phi$,
then $H$ is commensurable with $L$. Now Lemma
\ref{normaliserequalscommensuriser} tells us that $G(V)=Comm_{DG}(L)$. It
follows that $G(V)$ is also equal to $Comm_{DG}(H)$, and so $G(V)$ contains
$Comm_{G}(H)=G(v)$. Thus $G(v)$ is itself a $VPCn$--by--Fuchsian group, where
the normal $VPCn$ subgroup $H^{\prime}$ is commensurable with $H$. (But note
that at this stage it is still possible that $G(v)$ is not finitely
generated!) As we are assuming $H$ has large commensuriser, it follows that
the Fuchsian quotient group $\Theta=G(v)/H^{\prime}$ must be infinite. Further
as $G(v)$ is not $VPC(n+1)$, this quotient cannot be two-ended. This implies
that there are elements $\alpha$ and $\beta$ in $\Theta$ of infinite order
such that $\alpha$ and $\beta$ have non-zero geometric intersection number.
The pre-images in $G(v)$ of the infinite cyclic subgroups of $\Theta$
generated by $\alpha$ and $\beta$ are $VPC(n+1)$ subgroups $A^{\prime}$ and
$B^{\prime}$ of $G(v)$. Note that $A^{\prime}\cap B^{\prime}=H^{\prime}$.

If we regard $\alpha$ and $\beta$ as elements of $\Phi$, the pre-images in
$G(V)$ of the same infinite cyclic subgroups are $VPC(n+1)$ subgroups $A$ and
$B$ of $G(V)$. By replacing $\alpha$ and $\beta$ by their squares if needed,
we can ensure that they are orientable elements of $\Phi$, so that $A$ and $B$
will be orientable. As $\alpha$ and $\beta$ have non-zero geometric
intersection number, it follows that $A$ and $B$ are tori in $DG$ which cross.
As $A^{\prime}\ $and $B^{\prime}$ are subgroups of finite index in $A$ and $B$
respectively, they also are tori in $DG$ which cross. As $A^{\prime}$ and
$B^{\prime}$ are subgroups of $G(v)$ which is a subgroup of $G$, it follows
that $A^{\prime}$ and $B^{\prime}$ are tori in $(G,\partial G)$ which are
enclosed by $v$ and which cross.

Now we consider the torus decomposition $T_{n+1}(G,\partial G)$. Each of
$A^{\prime}$ and $B^{\prime}$ must be enclosed by some $V_{0}$--vertex of
$T_{n+1}(G,\partial G)$. As $A^{\prime}\ $and $B^{\prime}$ cross, they must
both be enclosed by a single $V_{0}$--vertex $u$ of $T_{n+1}(G,\partial G)$,
and neither is peripheral in $u$. Thus $u$ is not isolated, and so Theorem
\ref{torusdecompofapair} shows that $u$ must be of Seifert type adapted to
$\partial G$. Let $L^{\prime}$ denote the normal subgroup of $G(u)$ with
Fuchsian quotient. Recall that $A^{\prime}$ and $B^{\prime}$ are tori in
$(G,\partial G)$ which cross and that $A^{\prime}\cap B^{\prime}=H^{\prime}$.
Lemma \ref{torienclosedbydistinctverticesmustlieinone} implies that
$H^{\prime}$ is commensurable with $L^{\prime}$. Now Lemma
\ref{normaliserequalscommensuriser} implies that $G(u)$ equals $Comm_{G}%
(L^{\prime})$. It follows that $G(u)=Comm_{G}(L^{\prime})=Comm_{G}(H^{\prime
})=Comm_{G}(H)=G(v)$. Note that as $G(u)$ is finitely generated, it follows
that $G(v)$ must also be finitely generated. As each of $H^{\prime}$ and
$L^{\prime}$ is a normal $VPCn$ subgroup of $G(v)$ with Fuchsian quotient,
Lemma \ref{VPC-by-Fuchsiangrouphasuniquefibre} shows that $H^{\prime}$ and
$L^{\prime}$ must be equal.

Now Theorem \ref{FandF'havesameregnbhds} shows that each essential torus in
$(G,\partial G)$ which is enclosed by the vertex $u$ of $T_{n+1}(G,\partial
G)$ is also enclosed by some $V_{0}$--vertex $v$ of $\Gamma_{n,n+1}$. As the
tori which are enclosed by $u$ and are not peripheral in $u$ form a single
CCC, they must all be enclosed by $v$. It follows that the vertex $u$ of
$T_{n+1}(G,\partial G)$ is enclosed by the vertex $v$ of $\Gamma_{n,n+1}$.
Thus there is a refinement $\Gamma^{\prime}$ of $\Gamma_{n,n+1}$ and a vertex
$v^{\prime}$ of $\Gamma^{\prime}$, such that the projection map $p:\Gamma
^{\prime}\rightarrow\Gamma_{n,n+1}$ sends $v^{\prime}$ to $v$ and is an
isomorphism apart from the fact that certain edges incident to $v^{\prime}$
are mapped to $v$. Further the vertex $v^{\prime}$, like $u$, is of Seifert
type adapted to $\partial G$, and $G(v^{\prime})=G(u)$ maps isomorphically to
$G(v)$. This last fact implies that if $e$ is an edge of $\Gamma^{\prime}$
which is incident to $v^{\prime}$ and mapped to $v$, then the other vertex $w$
of $e$ has associated group equal to $G(e)$. As $\Gamma^{\prime}$ is minimal
there is at least one other edge incident to $w$, and each such edge must
carry a subgroup of the $VPC(n+1)$ group $G(w)=G(e)$. Thus each such edge
carries a $VPC(\leq n+1)$ subgroup of $G$. As usual, Lemma
\ref{edgesplittingsaredualtoannuliortori} tells us that the associated
splitting of $G$ must be dual to an essential annulus or torus in $(G,\partial
G)$. If one of these other edges incident to $w$ has associated splitting dual
to an essential torus, the edge group must equal $G(e)$, so that the vertex
$w$ of $\Gamma^{\prime}$ has two incident edges with associated splittings
over the same essential torus. In this case, as in the proof of Lemma
\ref{V0-vertexoftorustype}, Corollary 4.16 of \cite{SS2} implies that $w$ must
have valence $2$, and we modify $\Gamma^{\prime}$ by collapsing the edge $e$.
By repeating this process we will arrange that if $e$ is an edge of
$\Gamma^{\prime}$ which is incident to $v^{\prime}$ and mapped to $v$, with
the other vertex of $e$ being $w$, then each edge $e_{i}$, other than $e$,
which is incident to $w$ must carry a $VPCn$ subgroup $H_{i}$ of $G(w)=G(e)$,
so that the associated edge splitting is dual to an annulus $A_{i}$ in
$(G,\partial G)$. As usual we write $\partial H_{i}$ for the group carried by
$\partial A_{i}$. At this point we have proved that most of Definition
\ref{defnofSeiferttype} holds. It remains to show that $\partial H_{i}$ equals
$L^{\prime}$, for each $i$, where $L^{\prime}$ is the normal subgroup of
$G(u)=G(v)$ with Fuchsian quotient.

The image of $e_{i}$ in $\Gamma_{n,n+1}$ is an edge incident to $v$, so its
associated edge splitting is enclosed by $v$. The edge splitting of $G$
associated to the edge $e$ of $\Gamma^{\prime}$ is dual to a torus $T$ which
is also enclosed by $v$. We claim that the torus $T$ cannot be peripheral in
$v$. To see this we need to consider the universal covering $G$--trees
$T^{\prime}$ and $T_{n,n+1}$ of $\Gamma^{\prime}$ and $\Gamma_{n,n+1}$
respectively. Let $V^{\prime}$ be a vertex of $T^{\prime}$ above $v^{\prime}$
with stabiliser equal to $G(v^{\prime})$. Let $E$ be an edge of $T^{\prime}$
above $e$ which is incident to $V^{\prime}$ and has stabiliser equal to
$G(e)$, let $W$ denote the other vertex of $E$, and let $E_{i}$ denote an edge
of $T^{\prime}$ above $e_{i}$ which is incident to $W$ and has stabiliser
equal to $H_{i}$. Orient $E\ $and $E_{i}$ towards $V^{\prime}$, and let
$Z_{E}$ and $Z_{i}$ denote the associated almost invariant subsets of $G$. The
orientations imply that $Z_{i}^{\ast}\leq Z_{E}^{\ast}$. If $T$ were
peripheral in $v$, the fact that $Z_{i}$ is enclosed by $v$ would imply that
$Z_{E}^{\ast}\leq Z_{i}$ or $Z_{E}^{\ast}\leq Z_{i}^{\ast}$. The first
inequality would imply that $Z_{i}^{\ast}\leq Z_{i}$, which is obviously
impossible. The second inequality would imply that $Z_{i}$ and $Z_{E}$ were
equivalent. This is impossible as their stabilisers are not commensurable, as
$H_{i}$ is $VPCn$ and $G(e)$ is $VPC(n+1)$. This contradiction shows that the
torus $T$ cannot be peripheral in $v$, as claimed.

Now let $T_{i}$ denote the torus in $DG$ obtained by doubling the annulus
$A_{i}$ with fundamental group $H_{i}$. As $T$ is not peripheral in $v$, it
must cross some almost invariant subset of $G$ which is over a subgroup
commensurable with $H$ and belongs to the CCC of $\mathcal{F}_{n,n+1}$ which
gives rise to $v$. As in the proof of Lemma \ref{V0-vertexoftorustype}, Lemmas
\ref{XcrossesYimpliesitcrossesY'strongly} and \ref{FandF'havesamebetweenness}
imply that $T$ crosses some annulus enclosed by $v$, and so crosses the torus
obtained by doubling this annulus. Lemma
\ref{torienclosedbydistinctverticesmustlieinone} then implies that $T$ and
this torus are enclosed by the vertex $V$ of $T_{n+1}(DG)$ and are not
peripheral in $V$. As $H_{i}$ is contained in $G(e)$, Lemma
\ref{torienclosedbydistinctverticesmustlieinone} now shows that $T_{i}$ is
enclosed by $V$, and that $H_{i}$ is commensurable with $L$. (Recall that $L$
is the $VPCn$ normal subgroup of $G(V)$ with Fuchsian quotient $\Phi$.) Lemma
\ref{deltaHequalsL}\ then shows that $\partial H_{i}$ must equal $L$. In
particular, $G(e)$ contains $L$. As $G(e)$ is an edge torus of the $V_{0}%
$--vertex $u$ of $T_{n+1}(G,\partial G)$, it must contain $L^{\prime}$. As
$L\ $and $L^{\prime}$ are commensurable normal $VPCn$ subgroups of $G(e)$ each
of which has quotient isomorphic to $\mathbb{Z}$ or to $\mathbb{Z}_{2}%
\ast\mathbb{Z}_{2}$, Lemma \ref{normalsubgroupsofVPCgroupsofcolength1} implies
that $L^{\prime}$ equals $L$. We conclude that $\partial H_{i}$ equals
$L^{\prime}$, for each $i$. We have now established all the requirements in
Definition \ref{defnofSeiferttype}, so that $v$ must be of Seifert type in
$\Gamma_{n,n+1}$, as required. This completes the proof of Lemma
\ref{V0-vertexofSeiferttype}.
\end{proof}

We are finally ready to prove our main result. For the convenience of the
reader we restate it.

\begin{quote1}
(Main Result) Let $n\geq1$, and let $(G,\partial G)$ be an orientable
$PD(n+2)$ pair such that $G$ is not $VPC$. Let $\mathcal{F}_{n,n+1}$ denote
the family of equivalence classes of all nontrivial almost invariant subsets
of $G$ which are over a $VPCn$ subgroup, together with the equivalence classes
of all $n$--canonical almost invariant subsets of $G$ which are over a
$VPC(n+1)$ subgroup. Finally let $\Gamma_{n,n+1}$ denote the reduced algebraic
regular neighbourhood of $\mathcal{F}_{n,n+1}$ in $G$, and let $\Gamma
_{n,n+1}^{c}$ denote the completion of $\Gamma_{n,n+1}$. Thus $\Gamma_{n,n+1}$
and $\Gamma_{n,n+1}^{c}$ are bipartite graphs of groups structures for $G$,
with vertices of $V_{0}$--type and of $V_{1}$--type.

Then $\Gamma_{n,n+1}$ and $\Gamma_{n,n+1}^{c}$ have the following properties:

\begin{enumerate}
\item Each $V_{0}$--vertex $v$ of $\Gamma_{n,n+1}$ satisfies one of the
following conditions:

\begin{enumerate}
\item $v$ is isolated, and $G(v)$ is $VPC$ of length $n$ or $n+1$, and the
edge splittings associated to the two edges incident to $v$ are dual to
essential annuli or tori in $G$.

\item $v$ is of $VPC(n-1)$--by--Fuchsian type, and is of $I$--bundle type.
(See Definition \ref{defnofI-bundletype}.)

\item $v$ is of $VPCn$--by--Fuchsian type, and is of interior Seifert type.
(See Definition \ref{defnofinteriorSeiferttype}.)

\item $v$ is of commensuriser type. Further $v$ is of Seifert type (see
Definition \ref{defnofSeiferttype}), or of torus type (see Definition
\ref{defnoftorustype}) or of solid torus type (see Definition
\ref{defnofsolidtorustype}).
\end{enumerate}

\item The $V_{0}$--vertices of $\Gamma_{n,n+1}^{c}$ obtained by the completion
process are of special Seifert type (see Definition
\ref{defnofspecialSeiferttype}) or of special solid torus type (see Definition
\ref{defnofsolidtorustype}).

\item Each edge splitting of $\Gamma_{n,n+1}$ and of $\Gamma_{n,n+1}^{c}$ is
dual to an essential annulus or torus in $G$.

\item Any nontrivial almost invariant subset of $G$ over a $VPC(n+1)$ group
and adapted to $\partial G$ is enclosed by some $V_{0}$--vertex of
$\Gamma_{n,n+1}$, and also by some $V_{0}$--vertex of $\Gamma_{n,n+1}^{c}$.

\item If $H$ is a $VPC(n+1)$ subgroup of $G$ which is not conjugate into
$\partial G$, then $H$ is conjugate into a $V_{0}$--vertex group of
$\Gamma_{n,n+1}^{c}$.
\end{enumerate}
\end{quote1}

\begin{remark}
Part 3) follows immediately from parts 1) and 2), as the definitions of the
various types of $V_{0}$--vertex in the statements of parts 1) and 2) all
contain the requirement that the edge splittings be dual to an essential
annulus or torus.

Part 4) does not follow from the properties of an algebraic regular
neighbourhood as an almost invariant subset of $G$ over a $VPC(n+1)$ group
which is adapted to $\partial G$ need not be $n$--canonical, and so need not
lie in the family $\mathcal{F}_{n,n+1}$. Note that, from \cite{SS4}, we know
that there may be almost invariant subsets of $G$ over $VPC(n+1)$ subgroups
which are not adapted to $\partial G$.

Part 5) also does not follow from the properties of an algebraic regular
neighbourhood as a $VPC(n+1)$ subgroup $H$ of $G$ may be non-orientable.
\end{remark}

\begin{proof}
1) Recall that Theorem \ref{JSJforVPCoftwolengths} states that each $V_{0}%
$--vertex $v$ of $\Gamma_{n,n+1}$ satisfies one of the following conditions:

a) $v$ is isolated, and $G(v)$ is $VPC$ of length $n$ or $n+1$.

b) $v$ is of $VPCk$--by--Fuchsian type, where $k$ equals $n-1$ or $n$.

c) $v$ is of commensuriser type, and $G(v)$ is the full commensuriser
$Comm_{G}(H)$ for some $VPC$ subgroup $H$ of length $n$ or $n+1$, such that
$e(G,H)\geq2$.

We will consider each type of $V_{0}$--vertex of $\Gamma_{n,n+1}$ in turn.

a) Suppose that a $V_{0}$--vertex $v$ of $\Gamma_{n,n+1}$ is isolated, and
$G(v)$ is $VPC$ of length $n$ or $n+1$. Let $e$ denote either of the two edges
of $\Gamma_{n,n+1}$ incident to $v$. Then $G(e)$ equals $G(v)$ and so is $VPC$
of length $n$ or $n+1$. Now Lemma \ref{edgesplittingsaredualtoannuliortori}
implies that the edge splitting of $G$ associated to $e$ is dual to an
essential annulus or torus. It follows that we have case 1a) of Theorem
\ref{mainresult}.

b) Suppose that a $V_{0}$--vertex $v$ of $\Gamma_{n,n+1}$ is of $VPCk$%
--by--Fuchsian type, where $k$ equals $n-1$ or $n$. If $k=n$, part 1) of Lemma
\ref{V0vertexofVPCk-by-Fuchsiantype} shows that $v$ is of interior Seifert
type, so that we have case 1c) of Theorem \ref{mainresult}. If $k=n-1$, part
2) of Lemma \ref{V0vertexofVPCk-by-Fuchsiantype} shows that $v$ is of
$I$--bundle type, so that we have case 1b) of Theorem \ref{mainresult}.

c) Suppose that $v$ is a $V_{0}$--vertex of $\Gamma_{n,n+1}$ which is of
commensuriser type, and that $G(v)$ is the full commensuriser $Comm_{G}(H)$
for some $VPC$ subgroup $H$ of length $n$ or $n+1$, such that $e(G,H)\geq2$.

If $H$ has length $n+1$, we recall from Theorem \ref{JSJforVPCoftwolengths}
and Remark \ref{remarksonJSJ} that $v$ encloses elements $X$ and $Y$ of
$\mathcal{F}_{n,n+1}$ which are over a subgroup $H^{\prime}$ of finite index
in $H$, and which cross weakly. Now Corollary \ref{FiscontainedinF'} and
Remark \ref{adaptedtodGimpliesorientablegroup} imply that as $X$ and $Y$ are
$n$--canonical, they must be adapted to $\partial G$, and $H^{\prime}$ must be
orientable. Thus $X$ and $Y$ are dual to essential tori in $(G,\partial G)$.
Now Lemma \ref{propertiesofadapted} tells us that there are $H^{\prime}%
$--almost invariant subsets $\overline{X}$ and $\overline{Y}$ of $DG$ such
that $\overline{X}\cap G$ equals $X$ and $\overline{Y}\cap G$ equals $Y$. As
$X$ and $Y$ cross weakly, it follows that $\overline{X}$ and $\overline{Y}$
also cross weakly. But Proposition 7.4 of \cite{SS2} implies that no almost
invariant subset of $DG$ can cross $\overline{X}$ weakly, as $H^{\prime}$ has
only $2$ coends in $DG.$ This contradiction shows that $v$ cannot be of
commensuriser type.

If $H$ has length $n$, there are three cases depending on the index of $H$ in
$Comm_{G}(H)$. In all three cases, we have case 1d) of Theorem
\ref{mainresult}.

If this index is finite, Lemma \ref{V0-vertexofsolidtorustype} shows that $v$
is of solid torus type.

If this index is infinite and $Comm_{G}(H)$ is $VPC(n+1)$, Lemma
\ref{V0-vertexoftorustype} shows that $v$ is of torus type.

In the remaining case, Lemma \ref{V0-vertexofSeiferttype} shows that $v$ is of
Seifert type.

We have now shown that each $V_{0}$--vertex of the uncompleted graph of groups
$\Gamma_{n,n+1}$ satisfies part 1) of the theorem.

2) A $V_{0}$--vertex of $\Gamma_{n,n+1}^{c}$ which is obtained by the
completion process arises from a $V_{1}$--vertex $w$ of $\Gamma_{n,n+1}$. It
suffices to show that $w$ is of special Seifert type or of special solid torus
type. Recall that either $G(w)$ is $VPC(n+1)$ and $w$ has a single incident
edge $e$ with $G(e)$ of index $2$ in $G(w)$, or $G(w)$ is $VPCn$ and $w$ has a
single incident edge $e$ with $G(e)$ of index $2$ or $3$ in $G(w)$, or $G(w)$
is $VPCn$ and $w$ has three incident edges each carrying $G(w)$. We consider
each case separately.

\medskip

\textbf{Case: }$G(w)$\textbf{ is }$VPC(n+1)$,\textbf{ and }$w$\textbf{ has a
single incident edge }$e$\textbf{ with }$G(e)$\textbf{ of index }$2$\textbf{
in }$G(w)$\textbf{.}

\medskip

Lemma \ref{edgesplittingsaredualtoannuliortori} shows that the edge splitting
associated to $e$ is dual to an essential torus. It follows that $w$ is of
special Seifert type.

In the remaining cases, each edge incident to $w$ carries a $VPCn$ group, so
that Lemma \ref{edgesplittingsaredualtoannuliortori} shows that the associated
edge splitting is dual to an essential annulus.

\medskip

\textbf{Case: }$G(w)$\textbf{ is }$VPCn$,\textbf{ and }$w$\textbf{ has a
single incident edge }$e$\textbf{ with }$G(e)$\textbf{ of index }$2$\textbf{
or }$3$\textbf{ in }$G(w)$\textbf{.}

\medskip

Let $H$ denote $G(e)$, let $A$ denote the annulus associated to $e$, and let
$\partial H$ denote the subgroup of $H$ carried by $\partial A$. We claim that
$A$ is untwisted. Assuming this claim, it follows that $\partial H$ equals $H$
and so has index $2$ or $3$ in $G(w)$, which implies that $w$ is of special
solid torus type. It remains to prove the claim.

As usual we choose an aspherical space $M$ with fundamental group $G$ and with
aspherical subspaces corresponding to $\partial G$ whose union is denoted
$\partial M$. Let $M_{H}$ denote the cover of $M$ with fundamental group equal
to $H$. If $A$ is twisted, then its lift into $M_{H}$ has boundary in a single
component $\Sigma$ of $\partial M_{H}$ whose fundamental group must equal
$\partial H$. First suppose that $H$ is normal in $G(w)$, so that the quotient
$G(w)/H$ acts on $M_{H}$. Thus $\partial M_{H}$ has two or three distinct
components whose fundamental group equals $\partial H$, and each of these
components contains the boundary of an essential twisted annulus. The double
cover $M_{\partial H}$ must have two boundary components above each such
component of $\partial M_{H}$, each with fundamental group equal to $\partial
H$, giving a total of four or six such boundary components of $M_{\partial H}%
$. Part 3) of Proposition \ref{annuligenerate} now shows that $H$ has at least
$4$ coends in $G$. This implies that $\Gamma_{n,n+1}$ has a $V_{0}$--vertex
$v$ of commensuriser type with $G(v)=Comm_{G}(H)$. But $w$ is a $V_{1}%
$--vertex of $\Gamma_{n,n+1}$ and $G(w)$ commensurises $H$, so that
$G(w)\subset G(v)$. As $w$ has a single incident edge $e$ and $G(e)\neq G(w)$,
this is impossible. This contradiction shows that if $H$ is normal in $G(w)$,
then the annulus $A$ must be untwisted, as claimed. Now suppose that $H$ is
not normal in $G(w)$, so that $H$ has index $3$ in $G(w)$. Let $H_{1}$ denote
the intersection of the conjugates of $H$ in $G(w)$. Thus $H_{1}$ is a normal
subgroup of $G(w)$ of index some power of $3$. If $A_{1}$ denotes the cover of
$A$ with fundamental group $H_{1}$, the fact that this cover has odd degree
implies that $A_{1}$ is itself a twisted annulus in $M_{1}$. Now the preceding
argument yields a contradiction, showing that the annulus $A$ must be
untwisted as claimed. This completes the proof that $w$ must be of special
solid torus type when $G(w)$ is $VPCn$ and $w$ has a single incident edge $e$
with $G(e)$ of index $2$ or $3$ in $G(w)$.

\medskip

\textbf{Case: }$G(w)$\textbf{ is }$VPCn$,\textbf{ and }$w$\textbf{ has three
incident edges }$e_{1}$\textbf{, }$e_{2}$\textbf{ and }$e_{3}$\textbf{, each
carrying }$G(w)$\textbf{.}

\medskip

Let $K$ denote $G(w)$, and let $A_{i}$ denote the annulus associated to the
edge $e_{i}$. As in the preceding case, we claim that each $A_{i}$ is
untwisted. Assuming this claim, it follows that the boundary of each $A_{i}$
carries $K$, which implies that $w$ is of special solid torus type. It remains
to prove the claim.

As $K$ is a torsion free $VPCn$ group it is also $PDn$. Thus an annulus with
fundamental group $K$ is untwisted if and only if $K$ is orientable.

Now suppose that $K$ is non-orientable, so that each $A_{i}$ is twisted. As
usual we choose an aspherical space $M$ with fundamental group $G$ and with
aspherical subspaces corresponding to $\partial G$ whose union is denoted
$\partial M$. Let $M_{K}$ denote the cover of $M$ with fundamental group equal
to $K$, and let $M_{0}$ denote the cover of $M$ with fundamental group equal
to $K_{0}$, the orientable subgroup of $K$ of index $2$. The lift of each
annulus $A_{i}$ into $M_{K}$ has boundary in a single component $\Sigma_{i}$
of $\partial M_{K}$ whose fundamental group must equal $K_{0}$. As the $A_{i}%
$'s determine non-conjugate splittings of $G$, the $\Sigma_{i}$'s must be
distinct. As above, the double cover $M_{0}$ of $M_{K}$ has two boundary
components above each $\Sigma_{i}$ each with fundamental group equal to
$K_{0}$. It follows that $K$ has at least $6$ coends in $G$, which implies
that $\Gamma_{n,n+1}$ has a $V_{0}$--vertex $v$ of commensuriser type with
$G(v)=Comm_{G}(K)$. Further any almost invariant subset of $G$ over a subgroup
commensurable with $K$ must be enclosed by $v$. In particular, for each $i$,
the $K$--almost invariant subset of $G$ associated to $e_{i}$ is enclosed by
$v$. This implies that there is a path $\lambda_{i}$ in $\Gamma_{n,n+1}$ with
$v$ at one end and with the edge $e_{i}$ at the other end, and each interior
vertex of $\lambda_{i}$ is isolated. As $\Gamma_{n,n+1}$ is reduced bipartite,
each $\lambda_{i}$ contains at most two edges. As $v$ is a $V_{0}$--vertex and
$w$ is a $V_{1}$--vertex, it follows that each $\lambda_{i}$ consists of the
single edge $e_{i}$. Thus each $e_{i}$ has $v$ and $w$ as its endpoints. Let
$\Gamma_{K}$ denote the subgraph of $\Gamma_{n,n+1}$ which consists of the
union of the $e_{i}$'s. As each $e_{i}$ and $w$ has associated group $K$, and
$G(v)$ commensurises $K$, it follows that the group carried by $\Gamma_{K}$
also commensurises $K$. But as $G(v)$ is the full commensuriser of $K$, this
is impossible. This contradiction shows that $K$ must be orientable, so that
each $A_{i}$ is untwisted, as claimed. This completes the proof that $w$ must
be of special solid torus type when $G(w)$ is $VPCn$ and $w$ has three
incident edges each carrying $G(w)$, and so completes the proof of part 2) of
the theorem.

3) This follows immediately from parts 1) and 2), as the definitions of the
various types of $V_{0}$--vertex in the statements of parts 1) and 2) all
contain the requirement that the edge splittings be dual to an essential
annulus or torus.

4) Consider a nontrivial almost invariant subset of $G$ over a $VPC(n+1)$
group and adapted to $\partial G$. We need to show that this set is enclosed
by some $V_{0}$--vertex of $\Gamma_{n,n+1}$, and by some $V_{0}$--vertex of
$\Gamma_{n,n+1}^{c}$. By definition, any such set lies in the family
$\mathcal{F}^{\prime}$. Now Theorem \ref{FandF'havesameregnbhds} shows that
$\Gamma_{n,n+1}$ equals the algebraic regular neighbourhood of the family
$\mathcal{F}^{\prime}$. Thus any element of $\mathcal{F}^{\prime}$ is enclosed
by some $V_{0}$--vertex of $\Gamma_{n,n+1}$. The construction of the
completion $\Gamma_{n,n+1}^{c}$ of $\Gamma_{n,n+1}$ shows that any element of
$\mathcal{F}^{\prime}$ is also enclosed by some $V_{0}$--vertex of
$\Gamma_{n,n+1}^{c}$. This completes the proof of part 4).

5) Let $H$ be a $VPC(n+1)$ subgroup of $G$ which is not conjugate into
$\partial G$. We need to show that $H$ is conjugate into a $V_{0}$--vertex
group of $\Gamma_{n,n+1}^{c}$. First note that as $G$ is torsion free, so is
$H$. Thus $H$ must be $PD(n+1)$.

Suppose that $H$ is orientable. The hypothesis that $H$ is not conjugate into
$\partial G$ implies that $H$ is an essential torus in $(G,\partial G)$, so
that there is a nontrivial $H$--almost invariant subset $X$ of $G$ dual to
$H$. Now Theorem \ref{FandF'havesameregnbhds} shows that $X$ is enclosed by a
$V_{0}$--vertex of $\Gamma_{n,n+1}$, so that $H$ is conjugate into a $V_{0}%
$--vertex group of $\Gamma_{n,n+1}$, and hence is also conjugate into a
$V_{0}$--vertex group of $\Gamma_{n,n+1}^{c}$ as required.

Now we will suppose that $H$ is non-orientable and is not conjugate into any
$V_{0}$--vertex group of $\Gamma_{n,n+1}$. The proof of part 4) of Theorem
\ref{torusdecompofapair} shows that the commensuriser $K$ of $H$ in $G$ is
itself a non-orientable $PD(n+1)$ subgroup of $G$. Let $K_{0}$ denote the
orientation subgroup of $K$. Then $K_{0}$ is a maximal torus subgroup of $G$,
and the proof of part 4) of Theorem \ref{torusdecompofapair} shows that it is
not conjugate into $\partial G$. As in the preceding paragraph, it follows
that there is a $V_{0}$--vertex $v$ of $\Gamma_{n,n+1}$, so that $K_{0}$ is
conjugate into $G(v)$. Our assumption that $H$ is not conjugate into any
$V_{0}$--vertex group of $\Gamma_{n,n+1}$ implies that the same is true for
$K$. As $K$ contains $K_{0}$ with finite index, there is a vertex $w$ of
$T_{n+1}(G,\partial G)$ such that $K$ is conjugate into $G(w)$, and $w$ must
be a $V_{1}$--vertex. Hence there is an edge $e$ of $\Gamma_{n,n+1}$ which is
incident to $w$ such that $G(e)$ contains a conjugate of $K_{0}$. As all the
edge groups of $\Gamma_{n,n+1}$ are annulus or torus groups, the group $G(e)$
must equal this conjugate of $K_{0}$. Let $\sigma$ denote the splitting of $G$
over $K_{0}$ determined by the edge $e$, so that $\sigma$ is dual to an
essential torus in $G$, and let $L$ denote the vertex group of $\sigma$ which
contains $G(w)$. If $\partial G$ is empty, Theorem 8.1 of \cite{B-E} shows
that the pair formed by $L$ and one or two copies of $K_{0}$ is $PD(n+2)$,
where there will be two copies of $K_{0}$ if $\sigma$ is a HNN extension and
only one copy otherwise. In general, as discussed just before Definition
\ref{defnofatoroidal}, the pair becomes $PD(n+2)$ when some family of groups
in $\partial G$ is added to the copies of $K_{0}$. As $L$ contains a conjugate
of $K$, one of the copies of $K_{0}$ in $\partial L$ is not equal to its own
commensuriser in $L$. Thus Lemma \ref{boundarygroupismaximal} implies that $L$
contains this copy of $K_{0}$ with index $2$, and $\partial L$ consists only
of $K_{0}$. In particular it follows that $\sigma$ must be an amalgamated free
product and not a HNN extension. As $L$ contains a conjugate of $K$, and both
$K\ $and $L$ contain $K_{0}$ with index $2$, it follows that $G(w)$ is equal
to $L$ and must be a conjugate of $K$. Hence the edge $e$ of $\Gamma_{n,n+1}$
is separating, and if we remove the interior of $e$ from $\Gamma_{n,n+1}$,
then the component $\Gamma_{w}$ of the resulting subgraph which contains $w$
must carry the group $G(w)$. Now the minimality of $\Gamma_{n,n+1}$ implies
that $\Gamma_{w}$ consists solely of $w$, so that $w$ has valence $1$ in
$\Gamma_{n,n+1}$. It follows that $w$ becomes a $V_{0}$--vertex in the
completion $\Gamma_{n,n+1}^{c}$, so that $K$, and hence $H$, is conjugate into
a $V_{0}$--vertex group of $\Gamma_{n,n+1}^{c}$, as required.
\end{proof}

We observe the following result which again shows the similarity between the
algebra in this paper and the topology of $3$-manifolds.

\begin{lemma}
Let $(G,\partial G)$ be an orientable $PD(n+2)$ pair such that $G$ is not
$VPC$. Suppose that $X$ is a $n$--canonical almost invariant subset of $G$
which is over a $VPC(n+1)$ group $H$, and that $H$ intersects some group in
$\partial G$ in a $VPCn$ group $L$. Then $X$ is isolated in $\mathcal{F}%
_{n,n+1}$.
\end{lemma}

\begin{remark}
If $n=1$ and $M$ is an orientable Haken $3$--manifold, the corresponding
result holds. For $X$ corresponds to an essential torus $T$ in $M$ which
crosses no essential annulus in $M$, such that $\pi_{1}(T)$ intersects the
fundamental group of some boundary component of $M$ in an infinite cyclic
subgroup. This second condition implies that $T$ must be homotopic into a
component $W$ of $T(M)$ which meets $\partial M$, and the fact that $T$
crosses no essential annulus in $M$ now implies that $T$ must be homotopic
into a torus component of the frontier of $W$ in $M$. In particular, $T$
crosses no essential annulus or torus in $M$.
\end{remark}

\begin{proof}
As usual, Lemma 13.1 of \cite{SS2} tells us that $H$ has a subgroup
$H^{\prime}$ of finite index which normalises a subgroup $L^{\prime}$ of
finite index in $L$. Also as usual, we choose an aspherical space $M$ with
fundamental group $G$ and with aspherical subspaces corresponding to $\partial
G$ whose union is denoted $\partial M$. Let $M^{\prime}$ denote the covering
space of $M$ with fundamental group $L^{\prime}$. The hypothesis that $H$
intersects some group in $\partial G$ in the group $L$ implies that there is a
component of $\partial M^{\prime}$ with fundamental group $L^{\prime}$. The
action of $H^{\prime}/L^{\prime}$ on $M^{\prime}$ yields an infinite family of
distinct such boundary components. Now the proof of part 3) of Lemma
\ref{annuligenerate} shows that $e(G,L^{\prime})$ is infinite.

As $H^{\prime}$ normalises $L^{\prime}$, it also follows that $L^{\prime}%
\ $has large commensuriser. Thus there is a $V_{0}$--vertex $v$ of
$\Gamma_{n,n+1}$ of commensuriser type such that $G(v)=Comm_{G}(L^{\prime})$,
so that $G(v)$ must contain $H^{\prime}$. Also $\Gamma_{n,n+1}$ has a $V_{0}%
$--vertex $w$ which encloses $X$.

First we suppose that $v$ and $w$ are distinct. It follows that there is a
path joining $v$ and $w$ such that $H^{\prime}$ lies in the edge group of each
edge on the path. Let $e$ be an edge on this path. Recall that all edge groups
of $\ \Gamma_{n,n+1}$ are $VPCn$ or $VPC(n+1)$. Thus $G(e)$ must be $VPC(n+1)$
and must contain $H^{\prime}$ with finite index. Let $Y$ denote the almost
invariant subset of $G$ associated to the edge splitting of $G$ given by $e$.
Lemma \ref{edgesplittingsaredualtoannuliortori} tells us that $Y$ is adapted
to $\partial G$. Now Corollary \ref{FiscontainedinF'} tells us that $X$ is
also adapted to $\partial G$. But a $H^{\prime}$--almost invariant subset of
$G$ which is adapted to $\partial G$ is unique up to equivalence and
complementation, so it follows that $X$ must be equivalent to $Y$ or $Y^{\ast
}$. As $Y$ is isolated in $\mathcal{F}_{n,n+1}$, it follows that $X$ is also
isolated in $\mathcal{F}_{n,n+1}$.

Now suppose that $v=w$. Recall from Theorem \ref{mainresult} that any $V_{0}%
$--vertex of $\Gamma_{n,n+1}$ of commensuriser type must be of Seifert type,
of torus type, or of solid torus type. The last case cannot occur here, as
$G(v)$ contains the $VPC(n+1)$ group $H^{\prime}$. Thus $v$ is of Seifert
type, or of torus type. In the second case, $G(v)$ is $VPC(n+1)$ and so must
contain $H^{\prime}$ with finite index. Thus any almost invariant set enclosed
by $v$ and corresponding to a torus must be equivalent to $X$ or to $X^{\ast}%
$. Now the hypothesis that $X$ crosses no essential annulus in $M$ implies
that $X$ crosses no element of the CCC of $\mathcal{F}_{n,n+1}$ which gives
rise to $v$. As $X$ is enclosed by $v$, it must be associated to an edge
splitting of $\Gamma_{n,n+1}$ associated to an edge incident to $v$, and so
$X$ must be isolated in $\mathcal{F}_{n,n+1}$. Finally if $v$ is of Seifert
type, we use the facts that $v$ is of $VPCn$--by--Fuchsian type, and that tori
enclosed by $v$ correspond to loops in the base orbifold $X_{v}$, and annuli
enclosed by $v$ correspond to arcs in $X_{v}$. Now a loop in $X_{v}$ which
crosses no arc must be peripheral in $X_{v}$. It follows again that $X$ must
be associated to an edge splitting of $\Gamma_{n,n+1}$ associated to an edge
incident to $v$, and so $X$ must be isolated in $\mathcal{F}_{n,n+1}$. This
completes the proof of the lemma.
\end{proof}

\section{Comparing the decompositions of a $PD(n+2)$ pair and its
double\label{comparingdecompositions}}

In this section, we will apply Theorem \ref{mainresult} to understand the
effect of doubling on our decompositions $\Gamma_{n,n+1}(G)$ and
$\Gamma_{n,n+1}^{c}(G)$ of a $PD(n+2)$ pair $(G,\partial G)$.

First we need to improve our description of these decompositions. The same
description suffices for both decompositions, so we will only work with
$\Gamma_{n,n+1}^{c}(G)$. In section \ref{essentialannuli}, we used aspherical
spaces to clarify the concept of an essential annulus in a $PD(n+2)$ pair
$(G,\partial G)$. Now we need to greatly refine those ideas in order to
clarify "how vertices of $\Gamma_{n,n+1}^{c}(G)$ meet $\partial G$". Let
$\partial G=\{S_{1},...,S_{m}\}$, and recall that in section
\ref{essentialannuli}, we used a mapping cylinder construction to make a
$K(G,1)$ with $K(S_{i},1)$'s as disjoint subspaces.

Given any group $G$ and a graph of groups decomposition $\Gamma(G)$ of $G$,
there is a general construction of an aspherical space $X$ with $\pi_{1}(X)=G$
whose structure mimics the graph of groups $\Gamma(G)$. This is called a graph
of spaces. For each vertex group $G(v)$ of $\Gamma(G)$, we choose a
corresponding aspherical space $K(G(v),1)$, and for each edge group $G(e)$ we
choose a $K(G(e),1)$ and then take its product with the unit interval $I$. We
construct $X$ from the disjoint union of all the $K(G(v),1)$ and
$K(G(e),1)\times I$ by gluing each end of each $K(G(e),1)\times I$ to the
appropriate $K(G(v),1),$ by a map inducing the appropriate inclusion of
fundamental groups. Thus there is a natural map from $X$ to $\Gamma(G)$
defined by collapsing each $K(G(v),1)$ and each $K(G(e),1)$ to a point. We
will apply this construction to the decomposition $\Gamma_{n,n+1}^{c}(G)$.
Further, for each edge $e$ of $\Gamma_{n,n+1}^{c}(G)$ whose associated
splitting is dual to an annulus $\Lambda$ with fundamental group $H$, we will
choose our $K(G(e),1)$ to be $\Lambda$, which is an $I$--bundle over a
$K(H,1)$, as discussed in section \ref{essentialannuli}. Thus $\partial
\Lambda$, the boundary of the annulus $\Lambda$, is the corresponding
$\partial I$--bundle over $K(H,1)$.

Next we will carry out a similar construction of an aspherical space $Y$ which
represents $\partial G$. Each boundary component of each annulus associated to
an edge splitting of $\Gamma_{n,n+1}^{c}(G)$ induces a splitting of one of the
$S_{i}$'s. Thus each $S_{i}$ can be decomposed as a graph of groups structure
with edges corresponding to boundary components of these annuli. For each
$S_{i}$, we make a corresponding graph of spaces construction of a space
$Y_{i}$, with $\pi_{1}(Y_{i})=S_{i}$. Further, if an edge $e$ of this graph of
groups decomposition of $S_{i}$ has associated splitting over a boundary
component of the annulus $\Lambda$, we will choose the corresponding
$K(G(e),1)$ to be homeomorphic to that component of $\partial\Lambda$. Let $Y$
denote the disjoint union of the $Y_{i}$'s. Note that as $S_{i}$ is a
$PD(n+1)$ group and these splittings are over $PDn$ groups, it follows from
\cite{B-E} that each vertex space of $Y$ is naturally a $PD(n+1)$ pair.

Finally we combine the above constructions using a mapping cylinder
construction as follows. We take the disjoint union of the space $X$
constructed above from $\Gamma_{n,n+1}^{c}(G)$ with the product $Y\times I$,
and glue $Y\times\{0\}$ to $X$ so that each edge space glues by a
homeomorphism to a boundary component of the appropriate annulus $\Lambda$,
and each vertex space is glued to the appropriate vertex space of $X$. Denote
the resulting space by $M$, and denote $Y\times\{1\}$ by $\partial M$. We now
have a picture which mimics the topology of the $JSJ$ decomposition of a
$3$--manifold, as each edge splitting of $\Gamma_{n,n+1}^{c}(G)$ which is dual
to an annulus $\Lambda$ is represented by an embedding of $\Lambda$ in $M$
with $\partial\Lambda$ embedded in $\partial M$. Note that $M\ $has a natural
projection $p$ to $\Gamma_{n,n+1}^{c}(G)$ obtained by collapsing the
constituent edge and vertex spaces to a point.

In order to complete the analogy with $3$--manifold topology, we proceed as
follows. We subdivide $\Gamma_{n,n+1}^{c}(G)$ by adding a new vertex at the
middle of each edge $e$, and then for each vertex $v$ of $\Gamma_{n,n+1}%
^{c}(G)$, we let $M_{v}$ denote the pre-image under $p$ of the star of $v$ in
this subdivision. Let $\partial_{1}M_{v}$ denote the intersection of $M_{v}$
with the pre-image under $p$ of the new vertices, let $\partial_{0}M_{v}$
denote $M_{v}\cap\partial M$, and let $\partial M_{v}$ denote the union
$\partial_{0}M_{v}\cup\partial_{1}M_{v}$. Each component of $\partial_{1}%
M_{v}$ is an annulus or torus in $M$. Each component of $\partial_{0}M_{v}$ is
either a component of $\partial M$, or is naturally a $PD(n+1)$ pair with
boundary equal to the boundary of some annuli in $M$. Thus we have finally
assigned meaning to the "intersection with $\partial G$ of a vertex of
$\Gamma_{n,n+1}^{c}(G)$".

If $v$ is a $V_{0}$--vertex of $\Gamma_{n,n+1}^{c}(G)$ of interior Seifert
type, of special Seifert type, or is isolated with torus group, then
$\partial_{1}M_{v}$ consists entirely of essential tori in $(G,\partial G)$,
and $\partial_{0}M_{v}$ is empty. Thus the "intersection of $v$ with $\partial
G$" is empty.

If $v$ is a $V_{0}$--vertex of $\Gamma_{n,n+1}^{c}(G)$ of $I$--bundle type,
then $\partial_{1}M_{v}$ consists entirely of essential annuli in $(G,\partial
G)$, and $\partial_{0}M_{v}$ consists of one or two $PD(n+1)$ pairs the union
of whose boundary components is equal to the boundary of the annuli forming
$\partial_{1}M_{v}$. If $\partial_{0}M_{v}$ consists of two $PD(n+1)$ pairs,
each includes into $M_{v}$ by an isomorphism of fundamental groups, so that
the $I$--bundle is trivial. Otherwise, $\partial_{0}M_{v}$ consists of a
single $PD(n+1)$ pair such that the image of $\pi_{1}(\partial_{0}M_{v})$ in
$\pi_{1}(M_{v})$ has index $2$, so that the $I$--bundle is twisted. Thus the
"intersection of $v$ with $\partial G$" is the $\partial I$--bundle associated
to the $I$--bundle.

We note that if $\Sigma$ is a torus in $\partial G$, then it must be conjugate
into some vertex group of $\Gamma_{n,n+1}(G)$ or of $\Gamma_{n,n+1}^{c}(G)$.
Otherwise there is an essential annulus in $(G,\partial G)$ with a boundary
component in $\Sigma$, but then the $V_{0}$--vertex $v$ of $\Gamma_{n,n+1}%
^{c}(G)$ which encloses that annulus must be of commensuriser type so that
$G(v)$ contains $\Sigma$.

Now we can describe exactly what we mean by doubling $\Gamma_{n,n+1}(G)$ or
$\Gamma_{n,n+1}^{c}(G)$. Again the same description suffices for both. Let
$DM$ denote the space obtained by doubling $M$ along $\partial M$, so that
$DM$ is a $K(DG,1)$. We denote the two copies of $M$ in $DM$ by $M\ $and
$\overline{M}$, and let $\tau$ denote the involution of $DM$ which
interchanges $M$ and $\overline{M}$. We will describe a family of disjoint
tori in $DM$ which determines the decomposition of $DG\ $which we want. The
annuli and tori in $M$ which correspond to the edges of $\Gamma_{n,n+1}%
^{c}(G)$ determine tori in $DM$, as follows. A torus $T$ in $M$ yields two
tori $T\ $and $\tau T$ in $DM$, and an annulus $\Lambda$ in $M\ $yields a
torus $D\Lambda=\Lambda\cup\tau\Lambda$ in $DM$. In addition, for each torus
$\Sigma$ in $\partial G$ which is enclosed by a $V_{1}$--vertex of
$\Gamma_{n,n+1}^{c}(G)$, we add two parallel copies of the corresponding
component of $\partial M$, one copy in $M$ and the other in $\overline{M}$.
Clearly the tori in this family are all disjoint. Now this family of disjoint
tori in $DM$ determines a graph of groups structure of $DG$, which we denote
by $D\Gamma_{n,n+1}^{c}$. Thus all the edges of $D\Gamma_{n,n+1}^{c}$ have
associated splittings dual to tori in $DG$.

There is a natural map from $DM$ to $D\Gamma_{n,n+1}^{c}$, and it is easy to
describe the vertex spaces of $DM$. If $\Sigma$ is a torus in $\partial G$
which is conjugate into a $V_{1}$--vertex group of $\Gamma_{n,n+1}^{c}(G)$,
the two parallel copies of $\Sigma$ in $DM$ together bound a copy of
$\Sigma\times I$, which corresponds to an isolated vertex of $D\Gamma
_{n,n+1}^{c}$. We label such a vertex as a $V_{0}$--vertex of $D\Gamma
_{n,n+1}^{c}$. The other vertex spaces of $DM$ arise from vertex spaces of
$M$. If $M_{v}$ is disjoint from $\partial M$, so that $\partial_{0}M_{v}$ is
empty, there are two corresponding vertex spaces $M_{v}$ and $\tau M_{v}$ of
$DM$ each homeomorphic to $M_{v}$. If $v$ is a $V_{1}$--vertex of
$\Gamma_{n,n+1}^{c}(G)$ such that $\partial_{0}M_{v}$ consists only of torus
components of $\partial M$, there are again two corresponding vertex spaces,
each homeomorphic to $M_{v}$. If $v$ is a $V_{1}$--vertex of $\Gamma
_{n,n+1}^{c}(G)$ such that $\partial_{0}M_{v}$ is not empty and does not
consist only of torus components of $\partial M$, there is one corresponding
vertex space obtained by doubling $M_{v}$ along the non-torus components of
$\partial_{0}M_{v}$. Finally, if $v$ is a $V_{0}$--vertex of $\Gamma
_{n,n+1}^{c}(G)$ such that $\partial_{0}M_{v}$ is non-empty, there is one
corresponding vertex space obtained by doubling $M_{v}$ along $\partial
_{0}M_{v}$. We define each of the corresponding vertices of $D\Gamma
_{n,n+1}^{c}$ to be of type $V_{0}$ or $V_{1}$ so as to be of the same type as
$v$. With this labelling, $D\Gamma_{n,n+1}^{c}$ is bipartite. Note that an
isolated vertex of $D\Gamma_{n,n+1}^{c}$ either arises from a torus in
$\partial G$ which is conjugate into a $V_{1}$--vertex group of $\Gamma
_{n,n+1}^{c}(G)$, or it arises from an isolated annulus or torus vertex of
$\Gamma_{n,n+1}^{c}(G)$. As $\Gamma_{n,n+1}^{c}(G)$ is reduced, it follows
that $D\Gamma_{n,n+1}^{c}$ is also reduced.

Next we need to consider more detail about the structure of the $V_{0}%
$--vertex groups of $\Gamma_{n,n+1}^{c}(G)$. Recall that if $v$ is a $V_{0}%
$--vertex of $\Gamma_{n,n+1}^{c}(G)$ of interior Seifert type, then $v$ is of
$VPCn$--by--Fuchsian type. Let $L$ denote the $VPCn$ normal subgroup of
$G(v)$. Then the quotient group $G(v)/L$ is not virtually cyclic and is the
orbifold fundamental group of a compact $2$--orbifold $X_{v}$. Further there
is exactly one edge of $\Gamma_{n,n+1}^{c}(G)$ which is incident to $v$ for
each peripheral subgroup $K$ of $G(v)$, and this edge carries $K$. Thus there
is a natural projection of $M_{v}$ to $X_{v}$, in which $\partial
M_{v}=\partial_{1}M_{v}$ maps onto $\partial X_{v}$. This precisely mirrors
the picture in a $3$--manifold of an interior Seifert fibre space component of
the characteristic submanifold. We will show that a similar picture occurs for
any $V_{0}$--vertex $v$ of $\Gamma_{n,n+1}^{c}(G)$ of commensuriser type,
again mirroring the situation for $3$--manifolds. Thus, in all these cases,
there is a compact $2$--orbifold $X_{v}$ and a natural projection of $M_{v}$
to $X_{v}$ in which $\partial M_{v}$ maps onto $\partial X_{v}$. Note that any
$V_{0}$--vertex $v$ of $\Gamma_{n,n+1}^{c}(G)$ of commensuriser type is of
peripheral type, so that $\partial_{0}M_{v}$ is non-empty, which introduces
some new aspects to the discussion. Note also that if $v$ is of peripheral
Seifert type, then $G(v)$ is $VPCn$--by--Fuchsian and so $v$ has a natural
base orbifold $X_{v}$, but even this is not clear if $v$ is of torus type or
of solid torus type.

Consider a $V_{0}$--vertex of $\Gamma_{n,n+1}^{c}(G)$ of commensuriser type.
For a given such vertex $v$, there is a $VPCn$ subgroup $L$ of $G(v)$ such
that for each edge of $\Gamma_{n,n+1}^{c}(G)$ which is incident to $v$ and
associated to a splitting over an annulus $\Lambda$, the group carried by each
component of $\partial\Lambda$ is $L$. It follows that $\partial_{0}M_{v}$
consists of a disjoint union of torus components of $\partial M$ and of annuli
whose boundary components carry $L$. Thus each component of $\partial M_{v}$
is either a torus component of $\partial_{0}M_{v}$ or of $\partial_{1}M_{v}$,
or is a union of annuli in $\partial_{0}M_{v}$ and $\partial_{1}M_{v}$.

Let $T$ be a component of $\partial M_{v}$ which is a union of annuli in
$\partial_{0}M_{v}$ and $\partial_{1}M_{v}$. Either all the annuli in $T\ $are
untwisted and glued in a circular pattern, or there are precisely two twisted
annuli in $T$ separated by a chain of untwisted annuli. Note that all these
annuli carry subgroups of $G$, and so have torsion free $VPCn$ fundamental
group. It follows that $\pi_{1}(T)$ is a torsion free $VPC(n+1)$ group and so
is $PD(n+1)$. We claim that $\pi_{1}(T)$ is an orientable $PD(n+1)$ group. By
our construction of $M$, we know that $T$ has a neighbourhood homeomorphic to
an $I$--bundle over $T$. Further this $I$--bundle must be trivial as $T$ is a
boundary component of $M_{v}$. Now it follows that there is an excision
isomorphism between $H_{n+2}(T\times I,T\times\partial I)$ and $H_{n+2}%
(M,\partial M)\cong\mathbb{Z}$. As $H_{n+2}(T\times I)=0$, and $H_{n+1}%
(T\times\partial I)\cong H_{n+1}(T)\oplus H_{n+1}(T)$, it follows that
$H_{n+1}(T)$ is non-zero, which implies that $\pi_{1}(T)$ is an orientable
$PD(n+1)$ group, as required. We will abuse terminology and say that $T$ is a torus.

Note that an untwisted annulus has a natural projection to the unit interval
$I$, and a twisted annulus has a natural projection to the $1$--orbifold $Q$
which is the quotient of $I$ by a reflection involution. Hence if the torus
$T$ is a union of untwisted annuli, it has a natural map to the circle $S^{1}$
such that the restriction to each annulus is projection to a unit interval.
And if $T$ contains two twisted annuli, it has a natural map to the
$1$--orbifold $C$ which is the quotient of $S^{1}$ by a reflection involution
such that the restriction to each untwisted annulus is projection to a unit
interval, and the restriction to each twisted annulus is projection to a copy
of $Q$. Now recall that our aim is to show that if $v$ is a $V_{0}$--vertex of
$\Gamma_{n,n+1}^{c}(G)$ of commensuriser type, then there is a compact
$2$--orbifold $X_{v}$ and a natural projection of $M_{v}$ to $X_{v}$ in which
$\partial M_{v}$ maps onto $\partial X_{v}$. We will show further that this
projection can be chosen so that its restriction to each annulus component of
$\partial_{0}M_{v}$ and of $\partial_{1}M_{v}$ is the natural projection to a
$1$--suborbifold of $\partial X_{v}$ which is isomorphic to $I$ or $Q$, as
appropriate.$\pi$

We start by considering the case when $v$ is a $V_{0}$--vertex of
$\Gamma_{n,n+1}^{c}(G)$ of torus type, so that $G(v)$ is $VPC(n+1)$ and hence
$PD(n+1)$.

\begin{lemma}
\label{torustypevertexhasbaseorbifold}Let $v$ be a $V_{0}$--vertex of
$\Gamma_{n,n+1}^{c}(G)$ of torus type. Then the following hold:

\begin{enumerate}
\item If $G(v)$ is an orientable $PD(n+1)$ group, then $\partial M_{v}$
consists of two tori such that each includes into $M_{v}$ inducing an
isomorphism of fundamental groups.

\item If $G(v)$ is a non-orientable $PD(n+1)$ group, and if $G(v)_{0}$ denote
its orientation subgroup, then $\partial M_{v}$ consists of a single torus
whose inclusion into $M_{v}$ induces an injection of fundamental groups with
image $G(v)_{0}$.

\item In either case, there is a compact $2$--orbifold $X_{v}$ such that
$G(v)$ is $L$--by--$\pi_{1}^{orb}(X_{v})$, and a natural map from $M_{v}$ to
$X_{v}$ such that $\partial M_{v}$ maps onto $\partial X_{v}$. Further, each
annulus in $\partial_{0}M_{v}$ and $\partial_{1}M_{v}$ maps to a
$1$--suborbifold of $\partial X_{v}$ by the natural map. If the annulus is
untwisted, its image in $\partial X_{v}$ is isomorphic to the unit interval
$I$, and if the annulus is twisted, its image in $\partial X_{v}$ is
isomorphic to $Q$, the quotient of $I$ by a reflection involution.
\end{enumerate}
\end{lemma}

\begin{proof}
1) The above discussion shows that $\partial M_{v}$ consists of tori. There is
an excision isomorphism between $H_{n+2}(M_{v},\partial M_{v})$ and
$H_{n+2}(M,\partial M)\cong\mathbb{Z}$, and the fundamental class of $\partial
M_{v}$ maps to zero in $H_{n+1}(M_{v})$. Now consider the long exact homology
sequence of the pair $(M_{v},\partial M_{v})$. As $G(v)$ is $VPC(n+1)$ and
orientable, we know that $H_{n+2}(M_{v})=0$, and $H_{n+1}(M_{v})\cong%
\mathbb{Z}$. As $H_{n+1}(T)\cong\mathbb{Z}$ for any torus $T$, it follows that
$\partial M_{v}$ consists of at most two tori.

In case 1) of Definition \ref{defnoftorustype}, there is a torus $\Sigma$ in
$\partial_{0}M_{v}$ whose inclusion into $M_{v}$ induces an isomorphism of
fundamental groups.

In case 2) of Definition \ref{defnoftorustype}, there is a torus $\Sigma$ in
$\partial_{1}M_{v}$ whose inclusion into $M_{v}$ induces an isomorphism of
fundamental groups.

It follows that in either of these cases, $\partial M_{v}$ consists of two
tori and that the inclusion into $M_{v}$ of the second torus $T$ induces an
isomorphism $H_{n+1}(T)\rightarrow H_{n+1}(M_{v})$. Thus the map from $\pi
_{1}(T)$ to $\pi_{1}(M_{v})$ is onto. As $\pi_{1}(T)$ is $VPC(n+1)$, the
kernel of this map must be finite. Now the fact that $\pi_{1}(T)$ is torsion
free implies that this map is an isomorphism, as required.

In case 3) of Definition \ref{defnoftorustype}, there is a torus contained in
$M_{v}$ whose fundamental group equals $G(v)$. We could further refine our
above construction of $(M,\partial M)$ to arrange that $M_{v}$ contains such a
torus $\Sigma$. This must separate $M_{v}$ into two pieces each with
fundamental group $G(v)$. By considering each of these pieces separately as in
the preceding paragraph, we can again conclude that each component of
$\partial M_{v}$ includes into $M_{v}$ by an isomorphism of fundamental groups.

2) If $G(v)$ is a non-orientable $PD(n+1)$ group, then $H_{n+1}(M_{v})=0$. Now
the long exact homology sequence of the pair $(M_{v},\partial M_{v})$ shows
that $\partial M_{v}$ consists of a single torus $T$. As $T$ is orientable,
the image of $\pi_{1}(T)$ in $G(v)$ is contained in $G(v)_{0}$. Thus there is
an index $2$ subgroup $G^{\prime}$ of $G$ whose intersection with $G(v)$ is
$G(v)_{0}$, and which is naturally a $PD(n+2)$ pair $(G^{\prime},\partial
G^{\prime})$. Let $M^{\prime}$ denote the corresponding model space for this
pair, and apply part 1) of the lemma to the appropriate vertex space of
$M^{\prime}$. This will imply that the inclusion of $T=\partial M_{v}$ into
$M_{v}$ induces an injection of fundamental groups with image $G(v)_{0}$, as required.

3) If $G(v)$ is orientable and $T$ denotes a torus in $\partial M_{v}$, then
the pair $(M_{v},\partial M_{v})$ is homotopy equivalent to $(T\times
I,T\times\partial I)$. If $T$ consists of a circular chain of untwisted
annuli, then the quotient $\pi_{1}(T)/L$ is isomorphic to $\mathbb{Z}$, and we
take the orbifold $X_{v}$ to be the annulus. If $T$ contains two twisted
annuli, then the quotient $\pi_{1}(T)/L$ is isomorphic to $\mathbb{Z}_{2}%
\ast\mathbb{Z}_{2}$, and we take the orbifold $X_{v}$ to be the product
$C\times I$, where $C$ is the quotient of $S^{1}$ by a reflection involution.

If $G(v)$ is non-orientable, then $\partial M_{v}$ consists of a single torus
$T$ whose inclusion into $M_{v}$ induces an injection of fundamental groups
with image $G(v)_{0}$. If $T$ consists of a circular chain of untwisted
annuli, then the quotient $\pi_{1}(T)/L$ is isomorphic to $\mathbb{Z}$, and we
take the orbifold $X_{v}$ to be the Moebius band. If $T$ contains two twisted
annuli, then the quotient $\pi_{1}(T)/L$ is isomorphic to $\mathbb{Z}_{2}%
\ast\mathbb{Z}_{2}$, and we take the orbifold $X_{v}$ to be "a twisted
$I$--bundle" over $C$, which we describe as follows. Recall that the Moebius
band is a twisted $I$--bundle over the circle $S^{1}$. Let $\sigma$ denote the
reflection involution of $S^{1}$ with quotient $C$. Then $\sigma$ extends to
an involution $\overline{\sigma}$ of this $I$--bundle over $S^{1}$, and
$X_{v}$ is the quotient of this action. Note that $\overline{\sigma}$
preserves each fibre over the fixed points of $\sigma$, and fixes one of these
fibres pointwise, while reflecting the other one. Thus the fixed set of
$\overline{\sigma}$ consists of an interval, where the local picture of
$\overline{\sigma}$ is a reflection, and of an isolated point, where the local
picture is of a rotation through $\pi$. Thus $X_{v}$ has underlying surface a
disc $D$, the boundary of $D$ contains a single mirror interval whose
complement is thus a copy of $C$, and $X_{v}$ also has an order $2$ cone point
in the interior of $D$.
\end{proof}

Now we understand how "$v$ meets $\partial G$" when $v$ is of torus type, we
can apply this to the case where $v$ is a $V_{0}$--vertex of $\Gamma
_{n,n+1}^{c}(G)$ of peripheral Seifert type to obtain the following result.

\begin{corollary}
\label{Seiferttypevertexhasbaseorbifold}Let $v$ be a $V_{0}$--vertex of
$\Gamma_{n,n+1}^{c}(G)$ of peripheral Seifert type, and let $X_{v}$ denote the
base orbifold of $v$. Then there is a natural map from $M_{v}$ to $X_{v}$ such
that $\partial M_{v}$ maps onto $\partial X_{v}$. Further, each annulus in
$\partial_{0}M_{v}$ and $\partial_{1}M_{v}$ maps to a $1$--suborbifold of
$\partial X_{v}$ with non-empty boundary, by the natural map.
\end{corollary}

\begin{proof}
Recall from Definition \ref{defnofSeiferttype} that $\Gamma_{n,n+1}^{c}(G)$
can be refined by splitting at $v$ to a graph of groups structure
$\Gamma^{\prime}$ of $G$ with a vertex $v^{\prime}$ of $\Gamma^{\prime}$ such
that $G(v^{\prime})=G(v)$ and $v^{\prime}$ is of Seifert type adapted to
$\partial G$. The projection map $\Gamma^{\prime}\rightarrow\Gamma$ sends
$v^{\prime}$ to $v$ and is an isomorphism apart from the fact that certain
edges incident to $v^{\prime}$ are collapsed to $v$. Further if $e$ is an edge
of $\Gamma^{\prime}$ which is incident to $v^{\prime}$ and collapsed to $v$,
then the other vertex $w$ of $e$ is of torus type. As $G(e)$ is a torus, $w$
is of torus type 2) in Definition \ref{defnoftorustype}. We can
correspondingly refine $M$ to obtain the vertex space $M_{v^{\prime}}$. As
$v^{\prime}$ is of Seifert type adapted to $\partial G$, there is a natural
projection of $M_{v^{\prime}}$ to $X_{v^{\prime}}$ with all the required
properties. It also follows from Lemma \ref{torustypevertexhasbaseorbifold}
that the base orbifold of $w$ is an annulus or $C\times I$. In turn this
implies that $X_{v}$ and $X_{v^{\prime}}$ are isomorphic. Now it follows that
there is a natural projection of $M_{v}$ to $X_{v}$ in which $\partial M_{v}$
maps onto $\partial X_{v}$, and each annulus component of $\partial_{0}M_{v}$
and of $\partial_{1}M_{v}$ maps in the natural way to a $1$--suborbifold of
$\partial X_{v}$ with non-empty boundary, as required.
\end{proof}

We have shown that for any $V_{0}$--vertex $v$ of $\Gamma_{n,n+1}^{c}(G)$ of
commensuriser type, but not of solid torus type, there is a base orbifold
$X_{v}$ such that $\partial X_{v}$ is divided into the image of $\partial
_{0}M_{v}$, which we denote by $\partial_{0}X_{v}$, and the image of
$\partial_{1}M_{v}$, which we denote by $\partial_{1}X_{v}$. This precisely
mirrors the picture in a $3$--manifold of a peripheral Seifert type or torus
type component of the characteristic submanifold. Now we can continue our
discussion of the double $D\Gamma_{n,n+1}^{c}$. As discussed above, this gives
rise to a single vertex space $DM_{v}$ of $DM$ obtained by doubling $M_{v}$
along $\partial_{0}M_{v}$. It is clear that $DM_{v}$ has a natural map to the
$2$--orbifold $DX_{v}$ obtained by doubling $X_{v}$ along $\partial_{0}X_{v}$.
It is easy to see that $\pi_{1}^{orb}(DX_{v})$ is not virtually cyclic, so the
corresponding $V_{0}$--vertex $V$ of $D\Gamma_{n,n+1}^{c}$ is of
$VPCn$--by-Fuchsian type.

If $v$ is a $V_{0}$--vertex of $\Gamma_{n,n+1}^{c}(G)$ of $I$--bundle type, it
is of $VPC(n-1)$--by--Fuchsian type. Let $K$ denote the $VPC(n-1)$ normal
subgroup of $G(v)$. As before, we let $X_{v}$ denote the compact $2$--orbifold
whose orbifold fundamental group is $G(v)/K$, and whose boundary corresponds
to the edges of $\Gamma_{n,n+1}^{c}(G)$ which are incident to $v$. Thus
doubling $M_{v}$ along $\partial_{0}M_{v}$ yields a $V_{0}$--vertex $V$ of
$D\Gamma_{n,n+1}^{c}$ of $VPCn$--by--Fuchsian type with the same base
$2$--orbifold $X_{v}$.

The case of $V_{0}$--vertices of $\Gamma_{n,n+1}^{c}(G)$ of solid torus type
seems to be different, and we will need a separate and more subtle argument.
Let $v$ be a $V_{0}$--vertex of $\Gamma_{n,n+1}^{c}(G)$ of solid torus type.
Thus $G(v)$ is $VPCn$ and has a $VPCn$ subgroup $L$ such that all the annuli
in $\partial_{0}M_{v}$ and in $\partial_{1}M_{v}$ have boundary with
fundamental group $L$.

\begin{lemma}
\label{solidtorustypevertexhasbaseorbifold}Let $v$ be a $V_{0}$--vertex of
$\Gamma_{n,n+1}^{c}(G)$ of solid torus type, and let $L$ be as above. Then
there is a compact $2$--orbifold $X_{v}$, equal to a cone or the quotient of a
cone by a reflection, such that $G(v)$ is $L$--by--$\pi_{1}^{orb}(X_{v})$, and
a natural map from $M_{v}$ to $X_{v}$ such that $\partial M_{v}$ maps onto
$\partial X_{v}$. Further, each annulus in $\partial_{0}M_{v}$ and
$\partial_{1}M_{v}$ maps to a $1$--suborbifold of $\partial X_{v}$ with
non-empty boundary, by the natural map.
\end{lemma}

\begin{remark}
The orbifold fundamental group of $X_{v}$ must be finite in this case, as $L$
has finite index in $G(v)$. Note also that the definition of solid torus type,
Definition \ref{defnofsolidtorustype}, did not include the statement that
$L$\ is normal in $G(v)$. This seems to be quite nontrivial.
\end{remark}

\begin{proof}
As before there is an excision isomorphism between $H_{n+2}(M_{v},\partial
M_{v})$ and $H_{n+2}(M,\partial M)\cong\mathbb{Z}$. As $G(v)$ is $VPCn$, we
know that $H_{n+1}(M_{v})$ is zero. Now the long exact sequence of the pair
$(M_{v},\partial M_{v})$ shows that $\partial M_{v}$ consists of a single
torus $T$. This torus is a union of annuli in $\partial_{0}M_{v}$ and
$\partial_{1}M_{v}$, and all these annuli have boundary with fundamental group
$L$. Thus $L$ is a normal subgroup of $\pi_{1}(T)$ with quotient $\mathbb{Z}$
or $\mathbb{Z}_{2}\ast\mathbb{Z}_{2}$. However the inclusion of $T$ into
$M_{v}$ cannot induce an injective map of fundamental groups, as $G(v)$ is
$VPCn$, but $\pi_{1}(T)$ is $VPC(n+1)$. In this case, we need a more
complicated argument to find the base orbifold for $v$.

Let $M_{V}$ denote the double of $M_{v}$ along $\partial_{0}M_{v}$. The double
of each annulus in $\partial_{1}M_{v}$ is a torus component of $\partial
M_{V}$. Recall that any torus in $DG$ is enclosed by some $V_{0}$--vertex of
$\Gamma_{n+1}^{c}(DG)$. As pairs of these tori in $\partial M_{V}$ are joined
by an annulus in $\partial_{0}M_{v}$, it follows from the proof of Lemma
\ref{torienclosedbydistinctverticesmustlieinone} that all components of
$\partial M_{V}$ are enclosed by a single $V_{0}$--vertex $W$ of $\Gamma
_{n+1}^{c}(DG)$. Now let $A$ be an annulus in $\partial_{0}M_{v}$. As $A$ is
an annulus with ends in $W$, which is a $V_{0}$--vertex of $\Gamma_{n+1}%
^{c}(DG)$, it follows that $A$ cannot cross any component of $\partial M_{W}$.
Thus either $V$ itself is enclosed by $W$, or some component $T$ of $\partial
M_{W}$ is enclosed by $V$ and is not peripheral in $V$. Suppose there is such
a component $T$ of $\partial M_{W}$. As $T$ cannot cross any annulus in
$\partial_{0}M_{v}$, it follows that $T$ is enclosed by $v$ or $\tau v$, which
is impossible, as $G(v)$ is $VPCn$. We conclude that $V$ must be enclosed by
$W$. If $W$ is isolated, it follows that $V\ $is also isolated and hence that
$v$ must be isolated, contradicting the assumption that $v$ is of solid torus
type. If $W$ is of special Seifert type, it follows that $V$ must be isolated
or also of special Seifert type. We have just shown that $V$ cannot be
isolated, so it must be of special Seifert type. As $v$ is of solid torus
type, the only possibility is that it is of special solid torus type, with a
single incident edge carrying a subgroup of index $2$ in $G(v)$. In this case,
we can take $X_{v}$ to be a cone with cone point labeled $2$, and with
$\partial_{0}X_{v}$ and $\partial_{1}X_{v}$ each consisting of a single arc in
$\partial X_{v}$. If $W$ is $VPCn$--by--Fuchsian, it follows that $L$ is the
normal $VPCn$ subgroup of $G(W)$, and that $G(V)$ is the pre-image in $G(W)$
of a suborbifold $X$ of the base orbifold $X_{W}$ of $W$. Now we can choose a
projection map $p:M_{V}\rightarrow X$, so that $\partial M_{V}$ maps to
$\partial X$, and each annulus in $\partial_{1}M_{v}$, and its translate by
$\tau$, projects in the natural way to a $1$--orbifold contained in $\partial
X$. The involution $\tau$ of $DM$ induces an involution of $G(V)$ which is the
identity on $L$, and so induces a proper homotopy equivalence of $X$ of order
$2$. The Nielsen Realization Theorem \cite{Kerckhoff} implies that there is an
involution of $X$ in the given homotopy class, which we again denote by $\tau
$. Recall that each component of $\partial M_{V}$ is the double of an annulus
in $\partial_{1}M_{v}$. It follows that $\tau$ induces an involution on each
component of $\partial M_{V}$ which interchanges the two annuli. Hence the
involution $\tau$ on $X$ acts by a reflection on each component of $\partial
X$, and so fixes one or two points of each such component. We now need to
consider how $p:M_{V}\rightarrow X$ maps each annulus $\Lambda$ of
$\partial_{0}M_{v}$ into $X$. We already know that each component of
$\partial\Lambda$ is also a component of the boundary of an annulus in
$\partial_{1}M_{v}$, and so is mapped to a point of $\partial X$ fixed by
$\tau$.

If $\Lambda$ is untwisted, it follows that we can homotop $p$ restricted to
$\Lambda$ to factor through the natural projection of $\Lambda$ to the unit
interval. This yields a path $\lambda$ in $X$ joining two points of $\partial
X$ fixed by $\tau$. We now give $X$ a hyperbolic metric such that $\partial X$
consists of geodesics and $\tau$ is an isometry, and then homotop $\lambda$
rel $\partial\lambda$ to the unique geodesic in its homotopy class. As
$\tau\Lambda=\Lambda$, it follows that $\tau\lambda$ is homotopic rel
$\partial\lambda$ to $\lambda$. Now the uniqueness of hyperbolic geodesics in
a homotopy class implies that $\tau\lambda=\lambda$. It follows that $\lambda$
is contained in the fixed set of $\tau$, and that $\lambda$ must be a simple geodesic.

If $\Lambda$ is twisted, then $p$ maps $\partial\Lambda$ to a point $a$ of
$\partial X$ fixed by $\tau$. As $\pi_{1}(\Lambda)$ contains $\pi_{1}%
(\partial\Lambda)=L$ with index $2$, it follows that $p_{\ast}\pi_{1}%
(\Lambda)$ is a subgroup of $\pi_{1}^{orb}(X)$ of order $2$. Such a subgroup
of $\pi_{1}^{orb}(X)$ must be carried by a mirror $m$ of $X$ or by a cone
point $w$ of $X$ with even number attached. In the first case, we can homotope
the map from $\Lambda$ to $X$ to have image a path $\lambda$ joining $a$ to
$m$, and then further homotop $\lambda$ to be the shortest geodesic in its
homotopy class. As $\tau$ fixes $\Lambda$, it follows that $\tau$ must
preserve $m$, and that $\tau\lambda$ is homotopic to $\lambda$ fixing $a$ and
keeping the other end of $\lambda$ in $m$. Now the uniqueness of hyperbolic
geodesics in a homotopy class implies that $\tau\lambda=\lambda$. It follows
that $\lambda$ is contained in the fixed set of $\tau$, and that $\lambda$
must be a simple geodesic. In the second case, we can homotope the map from
$\Lambda$ to $X$ to have image a path $\lambda$ joining $a$ to $w$, and then
further homotop $\lambda$ to be the shortest geodesic in its homotopy class.
As $\tau$ fixes $\Lambda$, it follows that $\tau$ must fix $w$, and that
$\tau\lambda$ is homotopic to $\lambda$ rel $\partial\lambda$. Now the
uniqueness of hyperbolic geodesics in a homotopy class implies again that
$\tau\lambda=\lambda$, so that $\lambda$ must be contained in the fixed set of
$\tau$, and $\lambda$ must be a simple geodesic. In both cases, the image of
$\lambda$ is a $1$--suborbifold of $X$, isomorphic to $Q$, the quotient of the
unit interval by a reflection.

At this point, we have arranged that each annulus of $\partial_{0}M_{v}$ maps
to a "simple arc" in $X$ which is contained in the fixed set of $\tau$, so
that distinct such arcs cannot cross. It is conceivable that two of these arcs
coming from twisted annuli could share an endpoint at a cone point of $X$, but
that would imply the two twisted annuli in question carried the same subgroup
of $DG$, which is not possible. It follows that the image of $\partial M_{v}$
in $X$ is a connected $1$--orbifold $C$, without boundary, embedded in $X$.
Thus $C$ is a circle or the quotient of a circle by reflection. The image of
$\pi_{1}^{orb}(C)$ in $\pi_{1}^{orb}(X)$ must be finite, as it is contained in
$\pi_{1}(M_{v})/L$, which is the image of $\pi_{1}(M_{v})$. It follows that
$C$ bounds a suborbifold $Z$ of $X$ with finite orbifold fundamental group,
which can only be a cone or the quotient of a cone by a reflection. This is
the required base orbifold $X_{v}$ for $v$.
\end{proof}

Now we come to an important result about the $V_{1}$--vertices of
$D\Gamma_{n,n+1}^{c}$.

\begin{lemma}
\label{V1-vertexofDGammaisatoroidal}Let $(G,\partial G)$ be an orientable
$PD(n+2)$ pair, such that $G$ is not $VPC$. Let $V$ be a $V_{1}$--vertex of
$D\Gamma_{n,n+1}^{c}$, denote $G(V)$ by $K$, and let $\partial K$ denote the
family of subgroups of $K$ associated to the edges of $D\Gamma_{n,n+1}^{c}$
incident to $V$. Then $(K,\partial K)$ is an orientable atoroidal $PD(n+2)$ pair.
\end{lemma}

\begin{proof}
The fact that $(K,\partial K)$ is an orientable $PD(n+2)$ pair follows from
\cite{B-E}. Now suppose that $T$ is a torus in $(K,\partial K)$. We need to
show that $T$ is conjugate into a group in $\partial K$. Recall from the above
discussion of $DM\ $and $D\Gamma_{n,n+1}^{c}$, that $V$ must be obtained from
a $V_{1}$--vertex $v$ of $\Gamma_{n,n+1}^{c}(G)$. If $V$ is isolated, the
required result is immediate. If $V$ is equal to a $V_{1}$--vertex $v$ of
$\Gamma_{n,n+1}^{c}(G)$, then $T$ is a torus in $G$. Now any torus in $G$ is
conjugate into a group in $\partial G$ or is essential and so enclosed by some
$V_{0}$--vertex of $\Gamma_{n,n+1}^{c}(G)$. In either case, it must be
conjugate into a group in $\partial K$. If $V$ is obtained by doubling from
some $V_{1}$--vertex $v$ of $\Gamma_{n,n+1}^{c}(G)$, there are three possible
cases up to conjugacy, that $T\ $is an essential torus in $G$ or $\overline
{G}$, that $T$ is contained in a group in $\partial G$, or that $T$ is
decomposed into essential annuli lying in $(G,\partial G)$ or $(\overline
{G},\partial\overline{G})$ and enclosed by $v$ or by $\tau v$. The first case
is again not possible. In the second case, $T$ would be peripheral in $V$, as
required. Now we consider the third case in which $T$ is decomposed into
essential annuli lying in $(G,\partial G)$ or $(\overline{G},\partial
\overline{G})$ and enclosed by $v$ or by $\tau v$. We note that any essential
annulus in $(G,\partial G)$ is enclosed by some $V_{0}$--vertex of
$\Gamma_{n,n+1}^{c}(G)$. Thus if an essential annulus in $(G,\partial G)$ is
enclosed by a $V_{1}$--vertex $v$ of $\Gamma_{n,n+1}^{c}(G)$, it must be a
cover of an annulus associated to an edge of $\Gamma_{n,n+1}^{c}(G)$ incident
to $v$. Thus the annuli into which $T$ is decomposed are all covers of edge
annuli in $v$ or $\tau v$. If $v$ is a non-isolated $V_{1}$--vertex of
$\Gamma_{n,n+1}^{c}(G)$, it is not possible to have two distinct edge annuli
whose boundaries carry the same group $L$. It follows that $T$ is a subgroup
of the double of a single edge annulus of $v$, so that $T$ is a subgroup of an
edge torus of $V$, as required.
\end{proof}

Next we discuss more carefully the $V_{0}$--vertices of $D\Gamma_{n,n+1}$.
These are all obtained from $V_{0}$--vertices of $\Gamma_{n,n+1}(G)$ except
for those isolated vertices obtained from torus groups in $\partial G$
enclosed by a $V_{1}$--vertex of $\Gamma_{n,n+1}(G)$.

If $V$ is obtained from a $V_{0}$--vertex $v$ of $\Gamma_{n,n+1}(G)$ of
interior type, then $V$ is equal to $v$ or $\tau v$, so is isolated, of
special Seifert type or of interior Seifert type. In the last case, if $L$
denotes the $VPCn$ normal subgroup of $G(V)$ with Fuchsian quotient, then $L$
is contained in $G$ or $\overline{G}$, but is not contained in any group in
$\partial G$.

If $V$ is obtained from a $V_{0}$--vertex $v$ of $\Gamma_{n,n+1}(G)$ of
peripheral type, then $M_{V}$ is the double of $M_{v}$ along $\partial
_{0}M_{v}$. If $v$ is isolated, then $V$ is also isolated. If $v$ is of
commensuriser type, there is a $VPCn$ subgroup $L$ of $G$ such that
$G(v)=Comm_{G}(L)$ and all edge annuli of $v$ have boundary which carries $L$.
Thus $L$ is contained in $G$ and also contained in groups in $\partial G$.
Further Lemma \ref{torustypevertexhasbaseorbifold}, Corollary
\ref{Seiferttypevertexhasbaseorbifold} and Lemma
\ref{solidtorustypevertexhasbaseorbifold} together show that $G(v)$ is
$L$--by--$\pi_{1}^{orb}(X_{v})$, where $X_{v}$ is the base orbifold of $v$.
Hence $G(V)$ is $L$--by--$\pi_{1}^{orb}(DX_{v})$, where $DX_{v}$ is the double
of $X_{v}$ along $\partial_{0}X_{v}$. In almost all cases $\pi_{1}%
^{orb}(DX_{v})$ is not virtually cyclic, and $V$ is of $VPCn$--by--Fuchsian
type, where the normal $VPCn$ subgroup is again $L$. The only exception occurs
when $v$ is of special solid torus type, with a single incident edge carrying
a subgroup of $G(v)$ of index $2$. In this case, $V$ is of special Seifert type.

Finally if $V$ is obtained from a $V_{0}$--vertex $v$ of $\Gamma_{n,n+1}(G)$
of $I$--bundle type, then $V$ is of Seifert type, but this time the normal
$VPCn$ subgroup $L$ of $G(V)$ is not conjugate into $G$ or $\overline{G}$.

Now we can state our doubling result which is precisely analogous to the
situation in $3$--manifold theory.

\begin{theorem}
\label{doublingtheorem}Let $(G,\partial G)$ be an orientable $PD(n+2)$ pair,
such that $G$ is not $VPC$. Then the double $D\Gamma_{n,n+1}^{c}$ of
$\Gamma_{n,n+1}^{c}(G)$ is equal to the decomposition $\Gamma_{n+1}^{c}(DG)$
of $DG$.
\end{theorem}

\begin{remark}
\label{doublingremark}The double $D\Gamma_{n,n+1}$ of $\Gamma_{n,n+1}(G)$ need
not be equal to the decomposition $\Gamma_{n+1}(DG)$ of $DG$. For example,
consider the case when $G$ is the fundamental group of a $3$--manifold $M$,
and $\Gamma_{n,n+1}^{c}(G)$ has a $V_{0}$--vertex $v$ of special solid torus
type, with three incident edges each carrying $G(v)$. Thus the characteristic
submanifold of $M$ has a component $X$ which is a solid torus which meets
$\partial M$ in three annuli all of degree $1$ in $X$. Doubling $M$ yields a
component $DX$ of the characteristic submanifold of $DM$, which is the double
of $X$ along $X\cap\partial M$. Thus $DX$ is the product of a pair of pants
with the circle. Now $X$ corresponds to a $V_{1}$--vertex of $\Gamma
_{n,n+1}(G)$, but $DX$ corresponds to a $V_{0}$--vertex of $\Gamma_{n+1}(DG)$.
\end{remark}

\begin{proof}
The outline of our proof is to show that each $V_{0}$--vertex of
$D\Gamma_{n,n+1}^{c}$ is enclosed by some $V_{0}$--vertex of $\Gamma_{n+1}%
^{c}(DG)$, and that each $V_{0}$--vertex of $\Gamma_{n+1}^{c}(DG)D\Gamma
_{n,n+1}^{c}$ is enclosed by some $V_{0}$--vertex of $D\Gamma_{n,n+1}^{c}$.
Assuming these two facts, we can deduce the theorem as follows. Let $V$ be a
$V_{0}$--vertex of $D\Gamma_{n,n+1}^{c}$ enclosed by the $V_{0}$--vertex $W$
of $\Gamma_{n+1}^{c}(DG)$. As $W$ is enclosed by a $V_{0}$--vertex $V^{\prime
}$ of $D\Gamma_{n,n+1}^{c}$, it follows that $V$ is enclosed by $V^{\prime}$.
As $D\Gamma_{n,n+1}^{c}$ is reduced, it follows that $V=V^{\prime}$, and hence
that $G(V)=G(W)$. Let $\partial G(V)$ denote the family of subgroups of $G(V)$
associated to the edges of $D\Gamma_{n,n+1}^{c}$ incident to $V$, and
similarly for $W$. The facts that $V$ is enclosed by $W$, and $G(V)=G(W)$
implies that the $PD(n+2)$ pairs $(G(V),\partial G(V))$ and $(G(W),\partial
G(W))$ are isomorphic. Similarly if $W$ is a $V_{0}$--vertex of $\Gamma
_{n+1}^{c}(DG)$ enclosed by a $V_{0}$--vertex $V$ of $D\Gamma_{n,n+1}^{c}$, it
follows that the $PD(n+2)$ pairs $(G(V),\partial G(V))$ and $(G(W),\partial
G(W))$ are isomorphic. Together these facts imply that $D\Gamma_{n,n+1}^{c}$
is equal to $\Gamma_{n+1}^{c}(DG)$, as required.

First we will show that each $V_{0}$--vertex of $D\Gamma_{n,n+1}^{c}$ is
enclosed by some $V_{0}$--vertex of $\Gamma_{n+1}^{c}(DG)$.

Our discussion immediately before this theorem shows that if $V$ is a $V_{0}%
$--vertex of $D\Gamma_{n,n+1}^{c}$, it is isolated, of special Seifert type,
or of $VPCn$--by--Fuchsian type. Now any torus in $DG$ is enclosed by some
$V_{0}$--vertex of $\Gamma_{n+1}^{c}(DG)$, and crossing tori must be enclosed
by the same $V_{0}$--vertex of $\Gamma_{n+1}^{c}(DG)$. Thus if $V$ is not of
special Seifert type, it follows that $V$ is enclosed by some $V_{0}$--vertex
of $\Gamma_{n+1}^{c}(DG)$. If $V$ is of special Seifert type, we know that
$G(V)$ is conjugate into some $V_{0}$--vertex group of $\Gamma_{n+1}^{c}(DG)$,
so that again $V$ is enclosed by some $V_{0}$--vertex of $\Gamma_{n+1}%
^{c}(DG)$. Thus each $V_{0}$--vertex of $D\Gamma_{n,n+1}^{c}$ is enclosed by
some $V_{0}$--vertex of $\Gamma_{n+1}^{c}(DG)$, as required.

It remains to show that each $V_{0}$--vertex $W$ of $\Gamma_{n+1}%
^{c}(DG)D\Gamma_{n,n+1}^{c}$ is enclosed by some $V_{0}$--vertex of
$D\Gamma_{n,n+1}^{c}$.

If $W$ is isolated, let $T$ denote an edge torus of $W$. As $T$ is an edge
torus of $\Gamma_{n+1}^{c}(DG)$, it crosses no torus in $DG$, and so is
enclosed by some vertex of $D\Gamma_{n,n+1}^{c}$. Lemma
\ref{V1-vertexofDGammaisatoroidal} tells us that $V_{1}$--vertices of
$D\Gamma_{n,n+1}^{c}$ are atoroidal, so it follows that \ $T$, and hence $W$,
is enclosed by some $V_{0}$--vertex of $D\Gamma_{n,n+1}^{c}$, as required.

If $W$ is of special Seifert type, we let $T$ denote the edge torus of $W$.
Again $T$ is enclosed by some $V_{0}$--vertex $V$ of $D\Gamma_{n,n+1}^{c}$. As
$G(W)$ contains $T$ with index $2$, it follows that $W$ is enclosed by some
vertex $V^{\prime}$ of $D\Gamma_{n,n+1}^{c}$. If $V^{\prime}$ is a $V_{1}%
$--vertex, it follows that $T\ $is peripheral in both $V$ and in $V^{\prime}$.
Now Lemma \ref{boundarygroupismaximal} shows that $V^{\prime}$ must be of
special Seifert type. As no $V_{1}$--vertex of $\Gamma_{n,n+1}^{c}(G)$ can be
of special Seifert type, $V^{\prime}$ must be obtained by doubling from a
$V_{1}$--vertex $v$ of $\Gamma_{n,n+1}^{c}(G)$. This implies that $v$ must be
of special solid torus type. As no $V_{1}$--vertex of $\Gamma_{n,n+1}^{c}(G)$
can be of special solid torus type, it follows that $V^{\prime}$ is a $V_{0}%
$--vertex of $D\Gamma_{n,n+1}^{c}$ which encloses $W$, as required.

For the rest of this proof, we will assume that $W$ is of $VPCn$--by--Fuchsian
type with normal $VPCn$ subgroup $L$. Thus Lemma
\ref{normaliserequalscommensuriser} tells us that $G(W)=N_{DG}(L)$.

Let $T$ denote the universal covering $DG$--tree of the graph of groups
$\Delta$ determined by the doubling of $G$ along $\partial G$. Thus $\Delta$
has two vertices with associated groups $G$ and $\overline{G}$, and has edges
corresponding to the groups in $\partial G$. We consider the actions of $L$
and $G(W)=N_{DG}(L)$ on $T$, and the various cases which arise.

Case 1: $L$ fixes a vertex $z$ of $T$.

By a conjugation and possibly interchanging $G$ and $\overline{G}$, we can
assume that $G(z)=G$. In particular, $L\subset G$. Let $T^{\prime}$ denote the
subtree of $T$ consisting of all edges and vertices of $T$ fixed by $L$. The
action of $N_{DG}(L)$ on $T$ must preserve $T^{\prime}$.

Case 1a): $T^{\prime}=\{z\}$.

Thus $N_{DG}(L)$ fixes $z$. This implies that $N_{DG}(L)\subset G$. In
particular, any torus enclosed by $W$ is also enclosed by the vertex $U$ of
$\Delta$ with associated group $G$. As $W$ is filled by crossing tori, it
follows that $W$ itself is enclosed by $U$, and hence that $W$ is enclosed by
some $V_{0}$--vertex $Z$ of $\Gamma_{n,n+1}^{c}(G)$. Let $\overline{Z}$ denote
the corresponding $V_{0}$--vertex of $D\Gamma_{n,n+1}^{c}$, so that
$\overline{Z}$ is either equal to $Z$ or obtained by doubling $Z$. In either
case, $W$ is enclosed by $\overline{Z}$, as required.

Case 1b): $z$ has valence $1$ in $T^{\prime}$.

This implies that there are no essential annuli in $(G,\partial G)$ which
carry $L$, and the same holds for $\overline{G}$. Thus every vertex of
$T^{\prime}$ has valence $1$ in $T^{\prime}$, which implies that $T^{\prime}$
is equal to a single edge $e$ of $T$. As $N_{DG}(L)$ preserves $T^{\prime}$
and cannot interchange the ends of $e$, we deduce that $N_{DG}(L)$ fixes $e$.
Now the stabilizer of $e$ is a group $K$ in $\partial G$ and so is an
orientable $PD(n+1)$ group, and $N_{DG}(L)$ contains a torus subgroup $\Sigma$
which is also $PD(n+1)$. It follows that $\Sigma$ has finite index in $K$, so
that $K$ is $VPC(n+1)$, and hence also a torus. Hence $W$ must be isolated,
contradicting our assumption.

Case 1c): $z$ has valence $2$ in $T^{\prime}$.

This implies that all vertices of $T^{\prime}$ have valence $2$ in $T^{\prime
}$, so that $T^{\prime}$ is a line. As $N_{DG}(L)$ preserves this line, it
follows that there is a map from $N_{DG}(L)$ to $\mathbb{Z}$ or to
$\mathbb{Z}_{2}\ast\mathbb{Z}_{2}$ whose kernel $K$ fixes every point of
$T^{\prime}$. Note that $N_{DG}(L)\cap G=N_{G}(L)$. In the first case,
$N_{G}(L)$ equals $K$, and in the second case, $N_{G}(L)$ contains $K\ $with
index $2$. As $L$ fixes only two edges of $T$ incident to $z$, there is a
unique essential annulus in $(G,\partial G)$ which carries $L$. It follows
that there is a $V_{0}$--vertex $v$ of $\Gamma_{n,n+1}^{c}(G)$, such that
$G(v)=N_{G}(L)$, and that $v$ is isolated or of special solid torus type, with
a single incident edge carrying a subgroup of index $2$. Let $V$ denote the
$V_{0}$--vertex of $D\Gamma_{n,n+1}^{c}$ obtained by doubling $v$. Then $G(V)$
contains $N_{G}(L)$, and contains an element which acts on $T^{\prime}$ by a
translation of length $2$. Hence $G(V)$ is equal to $G(W)$. In particular, $W$
is enclosed by the $V_{0}$--vertex $V$ of $D\Gamma_{n,n+1}^{c}$, as required.

Case 1d): $z$ has valence at least $3$ in $T^{\prime}$.

As $L$ fixes three distinct edges of $T$ incident to $z$, it follows that
there are three distinct annuli in $(G,\partial G)$ whose boundaries carry
$L$. Hence there is a $V_{0}$--vertex $v$ of $\Gamma_{n,n+1}^{c}(G)$ of
special solid torus type or of commensuriser type, such that $G(v)=Comm_{G}%
(L)=$ $N_{G}(L)=N_{DG}(L)\cap G$. Hence for each edge $e$ of $T^{\prime}$
incident to $v$, with stabilizer $S$, we have $N_{DG}(L)\cap S=G(v)\cap S$.
Let $V$ denote the double of $v$, so that $G(V)\subset N_{DG}(L)=G(W)$. The
quotient of $T^{\prime}$ by $G(V)$ has two vertices, so the same holds for the
quotient of $T^{\prime}$ by $G(W)$. As both groups act on $T^{\prime}$ with
the same edge and vertex stabilizers, it follows that $G(V)=G(W)$. In
particular, $W$ is enclosed by the $V_{0}$--vertex $V$ of $D\Gamma_{n,n+1}%
^{c}$, as required.

Case 2: $L$ fixes no vertex of $T$.

In this case, there is a unique minimal $L$--subtree $\Lambda$ of $T$. As
$L\ $is $VPC$, Lemma \ref{splittingsofVPCgroups} tells us that $\Lambda$ is a
line. The uniqueness of $\Lambda$ implies that $N_{DG}(L)$ preserves $\Lambda
$. As in Case 1c), there is a map from $N_{DG}(L)$ to $\mathbb{Z}$ or to
$\mathbb{Z}_{2}\ast\mathbb{Z}_{2}$ whose kernel $K$ fixes every point of
$\Lambda$. And $N_{DG}(L)\cap G$ is either equal to $K$ or contains $K$ with
index $2$. As $G(W)=N_{DG}(L)$, we know that any torus $\Sigma$ enclosed by
$W$ must intersect $L$ in a subgroup of finite index in $L$. Hence the
$VPC(n+1)$ group $\Sigma$ acts on $\Lambda$, and $\Sigma\cap K$ is a $VPCn$
subgroup of $\Sigma$ which fixes the vertex $z$ and preserves the two edges of
$\Lambda$ incident to $z$. This determines an essential annulus in
$(G,\partial G)$ which carries $\Sigma\cap K$. Every torus enclosed by $W$
determines an essential annulus in $(G,\partial G)$ preserving the same two
edges of $\Lambda$ incident to $z$. As $W$ is not isolated or of special
Seifert type, there are infinitely many distinct such annuli, so that
$N_{DG}(L)\cap G=G(v)$, for some $V_{0}$--vertex $v$ of $\Gamma_{n,n+1}%
^{c}(G)$ of $I$--bundle type. Let $V$ be obtained by doubling $v$. As $G(V)$
contains $N_{DG}(L)\cap G$ and acts on $\Lambda$ by a translation of length
$2$, it follows that $G(V)=G(W)$. In particular, $W$ is enclosed by the
$V_{0}$--vertex $V$ of $D\Gamma_{n,n+1}^{c}$, as required.
\end{proof}

The machinery which we have just developed can be used to give simple proofs
of some other results. For example, we have the following result.

\begin{theorem}
Let $(G,\partial G)$ be an orientable $PD(n+2)$ pair, such that $G$ is not
$VPC$. Let $\mathcal{F}$ denote the family of all essential annuli in
$(G,\partial G)$, and let $\mathcal{F}_{n}$ denote the family of equivalence
classes of all nontrivial almost invariant subsets of $G$ which are over a
$VPCn$ subgroup. Then $\Gamma(\mathcal{F}:G)$, the algebraic regular
neighbourhood of $\mathcal{F}$ in $G$, is equal to $\Gamma_{n}(G)$, the
algebraic regular neighbourhood of $\mathcal{F}_{n}$ in $G$.
\end{theorem}

\begin{proof}
An essential annulus $A$ in $(G,\partial G)$ is enclosed by some $V_{0}%
$--vertex $v$ of $\Gamma_{n}(G)$. The discussion in this section shows that if
$v$ is of $I$--bundle type, then $A$ determines a loop in the base orbifold of
$v$, and if $v$ is of commensuriser type, then $A$ determines an arc in the
base orbifold of $v$. Further this yields a bijection between equivalence
classes of annuli in $(G,\partial G)$ and loops and arcs in base orbifolds of
$V_{0}$--vertices of $\Gamma_{n}(G)$. As any compact $2$--orbifold is filled
by essential (possibly singular) loops, and is also filled by essential
(possibly singular) arcs, it follows that $\Gamma(\mathcal{F}:G)$ is equal to
$\Gamma_{n}(G)$, as required.
\end{proof}

\section{Concluding remarks and problems\label{concludingremarks}}

In this section, we briefly discuss some problems which arise from our work in
this paper.

In the previous section, we showed, in Theorem \ref{doublingtheorem}, that if
$(G,\partial G)$ is an orientable $PD(n+2)$ pair, such that $G$ is not $VPC$,
then the double $D\Gamma_{n,n+1}^{c}$ of $\Gamma_{n,n+1}^{c}(G)$ is equal to
the decomposition $\Gamma_{n+1}^{c}(DG)$ of $DG$. This precisely mirrors the
situation in $3$--manifold theory. But in that theory, one can deduce the JSJ
decomposition theorem for $3$--manifolds with boundary from the corresponding
result for closed manifolds by "undoubling". This is a great simplification of
the direct proofs. Our first problem is to decide whether an analogous
argument is possible in the setting of $PD(n+2)$ pairs.

\begin{problem}
Can the main result of this paper, Theorem \ref{mainresult}, be deduced from
the properties of the decomposition $\Gamma_{n+1}^{c}(DG)$ of $DG$, by some
"undoubling" argument?
\end{problem}

Our main result in this paper, Theorem \ref{mainresult}, shows that the
situation for $PD(n+2)$ pairs is very similar to that for $3$--manifolds with
boundary. Now in $3$--manifold theory, the characteristic submanifold of a
Haken $3$--manifold has the enclosing property for Seifert pairs. Thus it is
reasonable to ask the following.

\begin{problem}
Do our results in this paper imply a result for $PD3$ pairs analogous to the
enclosing property for Seifert pairs in $3$--manifolds.
\end{problem}

This is Conjecture 10.4 of Wall's survey article \cite{Wall3}. However we are
unable to answer this question. Even a precise formulation of the statement
requires a theory of Poincar\'{e} duality triads, and/or of Poincar\'{e}
duality pairs with compressible boundary, neither of which has been developed
so far.

Since one of the main results in the JSJ\ theory of $3$--manifolds is
Johannson's Deformation Theorem \cite{JO}, it is natural to ask the following.

\begin{problem}
Is there a result for Poincar\'{e} duality pairs which is analogous to
Johannson's Deformation Theorem?
\end{problem}

Some of our discussion and definitions have been rather topological in order
to suit the setting of group pairs and duality. A more algebraic description
of the decomposition $\Gamma_{n,n+1}^{c}(G)$ should involve Poincar\'{e}
duality triads rather than pairs, since a vertex group of $\Gamma_{n,n+1}%
^{c}(G)$ has two distinct important types of subgroups, namely the edge groups
and the intersection groups with $\partial G$. We feel that a reformulation in
these terms and a strengthening of the statements would be necessary in order
to formulate the analogue of Johannson's Deformation Theorem. In \cite{Wall2},
Wall studied the notion of triads of Poincar\'{e} complexes. It is natural to
ask whether our algebraic decompositions can be realised using finite
aspherical Poincar\'{e} complexes. Suppose that $G$ is the fundamental group
of such a complex $X$ of dimension $n+2$ and that $G$ splits over a $VPCn$
group $H$. A basic question is whether $X$ splits over a subcomplex $Y$ with
$\pi_{1}(Y)=H$. One expects to have to change the complex $X$ to achieve this,
so it seems better to ask whether there is a finite complex $X^{\prime}$ of
the same simple homotopy type as $X$, such that $X^{\prime}$ splits over a
suitable subcomplex $Y$. This involves the study of some obstruction groups
most studied by Waldhausen \cite{Waldhausenalgtop}. For a general torsion free
$VPCn$ group $H$, these obstruction groups are not well understood, so the
problem seems difficult. But if we assume that $H$ is torsion free and
polycyclic, then Waldhausen's work shows that the obstruction groups are zero,
so that we can find such a complex $X^{\prime}$ as desired. However, even if
we try to restrict our attention to almost invariant subsets of $G$ over
torsion free polycyclic subgroups, it seems possible that the edge groups of
the decompositions we obtain may be $VPC$. If the edge groups are polycyclic,
as happens in the $3$--dimensional case, then the analogue of Johannson's
Deformation Theorem can be formulated. We expect this analogue to be correct
for $3$--dimensional Poincar\'{e} duality pairs, and finite aspherical
Poincar\'{e} $3$--complexes. In general it may be true only when the
peripheral groups are polycyclic.

\end{document}